\documentclass{article}

\usepackage{a4wide,cite,latexsym,amsfonts,amssymb,exscale,enumerate}
\usepackage{amsmath,amsthm}

\usepackage{pstricks}

\psset{linewidth=0.3pt,dimen=middle}
\psset{xunit=.70cm,yunit=0.70cm}

\usepackage[all]{xy}
\SelectTips{cm}{}

\usepackage{graphicx}

\newcommand{\E}[1]{\cal{E}^{( #1 )}}
\newcommand{\F}[1]{\cal{F}^{( #1 )}}
\newcommand{\onen}{{\mathbf 1}_{n}}
\newcommand{\onenn}[1]{{\mathbf 1}_{#1}}
\newcommand{\onenp}{{\mathbf 1}_{n'}}

\newcommand{\BNC}{\cal{N}\cal{H}}

\newcommand{\ep}{\underline{\epsilon}}

\newcommand{\Uq}{{\bf U}_q(\mathfrak{sl}_2)}
\newcommand{\U}{\dot{{\bf U}}}
\newcommand{\UA}{{_{\cal{A}}\dot{{\bf U}}}}
\newcommand{\Ucat}{\cal{U}}

\newcommand{\UcatD}{\dot{\cal{U}}}

\newcommand{\B}{\dot{\mathbb{B}}}

\newcommand{\Sq}{{\rm Sq}}
\newcommand{\und}[1]{\underline{#1}}

\newcommand{\Ucross}{\xy {\ar (2.5,-2.5)*{};(-2.5,2.5)*{}}; {\ar (-2.5,-2.5)*{};(2.5,2.5)*{} };
(4,0)*{};(-4,0)*{};\endxy}

\newcommand{\qbin}[2]{
\left[
 \begin{array}{c}
 #1 \\
 #2 \\
 \end{array}
 \right]
}

\newcommand{\xsum}[2]{
  \xy
  (0,.4)*{\sum};
  (0,3.7)*{\scs #2};
  (0,-2.9)*{\scs #1};
  \endxy
}

\newcommand{\refequal}[1]{\xy {\ar@{=}^{#1}
(-1,0)*{};(1,0)*{}};
\endxy}

\usepackage{fancyheadings}
\newcommand{\cat}[1]{\ensuremath{\mbox{\bfseries {\upshape {#1}}}}}

\newcommand{\To}{\Rightarrow}

\newcommand{\Hom}{{\rm Hom}}
\newcommand{\HOM}{{\rm HOM}}
\renewcommand{\to}{\rightarrow}
\newcommand{\maps}{\colon}

\newcommand{\bigb}[1]{
\begin{pspicture}(0,0)
 \rput(0,0){\psframebox[framearc=.5,fillstyle=solid]{\small $#1$}}
\end{pspicture}}

\newcommand{\END}{{\rm END}}

\newcommand{\im}{{\rm im\ }}

\newcommand{\rkq}{{\rm rk}_q}

\def\sgn{\mathop{\rm sgn}}
\def\det{\mathop{\rm det}}

\newcommand{\scs}{\scriptstyle}

\def\Id{\mathrm{Id}}

\def\mf{\mathfrak}
\def\shuffle{\,\raise 1pt\hbox{$\scriptscriptstyle\cup{\mskip
               -4mu}\cup$}\,}

\theoremstyle{definition}
\newtheorem{thm}{Theorem}[section]
\newtheorem{cor}[thm]{Corollary}

\newtheorem{lem}[thm]{Lemma}
\newtheorem{rem}[thm]{Remark}
\newtheorem{prop}[thm]{Proposition}
\newtheorem{defn}[thm]{Definition}

\numberwithin{equation}{section}

\def\emph#1{{\sl #1\/}}

\let\hat=\widehat
\let\tilde=\widetilde

\let\phi=\varphi
\let\theta=\vartheta

\usepackage{bbm}

\def\N{{\mathbbm N}}

\def\Z{{\mathbbm Z}}
\def\Q{{\mathbbm Q}}

\def\cal#1{\mathcal{#1}}%
\def\1{\mathbbm{1}}%
\def\nn{\notag}

\newcommand{\stccbub}[2]{\xybox{
 (0,0)*{\includegraphics[scale=0.5]{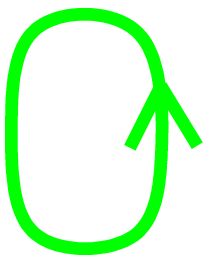}};
 (-5,0)*{\bigb{\pi_{#2}^{\spadesuit}}};(3,-7)*{#1};
}}

\newcommand{\stcbub}[2]{\xybox{
 (0,0)*{\reflectbox{\includegraphics[scale=0.5]{figs/ccbub.eps}}};
 (5,0)*{\bigb{\pi_{#2}^{\spadesuit}}};(-3,-7)*{#1};
}}

\newcommand{\lowrru}[1]{\xybox{%
  (-8,0)*{};
  (8,0)*{};
  (-6,-18)*{};(6,-9)*{} **\crv{(-6,-13) & (6,-15)} ?(1)*\dir{>};
  (6,-9)*{};(6,0)*{}  **\dir{-} ?(.3)*\dir{ }+(2,0)*{\scs {\bf j}};
}}

\newcommand{\lowllu}[1]{\xybox{%
  (-8,0)*{};
  (8,0)*{};
  (6,-18)*{};(-6,-9)*{} **\crv{(6,-13) & (-6,-15)} ?(1)*\dir{>};
  (-6,-9)*{};(-6,0)*{}  **\dir{-} ?(.3)*\dir{ }+(-2,0)*{\scs {\bf j}};
}}

\newcommand{\bbe}[1]{\xybox{%
  (-2,0)*{};
  (2,0)*{};
  (0,0);(0,-18) **\dir{-}; ?(.5)*\dir{<}+(2.3,0)*{\scriptstyle{#1}};
}}

\newcommand{\bbsid}{\xybox{%
  (-2,0)*{};
  (2,0)*{};
  (0,10);(0,4) **\dir{-};
}}
\newcommand{\bbpef}[1]{\xybox{%
  (-6,0)*{};
  (6,0)*{};
  (-4,0)*{}="t1";
  (4,0)*{}="t2";
  "t1";"t2" **\crv{(-4,-6) & (4,-6)}; ?(.15)*\dir{>} ?(.9)*\dir{>}
   ?(.5)*\dir{}+(0,-2)*{\scriptstyle{#1}};
}}
\newcommand{\bbpfe}[1]{\xybox{%
  (-6,0)*{};
  (6,0)*{};
  (-4,0)*{}="t1";
  (4,0)*{}="t2";
  "t2";"t1" **\crv{(4,-6) & (-4,-6)}; ?(.15)*\dir{>} ?(.9)*\dir{>}
  ?(.5)*\dir{}+(0,-2)*{\scriptstyle{#1}};
}}

\newcommand{\bbcfe}[1]{\xybox{%
  (-6,0)*{};
  (6,0)*{};
  (-4,0)*{}="t1";
  (4,0)*{}="t2";
  "t1";"t2" **\crv{(-4,6) & (4,6)}; ?(.15)*\dir{>} ?(.9)*\dir{>}
  ?(.5)*\dir{}+(0,2)*{\scriptstyle{#1}};
}}
\newcommand{\bbcef}[1]{\xybox{%
  (-6,0)*{};
  (6,0)*{};
  (-4,0)*{}="t1";
  (4,0)*{}="t2";
  "t2";"t1" **\crv{(4,6) & (-4,6)}; ?(.15)*\dir{>}
  ?(.9)*\dir{>} ?(.5)*\dir{}+(0,2)*{\scriptstyle{#1}};
}}

\newcommand{\ccbub}[2]{
\xybox{%
 (-6,0)*{};
  (6,0)*{};
  (-4,0)*{}="t1";
  (4,0)*{}="t2";
  "t2";"t1" **\crv{(4,6) & (-4,6)}; ?(.7)*\dir{}+(-2,0)*{\scs #2}
  ?(.05)*\dir{>} ?(1)*\dir{>};
  "t2";"t1" **\crv{(4,-6) & (-4,-6)};
   ?(.3)*\dir{}+(0,0)*{\bullet}+(0,-3)*{\scs {#1}};
}}
\newcommand{\cbub}[2]{
\xybox{%
 (-6,0)*{};
  (6,0)*{};
  (-4,0)*{}="t1";
  (4,0)*{}="t2";
  "t2";"t1" **\crv{(4,6) & (-4,6)};?(.7)*\dir{}+(-2,0)*{\scs #2};
   ?(0)*\dir{<} ?(.95)*\dir{<};
  "t2";"t1" **\crv{(4,-6) & (-4,-6)};
   ?(.3)*\dir{}+(0,0)*{\bullet}+(0,-3)*{\scs {#1}};
}}

\newcommand{\bbdl}[1]{\xybox{%
  (2,0);(0,-8) **\crv{(2,-2)&(0,-6)}; ?(.5)*\dir{>}
}}
\newcommand{\bbdlu}[1]{\xybox{%
  (2,0);(0,-8) **\crv{(2,-2)&(0,-6)}; ?(.5)*\dir{<}
}}
\newcommand{\bbdr}[1]{\xybox{%
  (-2,0);(0,-8) **\crv{(-2,-2)&(0,-6)}; ?(.5)*\dir{>}
}}
\newcommand{\bbdru}[1]{\xybox{%
  (-2,0);(0,-8) **\crv{(-2,-2)&(0,-6)}; ?(.5)*\dir{<}
}}

\title{Extended graphical calculus for categorified quantum sl(2)}
      \author{ Mikhail Khovanov, Aaron D. Lauda, Marco Mackaay, and Marko Sto\v si\'c}

\begin{document}

\date{June 14, 2010}

\maketitle

\begin{abstract}
A categorification of the Beilinson-Lusztig-MacPherson form of the quantum sl(2) was constructed in the paper arXiv:0803.3652  by the second author. Here we enhance the graphical calculus introduced and developed in that paper to include two-morphisms between divided
powers one-morphisms and their compositions. We obtain explicit
diagrammatical formulas for the decomposition of products of divided
powers one-morphisms as direct sums of indecomposable one-morphisms;
the latter are in a bijection with the Lusztig canonical basis
elements. These formulas have integral coefficients and imply that
one of the main results of Lauda's paper---identification of the
Grothendieck ring of his 2-category with the idempotented quantum
sl(2)---also holds when the 2-category is defined over the ring of
integers rather than over a field.
\end{abstract}

\tableofcontents

%
\section{Introduction}
%

An idempotented form $\U(\mf{sl}_n)$ of the quantum enveloping algebra of $\mf{sl}_n$ was introduced by Beilinson-Lusztig-MacPherson~\cite{BLM}, who also related it to the geometry of partial flag varieties.  This idempotented form was later extended by Lusztig to the idempotented version $\U(\mf{g})$ of the quantum group $\mathbf{U}_q(\mathfrak{g})$ associated to a Kac-Moody algebra $\mathfrak{g}$~\cite{Lus4}.  Lusztig also constructed a basis in $\U(\mf{g})$ generalizing the canonical basis~\cite{Kas2,Lus1,Lus2} in ${\bf U}^+(\mf{g})$.  Multiplication and comultiplication in $\U(\mf{sl}_2)$ have positive integral coefficients in the Lusztig canonical basis.

These results led to the visionary conjecture of Igor Frenkel that $\U(\mf{g})$
could be categorified at generic $q$ in a purely algebraic and combinatorial fashion.   Frenkel's key insight was that the canonical basis should play a fundamental role, lifting to become a collection of simple objects of some category. One of the motivations for this conjecture was Lusztig's realization~\cite{Lus1,Lus2,Lus3} of the canonical basis for ${\bf U}^+(\mf{g})$ via simple perverse sheaves on Lusztig quiver varieties.  Frenkel's proposal initiated the research program into categorification of quantum groups and their applications. Crane and Frenkel further conjectured the existence of a categorification of $\U(\mathfrak{sl}_2)$ when $q$ was a root of unity (this is still an open problem). Part of their conjecture included potential applications of categorified quantum groups to the construction of 4-dimensional topological quantum field theories\cite{CF}.

A categorification of $\U=\U(\mf{sl}_2)$ at generic $q$ was achieved by the second author in \cite{Lau1} demonstrating that the quantum enveloping algebra of $\mathfrak{sl}_2$ is just a shadow of a much richer algebraic
structure. A $\Bbbk$-linear 2-category $\UcatD$ was defined whose split Grothendieck ring was shown to be isomorphic to the integral version of $\U$. The 2-category $\UcatD$ is the idempotent completion, or Karoubi envelope, of a 2-category $\Ucat$ defined in terms of a graphical calculus.  The objects $n \in \Z$ of $\Ucat$ are parameterized by the weight lattice of $\mf{sl}_2$. For $\ep=\epsilon_1\dots\epsilon_m$ with $\epsilon_1,\dots,\epsilon_m \in \{+,-\}$ the 1-morphisms from $n$ to $n+2\sum_{i=1}^m \epsilon_i 1$ are given by directs sums of 1-morphisms $\cal{E}_{\ep}\onen\{t\} = \cal{E}_{\epsilon_1}\dots \cal{E}_{\epsilon_m}\onen\{t\}$ where $\cal{E}_{+}=\cal{E}$, $\cal{E_{-}}=\cal{F}$ and $t\in \Z$. For $\Bbbk$ a field, the 2-morphisms are given by $\Bbbk$-linear combinations of certain planar diagrams modulo local relations.  It was shown~\cite{Lau1} that the isomorphism classes of indecomposable 1-morphisms (up to grading shift) in $\UcatD$ bijectively correspond to elements in Lusztig's canonical basis.

The 2-category $\UcatD$ was extended to a 2-category $\UcatD(\mf{g})$ associated to an arbitrary root datum by the first two authors~\cite{KL3} generalizing their earlier work categorifying one-half of the quantum group associated to a symmetrizable Kac-Moody algebra~\cite{KL,KL2}.  For root datum associated to $\mathfrak{sl}_n$ it was shown that the 2-category $\UcatD$ categorifies the idempotented form of ${\bf U}_q(\mathfrak{sl}_n)$.

In this paper we extend the graphical calculus to the Karoubi envelope to allow explicit decompositions of 1-morphisms into indecomposable 1-morphisms in $\UcatD$ (Theorems~\ref{thm_cal_EaEb} and \ref{eq_cat_EaFb}). While the 2-category $\Ucat$ is naturally described in terms of a graphical calculus, there was previously no such purely diagrammatic description of its Karoubi envelope $\UcatD$.  As a consequence of the explicit decompositions in this article, we extend the main categorification result in \cite{Lau1} to the case where $\Bbbk=\Z$, see Corollary~\ref{cor_overZ}. We also obtain an explicit basis for the space of 2-morphisms between indecomposable 1-morphisms in Section~\ref{sec_basis}.

The extended graphical calculus reveals further surprising connections between symmetric functions and the combinatorics encoded in the relations of the 2-category $\Ucat$.  In \cite[Proposition 8.2]{Lau1} it was shown that there is an isomorphism
$$
\xymatrix{\HOM_{\Ucat}(\onen,\onen) \ar[r]^-{\cong} & \Lambda,}$$
where $\Lambda$ is the graded ring of symmetric functions in countably many variables.   A 2-morphism in $\Hom_{\Ucat}(\onen,\onen)$ is represented in the graphical calculus by a closed diagram.  Any such diagram can be reduced to a product of non-nested dotted bubbles with the same orientation. These dotted bubbles naturally correspond to the basis of $\Lambda$ given by the complete symmetric functions.  Proposition~\ref{prop_image_thickbub} identifies natural closed diagrams in the 2-category $\UcatD$ corresponding to Schur polynomials.  Reducing these closed diagrams using the relations in the 2-category $\UcatD$ one recovers the Jacobi-Trudy formula, or Giambelli formula, expressing the Schur polynomial in terms of complete symmetric functions.

After the categorification $\UcatD$ was defined in \cite{Lau1} a related but different approach to categorification appeared in the work of Rouquier~\cite{Rou2}, extending his earlier work with Chuang~\cite{CR}. Rouquier defines several 2-categories associated to a symmetrizable Kac-Moody algebra.  They are related to the 2-categories $\UcatD(\mf{g})$, however, already in the $\mf{sl}_2$ case these 2-categories appear to be different from $\UcatD$.  Both approaches have their benefits; for instance, the axiomatics in \cite{Lau1} ensures that the 2-category $\UcatD(\mf{sl}_2)$ has Grothendieck ring isomorphic to the integral idempotent form of $\mathbf{U}_q(\mf{sl}_2)$, while the less stringent axiomatics in \cite{Rou2} make it easier to check the existence of categorical actions in this sense~\cite{CKL,CKL2,CKL3} and lead to important results about 2-representations.
We do not know if counterparts of results in the present paper exist for Rouquier's 2-categories, since our results crucially depend on all the relations present in the 2-category $\UcatD$.

The explicit diagrammatic calculus developed in this article should be useful for enhancing and clarifying constructions of categorical quantum $\mathfrak{sl}_2$ actions.  Since the appearance of the more general 2-category $\UcatD(\mf{g})$, its 2-representation theory (see \cite[Section 3.4]{KL3}), is beginning to be studied.   Hill and Sussan~\cite{HS} have constructed an action of the 2-category $\UcatD$ (over field $\mathbb{F}_2$) associated to $\mathfrak{sl}_n$ on diagrammatic categories introduced in~\cite{HK} that categorify the irreducible highest weight representation of highest weight $2\omega_k$.  We expect that various categorical actions of $\U(\mf{g})$ extend to 2-representations of the 2-category $\UcatD(\mf{g})$, including the categorified representations studied in \cite{BS,BFK,FKS,LV,Rou2,Zheng}.

The formulas derived here for the 2-category $\UcatD$ might be useful for applications of categorified quantum groups to knot homology.   An action of the 2-category $\UcatD$ on categories of $\mathfrak{sl}_3$ foams arising in knot homology has been described by Mackaay~\cite{Mackaay}.  This action may be related to similar actions, constructed in \cite{Vaz,MV}, of the diagrammatic Soergel category \cite{EK,EKr}.  Categorified quantum groups appear in geometry in the work of Cautis-Kamnitzer-Licata~\cite{CKL,CKL2,CKL3} and Varagnolo-Vasserot~\cite{VV}, see also \cite{CK,Lau2}. The diagrammatic calculus developed here should be relevant to the striking new work of Webster~\cite{Web,Web2} on categorification of Reshetikhin-Turaev tangle invariants.

\bigskip
\noindent {\bf Acknowledgments:}
Arguably the most intriguing result of the paper is Theorem~\ref{thm_EaFb}, which gives an explicit graphical presentation for decomposing certain 1-morphisms in the 2-category $\UcatD$.  The discovery and proof of this theorem are due to the fourth author, M.S., and the three other authors would like to acknowledge this over M.S.'s objections.

We are grateful to Joshua Sussan for comments on an early version of this article. M.K. is grateful to the NSF for supporting him via grants DMS-0706924 and DMS-0739392.
A.L. was partially supported by the NSF grants  DMS-0739392 and DMS-0855713. M.K. and A.L. would like to thank the MSRI for support in Spring 2010. M.M. and M.S. were partially supported by the \textit{Programa Operacional Ci\^encia e Inova\c c\~ao 2010}, financed by FCT and cofinanced by the European Community fund FEDER, in part through the research project: New Geometry and Topology, PTDC/MAT/101503/2008. M.S. was also partially supported by the Ministry of Science of Serbia, project no. 144032.


%
\section{Thick calculus for the nilHecke ring}
\label{sec_nil_thick}
%

\subsection{The nilHecke ring and its diagrammatics} \label
{sec_nilHecke}

The nilHecke ring appeared in the study of the cohomology ring of
flag varieties~\cite{BGG,Dem} and is further related to the theory
of Schubert varieties, see ~\cite{KK,BilLak,Kum,Man}.

The nilHecke ring $\BNC_a$ is the unital ring of endomorphisms of the abelian
group $\Z[x_1, \dots , x_a]$ generated by the endomorphisms of
multiplication by $x_1, \dots, x_a$ and the divided difference
operators
$$ \partial_i(f(x)) = \frac{f(x)-s_i f(x)}{x_i -x_{i+1}}, \hspace{0.2in} 1\le i \le a-1,$$
where $s_i$ transposes $x_i$ and $x_{i+1}$ in the polynomial $f(x)$.
The defining relations are
\[
 \begin{array}{ll}
 x_i x_j =   x_j x_i , &  \\
   \partial_i x_j = x_j\partial_i \quad \text{if $|i-j|>1$}, &
   \partial_i\partial_j = \partial_j\partial_i \quad \text{if $|i-j|>1$}, \\
  \partial_i^2 = 0,  &
   \partial_i\partial_{i+1}\partial_i = \partial_{i+1}\partial_i\partial_{i+1},  \\
   x_i \partial_i - \partial_i x_{i+1}=1,  &   \partial_i x_i - x_{i+1} \partial_i =1.
  \end{array}
\]
In the above equations $x_i$ stands for the operator of
multiplication by $x_i$. We equip $\BNC_a$ with a grading such that
$\deg(\partial_i)=-2$ and $\deg(x_i)=2$.

To a permutation $w\in S_a$ assign $\partial_w =
\partial_{i_1}\dots \partial_{i_r},$ where $s_{i_1}\dots s_{i_r}$ is
a reduced expression of $w$. The element $\partial_w$ does not
depend on the choice of reduced expression.  Let $w_0$ be the
maximal length permutation in $S_a$.  We will denote
$\partial_{w_0}$ by $D_a$ throughout the paper.

Recall that
\begin{equation}
\partial_{w}\partial_{w'} \;\; = \;\;
\left\{
\begin{array}{ccl}
  \partial_{ww'} & \quad & \text{if $\ell(ww')=\ell(w)+\ell(w')$}, \\
  0 & \quad & \text{otherwise.}
\end{array}
\right.
\end{equation}
The following relations hold
\begin{equation} \label{eq_partial_Da_zero}
\partial_i D_a = D_a \partial_i = 0 \qquad \text{for all $1 \leq i\leq a-1$,}
\end{equation}
\begin{equation} \label{eq_DafDa}
 D_a f D_a  = D_a(f) D_a.
\end{equation}
$D_a(f) \in \Z[x_1,\dots, x_a]$ is the polynomial obtained by
applying the product $D_a$ of divided differences to $f$. If $f$ is
a monomial symmetric in some two variables, then $D_a(f)=0$ and
$D_afD_a=0$.

Let
\begin{eqnarray}
  \delta_a & =& x_1^{a-1}x_2^{a-2} \dots x_{a-1}, \\
  e_a & = & \delta_a D_a.
\end{eqnarray}
We have $D_a(\delta_a)=1$ and
\begin{equation} \label{eq_Daea}
  D_a e_a = D_a \delta_a D_a = D_a (\delta_a)D_a = D_a,
\end{equation}
implying that $e_a$ is an idempotent, $e_a^2=e_a$.

The center $Z(\BNC_a)$ of the nilHecke ring is the subring of
symmetric polynomials
\begin{equation}\label{eq_centerNH}
  Z(\BNC_a) = \Z[x_1,\dots, x_a]^{S_a} \subset \Z[x_1,\dots,x_a] \subset \BNC_a,
\end{equation}
and $\BNC_a$ is
isomorphic to the ring of $a!\times a!$ matrices with coefficients
in $Z(\BNC_a)$, see~\cite{Lau1}. In the diagram below
\begin{equation} \label{eq_nilHecke_inc}
  \xy
 (-45,0)*+{\BNC_a}="1";
 (0,0)*+{Z(\BNC_a) \cong Z(\BNC_a) e_a = e_a\BNC_a e_a}="2";
 (45,0)*+{\BNC_a}="4";
 {\ar@{_{(}->}_-{} "2";"1"};  {\ar@{^{(}->}^-{} "2";"4"};
 \endxy
\end{equation}
the leftmost arrow is the inclusion \eqref{eq_centerNH}, while the rightmost arrow is the non-unital inclusion of the
subring associated with the idempotent $e_a$ into $\BNC_a$.  The
isomorphism in the middle left takes a symmetric polynomial $y$ to
$ye_a$. The idempotent $e_a$ is minimal. We will frequently use the canonical isomorphism
\begin{equation} \label{eq_iso_center}
  \Z[x_1,\dots,x_a]^{S_a} \cong Z(\BNC_a) e_a
\end{equation}
obtained by combining the equality in \eqref{eq_centerNH} with the middle left isomorphism in \eqref{eq_nilHecke_inc}.

We find it convenient to use a graphical calculus to represent
elements in $\BNC_a$. We write
\begin{eqnarray}
   \xy  (0,0)*{\includegraphics[scale=0.4]{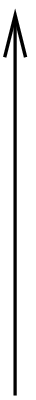}};  \endxy
   \;\;  \dots \;\;
   \xy  (0,0)*{\includegraphics[scale=0.4]{figs/thinup.eps}};  \endxy
   \;\;  \dots \;\;
   \xy  (0,0)*{\includegraphics[scale=0.4]{figs/thinup.eps}};  \endxy \quad := \quad
 1 \in \BNC_a
\end{eqnarray}
with a total of $a$ strands. The polynomial generators can be
written as
\begin{eqnarray}
  \xy  (0,0)*{\includegraphics[scale=0.4]{figs/thinup.eps}};  \endxy
   \;\;  \dots \;\;
   \xy  (0,0)*{\includegraphics[scale=0.4]{figs/thinup.eps}}; (0,0)*{\bullet}; \endxy
   \;\;  \dots \;\;
   \xy  (0,0)*{\includegraphics[scale=0.4]{figs/thinup.eps}};  \endxy\quad := \quad
   x_r
\end{eqnarray}
with the dot positioned on the $r$-th strand counting from the left,
and
\begin{eqnarray}
  \xy  (0,0)*{\includegraphics[scale=0.4]{figs/thinup.eps}};  \endxy
    \;\;  \dots \;\;
   \xy  (0,0)*{\includegraphics[scale=0.4]{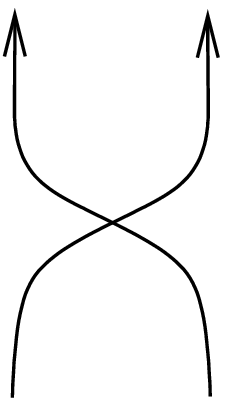}};  \endxy
   \;\;  \dots \;\;
   \xy  (0,0)*{\includegraphics[scale=0.4]{figs/thinup.eps}};  \endxy\quad := \quad
  \partial_r
\end{eqnarray}
with the crossing interchanging the $r$th and $(r+1)$st strands.

In the diagrammatic notation multiplication is given by stacking
diagrams on top of each other from bottom to top.  Relations in the nilHecke ring acquire a graphical interpretation. For example, the equalities $\partial_rx_r-x_{r+1}\partial_r = 1 =
x_r \partial_r - \partial_r x_{r+1}$ become diagrammatic
identities:
\begin{eqnarray}        \label{new_eq_iislide}
  \xy  (0,0)*{\includegraphics[scale=0.4]{figs/thincross.eps}};
  (-3.5,-4)*{\bullet};
  \endxy
    \quad - \quad
   \xy  (0,0)*{\includegraphics[scale=0.4]{figs/thincross.eps}};
   (3.9,3)*{\bullet};
  \endxy
  & = &
   \xy  (-3,0)*{\includegraphics[scale=0.4]{figs/thinup.eps}};
        (3,0)*{\includegraphics[scale=0.4]{figs/thinup.eps}};
  \endxy
  \quad = \quad
  \xy  (0,0)*{\includegraphics[scale=0.4]{figs/thincross.eps}};
  (-4,3)*{\bullet};
  \endxy
    \quad - \quad
   \xy  (0,0)*{\includegraphics[scale=0.4]{figs/thincross.eps}};
   (3.5,-4)*{\bullet};
  \endxy
\end{eqnarray}
and the relation $\partial_r \partial_r= 0$ becomes
\begin{eqnarray}
\xy (0,0)*{\includegraphics[scale=0.4]{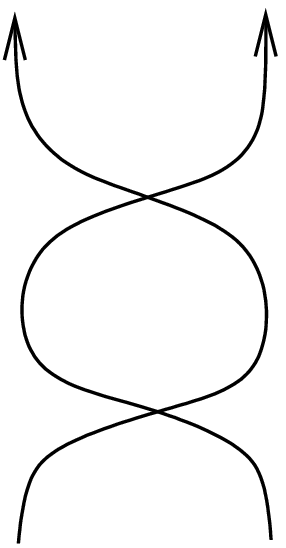}}; \endxy
  &=&
0
\end{eqnarray}
The relation $\partial_r\partial_{r+1}\partial_r =
\partial_{r+1}\partial_r \partial_{r+1}$ is depicted as
\begin{eqnarray}      \label{new_eq_r3_easy}
\xy (0,0)*{\includegraphics[scale=0.4]{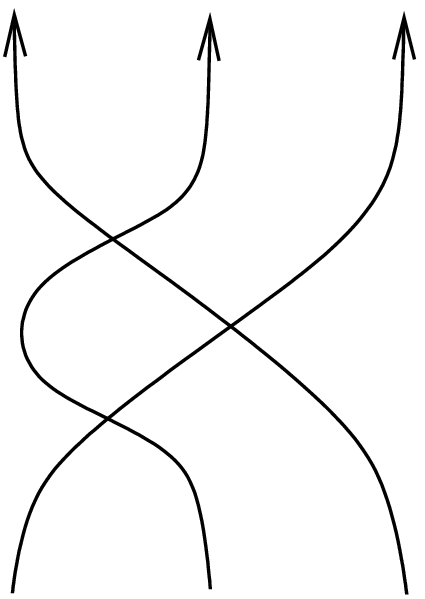}}; \endxy
  &=&
 \xy (0,0)*{\reflectbox{\includegraphics[scale=0.4]{figs/thinR3L.eps}}}; \endxy
\end{eqnarray}

The remaining relations in the nilHecke ring can be encoded by
the requirement that the diagrams are invariant under braid-like
isotopies.
\begin{equation}
 \xy
 (0,0)*{\includegraphics[scale=0.5]{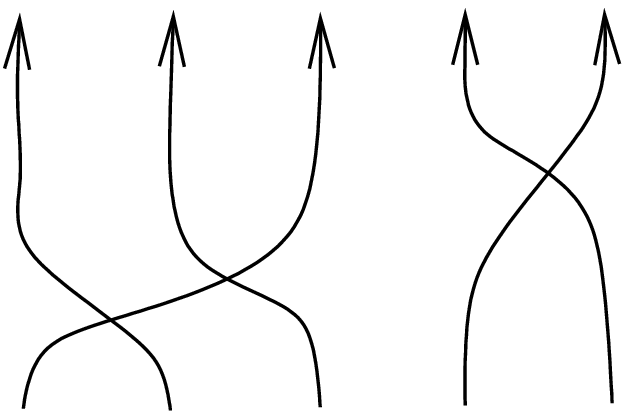}};
 (-15 ,4)*{\bullet};  (-15 ,0)*{\bullet};
 (7.8 ,-7)*{\bullet};  (14.3 ,5)*{\bullet};
  \endxy \qquad = \qquad
   \xy
 (0,0)*{\includegraphics[scale=0.5]{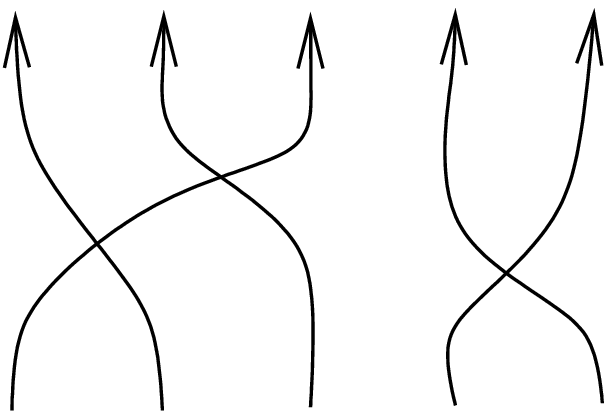}};
 (-14 ,4)*{\bullet};  (-12 ,0)*{\bullet};
 (7.2 ,-7)*{\bullet};  (14.2 ,5)*{\bullet};
  \endxy
\end{equation}

Inductively applying \eqref{new_eq_iislide} gives the relation
\begin{eqnarray}        \label{new_eq_IND_iislide}
  \xy  (0,0)*{\includegraphics[scale=0.4]{figs/thincross.eps}};
  (-3.5,-4)*{\bullet}+(-3.5,0)*{\scs m};
  \endxy
    \quad - \quad
   \xy  (0,0)*{\includegraphics[scale=0.4]{figs/thincross.eps}};
   (4,3)*{\bullet}+(3.5,0)*{\scs m};
  \endxy
  \quad = \quad
  \xy  (0,0)*{\includegraphics[scale=0.4]{figs/thincross.eps}};
  (-4,3)*{\bullet}+(-3.5,0)*{\scs m};
  \endxy
    \quad - \quad
   \xy  (0,0)*{\includegraphics[scale=0.4]{figs/thincross.eps}};
   (3.5,-4)*{\bullet}+(3.5,0)*{\scs m};
  \endxy
 \quad=\quad  \sum_{\ell_1+\ell_2=m-1}\;
   \xy  (-3,0)*{\includegraphics[scale=0.4]{figs/thinup.eps}};
        (3,0)*{\includegraphics[scale=0.4]{figs/thinup.eps}};
         (3,-4)*{\bullet}+(3,0)*{\scs \ell_2};
         (-3,-4)*{\bullet}+(-3,0)*{\scs \ell_1};
  \endxy
\end{eqnarray}

\subsection{Boxes, thick lines, and splitters}
\label{sec_box-line-split}

Here we develop basic diagrammatics for computations in the nilHecke
ring. We denote $D_a$ by
\begin{equation}
  \xy
 (0,0)*{\includegraphics[scale=0.5]{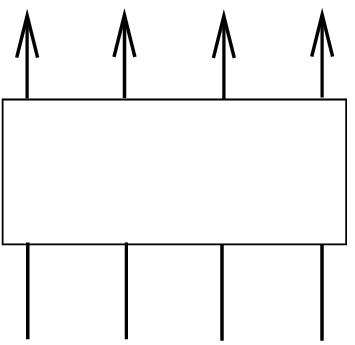}};
 (0,0)*{D_a};(0,-11)*{\underbrace{\hspace{0.7in}}};
 (0,-14)*{a};
  \endxy
 \quad := \;\;
   \xy
 (0,0)*{\includegraphics[scale=0.5]{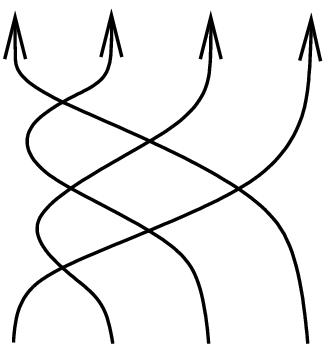}};
  \endxy  \nn
\end{equation}

To simplify diagrams, write
\[ 
  \xy
 (0,0)*{\includegraphics[scale=0.5]{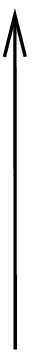}};
 (-3,-3)*{a};
  \endxy
  \quad : = \quad
  \xy
 (0,0)*{\includegraphics[scale=0.5]{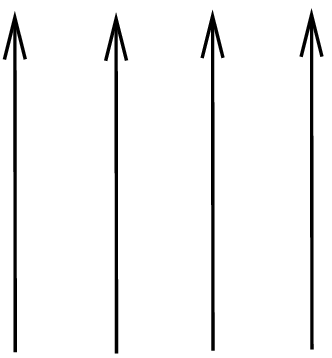}};
 (0,-11)*{\underbrace{\hspace{0.7in}}};  (0,-14)*{a};
  \endxy
\]
where we will omit the label $a$ if it appears in a coupon as below.
Next, let
\[ 
 \delta_a  \;\; = \;\;
     \xy
 (0,0)*{\includegraphics[scale=0.5]{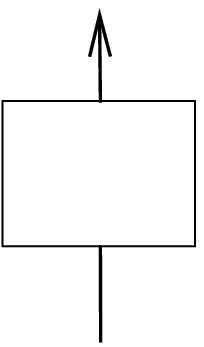}};
 (0,0)*{\delta_a};
  \endxy
 \;\;:= \;\;
 \xy
 (-2.6,0)*{\includegraphics[scale=0.5]{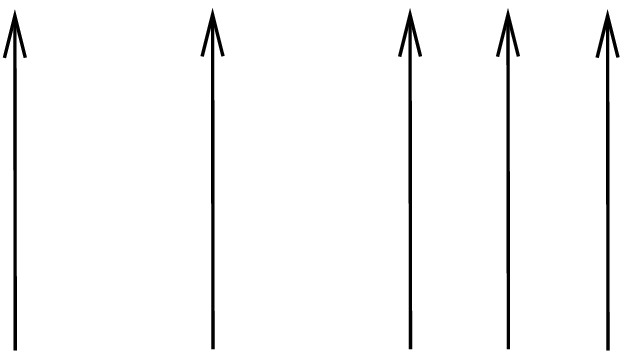}};
 (-17.7,0)*{\bullet}+(-4,1)*{\scs a-1};
 (-7.7,0)*{\bullet}+(-4,1)*{\scs a-2};
 (-3,-3)*{\cdots};
 (2.3,0)*{\bullet}+(-2,1)*{\scs 2};
 (7.3,0)*{\bullet};
  \endxy
\]
A box labelled $e_a$ denotes the idempotent $ e_a = \delta_a D_a$:
\begin{equation} \label{eq_def_ea}
   \xy
 (0,0)*{\includegraphics[scale=0.5]{figs/box-up.eps}};
 (0,0)*{e_a};
  \endxy
  \quad = \quad
  \xy
 (0,0)*{\includegraphics[scale=0.5]{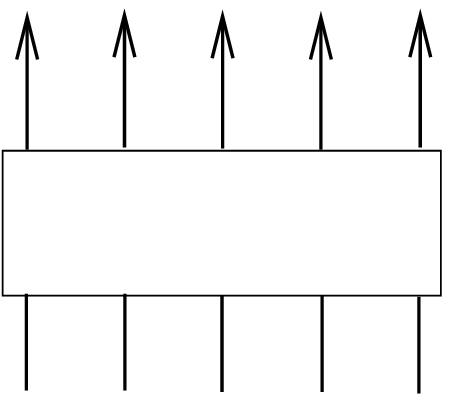}};
 (0,-1.5)*{e_{a}};
  \endxy
\quad := \quad
  \xy
 (0,0)*{\includegraphics[scale=0.5]{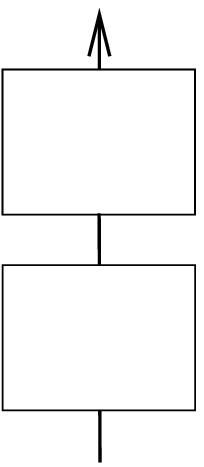}};
 (0,-5.5)*{D_a};(0,4.5)*{\delta_a};
  \endxy
  \quad = \quad
 \xy
 (0,0)*{\includegraphics[scale=0.5]{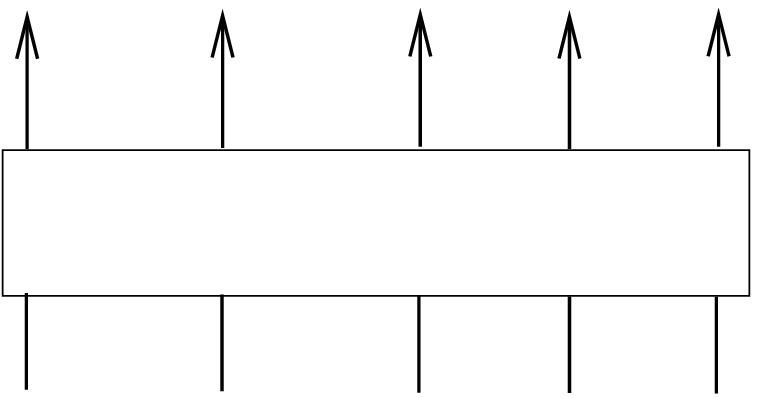}};
 (-17.7,5)*{\bullet}+(-4,1)*{\scs a-1};
 (-7.7,5)*{\bullet}+(-4,1)*{\scs a-2};
 (-3,4)*{\cdots};(-3,-7)*{\cdots};
 (2.3,5)*{\bullet}+(-2,1)*{\scs 2};
 (9.9,5)*{\bullet};(0,-1.5)*{D_a};
  \endxy
\end{equation}
More generally, we can use a box labelled $y$
\begin{equation} \label{eq_ybox}
  \xy
 (0,0)*{\includegraphics[scale=0.5]{figs/box-up.eps}};
 (0,0)*{y};(2.5,-8)*{a};
  \endxy
\end{equation}
to denote an element $y \in \BNC_a$.

Define a crossing by
\[ 
  \xy
 (0,0)*{\includegraphics[scale=0.5]{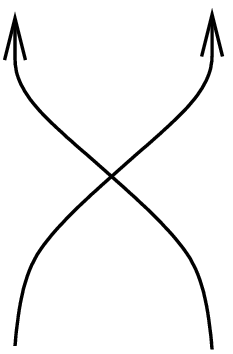}};
 (-7,-6)*{a};(7,-6)*{b};
  \endxy
  \quad := \quad
  \xy
 (0,0)*{\includegraphics[scale=0.5]{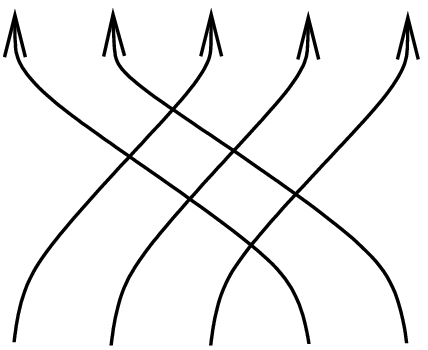}};
 (-5.5,-11)*{\underbrace{\hspace{0.45in}}};  (-5.5,-14)*{a};
 (7.5,-11)*{\underbrace{\hspace{0.25in}}};  (7.5,-14)*{b};
  \endxy
\]

From the definition we get
\begin{equation} \label{eq_various_partial}
     \xy
 (0,0)*{\includegraphics[scale=0.5]{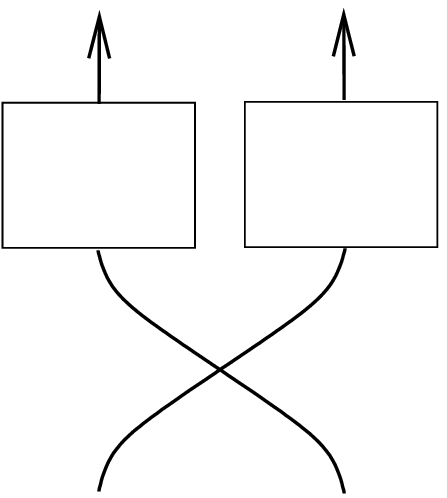}};
(-6,4)*{D_a};(6,4)*{D_b};
  \endxy
 \quad
 =
 \quad
      \xy
 (0,0)*{\includegraphics[scale=0.5]{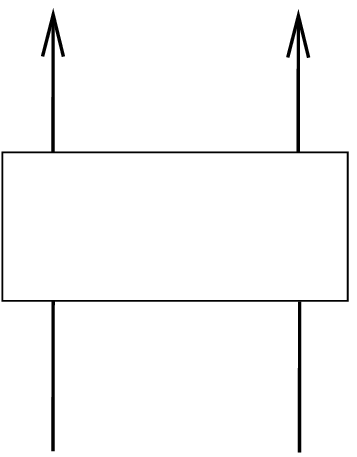}};
 (0,0)*{D_{a+b}};
  \endxy
\end{equation}

Equation \eqref{eq_DafDa} becomes
\begin{equation}
    \xy
 (0,0)*{\includegraphics[scale=0.5]{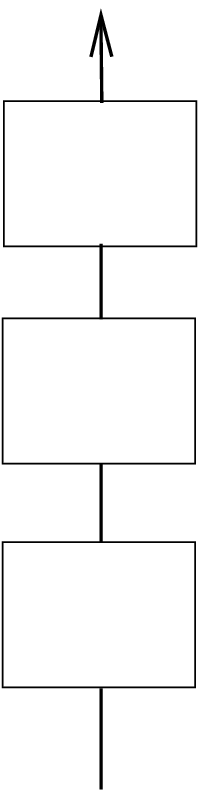}};
 (0,11)*{D_a};(0,0)*{f};(0,-11)*{D_a};
  \endxy
  \quad = \quad
    \xy
 (0,0)*{\includegraphics[scale=0.5]{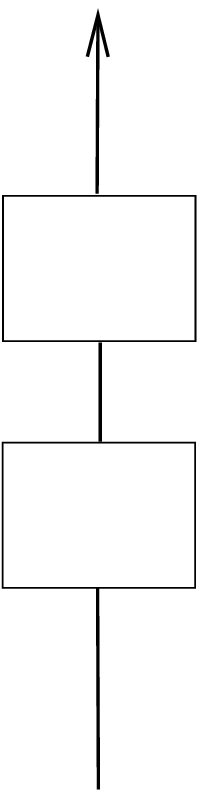}};
 (0,6)*{D_a(f)}; (0,-6)*{D_a};
  \endxy
\end{equation}
As a special case, if the degree of $f$ is less than minus the degree of $D_a$, the element $D_a f D_a$ is equal to zero
\begin{equation}
    \xy
 (0,0)*{\includegraphics[scale=0.5]{figs/three-box.eps}};
 (0,11)*{D_a};(0,0)*{f};(0,-11)*{D_a};
  \endxy
  \quad = \quad 0 \quad \text{if $\deg f < a(a-1)$ (recall that $\deg(x_i)=2$).}
\end{equation}

\begin{prop} The following identities hold
\begin{itemize}
\item
\begin{equation} \label{eq_partial-en}
  \xy
 (0,0)*{\includegraphics[scale=0.5]{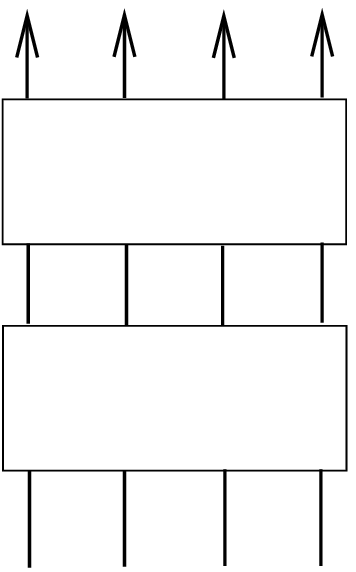}};
 (0,-5.5)*{e_{a}};(0,5.5)*{D_a};
  \endxy
 \quad = \quad
  \xy
 (0,0)*{\includegraphics[scale=0.5]{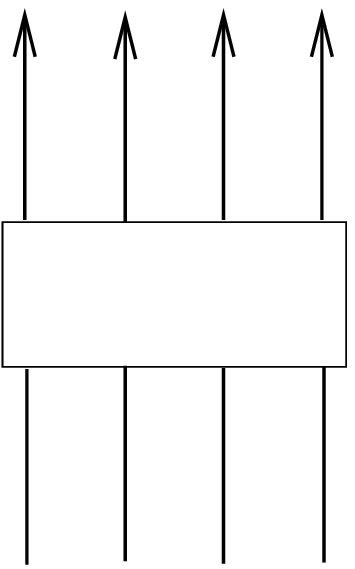}};
 (0,-0.5)*{D_a};
  \endxy
\end{equation}
\item
\begin{equation} \label{eq-boxes-1}
 \xy
 (1,-2)*{\includegraphics[scale=0.5]{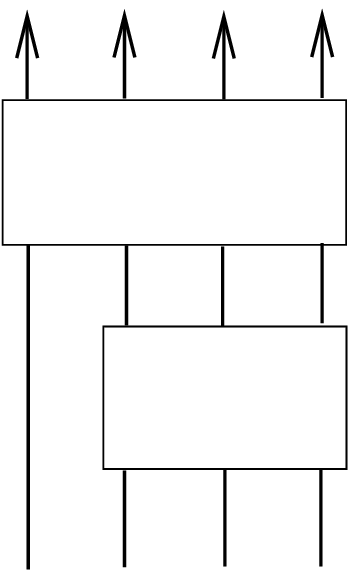}};
 (1,4)*{e_a};(4,-8)*{e_{a-1}};
\endxy
 \quad =\quad
  \xy
 (0,-2)*{\includegraphics[scale=0.5]{figs/c4-2.eps}};
 (0,-3)*{e_a};
  \endxy
\qquad \quad
 \xy
 (1,-2)*{\includegraphics[scale=0.5]{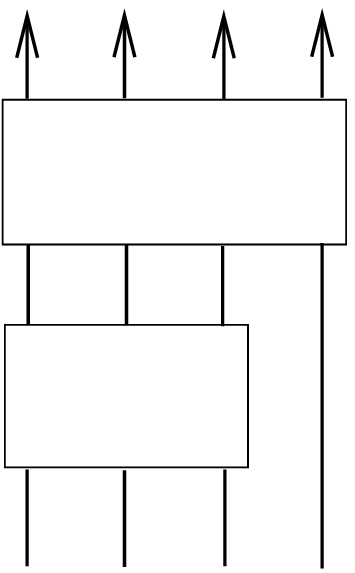}};
 (1,4)*{e_a};(-1,-8)*{e_{a-1}};
\endxy
 \quad =\quad
  \xy
 (0,-2)*{\includegraphics[scale=0.5]{figs/c4-2.eps}};
 (0,-3)*{e_a};
  \endxy
\end{equation}
More generally we have
\begin{equation} \label{eq_projector_absorb}
 \xy
 (1,-2)*{\includegraphics[scale=0.5]{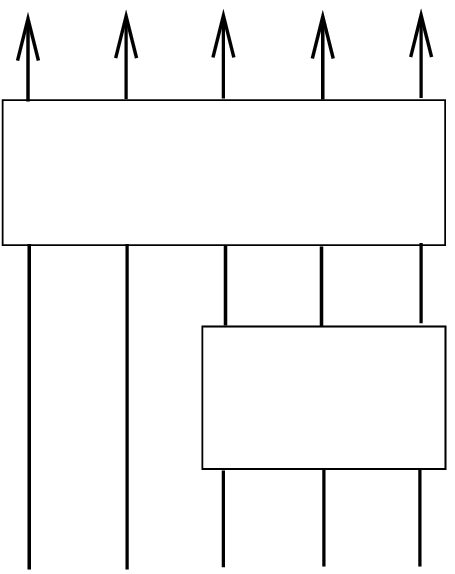}};
 (1,4)*{e_a};(6,-8)*{e_{a-k}};
\endxy
 \quad =\quad
  \xy
 (0,-2)*{\includegraphics[scale=0.5]{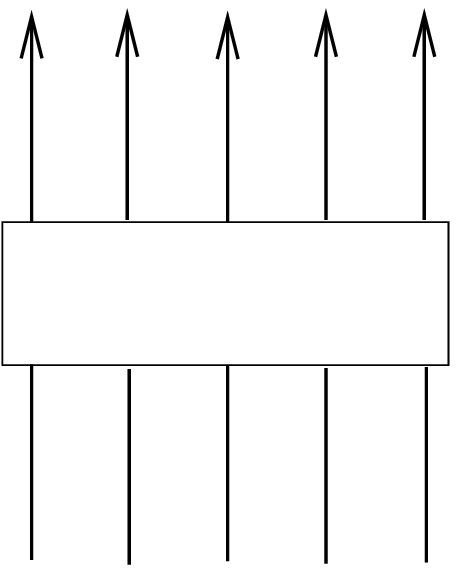}};
 (0,-3)*{e_a};
  \endxy
\qquad \quad
 \xy
 (1,-2)*{\reflectbox{\includegraphics[scale=0.5]{figs/c5ii.eps}}};
 (1,4)*{e_a};(-4,-8)*{e_{a-k}};
\endxy
 \quad =\quad
  \xy
 (0,-2)*{\includegraphics[scale=0.5]{figs/c4-2ii.eps}};
 (0,-3)*{e_a};
  \endxy
\end{equation}
for any $ 0 \leq k \leq a$.
\item
\begin{equation} \label{eq-boxes-2}
 \xy
 (-1,-2)*{\includegraphics[scale=0.5]{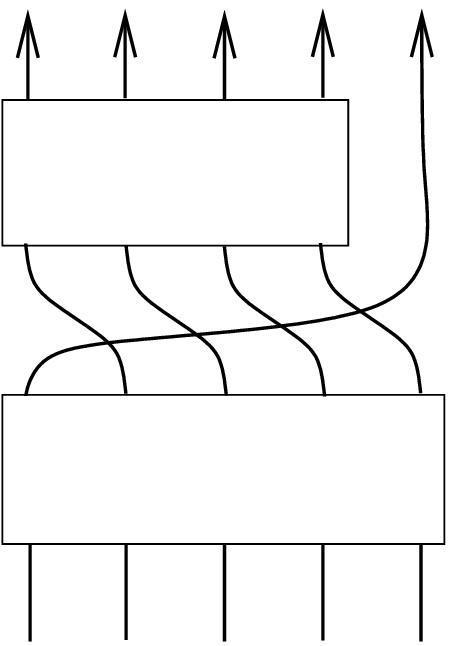}};
 (-3,5)*{e_{a-1}};(0,-10)*{e_{a}};
\endxy
 \;\; =\;\;
  \xy
 (-3,-2)*{\includegraphics[scale=0.5]{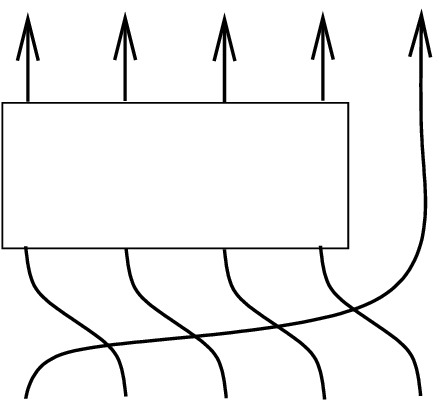}};
 (-4,-1)*{e_{a-1}};
  \endxy
\qquad \quad
 \xy
 (-1,-2)*{\reflectbox{\includegraphics[scale=0.5]{figs/c7-1.eps}}};
 (2,5)*{e_{a-1}};(0,-10)*{e_{a}};
\endxy
 \;\; =\;\;
  \xy
 (-3,-2)*{\reflectbox{\includegraphics[scale=0.5]{figs/c7-2.eps}}};
 (-1,-1)*{e_{a-1}};
  \endxy
\end{equation}

\item
\begin{eqnarray}
 \xy
 (1,-2)*{\includegraphics[scale=0.5]{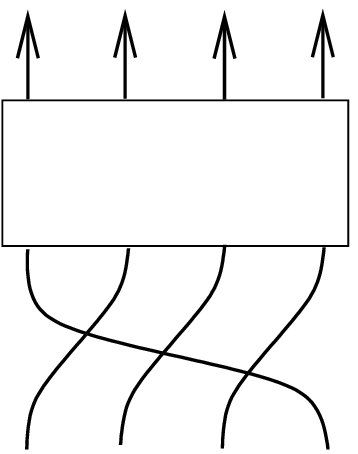}};
 (1,0)*{e_a};(-6,-6)*{\bullet}+(-2.2,0)*{\scs b};
\endxy
 & =& \left\{
\begin{array}{ccl}
  0 & &\text{if $b<a-1$,} \\ \\
    \xy
 (0,-2)*{\includegraphics[scale=0.5]{figs/c1-1.eps}};
 (0,-2)*{e_a};
  \endxy & &\text{if $b=a-1$,}
\end{array}
 \right. \label{eq-boxes_undercross1}
\\
 \xy
 (1,-2)*{\reflectbox{ \includegraphics[scale=0.5,]{figs/c8.eps}}};
 (1,0)*{e_a};(7.7,-6)*{\bullet}+(2.2,0)*{\scs b};
\endxy
 & =&  \left\{
\begin{array}{ccl}
  0 & &\text{if $b<a-1$,} \\ \\
   (-1)^{a-1} \;\;\xy
 (0,-2)*{\includegraphics[scale=0.5]{figs/c1-1.eps}};
 (0,-2)*{e_a};
  \endxy & &\text{if $b=a-1$.}
\end{array}
 \right. \label{eq-boxes_undercross2}
\end{eqnarray}
\end{itemize}
\end{prop}

\begin{proof}
These relations originally appeared in \cite{KL2}.  Relation
\eqref{eq_partial-en} is \eqref{eq_Daea}.
\end{proof}

We also observe that
\begin{equation} \label{eq_eaebcross}
  \xy
 (0,0)*{\includegraphics[scale=0.5]{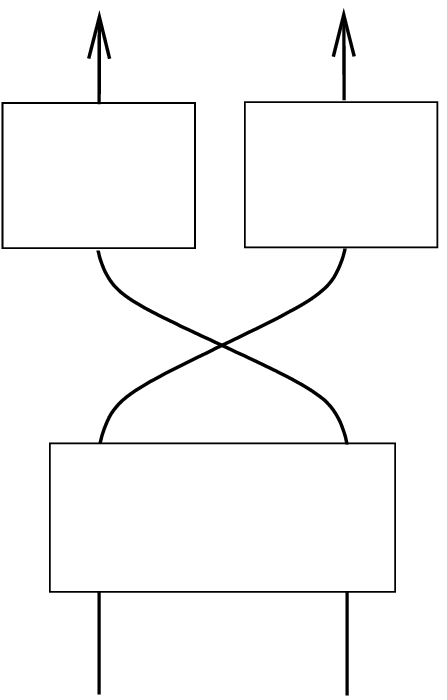}};
 (0,-8.5)*{e_{a+b}};(-6,8.5)*{e_{a}};(6,8.5)*{e_{b}};
  \endxy
\quad \refequal{\eqref{eq_def_ea}} \quad
  \xy
 (0,0)*{\includegraphics[scale=0.5]{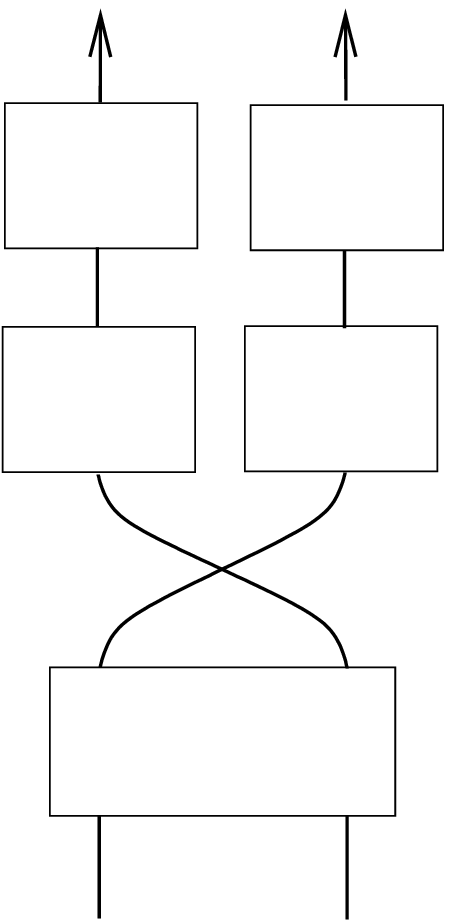}};
 (0,-15)*{e_{a+b}};(-6,14.5)*{\delta_a};(6,14.5)*{\delta_{b}};
 (-6,3)*{D_a};(6,3)*{D_b};
  \endxy
 \quad \refequal{\eqref{eq_various_partial}} \quad
  \xy
 (0,0)*{\includegraphics[scale=0.5]{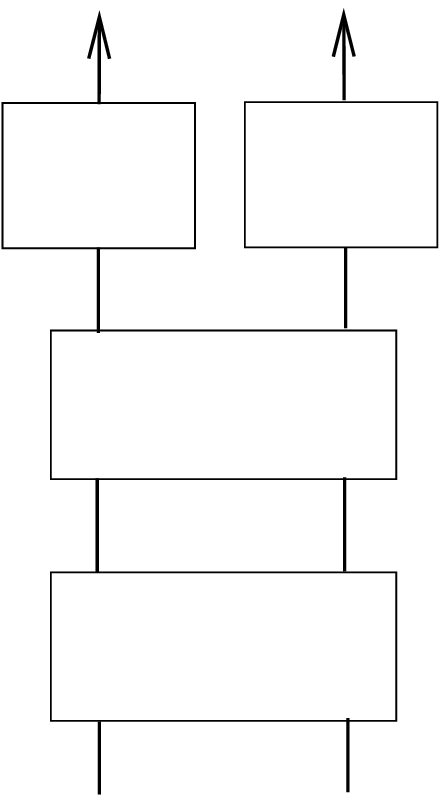}};
 (0,-13)*{e_{a+b}};(-6,11)*{\delta_a};(6,11)*{\delta_{b}};
 (0,-1)*{D_{a+b}};
  \endxy
 \end{equation}
\begin{equation} \qquad \nn
  \quad  \refequal{\eqref{eq_partial-en}} \quad
  \xy
 (0,0)*{\includegraphics[scale=0.5]{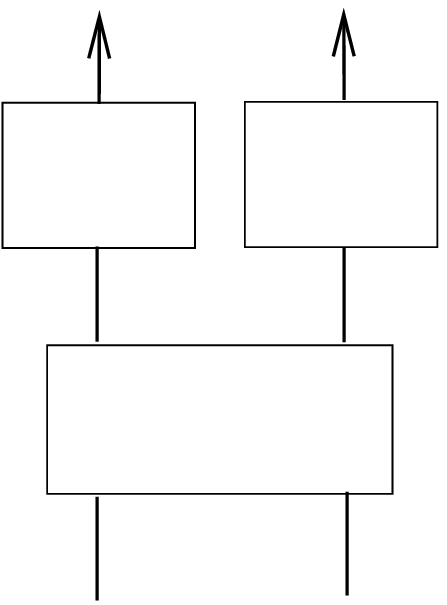}};
 (0,-6.5)*{D_{a+b}};(-6,6)*{\delta_a};(6,6)*{\delta_{b}};
  \endxy
   \quad \refequal{\eqref{eq_various_partial}}
 \quad
     \xy
 (0,0)*{\includegraphics[scale=0.5]{figs/def-tsplitu.eps}};
 (-7,-10)*{b};(7,-10)*{a};(-6,4)*{e_a};(6,4)*{e_b};
  \endxy
\end{equation}

\begin{lem}
  \begin{equation} \label{eq_lem-twobox}
                 \xy
 (0,0)*{\reflectbox{\includegraphics[scale=0.45]{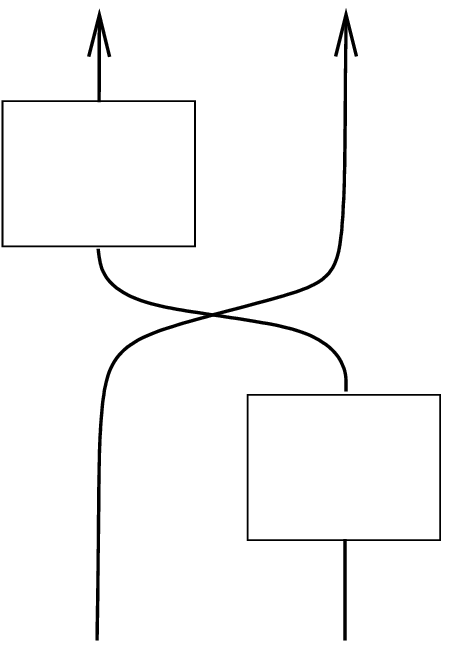}}};
(7,-7.5)*{b};(-6,-7)*{e_a};(6,7)*{e_a};(-7.5,7.5)*{b};
  \endxy
\quad = \quad
         \xy
 (2,0)*{\reflectbox{\includegraphics[scale=0.45]{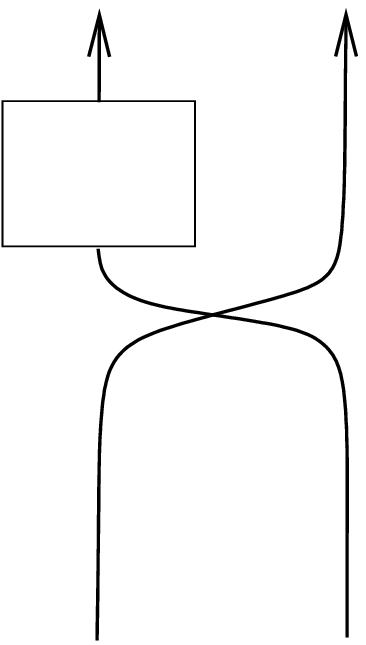}}};
(7,-7.5)*{b};(-7,-7.5)*{a};(6,7)*{e_a};(-7.5,7.5)*{b};
  \endxy
\qquad \qquad
 \xy
 (0,0)*{\includegraphics[scale=0.45]{figs/two-box-uldr.eps}};
(-7,-7.5)*{a};(6,-7)*{e_b};(-6,7)*{e_b};(7.5,7.5)*{a};
  \endxy
\quad = \quad
         \xy
 (-2,0)*{\includegraphics[scale=0.45]{figs/box-ul.eps}};
(-7,-7.5)*{a};(7,-7.5)*{b};(-6,7)*{e_b};(7.5,7.5)*{a};
  \endxy
  \end{equation}
\end{lem}

\begin{proof}
We prove the first identity; the second one is proven similarly.
The proof is by induction on $b$.  For $b=1$ we have
\begin{equation}
 \xy
 (0,0)*{\reflectbox{\includegraphics[scale=0.45]{figs/two-box-uldr.eps}}};
(7,-7.5)*{1};(-6,-7)*{e_a};(6,7)*{e_a};(-7.5,7.5)*{1};
  \endxy
   \quad = \quad
 \xy
 (4.8,0)*{\includegraphics[scale=0.45]{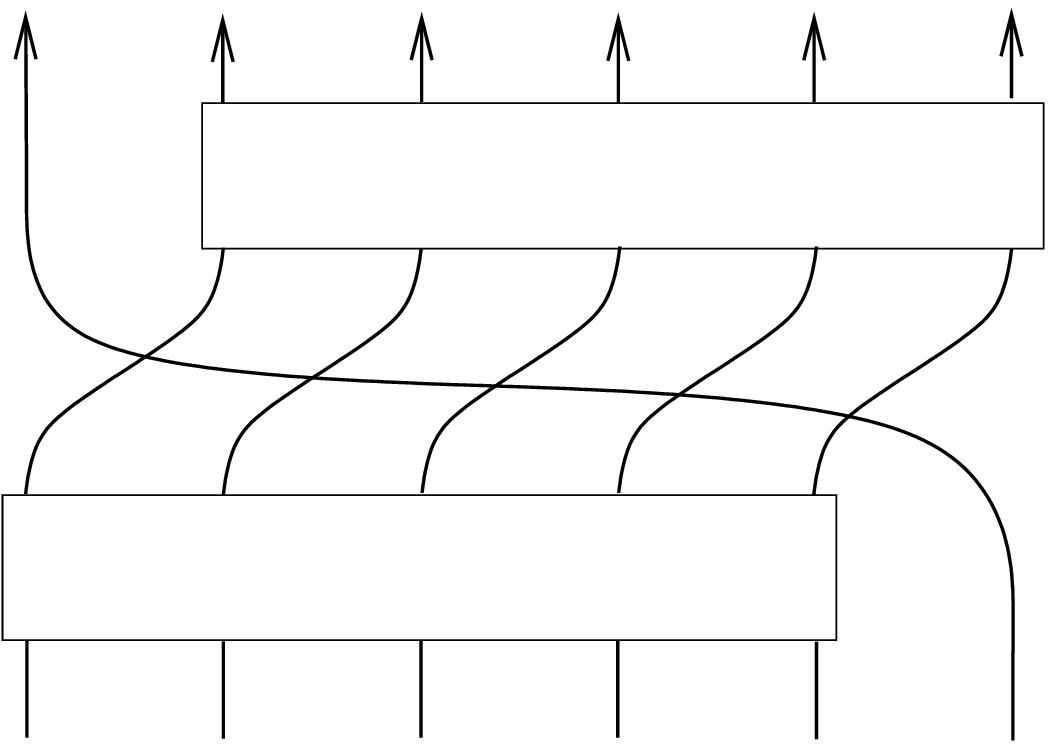}};
(30,-12)*{1};(1,-9)*{D_a};(9,9)*{e_a};(-20,10)*{1};
 (-17.3,-3)*{\bullet}+(-4,1)*{\scs a-1};
 (-8.3,-3)*{\bullet}+(-4,1)*{\scs a-2};
 (-4,-4)*{\cdots};
 (.7,-3)*{\bullet}+(-2.5,1)*{\scs 2};
 (9.8,-3)*{\bullet};
  \endxy
\end{equation}
Starting from the right, we slide the first dot across the line
using the nilHecke relation \eqref{new_eq_iislide} to get
\begin{equation} \label{eq_lem_slidedots}
= \quad
 \xy
 (4.8,0)*{\includegraphics[scale=0.45]{figs/lem1.eps}};
(30,-12)*{1};(1,-9)*{D_a};(9,9)*{e_a};(-20,10)*{1};
 (-17.3,-3)*{\bullet}+(-4,1)*{\scs a-1};
 (-8.3,-3)*{\bullet}+(-4,1)*{\scs a-2};
 (-4,-4)*{\cdots};
 (.7,-3)*{\bullet}+(-2.5,1)*{\scs 2};
 (16.5,2)*{\bullet};
  \endxy
\;\; + \;\;
 \xy
 (4.8,0)*{\includegraphics[scale=0.45]{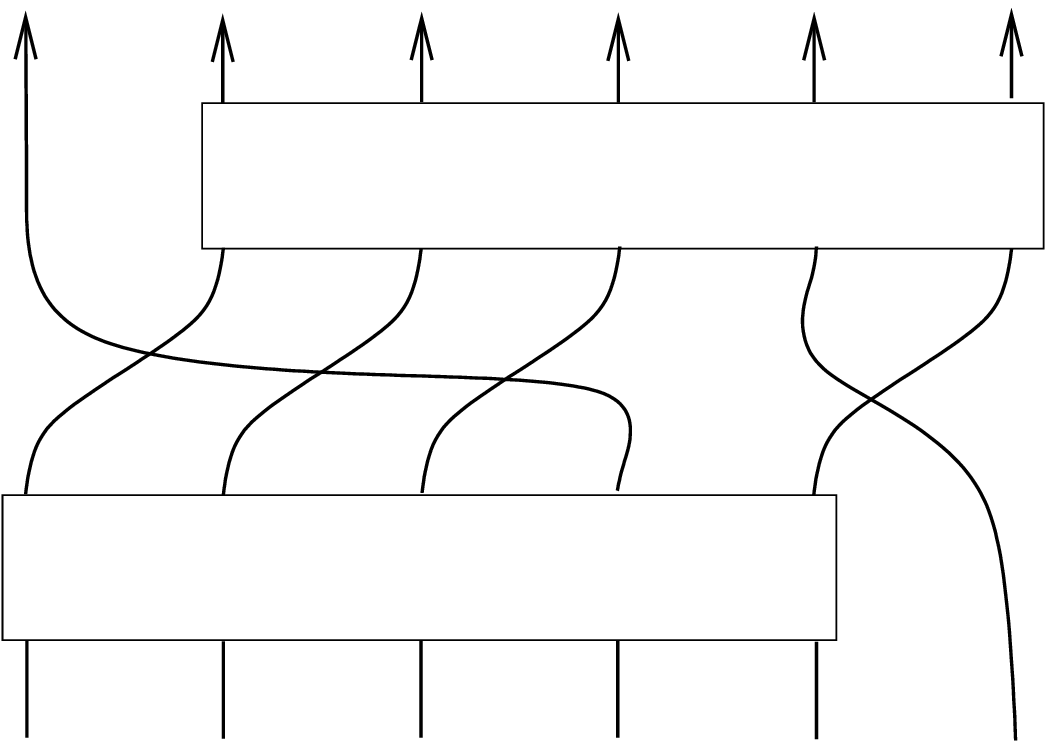}};
(30,-12)*{1};(1,-9)*{D_a};(9,9)*{e_a};(-20,10)*{1};
 (-17.3,-3)*{\bullet}+(-4,1)*{\scs a-1};
 (-8.3,-3)*{\bullet}+(-4,1)*{\scs a-2};
 (-4,-4)*{\cdots};
 (.7,-3)*{\bullet}+(-2.5,1)*{\scs 2};
  \endxy
\end{equation}
But the second term is zero since
\begin{equation}
  \xy
 (0,-2)*{\includegraphics[scale=0.5]{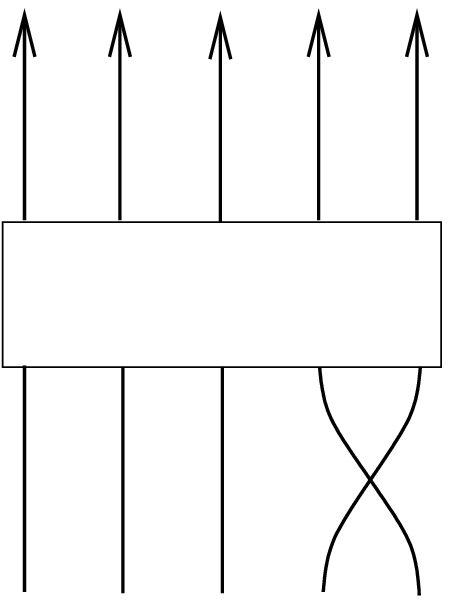}};
 (0,-2)*{e_a};
  \endxy
\;\; = \;\;
      \xy
 (0,-2)*{\includegraphics[scale=0.5]{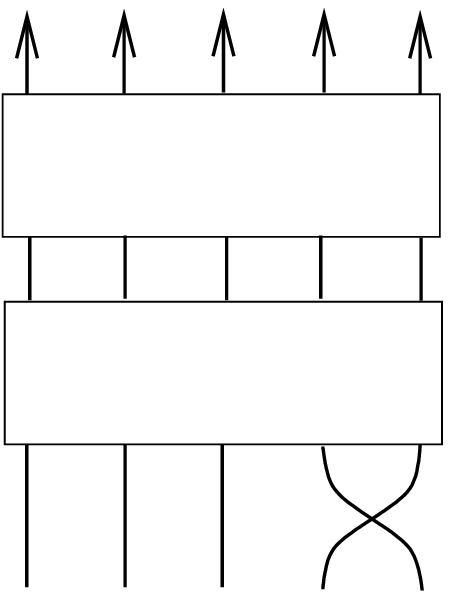}};
 (0,5)*{\delta_a};(0,-6)*{D_a};
  \endxy
 \;\; \refequal{\eqref{eq_partial_Da_zero}} \;\; 0.
\end{equation}
We continue in this way, going from right to left, sliding all dots
across the line using the inductive dot slide formula
\eqref{new_eq_IND_iislide}.  All terms will vanish except for the term
where all dots make it across the line and no crossing is resolved.
To see this suppose we slide all dots up to the strand with $k$
dots.  Then we have
\begin{equation}
 \xy
 (4.8,0)*{\includegraphics[scale=0.45]{figs/lem1.eps}};
(30,-12)*{1};(1,-9)*{D_a};(9,9)*{e_a};(-20,10)*{1};
 (-17.3,-3)*{\bullet}+(-4,1)*{\scs a-1}; (-13.5,-4)*{\cdots};
 (-8.3,-3)*{\bullet}+(-2.5,1)*{\scs k};
 (-4,-4)*{\cdots};
 (7.8,2)*{\bullet}+(-4,1)*{\scs k-1}; (12.5,3)*{\cdots};
 (16.5,2)*{\bullet};
  \endxy
\end{equation}
which by the inductive dot slide formula is
\begin{equation} \label{eq_lem_slidedots2}
\;\; \refequal{\eqref{new_eq_IND_iislide}}\;\;
 \xy
 (4.8,0)*{\includegraphics[scale=0.45]{figs/lem1.eps}};
(30,-12)*{1};(1,-9)*{D_a};(9,9)*{e_a};(-20,10)*{1};
 (-17.3,-3)*{\bullet}+(-4,1)*{\scs a-1}; (-13.5,-4)*{\cdots};
 (-1.3,2)*{\bullet}+(-2.5,1)*{\scs k};
 (-4,-4)*{\cdots};
 (7.8,2)*{\bullet}+(-4,1)*{\scs k-1}; (12.5,3)*{\cdots};
 (16.5,2)*{\bullet};
  \endxy
 \;\; + \;\; \sum_{b=0}^{k-1}\;\;
 \xy
 (4.8,0)*{\includegraphics[scale=0.45]{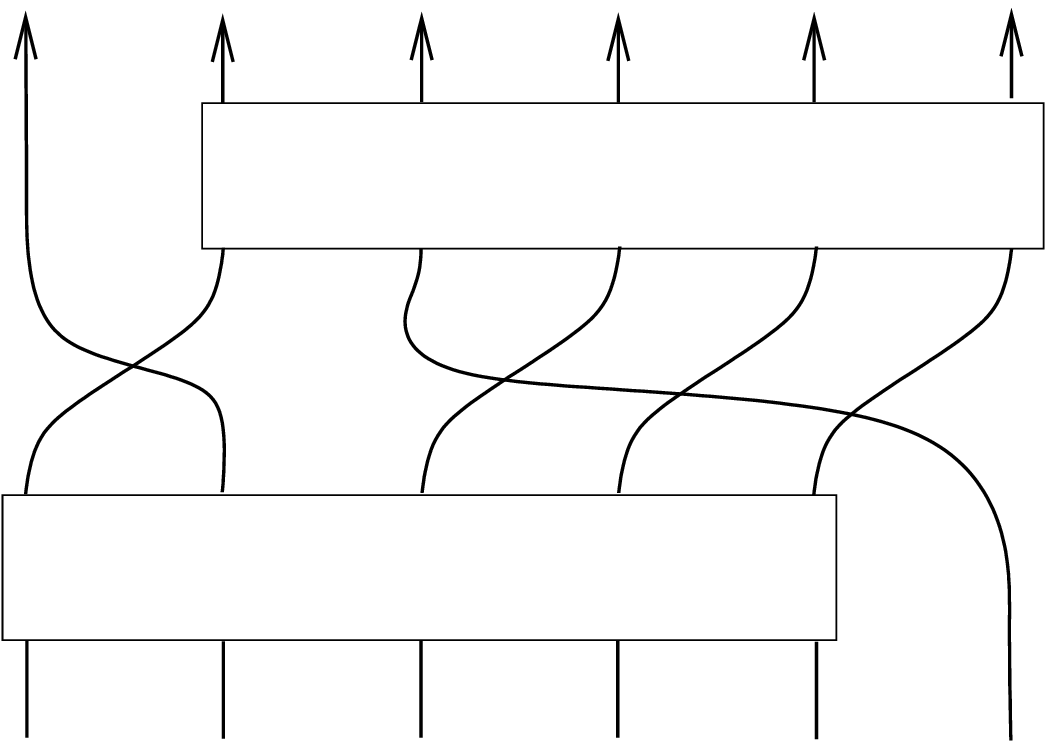}};
(30,-12)*{1};(1,-9)*{D_a};(9,9)*{e_a};(-20,10)*{1};
 (-17.3,-3)*{\bullet}+(-4,1)*{\scs a-1}; (-13.5,-4)*{\cdots};
 (-9,-3)*{\bullet}+(5.1,2)*{\scs k-1-b};
 (-.5,2)*{\bullet}+(-2,1)*{\scs b};
 (-4,-4)*{\cdots};
 (7.8,2)*{\bullet}+(-4,1)*{\scs k-1}; (12.5,3)*{\cdots};
 (16.5,2)*{\bullet};
  \endxy
\end{equation}
The terms in the summation can be rewritten using
\eqref{eq_various_partial} and other relations as
\begin{equation}
\xy
 (0,0)*{\includegraphics[scale=0.45]{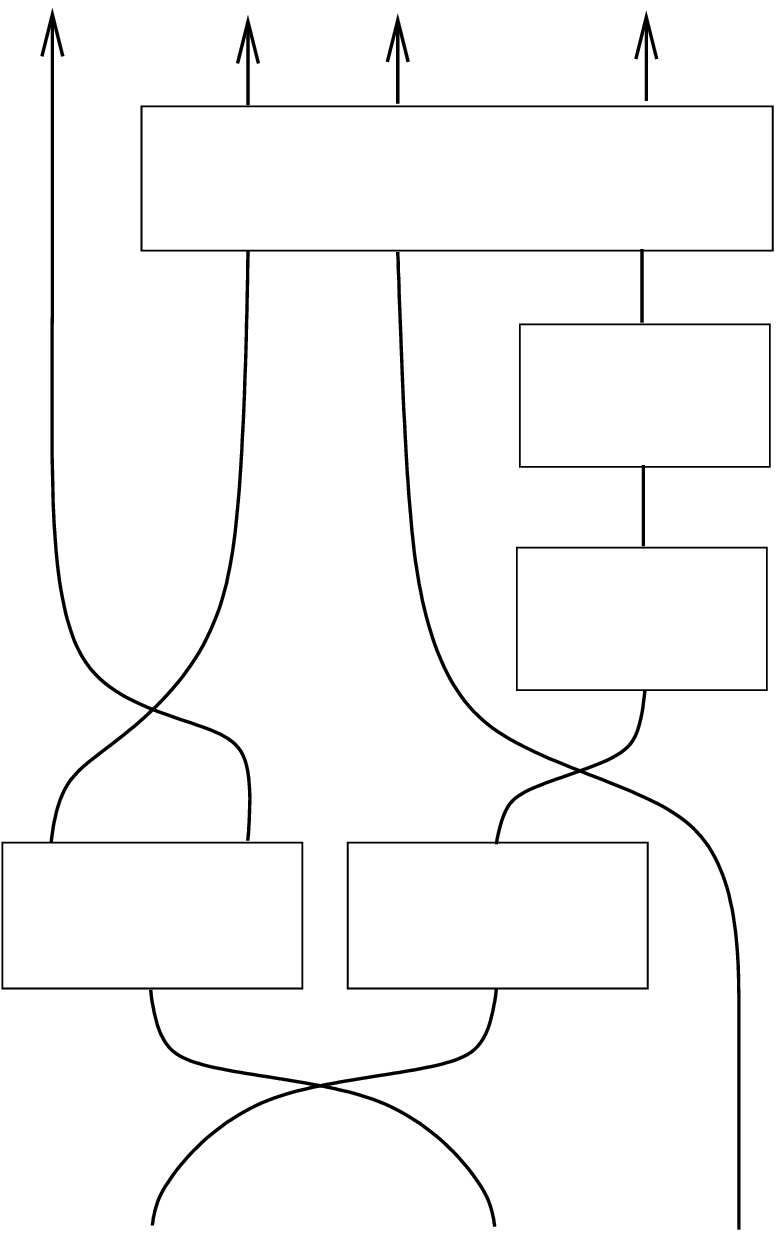}};
(18,-22)*{1};(-18,18)*{1};(5,-14)*{D_k};(-10.9,-14)*{D_{a-k}};
(11,10)*{e_k};(11,-0.5)*{\delta_k};
 (-15,-8)*{\bullet}+(-4,1)*{\scs a-1}; (-11,-9)*{\cdots};
 (-6.5,-8)*{\bullet}+(5.2,1)*{\scs k-1-b};
 (2.5,-3)*{\bullet}+(-2.5,1)*{\scs b};
 (1,20.5)*{e_a};
  \endxy
  \;\; = \;\;
\xy
 (0,0)*{\includegraphics[scale=0.45]{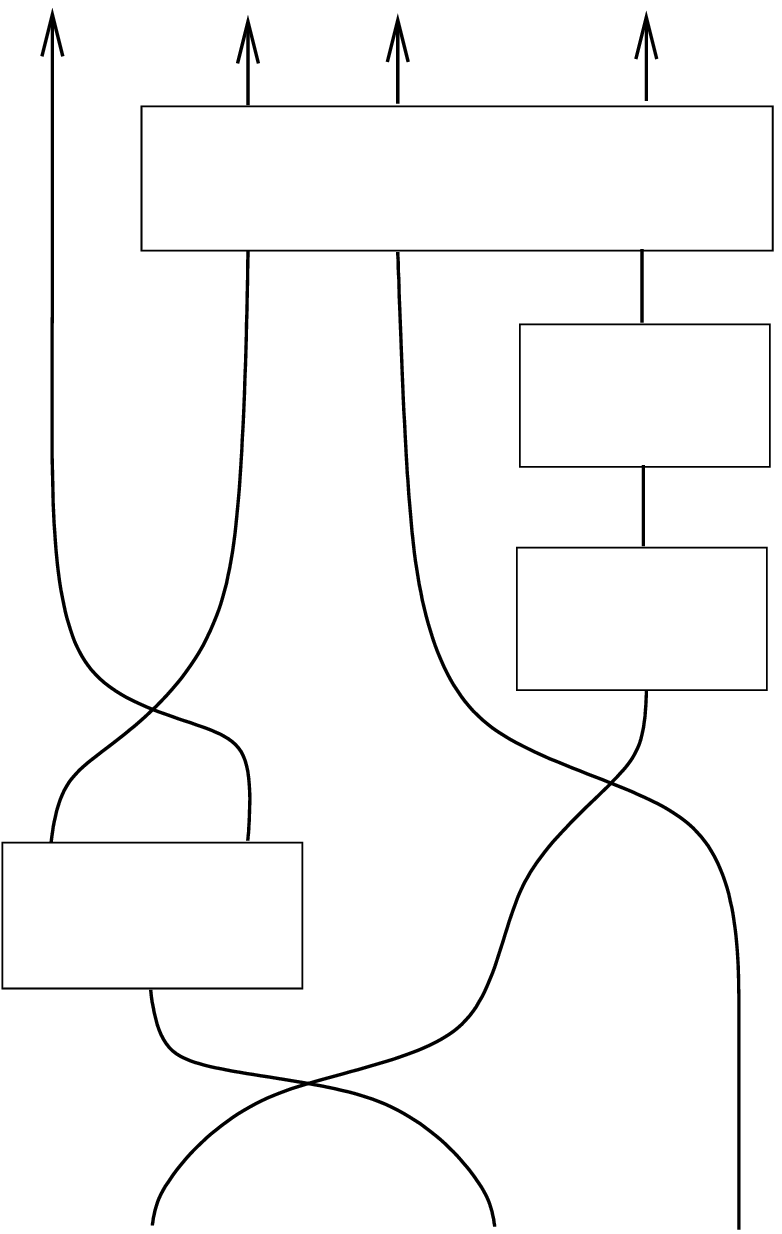}};
(18,-22)*{1};(-18,18)*{1};(-10.9,-14)*{D_{a-k}};
(11,10)*{e_k};(11,-0.5)*{e_k};
 (-15,-8)*{\bullet}+(-4,1)*{\scs a-1}; (-11,-9)*{\cdots};
 (-6.5,-8)*{\bullet}+(5.2,1)*{\scs k-1-b};
 (2.5,-3)*{\bullet}+(-2.5,1)*{\scs b};
 (1,20.5)*{e_a};
  \endxy
   \;\; \refequal{\eqref{eq_projector_absorb}} \;\;
 \xy
 (-3.8,0)*{\includegraphics[scale=0.45]{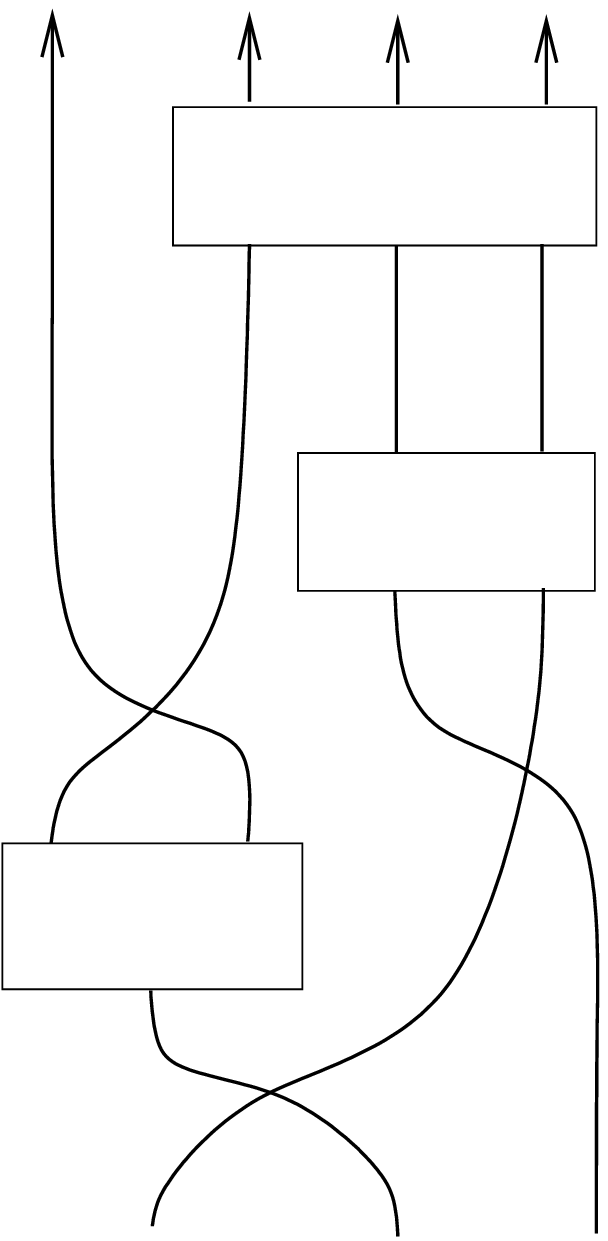}};
(12,-22)*{1};(-16,18)*{1};(-10.9,-14)*{D_{a-k}}; (4,4)*{e_{k+1}};
 (-15,-8)*{\bullet}+(-4,1)*{\scs a-1}; (-11,-9)*{\cdots};
 (-6.5,-8)*{\bullet}+(5.2,1)*{\scs k-1-b};
 (1,-3)*{\bullet}+(-2.5,1)*{\scs b};
 (1,20.5)*{e_a};
  \endxy
\end{equation}
where in the first equality we have slid the $D_{k}$ across the line
and combined it with the $\delta_k$ to form an $e_k$.  But each of
the above terms is zero by \eqref{eq-boxes_undercross1} since $b
\leq k-1$.  Hence, \eqref{eq_lem_slidedots2} becomes
\begin{eqnarray}
=\;\;
 \xy
 (4.6,0)*{\includegraphics[scale=0.45]{figs/lem1.eps}};
 (-4,-9)*{D_a};(5,9)*{e_a};
 (-9.8,3)*{\bullet}+(-3.5,1)*{\scs a-1};
 (-.7,3)*{\bullet}+(-3.5,1)*{\scs a-2};
 (8.3,3)*{\bullet}+(-2.5,1)*{\scs 2};(3,2)*{\cdots};
 (17.5,3)*{\bullet}
  \endxy
\;\; = \;\;
 \xy
 (0,0)*{\includegraphics[scale=0.45]{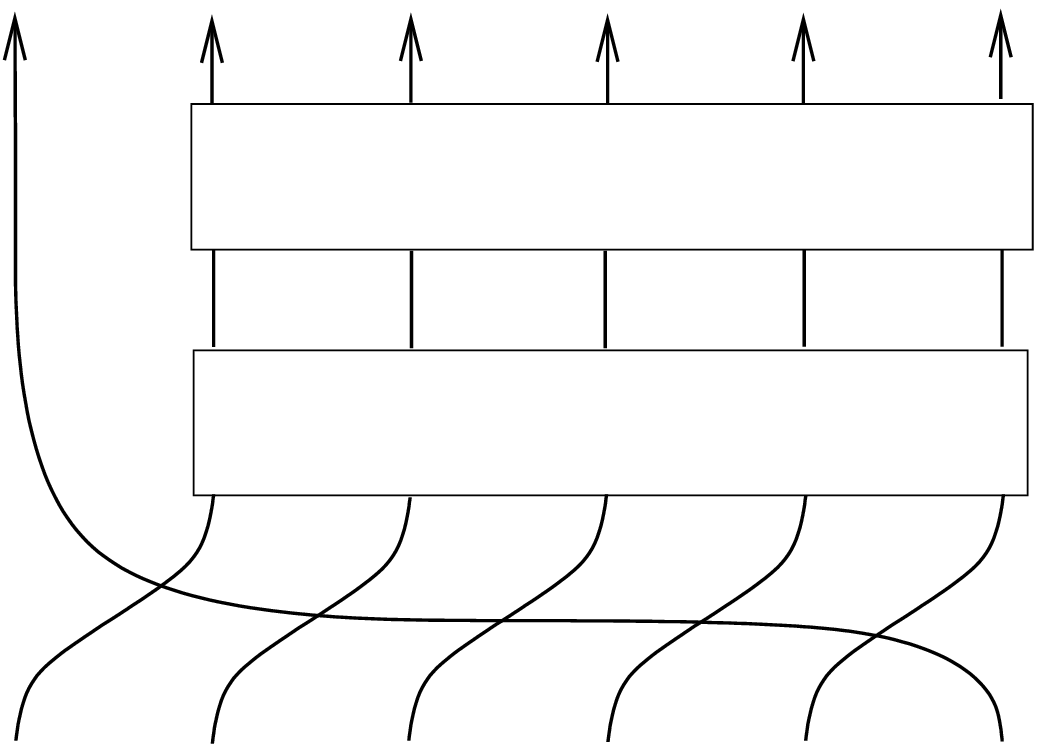}};
 (5,-2.5)*{e_a};(5,9)*{e_a};
  \endxy
 \;\; = \;\;
 \xy
 (2,0)*{\reflectbox{\includegraphics[scale=0.45]{figs/box-ul.eps}}};
 (-7,-7.5)*{a};(6,7)*{e_a};
 (-7.5,7.5)*{1};(7,-7.5)*{1};
  \endxy
 \nn
\end{eqnarray}
The case with generic $b$ follows immediately
\begin{equation}
 \xy
 (0,0)*{\reflectbox{\includegraphics[scale=0.45]{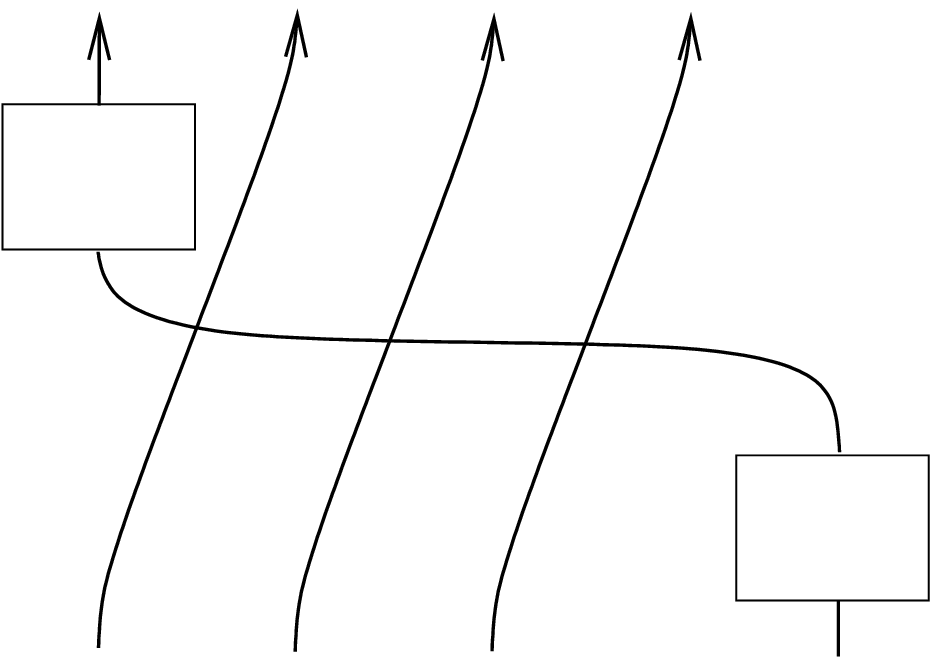}}};
 (-17,-9)*{e_a};(17,7)*{e_a};
 (5.5,-12)*{1};(-4,-12)*{1};(14,-12)*{1};
  \endxy
\;\; = \;\;
 \xy
 (0,0)*{\reflectbox{\includegraphics[scale=0.4]{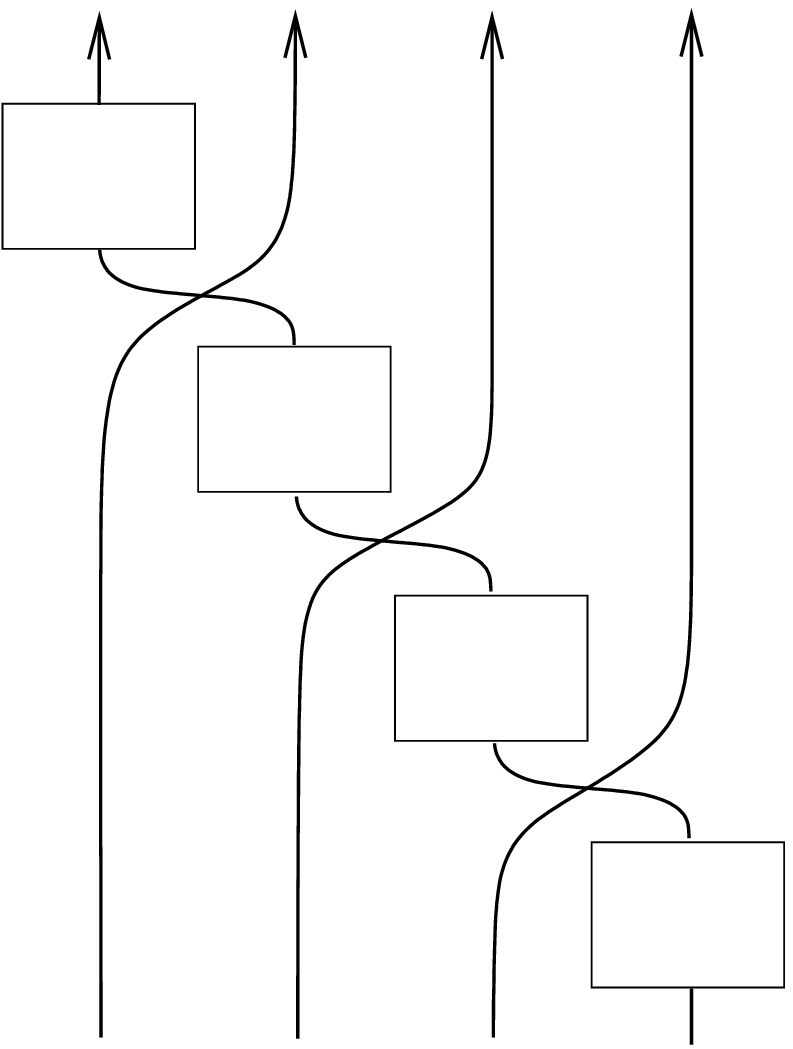}}};
(-11.5,-16)*{e_a};(-4,-6)*{e_a};(4,4.5)*{e_a};(12,14)*{e_a};
 (5.5,-16)*{1};(-1,-16)*{1};(14,-16)*{1};
  \endxy
 \;\; = \;\;
  \xy
 (1.8,0)*{\reflectbox{\includegraphics[scale=0.45]{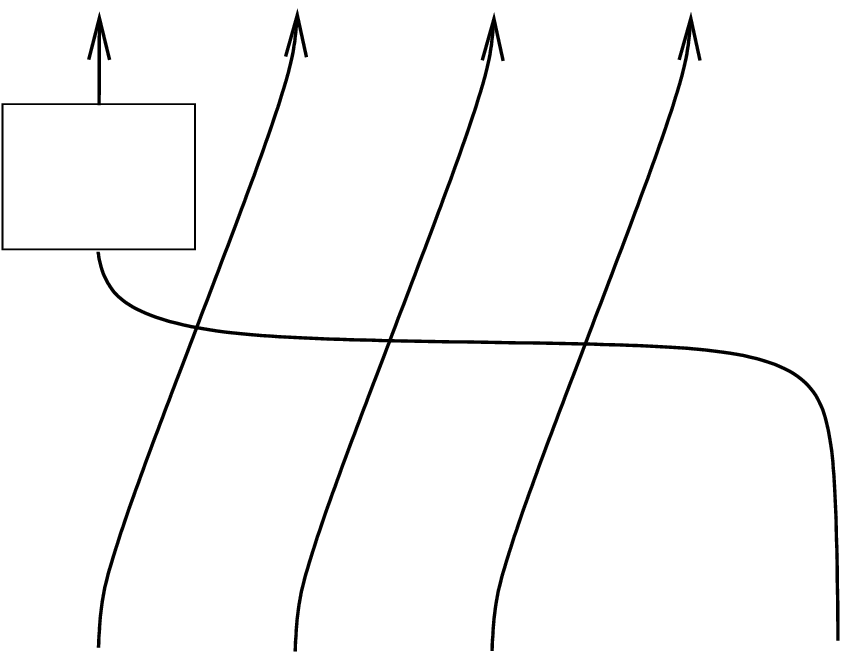}}};
 (17,7)*{e_a};
 (5.5,-12)*{1};(-4,-12)*{1};(14,-12)*{1};
  \endxy
\end{equation}
by repeatedly applying the $b=1$ case.
\end{proof}

\begin{lem} \label{lem_alpha_slide_up}
\begin{equation}
\vcenter{ \xy
  (-12,20)*{}="t1";(-4,20)*{}="t2";(4,20)*{}="t3"; (12,20)*{}="t4";
  (-12,4)*{}="m1";(-4,4)*{}="m2";(4,4)*{}="m3"; (12,4)*{}="m4";
  "t1";"m4" **\crv{(-12,10) & (12,10)}?(0)*\dir{<}?(.88)*\dir{}+(0,0)*{\bullet}+(2,2)*{\scs b};
  "t2";"m1" **\crv{(-4,10) & (-12,10)}?(0)*\dir{<};
  "t3";"m2" **\crv{(4,10) & (-4,10)}?(0)*\dir{<};"t4";"m3" **\crv{(12,10) & (4,10)}?(0)*\dir{<};
  (-15,4);(15,4) **\dir{-}; (15,4);(15,-4) **\dir{-}; (15,-4);(-15,-4) **\dir{-};
  (-15,-4); (-15,4)*{}**\dir{-};
  (-12,-4);(-12,-10) **\dir{-};(-4,-4);(-4,-10) **\dir{-};(4,-4);(4,-10) **\dir{-};(12,-4);(12,-10) **\dir{-};
  (0,0)*{D_a};
 \endxy}
\qquad = \qquad (-1)^{a-1} \sum_{\ell_1,\ell_2, \dots, \ell_a}
\;\;\vcenter{ \xy
  (-12,20)*{}="t1";(-4,20)*{}="t2";(4,20)*{}="t3"; (12,20)*{}="t4";
  (-12,4)*{}="m1";(-4,4)*{}="m2";(4,4)*{}="m3"; (12,4)*{}="m4";
  "t1";"m1" **\dir{-}?(0)*\dir{<}?(.4)*\dir{}+(0,0)*{\bullet}+(3,1)*{\scs \ell_1};
  "t2";"m2" **\dir{-}?(0)*\dir{<}?(.4)*\dir{}+(0,0)*{\bullet}+(5,1)*{\scs \ell_{2}};
  "t4";"m4" **\dir{-}?(0)*\dir{<}?(.4)*\dir{}+(0,0)*{\bullet}+(3,1)*{\scs \ell_a};
  (-15,4);(15,4) **\dir{-}; (15,4);(15,-4) **\dir{-}; (15,-4);(-15,-4) **\dir{-};
  (-15,-4); (-15,4)*{}**\dir{-};
  (-12,-4);(-12,-10) **\dir{-};(-4,-4);(-4,-10) **\dir{-};(4,-4);(4,-10) **\dir{-};(12,-4);(12,-10) **\dir{-};
  (0,0)*{D_a}; (4,10)*{\cdots};
 \endxy}
\end{equation}
where the sum is over all $\ell_1,\ell_2,\dots,\ell_a \geq 0$ such
that $\sum_{j=1}^{a}\ell_{j}=b+1-a$.
\end{lem}

\begin{proof}
Easy.
\end{proof}

We now introduce a new type of lines called {\em thick lines}. Thick
lines will have two interpretations.  In the context of graphical
calculus for the nilHecke ring, a thick line labelled $a$ is just
a thin line labelled $a$ with a projector $e_a$ on it.  Later, thick
lines will be interpreted as 2-morphisms in a suitable Karoubi
envelope 2-category, and their thick endpoints will correspond to
certain 1-morphisms.  For now, we denote
\[ 
  \xy
 (0,0)*{\includegraphics[scale=0.5]{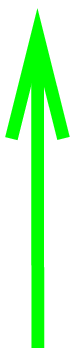}};
 (-2.5,-6)*{a}; 
  \endxy
  \quad : = \quad
 \xy
 (0,0)*{\includegraphics[scale=0.5]{figs/c1-1.eps}};
 (0,0)*{e_a}; 
  \endxy
  \quad
  = \quad
 \xy
 (0,0)*{\includegraphics[scale=0.5]{figs/box-up.eps}};
 (0,0)*{e_a}; 
  \endxy
\]
The notation is consistent, since $e_a$ is an idempotent, so that
cutting a thick line into two pieces and converting each of them
into an $e_a$ box results in the same element of the nilHecke ring.

The following diagrams will be referred to as {\em splitters} or
splitter diagrams.
\[ 
    \xy
 (0,0)*{\includegraphics[scale=0.5]{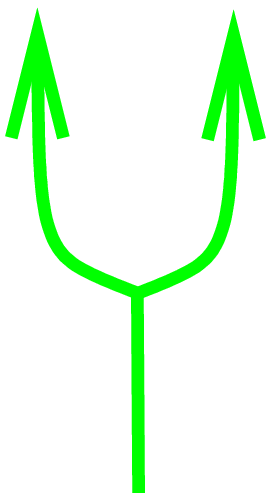}};
 (-5,-10)*{a+b};(-8,4)*{a};(8,4)*{b};
  \endxy
 \;\; :=\;\;
     \xy
 (0,0)*{\includegraphics[scale=0.5]{figs/def-tsplitu.eps}};
 (-7,-10)*{b};(7,-10)*{a};(-6,4)*{e_a};(6,4)*{e_b}; 
  \endxy
\qquad  \qquad 
    \xy
 (0,0)*{\includegraphics[scale=0.5,angle=180]{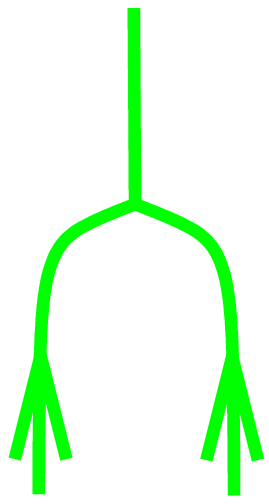}};
 (-5,10)*{a+b};(-8,-4)*{a};(8,-4)*{b}; 
  \endxy
  \;\; :=\;\;
     \xy
 (0,0)*{\includegraphics[scale=0.5,angle=180]{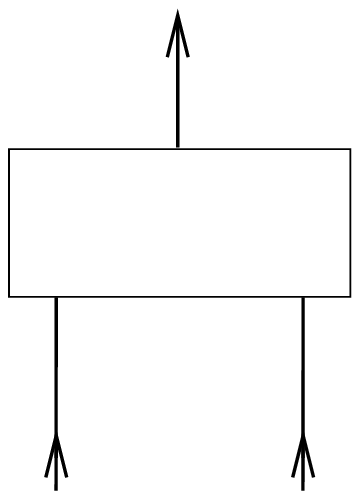}};
 (-8,-7)*{a};(8,-7)*{b};(0,1)*{e_{a+b}};(-5,10)*{a+b}; 
  \endxy
\]
The first definition is consistent in view of the relation
\eqref{eq_eaebcross}. Equation \eqref{eq_projector_absorb} ensures the
consistency of the second definition.

\begin{prop}[Associativity of splitters]
\begin{equation} 
  \xy
 (0,0)*{\includegraphics[scale=0.5]{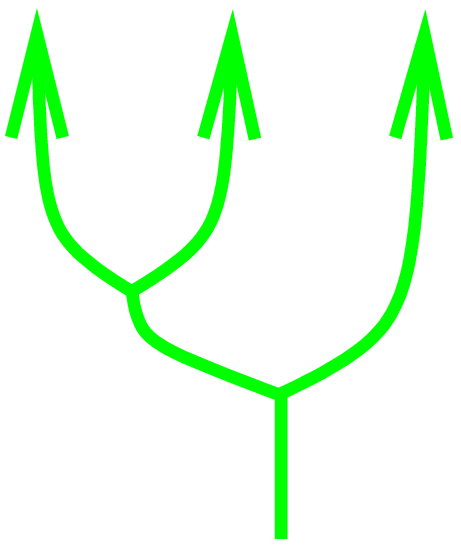}};
 (-5.5,-12)*{a+b+c};(-12,4)*{a};(-2.5,4)*{b};(12,4)*{c};
  \endxy
  \quad = \quad
   \xy
 (0,0)*{\reflectbox{\includegraphics[scale=0.5]{figs/uassoc.eps}}};
 (-10.5,-12)*{a+b+c};(-12,4)*{a};(-2.5,4)*{b};(12,4)*{c}; 
  \endxy
\qquad \qquad
   \xy
 (0,0)*{\reflectbox{\includegraphics[scale=0.5, angle=180]{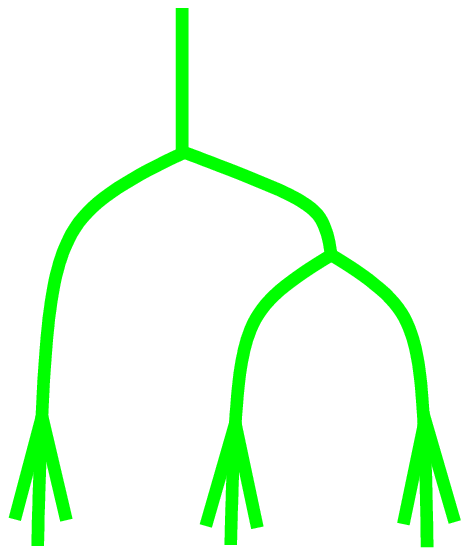}}};
 (-5.5,10)*{a+b+c};(-12,-4)*{a};(-2.5,-4)*{b};(12,-4)*{c};
  \endxy
  \quad = \quad
 \xy
 (0,0)*{\includegraphics[scale=0.5,angle=180]{figs/dassoc.eps}};
 (-10.5,10)*{a+b+c};(-12,-4)*{a};(-2.5,-4)*{b};(12,-4)*{c}; 
  \endxy \label{eq_split_assoc}
 \end{equation}
\end{prop}

\begin{proof}
\begin{eqnarray}
  \xy
 (0,0)*{\includegraphics[scale=0.5]{figs/uassoc.eps}};
 (-5.5,-12)*{a+b+c};(-12,4)*{a};(-2.5,4)*{b};(12,4)*{c};
  \endxy &  = &
  \xy
 (0,0)*{\includegraphics[scale=0.4]{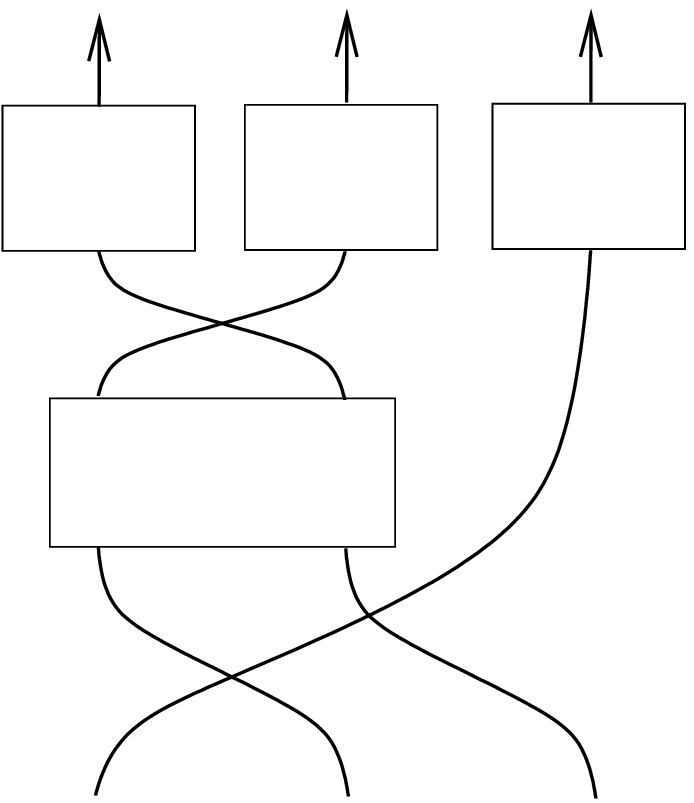}};
 (-5,-3.5)*{e_{a+b}};(-10,9)*{e_a};(0,9)*{e_b};(10.5,9)*{e_c};
  \endxy
 \;\; \refequal{\eqref{eq_eaebcross}} \;\;
   \xy
 (0,0)*{\includegraphics[scale=0.4]{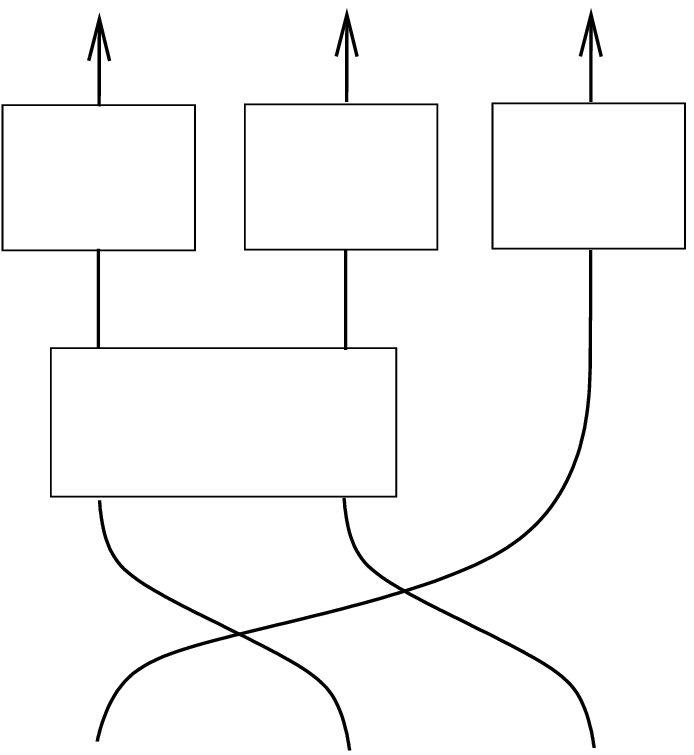}};
 (-5,-2.5)*{D_{a+b}};(-10,8)*{\delta_a};(0,8)*{\delta_{b}};(10.5,8)*{e_c};
  \endxy \nn \hspace{1.5in}
\end{eqnarray}
\begin{eqnarray}
& \refequal{\eqref{eq_def_ea}} &
   \xy
 (0,0)*{\includegraphics[scale=0.4]{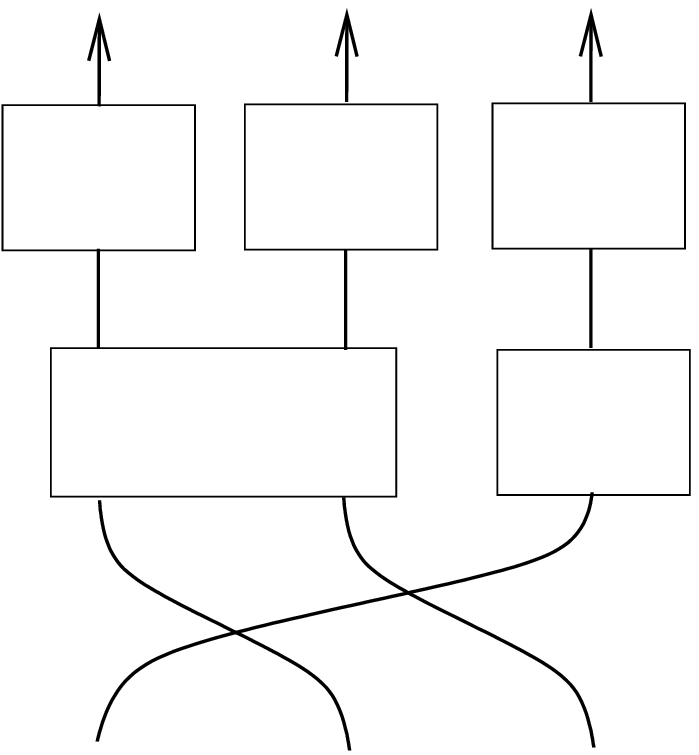}};
 (-5,-2.5)*{D_{a+b}};(-10,8)*{\delta_a};(0,8)*{\delta_{b}};(10.5,8)*{\delta_c};
 (10.5,-2.5)*{D_{c}}; 
  \endxy
\;\;  \refequal{\eqref{eq_various_partial}}  \;\;
    \xy
 (0,0)*{\includegraphics[scale=0.45]{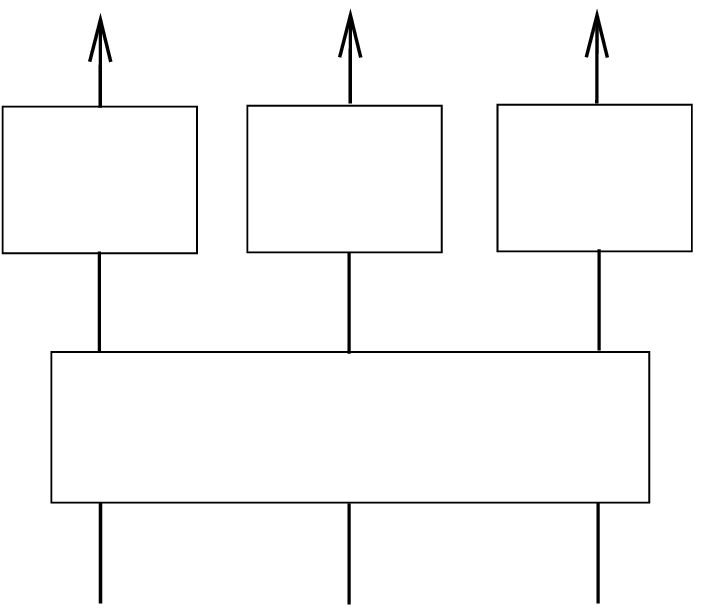}};
 (0,-6)*{D_{a+b+c}};(-11,5.5)*{\delta_a};(0,5.5)*{\delta_{b}};(11,5.5)*{\delta_{c}};
  \endxy
 \;\; \refequal{\eqref{eq_various_partial}} \;\; \; \xy
 (0,0)*{\reflectbox{\includegraphics[scale=0.4]{figs/assoc2.eps}}};
 (5,-2.5)*{D_{a+b}};(-10,8)*{e_a};(0,8)*{\delta_{b}};(10.5,8)*{\delta_{c}};
  \endxy \nn\\
   & \refequal{\eqref{eq_eaebcross}} &
    \xy
 (0,0)*{\reflectbox{\includegraphics[scale=0.4]{figs/assoc1.eps}}};
 (5,-3.5)*{e_{a+b}};(-10,9)*{e_a};(0,9)*{e_b};(10.5,9)*{e_c}; 
  \endxy
 \;\; = \;\;
   \xy
 (0,0)*{\reflectbox{\includegraphics[scale=0.5]{figs/uassoc.eps}}};
 (5.5,-12)*{a+b+c};(-12,4)*{a};(-2.5,4)*{b};(12,4)*{c}; 
  \endxy \nn
\end{eqnarray}
For the second part of the lemma we have
\begin{eqnarray}
    \xy
 (0,0)*{\reflectbox{\includegraphics[scale=0.5, angle=180]{figs/dassoc.eps}}};
 (-5.5,12)*{a+b+c};(-12,-4)*{a};(-2.5,-4)*{b};(12,-4)*{c};
  \endxy
& = &
  \xy
 (0,0)*{\includegraphics[scale=0.5]{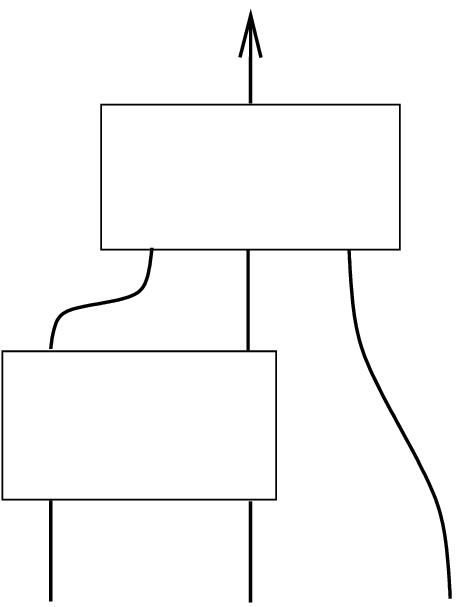}};
 (1,7)*{e_{a+b+c}};(-4,-6)*{e_{a+b}}; 
  \endxy
 \;\; \refequal{\eqref{eq_projector_absorb}} \;\;
  \xy
 (0,0)*{\includegraphics[scale=0.5]{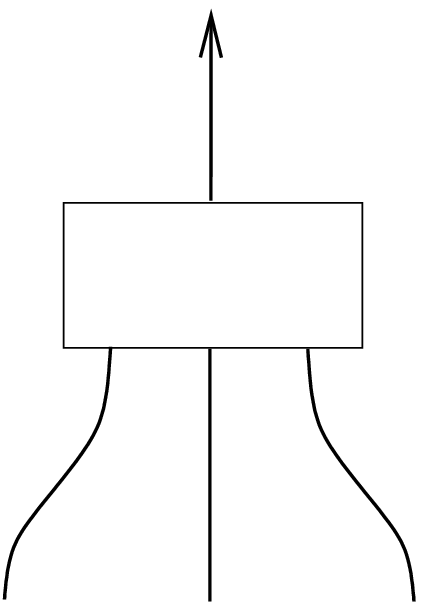}};
 (1,1)*{e_{a+b+c}}; 
  \endxy \\ \nn
 & \refequal{\eqref{eq_projector_absorb}}&
  \xy
 (0,0)*{\reflectbox{\includegraphics[scale=0.5]{figs/dassoc1.eps}}};
 (-1,7)*{e_{a+b+c}};(4,-6)*{e_{b+c}}; 
  \endxy
  \;\; = \;\;   \xy
 (0,0)*{\includegraphics[scale=0.5,angle=180]{figs/dassoc.eps}};
 (5.5,12)*{a+b+c};(-12,-4)*{a};(-2.5,-4)*{b};(12,-4)*{c}; 
  \endxy \nn
\end{eqnarray}
since smaller projectors can be absorbed into larger projectors via
\eqref{eq_projector_absorb}.
\end{proof}

Using the above relations we write
\begin{equation}
  \xy
 (0,0)*{\reflectbox{\includegraphics[scale=0.6]{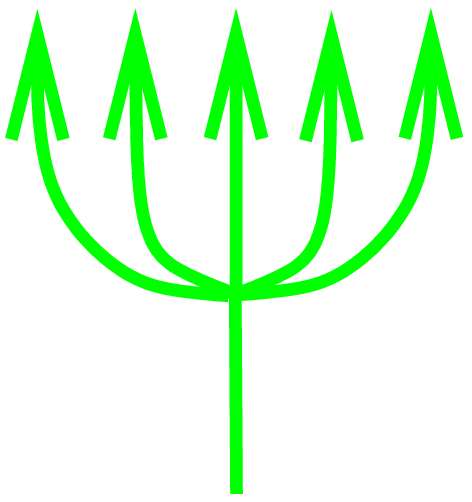}}};
 (-14,4)*{\scs \alpha_1};(14,4)*{\scs \alpha_n};(-8,4)*{\scs \alpha_2};(5,-9)*{\sum \alpha_i};
  \endxy
\qquad \qquad
  \xy
 (0,0)*{\reflectbox{\includegraphics[scale=0.6]{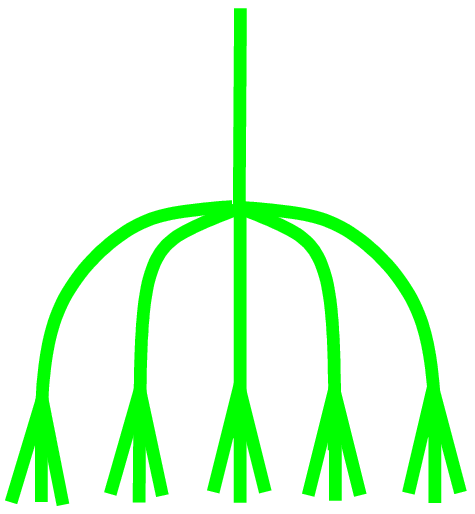}}};
 (-14,-4)*{\scs \alpha_1};(14,-4)*{\scs \alpha_n};(-8,-4)*{\scs \alpha_2};(5,9)*{\sum \alpha_i};
  \endxy
\end{equation}
to represent any of the equivalent ways of merging $n$ lines of
thicknesses $\alpha_i$ into a single line of thickness
$\sum_i\alpha_i$.  Define a thick crossing by

\begin{eqnarray} \label{eq_defn_thick_cross}
\xy
 (0,0)*{\includegraphics[scale=0.5]{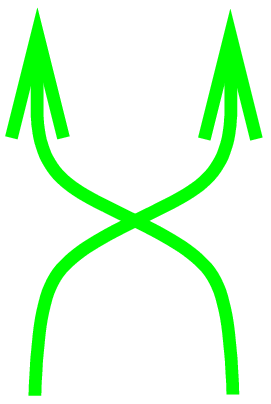}};
(-7,-7.5)*{a};(7,-7.5)*{b};
  \endxy \quad := \quad  \xy
 (0,0)*{\includegraphics[scale=0.45]{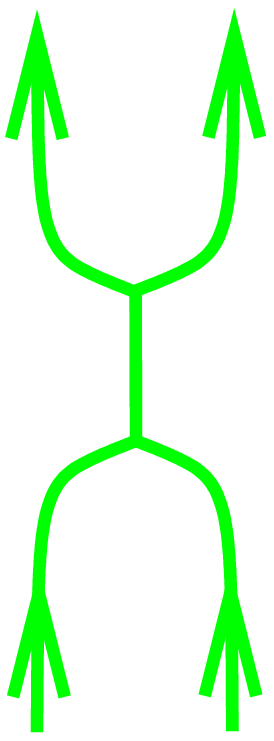}};
(-8,-11)*{a};(8,-11)*{b};(-8,11)*{b};(8,11)*{a};
  \endxy
 \quad = \quad
  \xy
 (0,0)*{\includegraphics[scale=0.5]{figs/pp1.eps}};
 (0,-8.5)*{e_{a+b}};(-6,8.5)*{e_{b}};(6,8.5)*{e_{a}};
  \endxy
   \quad \refequal{\eqref{eq_eaebcross}} \quad
     \xy
 (0,0)*{\includegraphics[scale=0.5]{figs/def-tsplitu.eps}};
 (-7,-10)*{a};(7,-10)*{b};(-6,4)*{e_b};(6.5,4)*{e_a};
  \endxy
\end{eqnarray}

\begin{prop}\label{prop_almostRthree}
For $a,b,c \geq 0$
\begin{equation} \label{eq_triangle}
 \xy
 (0,0)*{\includegraphics[scale=0.5]{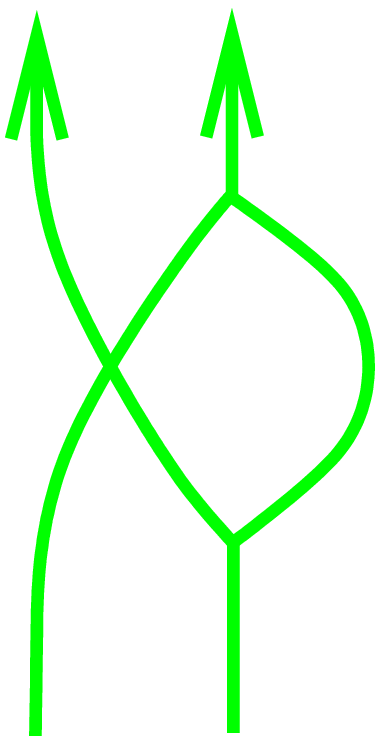}};
 (-10,-13)*{c};(8,-13)*{a+b};(-10,6)*{a};(11,4)*{b};
  \endxy
 \quad = \quad
         \xy
 (0,0)*{\includegraphics[scale=0.5]{figs/tH.eps}};
(-8,-11)*{c};(12,-11)*{a+b};(-8,11)*{a};(12,11)*{b+c};
  \endxy
\qquad \qquad \xy
 (0,0)*{\reflectbox{\includegraphics[scale=0.5]{figs/triangle.eps}}};
 (10,-13)*{c};(-8,-13)*{a+b};(10,6)*{a};(-11,4)*{b};
  \endxy
 \quad = \quad
         \xy
 (0,0)*{\includegraphics[scale=0.5]{figs/tH.eps}};
(8,-11)*{c};(-12,-11)*{a+b};(8,11)*{a};(-12,11)*{b+c};
  \endxy
\end{equation}
\end{prop}

\begin{proof}  The first identity is
\begin{eqnarray}
   \xy
 (0,0)*{\includegraphics[scale=0.5]{figs/triangle.eps}};
 (-10,-13)*{c};(8,-13)*{a+b};(-10,6)*{a};(11,4)*{b};
  \endxy
 \quad = \quad
   \xy
 (0,0)*{\includegraphics[scale=0.5]{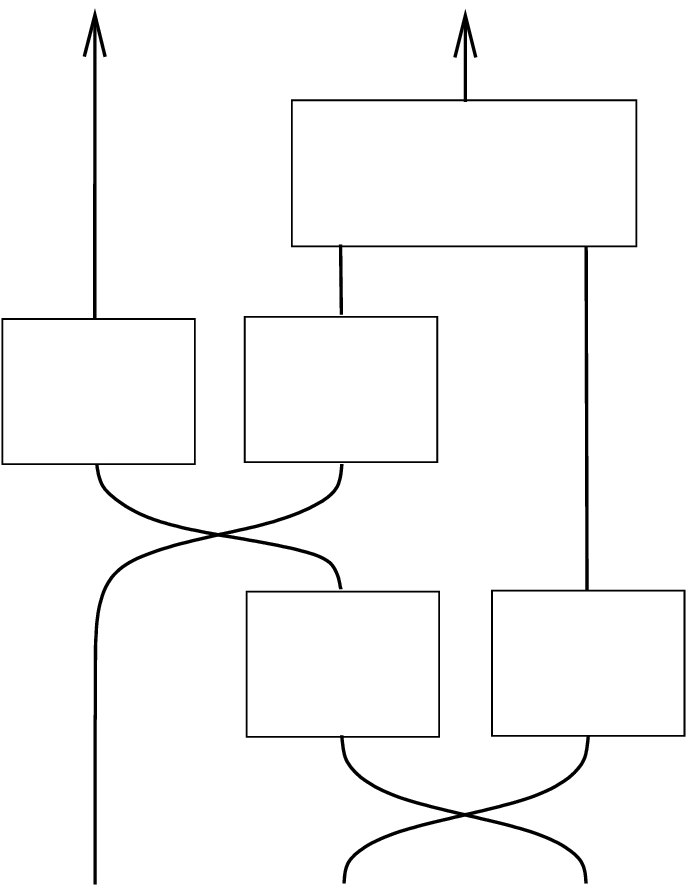}};
 (7,13.5)*{e_{b+c}};(0,-12)*{e_a};(13,-12)*{e_b};(-12.5,2)*{e_a};(0,2)*{e_c};
  \endxy
 \quad \refequal{\eqref{eq_lem-twobox},\eqref{eq_projector_absorb}} \quad
    \xy
 (0,0)*{\includegraphics[scale=0.5]{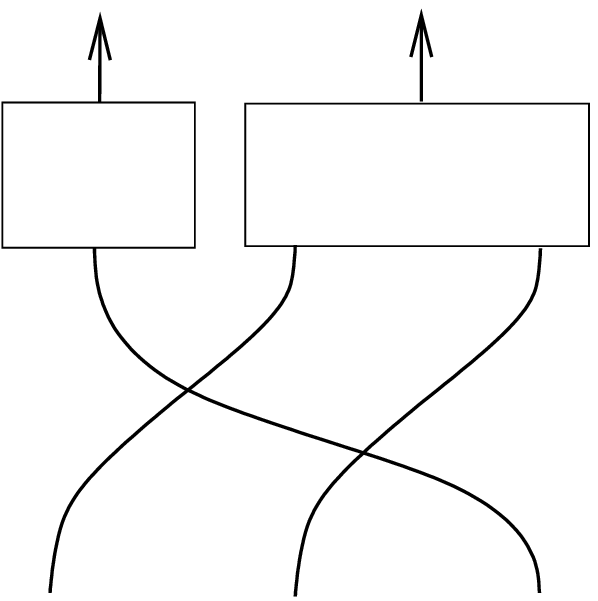}};
 (7,6)*{e_{b+c}};(-10,6)*{e_a};
  \endxy \nn \\
  \quad = \quad
   \xy
 (0,0)*{\includegraphics[scale=0.5]{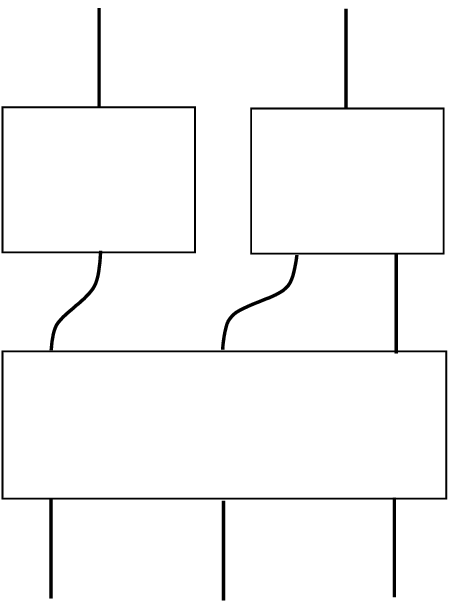}};
 (6.5,6)*{\delta_{b+c}};(0.5,-6.5)*{D_{a+b+c}};(-6.5,6)*{\delta_a};
  \endxy
   \quad \refequal{\eqref{eq_partial-en}} \quad
   \xy
 (0,0)*{\includegraphics[scale=0.5]{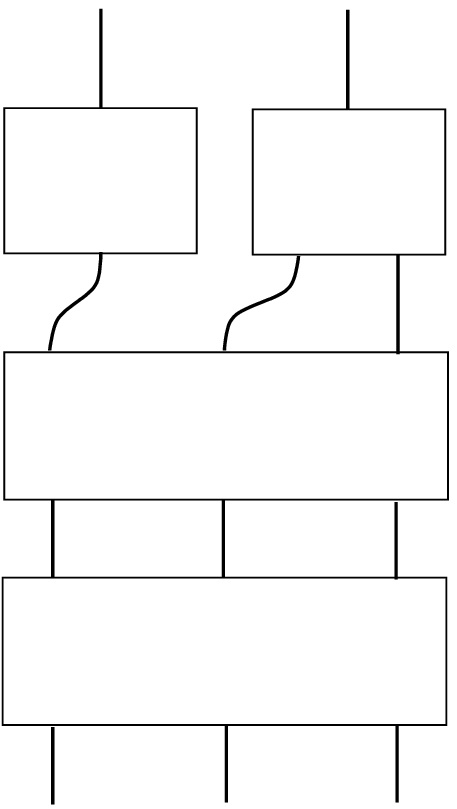}};
 (6.5,11.5)*{\delta_{b+c}};(0.5,-1)*{D_{a+b+c}};(-6.5,11.5)*{\delta_a};
 (0.5,-12.5)*{e_{a+b+c}};
  \endxy
 \quad = \quad
\xy
 (0,0)*{\includegraphics[scale=0.5]{figs/pp1.eps}};
 (0,-8.5)*{e_{a+b+c}};(-6,8.5)*{e_{a}};(6,8.5)*{e_{b+c}};
  \endxy
  \quad = \quad
         \xy
 (0,0)*{\includegraphics[scale=0.5]{figs/tH.eps}};
(-8,-11)*{c};(12,-11)*{a+b};(-8,11)*{a};(12,11)*{b+c};
  \endxy \nn
\end{eqnarray}
The second identity is proven similarly.
\end{proof}

\begin{cor}
\begin{equation}
 \xy
  (0,0)*{\includegraphics[scale=0.5]{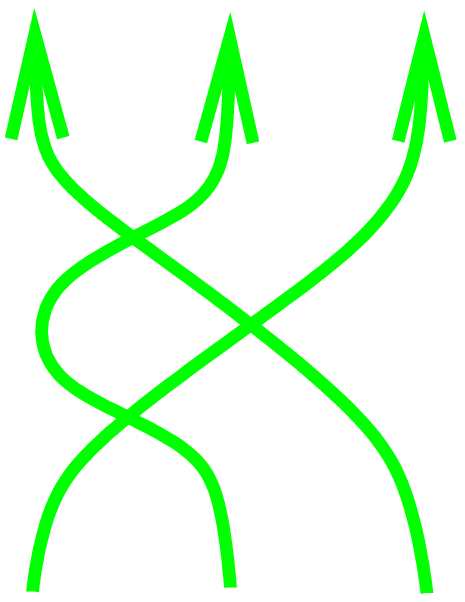}};
  (-13,-13)*{a};   (3,-13)*{b};    (13,-13)*{c};
  (-13,10)*{c};   (3,10)*{b};    (13,10)*{a};
 \endxy \quad = \quad
 \xy
  (0,0)*{\includegraphics[scale=0.5]{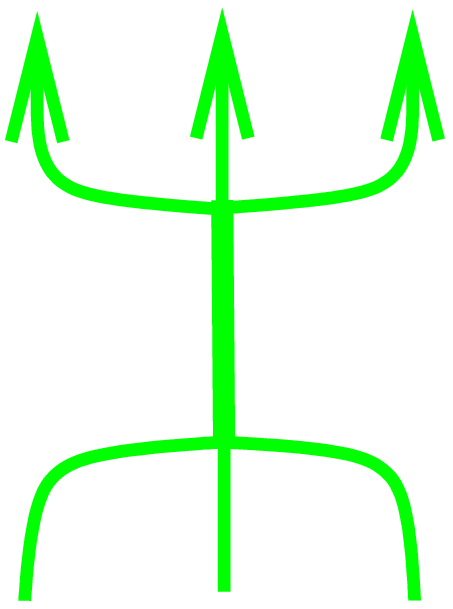}};
  (-13,-13)*{a};   (3,-13)*{b};    (13,-13)*{c};
  (-13,10)*{c};   (3,10)*{b};    (13,10)*{a};
 \endxy
 \quad = \quad
  \xy
  (0,0)*{\reflectbox{\includegraphics[scale=0.5]{figs/r3.eps}}};
    (-13,-13)*{a};   (3,-13)*{b};    (13,-13)*{c};
  (-13,10)*{c};   (3,10)*{b};    (13,10)*{a};
 \endxy
\end{equation}
\end{cor}
\begin{proof}
This follows from Proposition~\ref{prop_almostRthree}.
\end{proof}

%
\subsection{Partitions, symmetric functions, and their diagrammatics}
%

By a partition $\alpha=(\alpha_1,\alpha_2,\ldots,\alpha_a)$ we mean
a nonincreasing sequence of nonnegative integers, i.e.
$\alpha_1\ge\alpha_2\ge\ldots\ge\alpha_a\ge 0$. By convention,
$\alpha_i=0$ for $i\ge a+1$. We define $|\alpha|=\sum_{i=1}^a
{\alpha_i}$.

By $P(a,b)$ we denote the set of all partitions $\alpha$ with at
most $a$ parts (i.e. with $\alpha_{a+1}=0$) such that $\alpha_1\le
b$. That is, $P(a,b)$ consists of partitions that fit into a box of
size $a \times b$. Moreover, the set of all partitions with at most
$a$ parts (i.e. the set $P(a,\infty)$) we denote simply by $P(a)$.  $P(0)=\{ \emptyset\}$ is the set of all partitions with at most $0$ parts, so that $P(0)$ contains only the empty partition.

The cardinality of $P(a,b)$ is ${a+b \choose a}$. We call
\begin{equation} \label{eq_Pab_card}
  |P(a,b)|_q:=\sum_{\alpha\in P(a,b)} {q^{2|\alpha|-ab}}=\left[{a+b
\atop a}\right]
\end{equation}
the $q$-cardinality of $P(a,b)$. The last term in
the above equations is the balanced q-binomial.

The dual (conjugate) partition of  $\alpha$ is the partition
$\bar{\alpha}=(\bar{\alpha}_1,\bar{\alpha}_2,\ldots)$ with
$\bar{\alpha}_j=\sharp\{i|\alpha_i\ge j\}$ and given by reflecting
the Young diagram of $\alpha$ along the diagonal. For a partition
$\alpha\in P(a,b)$ we define the {\it complementary partition}
$\alpha^c=(b-\alpha_a,\ldots,b-\alpha_2,b-\alpha_1)$. To be more explicit we let $K=(b^a)\in P(a,b)$ and write $K-\alpha$ in place of $\alpha^c$ to emphasize that $\alpha \in P(a,b)$. Finally we
define $\hat{\alpha}:=\overline{\alpha^c}$. Note that $\bar{\alpha}$ and
$\hat{\alpha}$ belong to $P(b,a)$ and $\alpha^c $ to $P(a,b)$.

Schur polynomials $\pi_{\alpha}$, taken over all partitions $\alpha$ with at most $a$ parts, form an additive basis of the graded ring of symmetric functions $\Z[x_1,\dots,x_a]^{S_a}$.   Here we review basic definitions and properties of Schur polynomials -- for more details see \cite{Fulton,Mac,Man}. Recall from Section~\ref{sec_nilHecke} that we set $\deg(x_i)=2$.

For every partition $\alpha$, the Schur polynomial
$\pi_{\alpha}=\pi_{\alpha_1,\alpha_2,\ldots,\alpha_a}$ is given by
the formula:
$$\pi_{\alpha}(x_1,x_2,\ldots,x_a)=\frac{|x_i^{\alpha_j+a-j}|}{\Delta},$$
where $\Delta=\prod_{1\le r < s \le a} (x_r-x_s)$, and
$|x_i^{\alpha_j+a-j}|$ is the determinant of the $a\times a$ matrix
whose $(i,j)$ entry is $x_i^{\alpha_j+a-j}$. We let
$\pi_{\alpha}(x_1,\ldots,x_a)=0$ if some entry of $\alpha$ is
negative ($\alpha$ is not a partition then), or if $\alpha_{a+1}>0$.

The elementary symmetric polynomials $\varepsilon_{m}(x_1,\ldots,x_a)$,
$m=0,\ldots,a$ and the complete symmetric polynomials
$h_m(x_1,\ldots,x_a)$ are special Schur polynomials:
\begin{eqnarray}
h_m(x_1,\ldots ,x_a)&=&\pi_{(m)}(x_1,\ldots ,x_a),\\
\varepsilon_m(x_1,\ldots ,x_a)&=&\pi_{(1^m)}(x_1,\ldots ,x_a).
\end{eqnarray}
These two families of functions are related by
\begin{equation} \label{eq_eh_rel}
  \sum_{r=0}^m(-1)^r\varepsilon_r h_{m-r} =0
\end{equation}
for all $m \geq 1$.  By convention, $\varepsilon_m =h_m=0$ for $m<0$ and
$m>a$.

Every symmetric function can be written as a polynomial in
elementary symmetric functions, and as a polynomial in complete
symmetric functions. Explicit formulas for the Schur functions are
the Jacobi-Trudy formulas, also known as Giambelli formulas. For $\alpha \in
P(a)$, the Jacobi-Trudy formula is given by the determinant:
\begin{equation}
\pi_{\alpha}=\left|
\begin{array}{cccc}
h_{\alpha_1} & h_{\alpha_1+1} & \cdots & h_{\alpha_1+a-1} \\
h_{\alpha_2-1} & h_{\alpha_2} & \cdots & h_{\alpha_2+a-2} \\
\vdots & \vdots & \ddots & \vdots \\
h_{\alpha_a-(a-1)} & h_{\alpha_a-(a-2)} & \cdots & h_{\alpha_a}
\end{array}
\right|= |h_{\alpha_i+j-i}|,
\end{equation}
and if $\bar{\alpha}\in P(b)$, then the dual formula is given by the determinant
\begin{equation} \label{eq_Schur_e}
\pi_{\alpha}=\left|
\begin{array}{cccc}
\varepsilon_{\bar{\alpha}_1} & \varepsilon_{\bar{\alpha}_1+1} & \cdots & \varepsilon_{\bar{\alpha}_1+b-1} \\
\varepsilon_{\bar{\alpha}_2-1} & \varepsilon_{\bar{\alpha}_2} & \cdots & \varepsilon_{\bar{\alpha}_2+b-2} \\
\vdots & \vdots & \ddots & \vdots \\
\varepsilon_{\bar{\alpha}_b-(b-1)} & \varepsilon_{\bar{\alpha}_b-(b-2)} & \cdots & \varepsilon_{\bar{\alpha}_b} \\
\end{array}\right|=
|\varepsilon_{\bar{\alpha}_i+j-i}|.
\end{equation}

Throughout the paper we shall use the Littlewood-Richardson
coefficients. They are the coefficients in the expression of the
product of two Schur polynomials as a linear combination of the
Schur polynomials. The Littlewood-Richardson coefficients
$c_{\alpha,\beta}^{\gamma}$ are given by:
\begin{equation} \label{eq_prod_schur}
\pi_{\alpha}\pi_{\beta}=\sum_{\gamma}
{c_{\alpha,\beta}^{\gamma}\pi_{\gamma}}.
\end{equation}
The coefficients $c_{\alpha,\beta}^{\gamma}$ are nonnegative integers that can be nonzero only when $|\gamma|=|\alpha|+|\beta|$.  Note that for all $\alpha,\beta,\gamma \in P(a)$ the coefficients $c_{\alpha,\beta}^{\gamma}$ are independent of the choice of $a$.

For any three partitions $\alpha,\beta,\gamma$, we have
\begin{equation}
c_{\alpha,\beta}^{\gamma}=c_{\beta,\alpha}^{\gamma}=c_{\bar{\alpha},\bar{\beta}}^{\bar{\gamma}}.
\end{equation}
Moreover, if $\underline{x}=(x_1,x_2,\ldots,x_a)$ and
$\underline{y}=(y_1,y_2,\ldots,x_b)$ are two sets of
variables, then
\begin{equation} \label{eq_prod_schur_twovar}
\pi_{\gamma}(\underline{x},\underline{y})=\sum_{\alpha,\beta}
{c_{\alpha,\beta}^{\gamma} \pi_{\alpha}(\underline{x})
\pi_{\beta}(\underline{y})}.
\end{equation}
Also we shall use iterated Littlewood-Richardson coefficients. For
any $k\ge 2$ and partitions $\alpha_1,\ldots,\alpha_k$ and $\beta$
define $c_{\alpha_1,\ldots,\alpha_k}^{\beta}$ by
\begin{equation}
\pi_{\alpha_1}\pi_{\alpha_2}\ldots\pi_{\alpha_k}=\sum_{\beta}{c_{\alpha_1,\ldots,\alpha_k}^{\beta}\pi_{\beta}}.
\end{equation}
These coefficients can be expressed in terms of the ``ordinary"
Littlewood-Richardson coefficients:
\begin{equation}
c_{\alpha_1,\ldots,\alpha_k}^{\beta}=\sum_{\varphi_1,\ldots,\varphi_{k-2}}
{c_{\alpha_1,\alpha_2}^{\varphi_1}c_{\varphi_1,\alpha_3}^{\varphi_2}\ldots
c_{\varphi_{k-3},\alpha_{k-1}}^{\varphi_{k-2}}c_{\varphi_{k-2},\alpha_k}^{\beta}}.
\end{equation}
Finally, if $\alpha,\beta,\gamma \in P(a,b)$ and $K=(b^a)$ then
$c_{\alpha,\beta}^K=\delta_{\alpha,\beta^c}$ and
$c_{\alpha,\beta,\gamma}^K=c_{\alpha,\beta}^{\gamma^c}$.

For a partition $\alpha=(\alpha_1,\alpha_2, \dots, \alpha_a)$ and $m\in \Z$ we write
\begin{equation} \label{eq_def_partition_add}
\alpha+m := \big(\alpha_1+m,\alpha_{2}+m, \dots ,
\alpha_a+m\big).
\end{equation}
Our conventions explained above imply that $\pi_{\alpha+m}=0$ if $\alpha+m$ is not a partition, i.e. if $\alpha_a+m<0$.

We also need the relation between the divided differences and the
Schur polynomials.  Let $\alpha=(\alpha_1,\alpha_2,\ldots,\alpha_a)$
be a partition. Then
\begin{equation} \label{eq_SchurDa}
D_a (x_1^{\alpha_1+a-1}x_2^{\alpha_2+a-2}\ldots
x_a^{\alpha_a})=\pi_{\alpha}(x_1,x_2,\ldots,x_a),
\end{equation}
where $D_a$ was defined in Section~\ref{sec_nilHecke}, see \cite[Proposition 2.3.2]{Man} or \cite[Chapter 10]{Fulton}.

For any nonnegative integers $b_1,b_2,\ldots,b_a$
\begin{equation} \label{eq_divided1}
D_a(x_1^{b_1}\dots x_r^{b_r}\dots x_s^{b_s} \dots x_a^{b_a})
=-D_a(x_1^{b_1}\dots x_r^{b_s}\dots x_s^{b_r} \dots x_a^{b_a}).
\end{equation}
In particular, if $b_r=b_s$ for some $r<s$ then
\begin{equation} \label{eq_divided2}
D_a(x_1^{b_1}\dots x_r^{b_r}\dots x_s^{b_s} \dots x_a^{b_a})=0.
\end{equation}



An element $y \in \BNC_a$ can be denoted by a box $y$ on a black line labelled $a$, see \eqref{eq_ybox}. In particular, the elementary symmetric polynomials are represented in the graphical calculus by
\begin{equation}
  \xy
 (0,0)*{\includegraphics[scale=0.5]{figs/single-up.eps}};
 (-2.5,-8)*{a};(0,0)*{\bigb{\varepsilon_s}};
  \endxy
  \quad : = \quad
   \varepsilon_s(x_1,x_2, \dots, x_a)
\end{equation}
and viewed as elements of $Z(\BNC_a)\cong \Z[x_1,\dots,x_a]^{S_a}$.
Note that our conventions imply that the above diagram is zero if $s>a$, or $s<0$.
For example, $\varepsilon_2(x_1,x_2,x_3) = x_1 x_2+ x_1 x_3+ x_2 x_3$ so that
\begin{align}
    \xy
 (0,0)*{\includegraphics[scale=0.5]{figs/single-up.eps}};
 (-2.5,-7)*{3};(0,0)*{\bigb{\varepsilon_2}};
  \endxy \qquad = \qquad  \xy
 (-5,0)*{\includegraphics[scale=0.5]{figs/single-up.eps}};
 (0,0)*{\includegraphics[scale=0.5]{figs/single-up.eps}};
 (5,0)*{\includegraphics[scale=0.5]{figs/single-up.eps}};
 (0,0)*{\bullet};  (-5,0)*{\bullet};
  \endxy \;\; + \;\;  \xy
 (-5,0)*{\includegraphics[scale=0.5]{figs/single-up.eps}};
 (0,0)*{\includegraphics[scale=0.5]{figs/single-up.eps}};
 (5,0)*{\includegraphics[scale=0.5]{figs/single-up.eps}};
 (5,0)*{\bullet}; (-5,0)*{\bullet};
  \endxy  \;\; +\;\;   \xy
 (-5,0)*{\includegraphics[scale=0.5]{figs/single-up.eps}};
 (0,0)*{\includegraphics[scale=0.5]{figs/single-up.eps}};
 (5,0)*{\includegraphics[scale=0.5]{figs/single-up.eps}};
 (0,0)*{\bullet}; (5,0)*{\bullet};
  \endxy\nn
\end{align}
We represent complete symmetric functions $h_s$ and more general Schur polynomials $\pi_{\alpha}$ using a similar notation
\begin{eqnarray}
  \xy
 (0,0)*{\includegraphics[scale=0.5]{figs/single-up.eps}};
 (-2.5,-8)*{a};(0,0)*{\bigb{h_s}};
  \endxy
   & =&
   \xy
  (-8,0)*{};
 (-6,0)*{\includegraphics[scale=0.5]{figs/single-up.eps}};
 (-2,0)*{\includegraphics[scale=0.5]{figs/single-up.eps}};
 (2,0)*{\includegraphics[scale=0.5]{figs/single-up.eps}};
 (6,0)*{\includegraphics[scale=0.5]{figs/single-up.eps}};
 (0,0)*{\bigb{\hspace{0.2in}h_s\hspace{0.2in}}};
  \endxy
  \quad : = \quad
   h_s(x_1,x_2, \dots, x_a),
  \\
    \xy
 (0,0)*{\includegraphics[scale=0.5]{figs/single-up.eps}};
 (-2.5,-8)*{a};(0,0)*{\bigb{\pi_{\alpha}}};
  \endxy
  &=&
  \xy
  (-8,0)*{};
 (-6,0)*{\includegraphics[scale=0.5]{figs/single-up.eps}};
 (-2,0)*{\includegraphics[scale=0.5]{figs/single-up.eps}};
 (2,0)*{\includegraphics[scale=0.5]{figs/single-up.eps}};
 (6,0)*{\includegraphics[scale=0.5]{figs/single-up.eps}};
 (0,0)*{\bigb{\hspace{0.2in}\pi_{\alpha}\hspace{0.2in}}};
  \endxy
  \quad : = \quad
   \pi_{\alpha}(x_1,x_2, \dots, x_a).
\end{eqnarray}

For a central element $y \in Z(\BNC_a)$, a thick line with a box labelled $y$
\begin{equation} \label{eq_ygreenbox}
  \xy
 (0,0)*{\includegraphics[scale=0.5]{figs/single-tup.eps}};
 (-2.5,-8)*{a};(0,-2)*{\bigb{y}};
  \endxy
\end{equation}
will denote the element $ye_a$.  This notation is consistent since $ye_a=e_ay=e_aye_a$.  Moreover, for central $y,z$
\begin{equation} \label{eq_yz_greenbox}
  \xy
 (0,0)*{\includegraphics[scale=0.5]{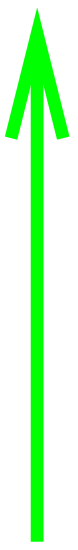}};
 (-2.5,-12)*{a};(0,2)*{\bigb{y}};(0,-5)*{\bigb{z}};
  \endxy
 \quad = \quad
   \xy
 (0,0)*{\includegraphics[scale=0.5]{figs/tlong-up.eps}};
 (-2.5,-12)*{a};(0,-1)*{\bigb{yz}};
  \endxy
\end{equation}
since $y e_a z e_a = yz e_a$.  Note that it would be unnecessary and perhaps confusing to use notation \eqref{eq_ygreenbox} to represent $e_aye_a$ for an arbitrary $y \in \BNC_a$.  First, \eqref{eq_yz_greenbox} would fail for general $y,z\in \BNC_a$.  Second, for any $y \in \BNC_a$ there exists a unique $y' \in Z(\BNC_a)$ such that $e_a y e_a = y' e_a$.

We
write
\begin{equation}
  \xy
 (0,0)*{\includegraphics[scale=0.5]{figs/single-tup.eps}};
 (-2.5,-8)*{a};(0,-2)*{\bigb{y}}+(4.5,2)*{\scs p};
  \endxy
  \quad = \quad
    \xy
 (0,0)*{\includegraphics[scale=0.5]{figs/single-tup.eps}};
 (-2.5,-8)*{a};(0,-2)*{\bigb{y^p}};
  \endxy
  \quad = \quad \left( \;\; \xy
 (0,0)*{\includegraphics[scale=0.5]{figs/single-tup.eps}};
 (-2.5,-8)*{a};(0,-2)*{\bigb{y}};
  \endxy \;\; \right)^p
\end{equation}
As special cases we have the following notations:
\begin{equation}
  \xy
 (0,0)*{\includegraphics[scale=0.5]{figs/single-tup.eps}};
 (-2.5,-8)*{a};(0,-2)*{\bigb{\varepsilon_s}};
  \endxy
  \quad := \quad
     \xy
 (0,0)*{\includegraphics[scale=0.5]{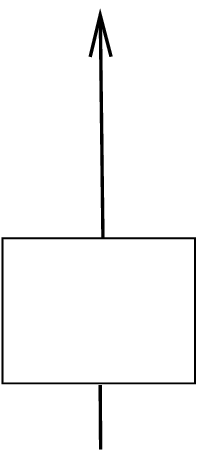}};
 (0,-4.5)*{e_a};(0,4.5)*{\bigb{\varepsilon_s}};
  \endxy
  \quad = \quad
 \xy
 (0,0)*{\includegraphics[scale=0.5]{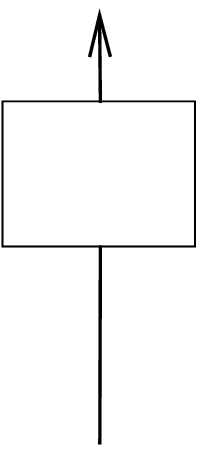}};
 (0,-6)*{\bigb{\varepsilon_s}};(0,2.5)*{e_a};
  \endxy
 \qquad \qquad
  \xy
 (0,0)*{\includegraphics[scale=0.5]{figs/single-tup.eps}};
 (-2.5,-8)*{a};(0,-2)*{\bigb{h_s}};
  \endxy
  \quad := \quad
   \xy
 (0,0)*{\includegraphics[scale=0.5]{figs/jj2.eps}};
 (0,-4.5)*{e_a};(0,4.5)*{\bigb{h_s}};
  \endxy
  \quad = \quad
  \xy
 (0,0)*{\includegraphics[scale=0.5]{figs/jj1.eps}};
 (0,-6)*{\bigb{h_s}};(0,2.5)*{e_a};
  \endxy
\end{equation}
since $\varepsilon_s$, $h_s$, are central in $\BNC_a$ under the isomorphism \eqref{eq_iso_center}.

One can check that the following identities are satisfied.
\begin{equation} \label{eq_thickbubslide}
    \xy
 (0,0)*{\includegraphics[scale=0.5]{figs/tsplit.eps}};
 (-5,-11)*{a+b};(-8,8)*{a};(8,8)*{b};
  (0,-6)*{\bigb{\varepsilon_s}};
  \endxy
  \quad =
 \quad
 \sum_{\ell=0}^{s} \;
    \xy
 (0,0)*{\includegraphics[scale=0.5]{figs/tsplit.eps}};
 (-5,-11)*{a+b};(-8,8)*{a};(8,8)*{b};
 (-5,2)*{\bigb{\varepsilon_{s-\ell}}};(5,2)*{\bigb{\varepsilon_{\ell}}};
  \endxy
\qquad \qquad
    \xy
 (0,0)*{\includegraphics[scale=0.5,angle=180]{figs/tsplitd.eps}};
 (-5,11)*{a+b};(-8,-8)*{a};(8,-8)*{b};
  (0,6)*{\bigb{\varepsilon_s}};
  \endxy
  \quad =
 \quad
 \sum_{\ell=0}^{s} \;
    \xy
 (0,0)*{\includegraphics[scale=0.5,angle=180]{figs/tsplitd.eps}};
 (-5,11)*{a+b};(-8,-8)*{a};(8,-8)*{b};
 (-5,-2)*{\bigb{\varepsilon_{s-\ell}}};
 (5,-2)*{\bigb{\varepsilon_{\ell}}};
  \endxy
\end{equation}
Equation \eqref{eq_eh_rel} becomes
\begin{equation}
 \sum_{r=0}^m (-1)^r\quad \xy
 (0,0)*{\includegraphics[scale=0.5]{figs/tlong-up.eps}};
 (-2.5,-12)*{a};(0,2)*{\bigb{\varepsilon_r}};(0,-5)*{\bigb{h_{m-r}}};
  \endxy
  \quad = \;\; 0.
\end{equation}

Sometimes it is convenient to split thick edges into thin edges.
Since $\delta_1=1$ and $e_1=1$ the definitions above imply
\begin{equation} 
  \xy
 (0,0)*{\includegraphics[scale=0.5]{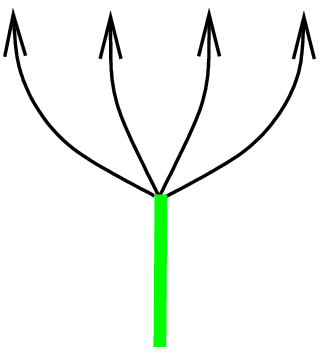}};
 (-3,-5)*{a};
  \endxy
  \quad = \quad
  \xy
 (0,0)*{\includegraphics[scale=0.5]{figs/c1-1.eps}};
 (0,0)*{D_a};
  \endxy
  \qquad \qquad
   \xy
 (0,0)*{\includegraphics[scale=0.5,angle=180]{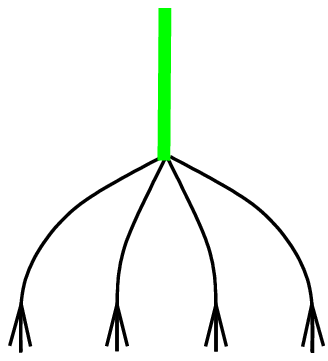}};
 (-3,5)*{a};
  \endxy
  \quad = \quad
  \xy
 (0,0)*{\includegraphics[scale=0.5]{figs/c1-1.eps}};
 (0,0)*{e_a};
  \endxy
\end{equation}
It follows from equations \eqref{eq_SchurDa}--\eqref{eq_divided2} above
that
\begin{equation} \label{eq_explode_rel}
  \xy
 (0,0)*{\includegraphics[scale=0.5]{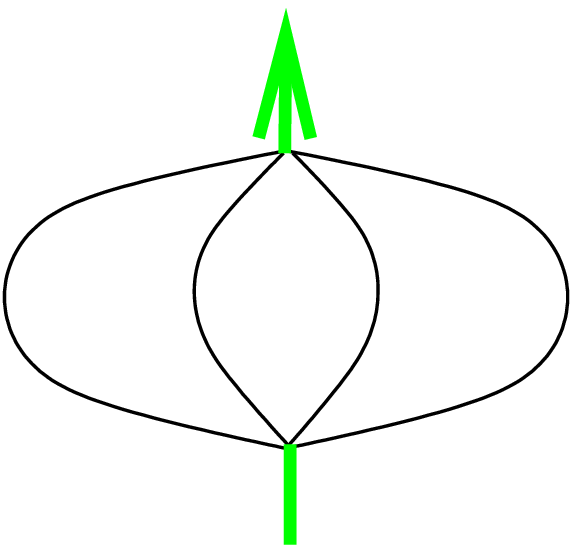}};
 (-3,-12)*{a};
  \endxy
  \quad = \quad
0 \quad \text{if $a>1$},
  \qquad \qquad
  \xy
 (0,0)*{\includegraphics[scale=0.5]{figs/texplode.eps}};
 (-3,-12)*{a};(-14,0)*{\bullet}+(-3.5,1)*{\scs a-1};
 (-4.5,0)*{\bullet}+(-3.5,1)*{\scs a-2};(4.5,0)*{\bullet}+(-2.5,1)*{\scs 1};
 (0,-2)*{\cdots};
  \endxy
  \quad = \quad
    \xy
 (0,0)*{\includegraphics[scale=0.5]{figs/tlong-up.eps}};
 (-3,-12)*{a};
  \endxy ,
\end{equation}
\begin{equation} \label{eq_split_bs}
\xy
 (0,0)*{\includegraphics[scale=0.5]{figs/texplode.eps}};
 (-14,0)*{\bullet}+(-2.5,1)*{\scs b_1};
 (-4.5,0)*{\bullet}+(-2.5,1)*{\scs b_r};
 (4.5,0)*{\bullet}+(2.5,1)*{\scs b_s};
 (14,0)*{\bullet}+(3,1)*{\scs b_{a}};
 (0,-2)*{\cdots}; (9,-2)*{\cdots};(-9,-2)*{\cdots}; (-3,-12)*{a};
  \endxy
 \;\; = \;\; - \; \;
\xy
 (0,0)*{\includegraphics[scale=0.5]{figs/texplode.eps}};
 (-14,0)*{\bullet}+(-2.5,1)*{\scs b_1};
 (-4.5,0)*{\bullet}+(-2.5,1)*{\scs b_s};
 (4.5,0)*{\bullet}+(2.5,1)*{\scs b_r};
 (14,0)*{\bullet}+(3,1)*{\scs b_{a}};
 (0,-2)*{\cdots}; (9,-2)*{\cdots};(-9,-2)*{\cdots}; (-3,-12)*{a};
  \endxy,
 \qquad \quad
\xy
 (0,0)*{\includegraphics[scale=0.5]{figs/texplode.eps}};
 (-14,0)*{\bullet}+(-2.5,1)*{\scs b_1};
 (-4.5,0)*{\bullet}+(-2.5,1)*{\scs b_r};
 (4.5,0)*{\bullet}+(2.5,1)*{\scs b_r};
 (14,0)*{\bullet}+(3,1)*{\scs b_{a}};
 (0,-2)*{\cdots}; (9,-2)*{\cdots};(-9,-2)*{\cdots}; (-3,-12)*{a};
  \endxy
  \;\; = \;\; 0.
\end{equation}

Equation \eqref{eq_SchurDa} implies that the Schur polynomial $\pi_{\alpha}(x_1,\dots,x_a)$, viewed as an element of $Z(\BNC_a)e_a$ via the isomorphism in equation \eqref{eq_iso_center}, has the following presentation
\begin{equation}\label{eq_schur_thin}
  \xy
 (0,0)*{\includegraphics[scale=0.5]{figs/tlong-up.eps}};
 (-2.5,-11)*{a};(0,-2)*{\bigb{\pi_{\alpha}}};
  \endxy
 \quad = \quad
 \xy
 (0,0)*{\includegraphics[scale=0.5]{figs/texplode.eps}};
 (-3,-12)*{a};(-14,0)*{\bullet}+(-4.5,1)*{\scs \alpha_1+}+(0,-3)*{\scs a-1};
 (-4.5,0)*{\bullet}+(-4.5,1)*{\scs \alpha_2+}+(0,-3)*{\scs a-2};
 (4.5,0)*{\bullet}+(4,1)*{\scs \alpha_{a-1}}+(0,-3)*{\scs + 1};
 (14,0)*{\bullet}+(3,1)*{\scs \alpha_{a}};
 (0,-2)*{\cdots};
  \endxy
 \quad =\quad
  \xy
 (0,0)*{\includegraphics[scale=0.5]{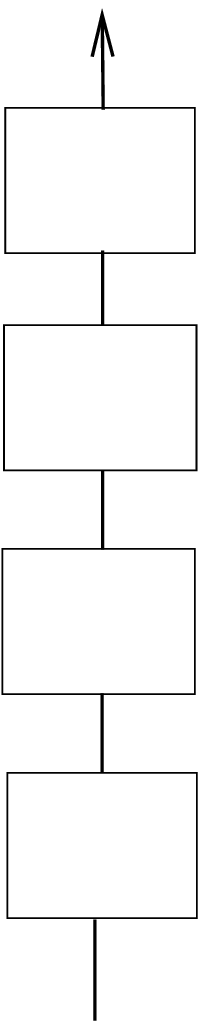}};
 (0,17)*{e_a};(0,6)*{x^{\alpha}};(0,-5.5)*{\delta_a};(0,-17)*{D_a};
  \endxy
\end{equation}
where $x^{\alpha}:=x_1^{\alpha_1}x_2^{\alpha_2}\dots x_a^{\alpha_a}$.
We call the diagrammatic operation given by the first equality above
{\em exploding} the Schur function.  More generally, to reduce the diagram
\begin{equation}
  \xy
 (0,0)*{\includegraphics[scale=0.5]{figs/texplode.eps}};
 (-14,0)*{\bullet}+(-2.5,1)*{\scs b_1};
 (-4.5,0)*{\bullet}+(-2.5,1)*{\scs b_r};
 (4.5,0)*{\bullet}+(2.5,1)*{\scs b_s};
 (14,0)*{\bullet}+(3,1)*{\scs b_{a}};
 (0,-2)*{\cdots}; (9,-2)*{\cdots};(-9,-2)*{\cdots}; (-3,-12)*{a};
  \endxy
\end{equation}
permute $b_i$'s to a non-increasing sequence $b'_1 \geq b'_2\geq
\dots \geq b'_a$.  If at least one of these inequalities is not
strict the diagram is zero by \eqref{eq_split_bs}.  Otherwise let
$b''_i=b_i'-a+i$.  The diagram is zero if $b''_a < 0$ and equal
to $\pm \pi_{(b''_1,\dots, b''_a)}e_a$ otherwise.

The following identities follow from (\ref{eq_thickbubslide}) and the properties of Schur functions in \eqref{eq_prod_schur_twovar}:
\begin{equation} \label{eq_schur_fork_slide}
    \xy
 (0,0)*{\includegraphics[scale=0.5]{figs/tsplit.eps}};
 (-5,-11)*{a+b};(-8,8)*{a};(8,8)*{b};
  (0,-6)*{\bigb{\pi_{\gamma}}};
  \endxy
  \quad =
 \quad
 \sum_{\alpha,\beta} \; c_{\alpha,\beta}^{\gamma} \;\;
    \xy
 (0,0)*{\includegraphics[scale=0.5]{figs/tsplit.eps}};
 (-5,-11)*{a+b};(-8,8)*{a};(8,8)*{b};
 (-5,2)*{\bigb{\pi_{\alpha}}};(5,2)*{\bigb{\pi_{\beta}}};
  \endxy
\qquad \qquad
    \xy
 (0,0)*{\includegraphics[scale=0.5,angle=180]{figs/tsplitd.eps}};
 (-5,11)*{a+b};(-8,-8)*{a};(8,-8)*{b};
  (0,6)*{\bigb{\pi_{\gamma}}};
  \endxy
  \quad =
 \quad
 \sum_{\alpha,\beta} \;c_{\alpha,\beta}^{\gamma}\;\;
    \xy
 (0,0)*{\includegraphics[scale=0.5,angle=180]{figs/tsplitd.eps}};
 (-5,11)*{a+b};(-8,-8)*{a};(8,-8)*{b};
 (-5,-2)*{\bigb{\pi_{\alpha}}};
 (5,-2)*{\bigb{\pi_{\beta}}};
  \endxy
\end{equation}
By \eqref{eq_prod_schur} we have the identity
\begin{equation} \label{eq_thickline_schur_mult}
 \xy
 (0,0)*{\includegraphics[scale=0.5]{figs/tlong-up.eps}};
 (-2.5,-12)*{a};(0,1)*{\bigb{\pi_{\alpha}}};(0,-6)*{\bigb{\pi_{\beta}}};
  \endxy
  \quad = \quad\sum_{\gamma \in P(a)} \; c_{\alpha,\beta}^{\gamma}\quad \xy
 (0,0)*{\includegraphics[scale=0.5]{figs/tlong-up.eps}};
 (-2.5,-12)*{a};(0,-2)*{\bigb{\pi_{\gamma}}};
  \endxy
\end{equation}

\begin{rem} \label{rem_diagram_conventions}
 For $\alpha \in P(a)$ the coefficients $c_{\beta,\gamma}^{\alpha}$ vanish whenever $\beta$ or $\gamma$ have more than $a$ parts.  That is, we must have $\beta,\gamma \in P(a)$ for $c_{\beta,\gamma}^{\alpha}\neq 0$.  This observation together with our conventions for representing Schur polynomials on thick lines imply that
\begin{equation}
c_{\beta,\gamma}^{\alpha} \;
   \xy
 (0,0)*{\includegraphics[scale=0.5]{figs/tlong-up.eps}};
 (-2.5,-11)*{b};(0,0)*{\bigb{\pi_{\beta}}};
 (9,9)*{n};
  \endxy
\end{equation}
is zero unless $\beta,\gamma \in P(a)$ and $\beta \in P(b)$.  We will often use these conventions to simplify summations when the summation is over all values of the indices for which a diagram is nonzero.  Thus we can write \eqref{eq_thickline_schur_mult} as
\begin{equation}
 \xy
 (0,0)*{\includegraphics[scale=0.5]{figs/tlong-up.eps}};
 (-2.5,-12)*{a};(0,1)*{\bigb{\pi_{\alpha}}};(0,-6)*{\bigb{\pi_{\beta}}};
  \endxy
  \quad = \quad\sum_{\gamma } \; c_{\alpha,\beta}^{\gamma}\quad \xy
 (0,0)*{\includegraphics[scale=0.5]{figs/tlong-up.eps}};
 (-2.5,-12)*{a};(0,-2)*{\bigb{\pi_{\gamma}}};
  \endxy
\end{equation}
with no confusion.
\end{rem}

\begin{prop} \label{prop_schur_splitter}
 \begin{equation}
   \xy
 (0,0)*{\includegraphics[scale=0.65]{figs/ufork-u.eps}};
 (0,-6)*{\bigb{\pi_{\alpha}}};
 (-3,-11)*{a};
  \endxy
  \quad = \quad
  \xy
 (0,0)*{\includegraphics[scale=0.65]{figs/ufork-u.eps}};
 (0,4)*{\bigb{\hspace{0.3in}\pi_{\alpha}\hspace{0.3in}}};
 (-3,-11)*{a};
  \endxy \qquad \qquad \qquad
    \xy
 (0,0)*{\includegraphics[scale=0.65,angle=180]{figs/ufork-d.eps}};
 (0,6)*{\bigb{\pi_{\alpha}}};
 (-3,11)*{a};
  \endxy
  \quad = \quad
  \xy
 (0,0)*{\includegraphics[scale=0.65,angle=180]{figs/ufork-d.eps}};
 (0,-4)*{\bigb{\hspace{0.3in}\pi_{\alpha}\hspace{0.3in}}};
 (-3,11)*{a};
 \endxy
 \end{equation}
\end{prop}

\begin{proof} The first equation is proven as follows:
\begin{equation}
    \xy
 (0,0)*{\includegraphics[scale=0.65]{figs/ufork-u.eps}};
 (0,-6)*{\bigb{\pi_{\alpha}}};
 (-3,-11)*{a};
  \endxy
  \quad = \quad  \xy
 (0,0)*{\includegraphics[scale=0.5]{figs/four-box.eps}};
 (0,17)*{D_a};(0,6)*{e_a};(0,-5.5)*{x^{\alpha}};(0,-17)*{e_a};
  \endxy
 \quad
 \refequal{\eqref{eq_partial-en}}
 \quad
     \xy
 (0,0)*{\includegraphics[scale=0.5]{figs/three-box.eps}};
 (0,11)*{D_a};(0,0)*{x^{\alpha}};(0,-11)*{e_a};
  \endxy
  \quad
 \refequal{\eqref{eq_DafDa}, \eqref{eq_SchurDa}}
 \quad
      \xy
 (0,0)*{\includegraphics[scale=0.5]{figs/two-box.eps}};
 (0,6)*{\pi_{\alpha}};(0,-5.5)*{D_a};
  \endxy
  \quad
  =:\quad
     \xy
 (0,0)*{\includegraphics[scale=0.65]{figs/ufork-u.eps}};
 (0,4)*{\bigb{\hspace{0.3in}\pi_{\alpha}\hspace{0.3in}}};
 (-3,-9)*{a};
  \endxy
\end{equation}
The second equation is proven similarly.
\end{proof}

 \subsection{Splitter equations, explosions, and idempotents } \label{sec_split_exp_idem}

Let $\alpha\in P(a)$ and $\beta\in P(b)$ be two partitions, and let
$w\in S_{a+b}$ be a permutation taking the sequence $(\alpha_1+a-1,
\dots, \alpha_{a-1}+1,\alpha_a, \beta_1+b-1, \dots, \beta_{b-1}+1,
\beta_b)$ to a nonincreasing sequence $\gamma'= (\gamma'_1, \dots,
\gamma'_{a+b})$, so that $\gamma' \in P(a+b)$.  The permutation $w$
is unique if and only if all parts of $\gamma'$ are distinct.  Let
$\gamma=(\gamma_1, \dots , \gamma_{a+b})$ with
$\gamma_k=\gamma'_k+k-a-b$.

\begin{prop} \label{prop_schur_leftright}
With the above notations
\begin{equation} \label{eq_schur_left}
  \xy
 (0,0)*{\includegraphics[scale=0.5]{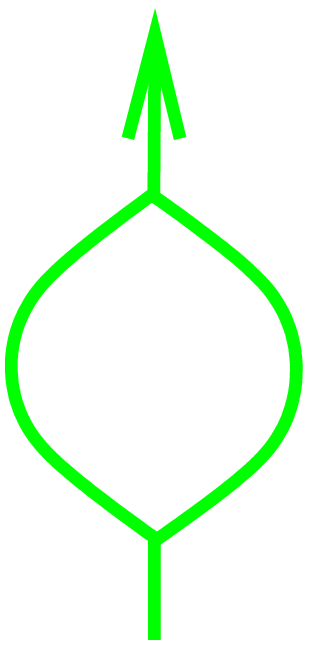}};
 (-6,5)*{a};  (6,5)*{b};(-6,-14)*{a+b};
 (-7,-2)*{\bigb{\pi_{\alpha}}};(7,-2)*{\bigb{\pi_{\beta}}};
  \endxy
  \qquad \quad = \quad \quad
      \left\{
\begin{array}{ccc}
 (-1)^{\ell(w)}   \left(\;\;\xy
 (0,0)*{\includegraphics[scale=0.5]{figs/single-tup.eps}};
 (-6 ,-8)*{a+b}; (0,-2)*{\bigb{\pi_{\gamma}}};
  \endxy\;\; \right) & \quad & \text{if $\gamma$ is a partition} \\ & & \\
  0 & \quad & \text{otherwise.}
\end{array}
  \right.
 \end{equation}
Note that $\gamma$ is a partition if and only if $\gamma'_1 > \gamma'_2> \dots >
\gamma'_{a+b}$.
\end{prop}

\begin{proof}
Using \eqref{eq_schur_thin} to explode the Schur polynomials
$\pi_{\alpha}$ and $\pi_{\beta}$  into thin lines we have
\begin{equation} \label{eq_proof_explode}
  \xy
 (0,0)*{\includegraphics[scale=0.5]{figs/tline-bubble.eps}};
 (-6,5)*{a};  (6,5)*{b};(-6,-14)*{a+b};
 (-7,-2)*{\bigb{\pi_{\alpha}}};(7,-2)*{\bigb{\pi_{\beta}}};
  \endxy \;\;\refequal{\eqref{eq_schur_thin}}\;\;
     \xy
 (0,0)*{\includegraphics[scale=0.5]{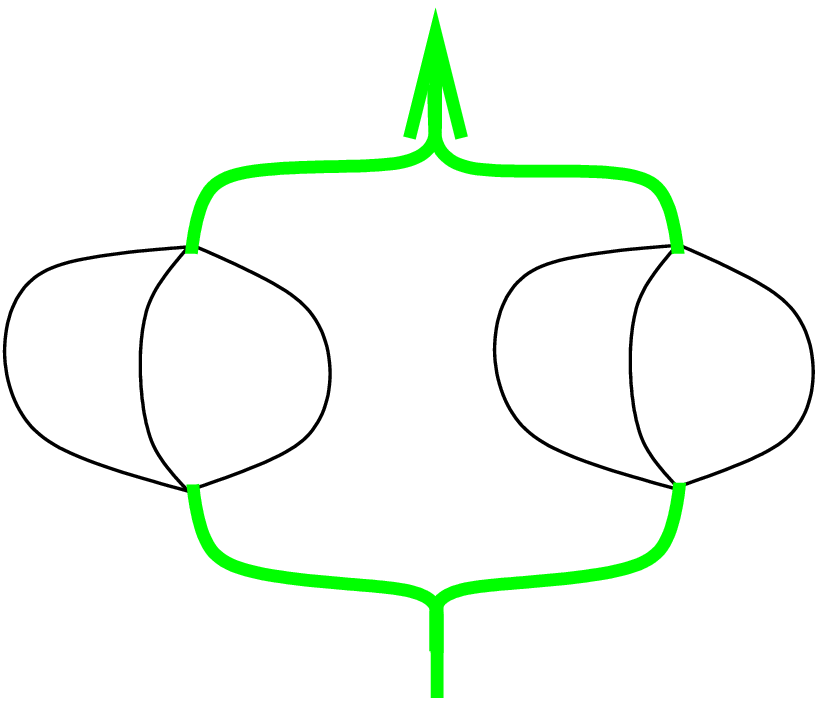}};
 (-7,-16)*{a+b};
 (-20.5,0)*{\bullet}+(-2.5,1)*{\scs s_1};
 (-14,0)*{\bullet}+(-2.5,1)*{\scs s_2};
 (-4,0)*{\bullet}+(-2.5,1)*{\scs s_a};
 (4,0)*{\bullet}+(-2.5,1)*{\scs r_1};
  (11,0)*{\bullet}+(-2.5,1)*{\scs r_2};
  (20.5,0)*{\bullet}+(-2.5,1)*{\scs r_b};(18,11)*{n};
 (15,-4)*{\cdots};(-9,-4)*{\cdots};
  \endxy
  \;\; = \;\;
       \xy
 (0,0)*{\includegraphics[scale=0.5]{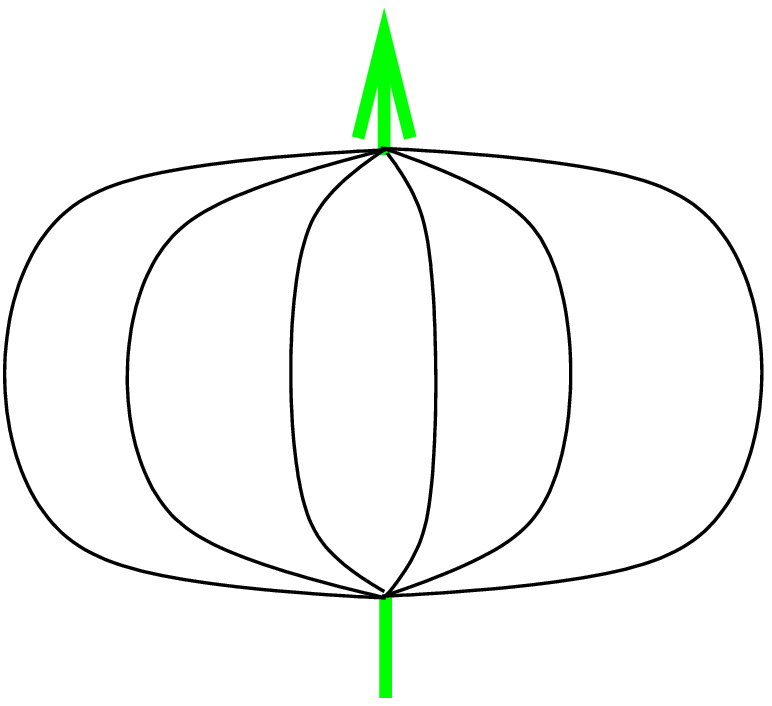}};
 (-7,-16)*{a+b};
 (-19,0)*{\bullet}+(-2.5,2)*{\scs s_1};
 (-13,0)*{\bullet}+(-2.5,2)*{\scs s_2};
 (-4.5,0)*{\bullet}+(-2.5,2)*{\scs s_a};
 (2.5,0)*{\bullet}+(-2.5,2)*{\scs r_1};
  (9.5,0)*{\bullet}+(-2.5,2)*{\scs r_2};
  (19,0)*{\bullet}+(-2.5,2)*{\scs r_b};(18,11)*{n};
 (15,-4)*{\cdots};(-9,-4)*{\cdots};
  \endxy
\end{equation}
where $s_j = a-j+\alpha_j$ and $r_j=b-j+b_j$.  Permuting the labels using \eqref{eq_split_bs} so that they are
nonincreasing from left to right we have
\begin{equation}
    \xy
 (0,0)*{\includegraphics[scale=0.5]{figs/tline-bubble.eps}};
 (-6,5)*{a};  (6,5)*{b};(-6,-14)*{a+b};
 (-7,-2)*{\bigb{\pi_{\alpha}}};(7,-2)*{\bigb{\pi_{\beta}}};
  \endxy  \quad = \quad    (-1)^{\ell(w)}  \;\; \xy
 (0,0)*{\includegraphics[scale=0.5]{figs/texplode2.eps}};
 (-7,-16)*{a+b};
 (-19,0)*{\bullet}+(-2.5,2)*{\scs \gamma'_1};
 (-13,0)*{\bullet}+(-2.5,2)*{\scs \gamma'_2};
 (-4.5,0)*{\bullet}+(-2.5,2)*{\scs \gamma'_a};
 (2.5,0)*{\bullet}+(3,2)*{\scs \gamma'_{a+1}};
  (9.5,0)*{\bullet}+(3,2)*{\scs \gamma'_{a+2}};
  (19,0)*{\bullet}+(3,2)*{\scs \gamma'_{a+b}};(18,11)*{n};
 (15,-4)*{\cdots};(-9,-4)*{\cdots};
  \endxy
\end{equation}
completing the proof.
\end{proof}

We have the following important special case of the above
Proposition.
\begin{cor} \label{cor_oval_small}
Let $\alpha\in P(a,b)$ and $\beta\in P(b,a)$ be two partitions. Then
\begin{equation} \label{eq_cor_oval_small}
  \xy
 (0,0)*{\includegraphics[scale=0.5]{figs/tline-bubble.eps}};
 (-6,5)*{a};  (6,5)*{b};(-6,-14)*{a+b};
 (-7,-2)*{\bigb{\pi_{\alpha}}};(7,-2)*{\bigb{\pi_{\beta}}};
  \endxy
\;\; = \;\;
    \left\{
\begin{array}{ccc}
 (-1)^{|\hat{\alpha}|}   \left(\;\;\xy
 (0,0)*{\includegraphics[scale=0.5]{figs/single-tup.eps}};
 (-6 ,-6)*{a+b};
  \endxy\;\; \right) & \quad & \text{if $\beta=\hat{\alpha}$} \\ & & \\
  0 & \quad & \text{otherwise.}
\end{array}
  \right.
 \end{equation}
\end{cor}

\begin{proof}
After exploding the Schur functions we obtain the diagram depicted on the right of \eqref{eq_proof_explode}, with $0 \leq s_i \leq a+b-1$ and $0 \leq r_j \leq a+b-1$ due to restrictions on $\alpha$ and $\beta$.  Since there are $a+b$ $s_i$'s and $r_j$'s, they must take each value between $0$ and $a+b-1$ exactly once.  Therefore, the $r_j$'s are determined by the $s_i$'s.  It is easy to check that $s_i$'s and $r_j$'s are all distinct exactly when $\beta=\hat{\alpha}$.  The minus sign in the formula follows from a straightforward computation.
\end{proof}

\begin{cor} \label{eq_nil_EaEb}For $\alpha$, $\beta$ as in Proposition~\ref{prop_schur_leftright}
\begin{equation}
    \xy
 (0,0)*{\includegraphics[scale=0.5]{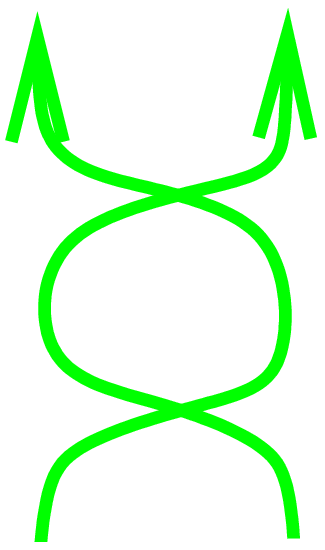}};
 (-9,-11)*{a};  (9,-11)*{b};
 (-7,-2)*{\bigb{\pi_{\beta}}};(7,-2)*{\bigb{\pi_{\alpha}}};
  \endxy
 \quad = \quad (-1)^{\ell(w)} \;\;
 \xy
 (0,0)*{\includegraphics[scale=0.5]{figs/tH.eps}};
(-8,-11)*{a};(8,-11)*{b}; (0,0)*{\bigb{\pi_{\gamma}}};
  \endxy
\end{equation}
Recall that our conventions are such that $\pi_{\gamma}=0$ if $\gamma$ is not a partition.
\end{cor}

\begin{proof}
This follows from Proposition~\ref{prop_schur_leftright} by
expanding the diagram on the left using the definition of the thick
crossing \eqref{eq_defn_thick_cross}.
\end{proof}

For every partition $\alpha\in P(a,b)$ define
\begin{eqnarray} \label{eq_def_sigmalambda_alpha}
  \sigma_{\alpha}:=
  \xy
 (0,0)*{\includegraphics[scale=0.5]{figs/tsplit.eps}};
 (-5,-11)*{a+b};(-8,8)*{a};(8,8)*{b};
 (-5,2)*{\bigb{\pi_{\alpha}}};
  \endxy, \qquad \quad
 \lambda_{\alpha}:= (-1)^{|\hat{\alpha}|}
     \xy
 (0,0)*{\includegraphics[scale=0.5,angle=180]{figs/tsplitd.eps}};
 (-5,11)*{a+b};(-8,-8)*{a};(8,-8)*{b};
 (5,-2)*{\bigb{\pi_{\hat{\alpha}}}};
  \endxy, \qquad \quad e_{\alpha}=\sigma_{\alpha}\lambda_{\alpha}.
\end{eqnarray}
In this section we view $\sigma_{\alpha}$, $\lambda_{\alpha}$, and $e_{\alpha}$ as elements of $\BNC_{a+b}$, $\deg(\sigma_{\alpha}) = -\deg(\lambda_{\alpha})=2|\alpha|-2ab$, and $\deg(e_{\alpha})=0$.
Corollary~\ref{cor_oval_small} says that
\begin{equation}
  \lambda_{\beta} \sigma_{\alpha} = \delta_{\alpha,\beta} e_{a+b}, \qquad \text{$\alpha,\beta \in P(a,b)$.}
\end{equation}
This implies
\begin{equation}
e_{\beta}e_{\alpha}=\delta_{\alpha,\beta}e_{\alpha}.
\end{equation}

Let
\begin{equation}
  e_{a,b} \quad = \quad \xy
 (4,0)*{\includegraphics[scale=0.5]{figs/tlong-up.eps}};
 (-4,0)*{\includegraphics[scale=0.5]{figs/tlong-up.eps}};
 (-7 ,-12)*{a};
 (7 ,-12)*{b};
  \endxy
\end{equation}
This is an idempotent in $\BNC_{a+b}$.

\begin{thm} \label{thm_nil-eaeb}
\begin{equation}
  e_{a,b} = \sum_{\alpha \in P(a,b)} e_{\alpha}
\end{equation}
or diagrammatically
\begin{equation}
  \xy
 (4,0)*{\includegraphics[scale=0.5]{figs/tlong-up.eps}};
 (-4,0)*{\includegraphics[scale=0.5]{figs/tlong-up.eps}};
 (-7 ,-12)*{a};
 (7 ,-12)*{b};
  \endxy
 \quad = \quad \sum_{\scriptscriptstyle
  \alpha\in P(a,b)
} (-1)^{|\hat{\alpha}|}  \xy
 (0,0)*{\includegraphics[scale=0.5]{figs/tH.eps}};
 (-8,-13)*{a};(8,-13)*{b};(-8,13)*{a};(8,13)*{b};
 (-5,7)*{\bigb{\pi_{\alpha} }};
 (4,-9)*{\bigb{\pi_{\hat{\alpha}} }};
  \endxy
\end{equation}
\end{thm}

First we prove the following Lemma.

\begin{lem} \label{lem_EaEone}
\begin{equation}
  \xy
 (4,0)*{\includegraphics[scale=0.5]{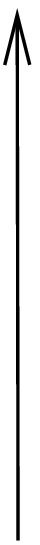}};
 (-4,0)*{\includegraphics[scale=0.5]{figs/tlong-up.eps}};
 (-7 ,-12)*{a};
 (7 ,-12)*{\scs 1};
  \endxy
 \quad = \quad \sum_{s=0}^a (-1)^{s}  \xy
 (0,0)*{\includegraphics[scale=0.4]{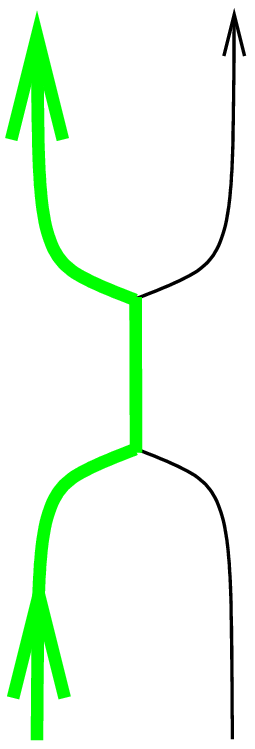}};
 (-7,-13)*{a};(8,-13)*{\scs 1};(-7,13)*{a};(7,13)*{\scs 1};
 (-5,7)*{\bigb{\varepsilon_{a-s} }};
 (4,-9)*{\bullet}+(3,1)*{s};
  \endxy
\end{equation}
\end{lem}

\begin{proof} Note that $\varepsilon_b=\pi_{1^b}$. W have
\begin{equation}\label{eq_lemEaEb}
    \xy
 (4,0)*{\includegraphics[scale=0.5]{figs/long-up.eps}};
 (-4,0)*{\includegraphics[scale=0.5]{figs/tlong-up.eps}};
 (-7 ,-12)*{a};
 (7 ,-12)*{\scs 1};
  \endxy
\;\; \refequal{\eqref{eq_cor_oval_small}} \;\;
  \xy
  (-14,0)*{\reflectbox{\includegraphics[scale=0.45]{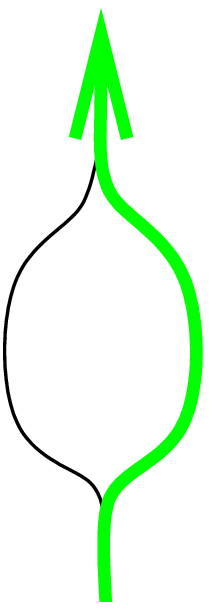}}};
 (-16.5,-12)*{a};(-19,-2)*{\bigb{\pi_{1^{a-1}}}}; (-22,0)*{};
 (-4,0)*{\includegraphics[scale=0.5]{figs/long-up.eps}}; (0 ,-12)*{\scs 1};
  \endxy
\;\; \refequal{\eqref{new_eq_iislide}} \;\;
     \xy
 (-5,0)*{\reflectbox{\includegraphics[scale=0.4]{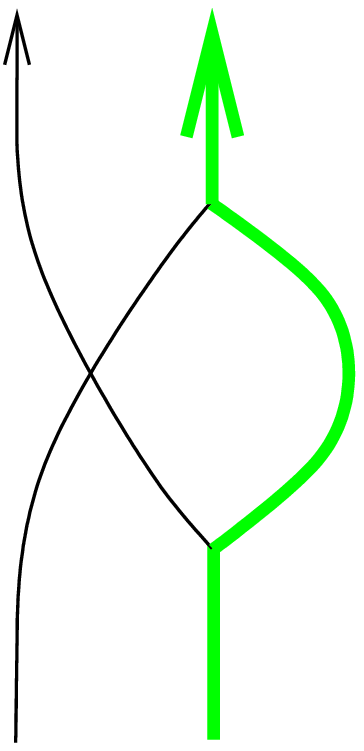}}};
 (-9,-13)*{a};(4,-13)*{\scs 1};
 (-12,0)*{\bigb{\pi_{1^{a-1}}}}; (-22,0)*{};
 (-4,4)*{\bullet};
  \endxy
  \quad - \quad
      \xy
 (-5,0)*{\reflectbox{\includegraphics[scale=0.4]{figs/triangle-2.eps}}};
 (-9,-13)*{a};(4,-13)*{\scs 1};
 (-12,0)*{\bigb{\pi_{1^{a-1}}}}; (-22,0)*{};
 (1,-5)*{\bullet};
  \endxy
\end{equation}
(Using \eqref{eq_schur_fork_slide} the Schur polynomials on the
thick line can be slid to the top)
\begin{equation}= \quad
\xy
 (-5,0)*{\reflectbox{\includegraphics[scale=0.4]{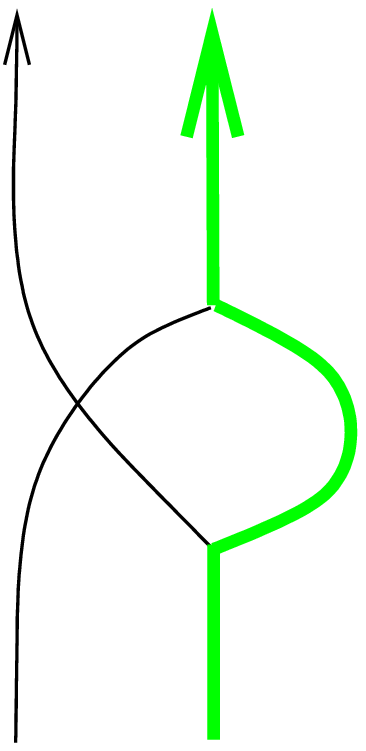}}};
 (-10,-13)*{a};(4,-13)*{\scs 1};
 (-6,6)*{\bigb{\pi_{1^a}}};
  \endxy
   \quad  - \quad  \sum_{r=0}^{a-1}(-1)^r    \xy
 (-5,0)*{\reflectbox{\includegraphics[scale=0.4]{figs/triangle-3.eps}}};
 (-10,-13)*{a};(4,-13)*{\scs 1};
 (-6,6)*{\bigb{\varepsilon_{a-1-r}}};
 (1,-5)*{\bullet};  (-3,1)*{\bullet}+(2,1)*{r};
  \endxy
\end{equation}
Sliding the $r$ dots down on the second term using
\eqref{new_eq_IND_iislide}, and using that
\begin{equation} \label{eq_split_right_thin}
  \xy
  (-14,0)*{\reflectbox{\includegraphics[scale=0.45]{figs/split-thinthick2.eps}}};
 (-16.5,-11)*{a};(-9.5,-2)*{\bullet}+(2,1)*{x};
  \endxy \quad = \quad 0, \qquad \quad \text{if $x< a-1$}
\end{equation}
 gives
\begin{equation}
\eqref{eq_lemEaEb} \quad = \quad\xy
 (-5,0)*{\reflectbox{\includegraphics[scale=0.4]{figs/triangle-3.eps}}};
 (-8,-13)*{a};(4,-13)*{\scs 1};
 (-6,6)*{\bigb{\varepsilon_{a}}};
  \endxy
   \quad  - \quad  \sum_{r=0}^{a-1}(-1)^r    \xy
 (-5,0)*{\reflectbox{\includegraphics[scale=0.4]{figs/triangle-3.eps}}};
 (-8,-13)*{a};(4,-13)*{\scs 1};
 (-6,6)*{\bigb{\varepsilon_{a-1-r}}};
 (1,-5)*{\bullet}+(5,1)*{r+1};
  \endxy
\end{equation}
completing the proof by Proposition~\ref{prop_almostRthree}.
\end{proof}

\begin{proof}[Proof of Theorem \ref{thm_nil-eaeb}]
By exploding the thick line of thickness $b$ and using Lemma
\ref{lem_EaEone} for the second equality we can write
\begin{equation}
  \xy
 (4,0)*{\includegraphics[scale=0.5]{figs/tlong-up.eps}};
 (-4,0)*{\includegraphics[scale=0.5]{figs/tlong-up.eps}};
 (-7 ,-12)*{a};
 (7 ,-12)*{b};
  \endxy \quad = \quad
  \xy
  (-21,0)*{\includegraphics[scale=0.5]{figs/tlong-up.eps}};
 (-24 ,-12)*{a};
 (0,0)*{\includegraphics[scale=0.5]{figs/texplode.eps}};
 (-14,0)*{\bullet}+(-2.5,2)*{\scs b-1};
 (-4.5,0)*{\bullet}+(-2.5,2)*{\scs b-2};
 (4.5,0)*{\bullet}+(2.5,1)*{\scs 1};
 (14,0)*{}+(3,1)*{\scs };
 (0,-2)*{\cdots};  (-3,-12)*{b};
  \endxy
 \quad = \quad
 \sum_{\beta_1=0}^a (-1)^{\beta_1}\;
   \xy
 (-24 ,-14)*{a}; (-2,-14)*{b};
 (-3.5,0)*{\includegraphics[scale=0.55]{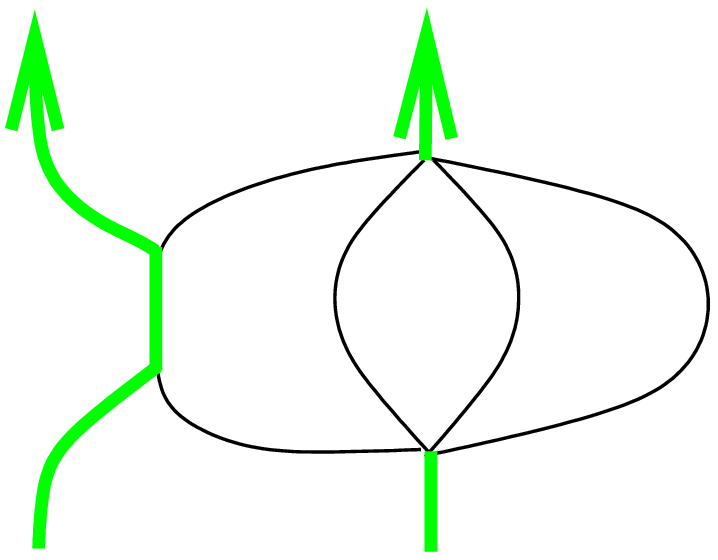}};
 (-8,-9.5)*{\bullet}+(0,2)*{\scs b-1+\beta_1};
 (-4.5,0)*{\bullet}+(-2.5,2)*{\scs b-2};
 (5.5,0)*{\bullet}+(2.5,1)*{\scs 1};
 (-21,5)*{\bigb{\varepsilon_{a-\beta_{1}} }};
 (14,0)*{}+(3,1)*{\scs };
 (0,-2)*{\cdots};
  \endxy
\end{equation}
Now repeatedly applying Lemma \ref{lem_EaEone} we get
\begin{equation}
\quad = \quad \sum_{\beta_1=0}^{a}\sum_{\beta_2=0}^{a+1} \dots
\sum_{\beta_b=0}^{a+b-1}(-1)^{\beta_1+\beta_2+\dots+\beta_b}
   \xy
 (-24 ,-29)*{a}; (8,-29)*{b};
 (-3.5,4)*{\includegraphics[scale=0.55]{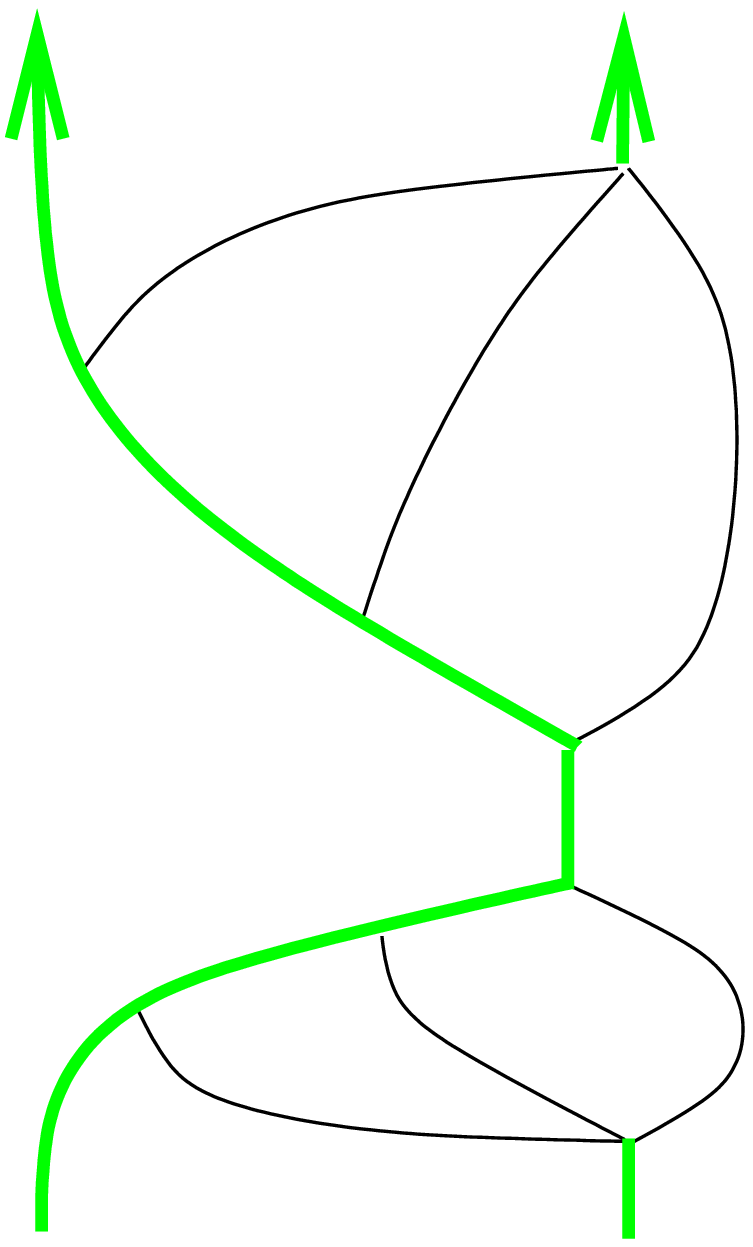}};
 (-9,-23.7)*{\bullet}+(-4,-2)*{\scs b-1+\beta_1};
 (-1.5,-18)*{\bullet}+(-6.5,0)*{\scs b-2+\beta_2};
 (17,-19)*{\bullet}+(2.5,1)*{\scs \beta_b};
 (-1,1)*{\bigb{\varepsilon_{a+b-1-\beta_b} }};
 (-15,9)*{\bigb{\varepsilon_{a+b-2-\beta_{b-1}} }};
 (-24,24)*{\bigb{\varepsilon_{a-\beta_{1}} }};
 (8,-17)*{\cdots}; (-13,17)*{\vdots};
  \endxy
 \end{equation}
 \begin{equation}\label{eq_bigugly}
 \quad = \quad \sum_{\beta_1=0}^{a}\sum_{\beta_2=0}^{a+1} \dots \sum_{\beta_b=0}^{a+b-1}(-1)^{\beta_1+\beta_2+\dots+\beta_b}
   \xy
 (-24 ,-29)*{a}; (8,-29)*{b};
 (0,4)*{\includegraphics[scale=0.55]{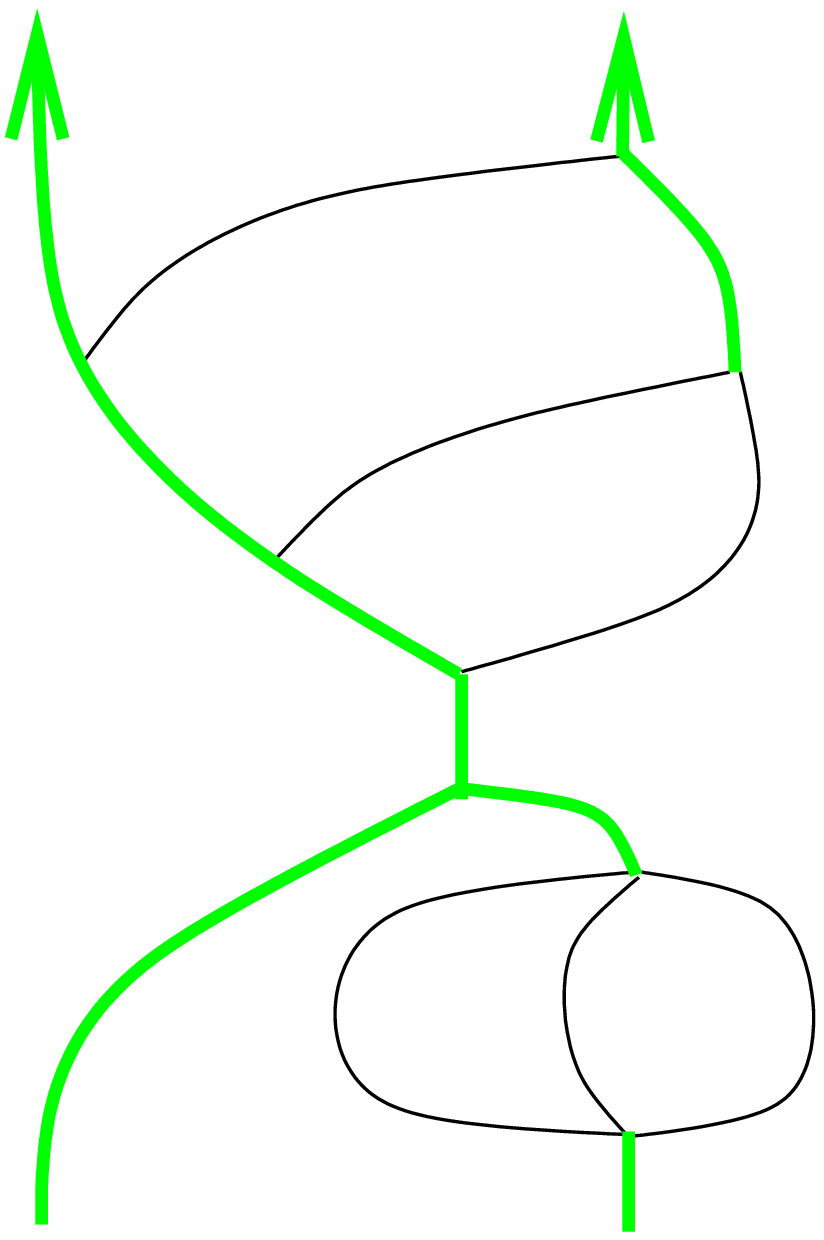}};
 (-4.2,-19)*{\bullet}+(-5.6,1.3)*{\scs b-1+\beta_1};
 (8.6,-19)*{\bullet}+(-5.6,1.3)*{\scs b-2+\beta_2};
 (22.5,-19)*{\bullet}+(2.5,1.5)*{\scs \beta_b};
 (-5,4)*{\bigb{\varepsilon_{a+b-1-\beta_b} }};
 (-18,11)*{\bigb{\varepsilon_{a+b-2-\beta_{b-1}} }};
 (-22,25)*{\bigb{\varepsilon_{a-\beta_{1}} }};
 (15,-18)*{\cdots}; (0,22)*{\vdots};
  \endxy
\end{equation}
On the upper half of the diagram we repeatedly make use of the
relation
\begin{equation}
   \xy
 (0,0)*{\includegraphics[scale=0.5]{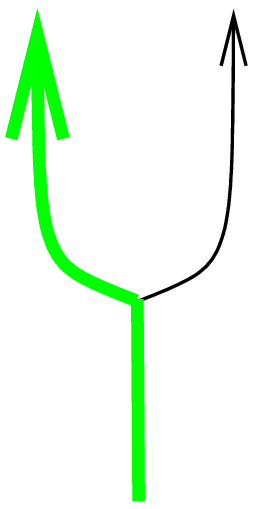}};
 (-3,-11)*{c};(-11,8)*{c-1};(8,8)*{};
  (0,-7)*{\bigb{\varepsilon_x}};
  \endxy
  \quad =
 \quad
    \xy
 (0,0)*{\includegraphics[scale=0.5]{figs/tonesplit.eps}};
 (-3,-11)*{c};(-11,8)*{c-1};(8,8)*{};
 (-5,2)*{\bigb{\varepsilon_{x}}};
  \endxy
  \quad + \quad
      \xy
 (0,0)*{\includegraphics[scale=0.5]{figs/tonesplit.eps}};
 (-3,-11)*{c};(-11,8)*{c-1};(8,8)*{};
 (-5,2)*{\bigb{\varepsilon_{x-1}}};(5,2)*{\bullet};
  \endxy
\end{equation}
with the appropriate thickness $c$ to slide all the elementary
symmetric functions to the top of the diagram.  In fact, only the
second term in the sum above where the dot is on the thin line will
be nonzero in \eqref{eq_bigugly} since all terms that are produced
sliding the elementary symmetric functions to the top have the form

\begin{equation}
        \xy
 (0,0)*{\includegraphics[scale=0.5]{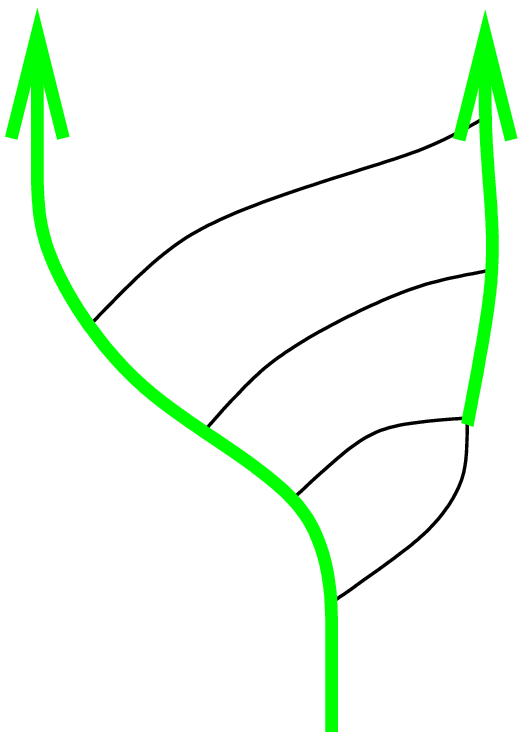}};
 (-2,-17)*{a+b};(-15,13)*{a};(15,13)*{b};
 (0,5)*{\vdots};
 (5,3)*{\bullet}+(0,2)*{\scs m_2};(5,-4)*{\bullet}+(0,2)*{\scs m_1};
 (5,10)*{\bullet}+(0,2)*{\scs m_{b-1}};
  \endxy
\end{equation}
where $0 \leq m_j \leq j$.  But using associativity of splitters we
have
\begin{equation}
        \xy
 (0,0)*{\includegraphics[scale=0.5]{figs/eaeb8.eps}};
 (-2,-17)*{a+b};(-15,13)*{a};(15,13)*{b};
 (0,5)*{\vdots};
 (5,2.8)*{\bullet}+(0,2)*{\scs m_2};(5,-4)*{\bullet}+(0,2)*{\scs m_1};
 (5,10)*{\bullet}+(0,2)*{\scs m_{b-1}};
  \endxy
  \quad = \quad
         \xy
 (5,0)*{\includegraphics[scale=0.5]{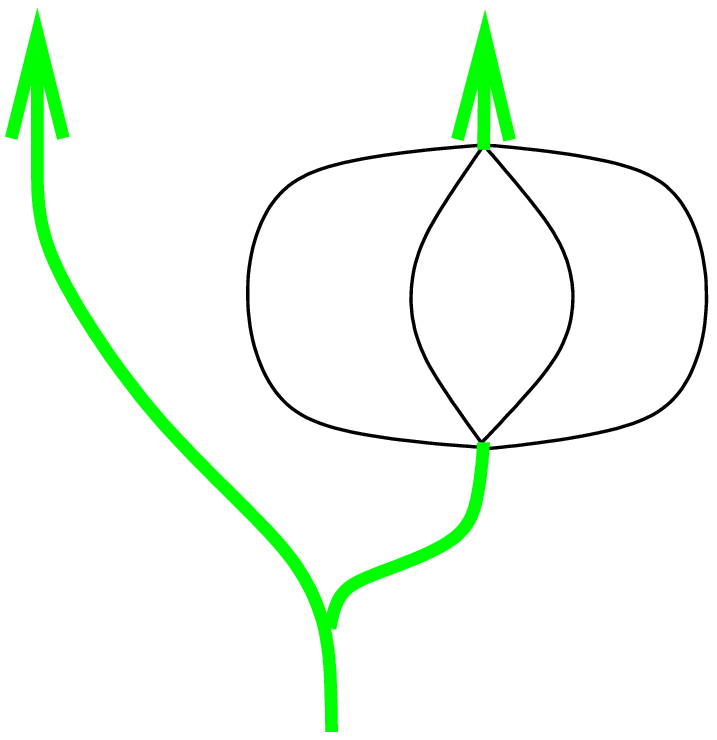}};
 (-2,-17)*{a+b};(-15,13)*{a};(15,13)*{b};
 (4,6)*{\dots};
 (-0.5,2)*{\bullet}+(-3.7,1)*{\scs m_{b-1}};
 (8,2)*{\bullet}+(-2.5,1)*{\scs m_2};
 (15.8,2)*{\bullet}+(-2.5,1)*{\scs m_1};
  \endxy
\end{equation}
so that the only term that is nonzero occurs when each $m_j=j$.

Adjusting the summations in \eqref{eq_bigugly} to include only the
nonzero terms we have
 \begin{equation}\label{eq_bigugly2}
 \quad = \quad \sum_{\beta_1=0}^{a}\sum_{\beta_2=0}^{a} \dots \sum_{\beta_b=0}^{a}(-1)^{\beta_1+\beta_2+\dots+\beta_b}
   \xy
 (-24 ,-29)*{a}; (8,-29)*{b};
 (0,4)*{\includegraphics[scale=0.55]{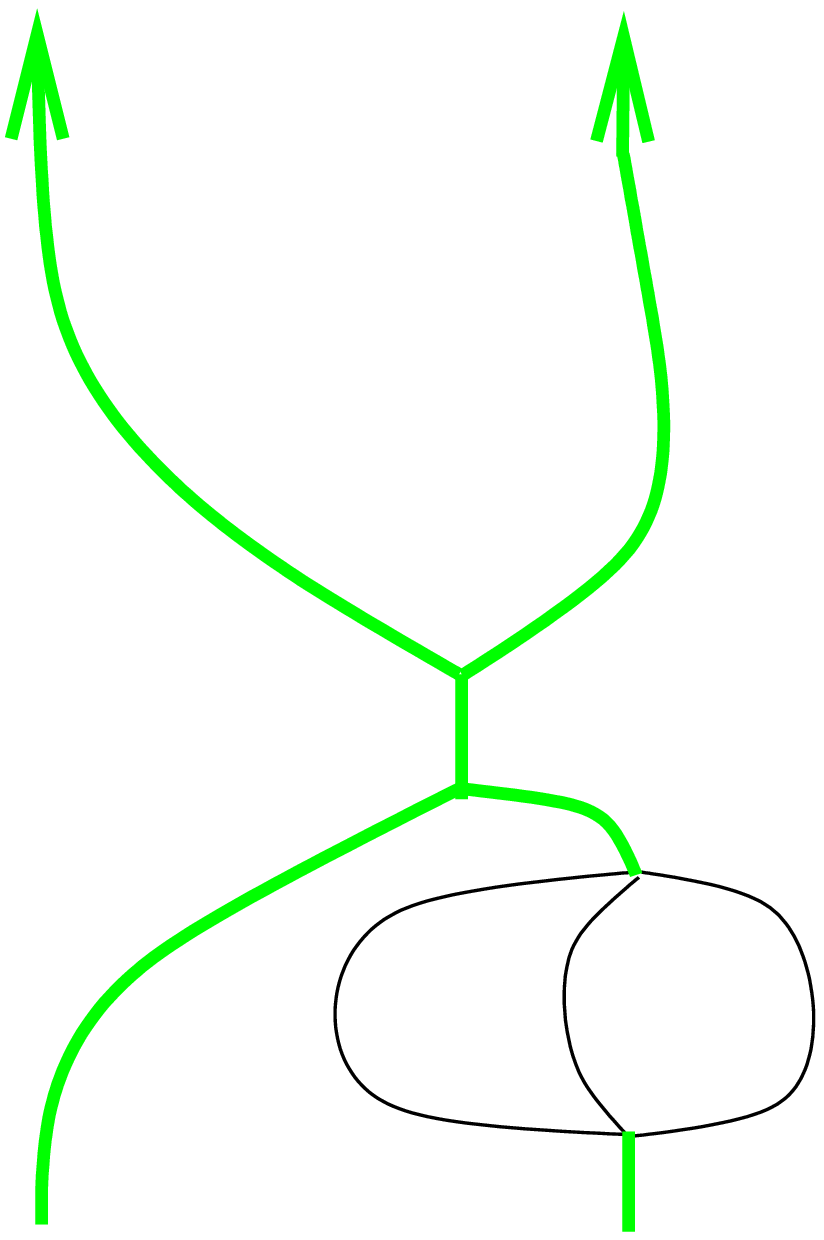}};
 (-4.5,-18)*{\bullet}+(-6.5,0)*{\scs b-1+\beta_1};
 (9,-18)*{\bullet}+(-6.5,0)*{\scs b-2+\beta_2};
 (22.5,-19)*{\bullet}+(2.5,1)*{\scs \beta_b};
 (-5,4)*{\bigb{\varepsilon_{a-\beta_b} }};
 (-18,11)*{\bigb{\varepsilon_{a-\beta_{b-1}} }};
 (-22,25)*{\bigb{\varepsilon_{a-\beta_{1}} }};
 (15,-18)*{\cdots}; (-22,19)*{\vdots};
  \endxy
\end{equation}
Using \eqref{eq_schur_thin} the bottom right subdiagram is clearly
an exploded Schur polynomial $\pi_{\beta}$ for some $\beta \in
P(b,a)$. Moreover, using our conventions that $\varepsilon_{a-\beta_i}$
is zero whenever $a-\beta_i<0$ or $a-\beta_i>a$ we can sum over all
$\beta_i \in \Z$.  Then the summation is over all $\gamma \in
P(b,a)$ where a given $\pi_{\gamma}$, with
$\gamma=(\gamma_1,\ldots,\gamma_b)$, appears $b!$ times in the sum
corresponding to the values of $\beta_i$ satisfying
\begin{equation}\label{for}
\beta_i+b-i=\gamma_{\sigma_i}+b-\sigma_i,\quad \text{for all
$i=1,\ldots,b$},
\end{equation}
where $\sigma \in S_b$ is a permutation, and we write $\sigma_i$ for
the image of $i$ under $\sigma$, that is $\sigma_i=\sigma(i)$. The
sign of $\pi_{\gamma}$ would then be $\sgn(\sigma)$, because of the
permutation of the labels on the dotted strands required to make
them non-increasing.

Then \eqref{eq_bigugly2} can be written as
\begin{equation}
  \xy
 (4,0)*{\includegraphics[scale=0.5]{figs/tlong-up.eps}};
 (-4,0)*{\includegraphics[scale=0.5]{figs/tlong-up.eps}};
 (-7 ,-12)*{a};
 (7 ,-12)*{b};
  \endxy
  \quad = \quad \sum_{\gamma \in P(b,a)}(-1)^{|\gamma|} \sum_{\sigma \in
  S_b} \sgn(\sigma)
    \xy
 (-22 ,-20)*{a}; (8,-20)*{b};(-23 ,28)*{a}; (7,28)*{b};
 (-4,6)*{\includegraphics[scale=0.55]{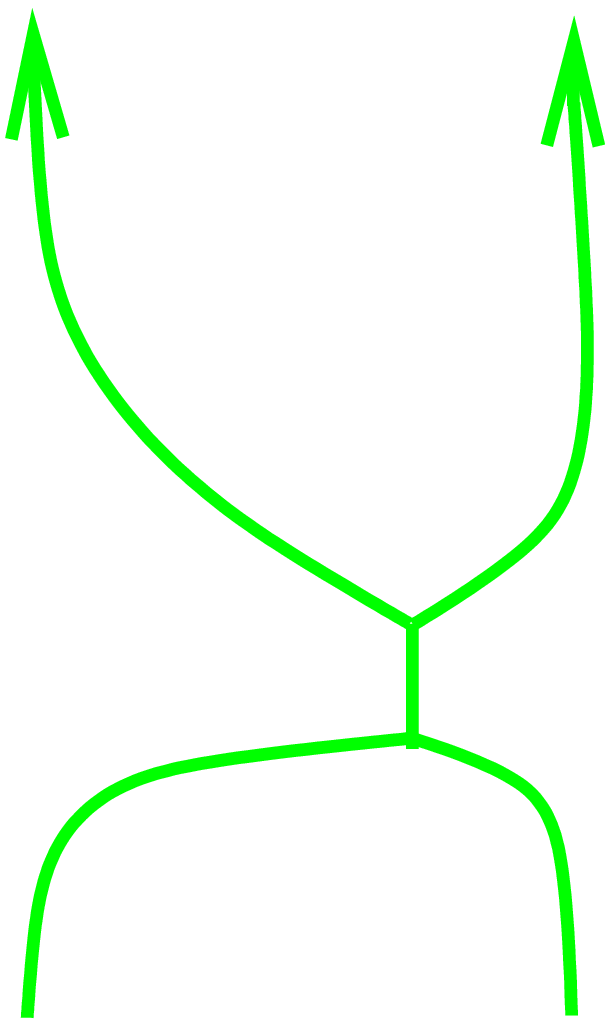}};
 (-8,4)*{\bigb{\varepsilon_{a-\gamma_{\sigma_1}+\sigma_1-1}}};
 (-14,11)*{\bigb{\varepsilon_{a-\gamma_{\sigma_2}+\sigma_2-2}} };
 (-18,23)*{\bigb{\varepsilon_{a-\gamma_{\sigma_b}+\sigma_b-b}} };
 (10,-14)*{\bigb{\pi_{\gamma} }};
 (-22,18)*{\vdots};
  \endxy
\end{equation}
Finally, we have
$$\sum_{\sigma\in S_b} {\sgn(\sigma)\varepsilon_{a-\gamma_{\sigma_1}+\sigma_1-1}
\varepsilon_{a-\gamma_{\sigma_2}+\sigma_2-2}\ldots\varepsilon_{a-\gamma_{\sigma_b}+\sigma_b-b}}
=\det{[\varepsilon_{a-\gamma_i+i-j}]}_{i,j=1}^b.$$ By interchanging the order of the rows and columns in the last $b\times b$ matrix (i.e.
by reordering $i\mapsto b+1-i$, $j\mapsto b+1-j$), and applying the Giambelli-Jacobi-Trudy
formula we see
$$\det{[\varepsilon_{a-\gamma_{b+1-i}+j-i}]}_{i,j=1}^b=\pi_{\overline{(a-\gamma_b,a-\gamma_{b-1},
\ldots,a-\gamma_{1})}}=\pi_{\hat{\gamma}},$$ as desired.
\end{proof}

\subsection{The nilHecke algebra as a matrix algebra over its center}

The set of sequences
\begin{equation}
  \Sq(a) := \{
  \und{\ell} = \ell_1 \dots \ell_{a-1} \mid 0 \leq \ell_{\nu} \leq \nu, \nu =  1,2, \dots a-1
  \}
\end{equation}
has size $|\Sq(a)|=a!$.  Let $|\und{\ell}|=\sum_{\nu} \ell_{\nu}$, $\hat{\ell_j}=j-\ell_j$ and
\begin{equation}
  \hat{\und{\ell}}=\hat{\ell}_1\dots\hat{\ell}_{n-1}= 1-\ell_1\; 2-\ell_2 \; \cdots \; a-1-\ell_{a-1}.
\end{equation}

Let $\varepsilon_r^{(a)}$ denote the $r$th elementary symmetric polynomial in $a$ variables.   The {\em standard elementary monomials} are given  by
\begin{equation}
  \varepsilon_{\und{\ell}} := \varepsilon_{\ell_1}^{(1)}\varepsilon_{\ell_2}^{(2)} \dots \varepsilon_{\ell_{a-1}}^{(a-1)}.
\end{equation}
In the graphical calculus for $\BNC_a$ we can write
\begin{equation}
  \varepsilon_{\und{\ell}} \quad = \quad
  \xy
 (-15,0)*{\includegraphics[scale=0.5]{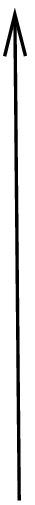}};
 (-10,0)*{\includegraphics[scale=0.5]{figs/single-up-long.eps}};
 (-5,0)*{\includegraphics[scale=0.5]{figs/single-up-long.eps}};
 (0,0)*{\includegraphics[scale=0.5]{figs/single-up-long.eps}};
 (5,0)*{\includegraphics[scale=0.5]{figs/single-up-long.eps}};
 (-15,8)*{\bullet}+(-2.5,0)*{\scs \ell_1};
 (-11.5,4)*{ \bigb{\;\varepsilon_{\ell_2} \;}};
  (-7.5,-8)*{ \bigb{\quad\; \varepsilon_{\ell_{a-1}} \;\quad}};
  (-7.5,-1)*{\vdots};
  \endxy
  \quad = \qquad
    \xy
 (-15,0)*{\includegraphics[scale=0.5]{figs/single-up-long.eps}};
 (-10,0)*{\includegraphics[scale=0.5]{figs/single-up-long.eps}};
 (-5,0)*{\includegraphics[scale=0.5]{figs/single-up-long.eps}};
 (0,0)*{\includegraphics[scale=0.5]{figs/single-up-long.eps}};
 (5,0)*{\includegraphics[scale=0.5]{figs/single-up-long.eps}};
  (-5,0)*{ \bigb{\quad\qquad \varepsilon_{\und{\ell}} \qquad\quad}};
  \endxy
\end{equation}
For further references related to standard elementary monomials see \cite{Win}.

Let
\begin{equation}
 \sigma_{\und{\ell}} \;\; = \;\;  \xy
 (0,0)*{\includegraphics[scale=0.7]{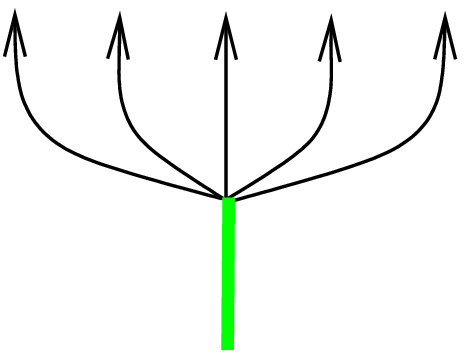}};
  (0,4)*{ \bigb{\qquad\qquad\; \varepsilon_{\und{\ell}} \qquad\qquad\;}};
 (-3,-10)*{a};
  \endxy
  \qquad  \quad \lambda_{\und{\ell}} \;\; = \;\; (-1)^{|\hat{\und{\ell}}|}
  \xy
 (0,0)*{\includegraphics[angle=180,scale=0.7]{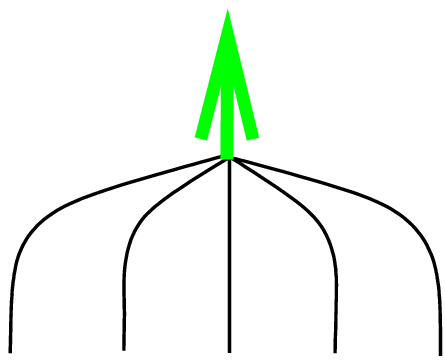}};
 (-7,-5)*{\bullet}+(-3,-1.5)*{\scs \hat{\ell}_{1}};;
 (.3,-5)*{\bullet}+(2.5,-1.5)*{\scs \hat{\ell}_{r}};;
 (7.5,-5)*{\bullet}+(4,-1.5)*{\scs \hat{\ell}_{a-2}};;
 (14.5,-5)*{\bullet}+(4.2,-1.5)*{\scs \hat{\ell}_{a-1}};
 (3,-3)*{\cdots}; (-3.3,-3)*{\cdots}; (-3,8)*{a};
  \endxy
\end{equation}

\begin{lem} \label{lem_nil_orth}
\begin{equation}
  \lambda_{\und{\ell'}}\sigma_{\und{\ell}} \;\; = \;\;
  \delta_{\underline{\ell},\underline{\ell'}}
 \xy
  (0,0)*{\includegraphics[scale=0.5]{figs/single-tup.eps}};
  (-2.5,-6)*{a};
 \endxy
 \;\; = \;\;
  \delta_{\und{\ell},\und{\ell'}} e_a
\end{equation}
for $\und{\ell}$, $\und{\ell'} \in \Sq(a)$.
\end{lem}

\begin{proof}
Using the (co)associativity of splitters, together with Proposition~\ref{prop_schur_splitter} for $\varepsilon_m(x_1,\ldots ,x_b) = \pi_{(1^m)}(x_1,\ldots ,x_b)$ shows that
\begin{equation}
 \lambda_{\und{\ell'}}\sigma_{\und{\ell}} \;\; = \;\; (-1)^{|\hat{\ell'}|}
 \vcenter{ \xy
   (0,0)*{\includegraphics[scale=0.7]{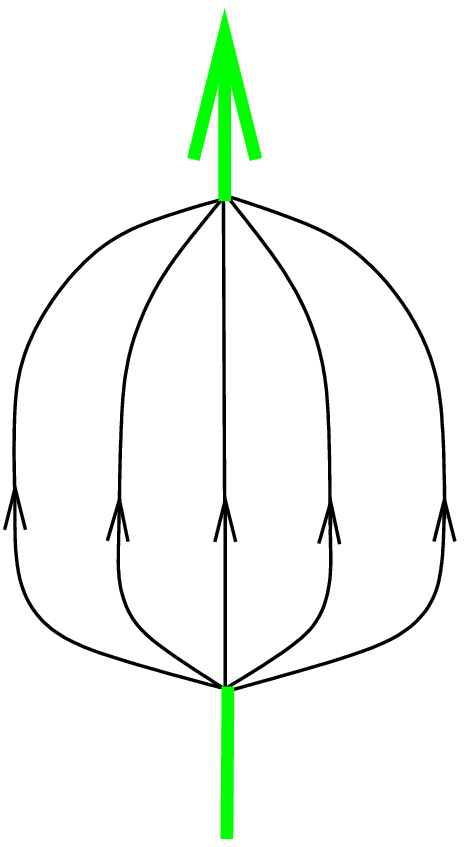}};
  (0,-13)*{ \bigb{\qquad\qquad\; \varepsilon_{\und{\ell}} \qquad\qquad\;}};
 (-3,-28)*{a};
 (-7.8,0)*{\bullet}+(-3,-1.5)*{\scs \hat{\ell'}_{1}};;
 (-.5,0)*{\bullet}+(2.5,-1.5)*{\scs \hat{\ell'}_{r}};;
 (7,0)*{\bullet}+(4,-1.5)*{\scs \hat{\ell'}_{a-2}};;
 (15,0)*{\bullet}+(4.3,-1.5)*{\scs \hat{\ell'}_{a-1}};
 (3,5)*{\cdots}; (-3.3,5)*{\cdots};
 (-5,20)*{a};
  \endxy}
  \;\; = \;\; (-1)^{|\hat{\und{\ell}'}|}
  \vcenter{\xy
   (2,0)*{\includegraphics[scale=0.5]{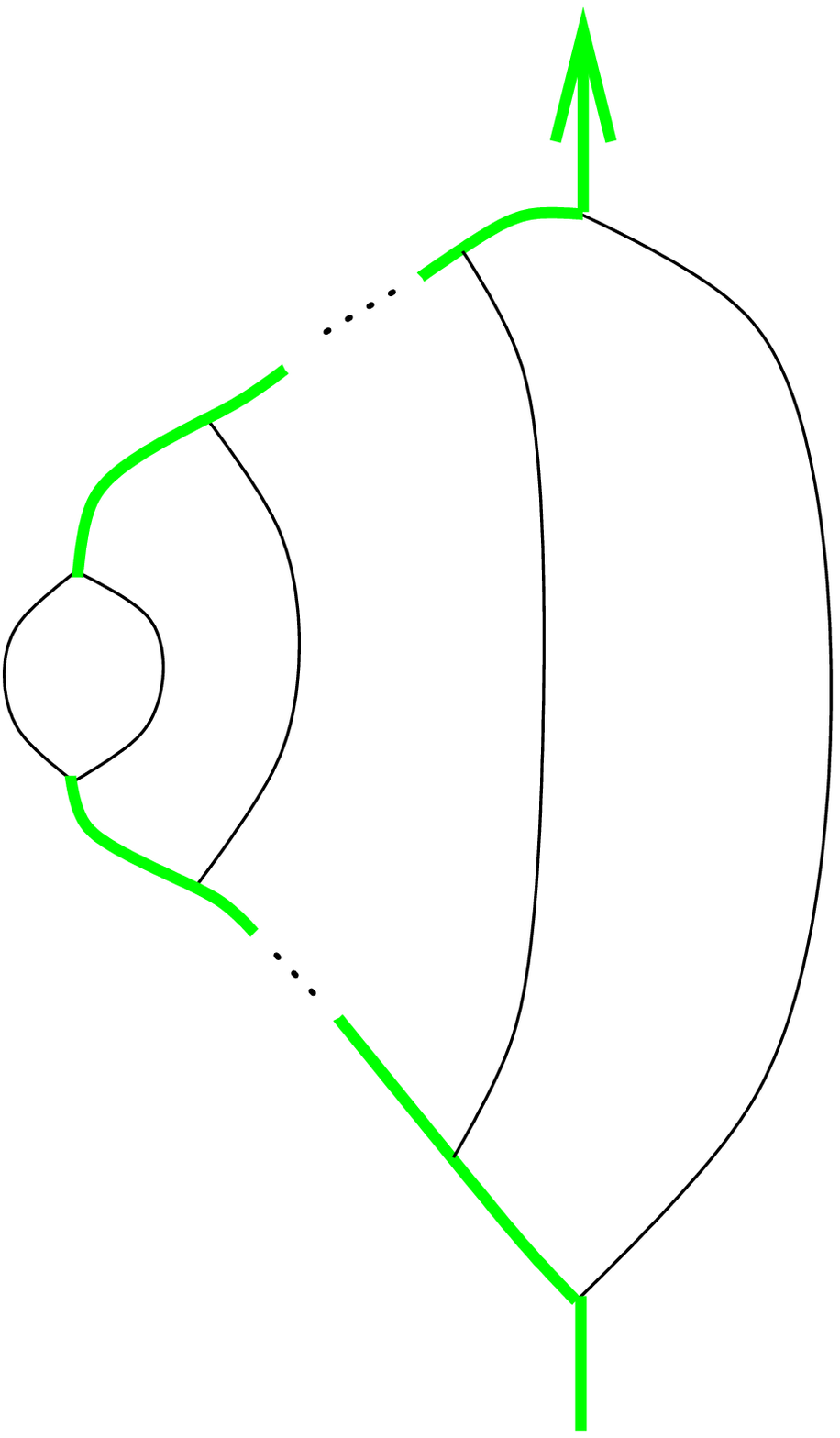}};
   (-16,-8)*{ \bigb{ \varepsilon_{\ell_{2}}}};
    (-1,-21)*{ \bigb{ \varepsilon_{\ell_{a-2}}}};
    (6,-29)*{ \bigb{ \varepsilon_{\ell_{a-1}}}};
    (-22,3)*{\bullet}+(-2.5,1.5)*{\scs \ell_{1}};
    (-12.5,3)*{\bullet}+(3.5,1.5)*{\scs 1-\ell_{1}'};
    (-5,3)*{\bullet}+(4,1.5)*{\scs 2-\ell_{2}'};
    (10,3)*{\bullet}+(7,1.5)*{\scs a-2-\ell_{a-2}'};
    (25.5,3)*{\bullet}+(7.5,1.5)*{\scs a-1-\ell_{a-1}'};
    (8,-39)*{a};
  \endxy}
\end{equation}
The lemma follows since
\begin{equation}
  \xy
  (0,0)*{\reflectbox{\includegraphics[scale=0.5]{figs/split-thinthick2.eps}}};
 (-5,-13)*{j+1};(5,-2)*{\bullet}+(6,1)*{j-\ell'};
  (-5,-2)*{ \bigb{ \varepsilon_{\ell}}};
  \endxy \quad = \quad (-1)^{j-\ell'}\delta_{\ell,\ell'}
  \xy
   (0,0)*{\includegraphics[scale=0.5]{figs/tlong-up.eps}};
    (-5,-12)*{j+1};
  \endxy
\end{equation}
for $0 \leq \ell, \ell' \leq j$ by Proposition~\ref{prop_schur_leftright}.
\end{proof}

\begin{rem}
There is a symmetric bilinear form
\begin{eqnarray}
 \langle , \rangle \maps \Z[x_1,\dots,x_a] \times \Z[x_1,\dots,x_a] &\to&
 \Z[x_1,\dots,x_a]^{S_a} \\
  (f,g) &\mapsto& \langle f,g\rangle := D_a(fg).
\end{eqnarray}
With respect to this bilinear form the Schubert polynomials \cite{Man} form an orthogonal basis for $\Z[x_1,\dots, x_a]$ as a module over $\Z[x_1,\dots, x_a]^{S_a}$.   The previous lemma shows that the dual basis to the standard elementary monomials is given by monomials $(-1)^{|\und{\hat{\ell}}|}x_2^{\hat{\ell}_1}\dots x_{a-1}^{\hat{\ell}_{a-2}}x_{a}^{\hat{\ell}_{a-1}}$.
\end{rem}

\begin{prop}
Let $e_{\und{\ell}} = \sigma_{\und{\ell}}\lambda_{\und{\ell}}$. The set $\{e_{\und{\ell}}\}_{\und{\ell} \in \Sq(a)}$ are
mutually-orthogonal idempotents that add up to $1\in \BNC_a$:
\begin{equation}
  e_{\und{\ell}}e_{\und{\ell'}} = \delta_{\und{\ell},\und{\ell'}} e_{\und{\ell}},
  \qquad \sum_{\und{\ell} \in \Sq(a)} e_{\und{\ell}} =1.
\end{equation}
\end{prop}

\begin{proof}
Lemma~\ref{lem_nil_orth} shows that $e_{\und{\ell}}$ are orthogonal idempotents.  We prove that the sum of the idempotents $e_{\und{\ell}}$ over $\und{\ell} \in \Sq(a)$ is $1 \in \BNC_a$ by induction on $a$.  The base case when $a=2$ is just the nilHecke relation.  For the induction step assume the result holds giving a decomposition of $1 \in \BNC_{a}$.  We show that it also holds for $1\in \BNC_{a+1}$.  Using the induction hypothesis we have
\begin{equation} \label{eq_thm_nil}
    \BNC_a \ni 1 \;\; = \;\; \vcenter{\xy
 (-6,-5)*{\includegraphics[scale=0.5]{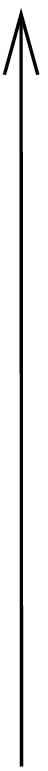}};
 (0,-5)*{\includegraphics[scale=0.5]{figs/single-up-long2.eps}};
 (12,-5)*{\includegraphics[scale=0.5]{figs/single-up-long2.eps}};
 (20,-5)*{\includegraphics[scale=0.5]{figs/single-up-long2.eps}};
 (6,0)*{\cdots};
  \endxy}
  \quad = \quad \sum_{\und{\ell} \in \Sq(a)} (-1)^{|\hat{\und{\ell}}|}
  \vcenter{\xy
  (-1.5,2)*{\xy
   (0,0)*{\includegraphics[scale=0.7]{figs/uexplode.eps}};
  (0,4)*{ \bigb{\qquad\qquad\; \varepsilon_{\und{\ell}} \qquad\qquad\;}};
  \endxy};
 (1,-12)*{ \xy
  (0,0)*{\includegraphics[angle=180,scale=0.7]{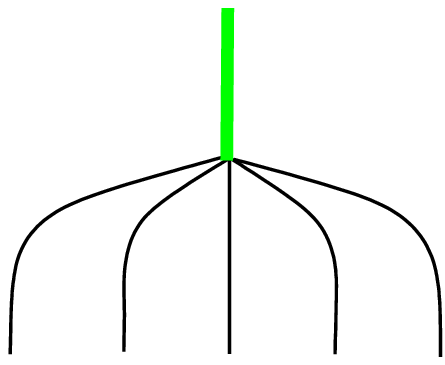}};
 (-7,-5)*{\bullet}+(-3,-1.5)*{\scs \hat{\ell}_{1}};;
 (.3,-5)*{\bullet}+(2.5,-1.5)*{\scs \hat{\ell}_{r}};;
 (7.5,-5)*{\bullet}+(4,-1.5)*{\scs \hat{\ell}_{a-2}};;
 (14.5,-5)*{\bullet}+(4.2,-1.5)*{\scs \hat{\ell}_{a-1}};
 (3,-3)*{\cdots}; (-3.3,-3)*{\cdots}; (-3,8)*{a}; \endxy};
 (24,-5)*{\includegraphics[scale=0.5]{figs/single-up-long2.eps}};
  \endxy}
\end{equation}
so that the Theorem follows by applying the Lemma~\ref{lem_EaEone} in the form
\begin{equation}
  \xy
 (4,0)*{\includegraphics[scale=0.5]{figs/long-up.eps}};
 (-4,0)*{\includegraphics[scale=0.5]{figs/tlong-up.eps}};
 (-7 ,-12)*{a};
 (7 ,-12)*{\scs 1};
  \endxy
 \quad = \quad \sum_{\ell_a=0}^a (-1)^{a-\ell_a}  \xy
 (0,0)*{\includegraphics[scale=0.4]{figs/tH-one.eps}};
 (-7,-13)*{a};(8,-13)*{\scs 1};(-7,13)*{a};(7,13)*{\scs 1};
 (-5,7)*{\bigb{\varepsilon_{\ell_a} }};
 (4,-9)*{\bullet}+(5,1)*{\scs a-\ell_a};
  \endxy
\end{equation}
to the middle of \eqref{eq_thm_nil}.
\end{proof}

Elements $\sigma_{\und{\ell}}$, $\lambda_{\und{\ell}}$ give an explicit
realization of $\BNC_a$ as the algebra of $a! \times a!$ matrices
over its center $Z(\BNC_a)= \Z[x_1,\dots,x_a]^{S_a}$.  Suppose that
rows and columns of $a!\times a!$ matrices are enumerated by
elements of $\Sq(a)$.  The isomorphism takes the matrix with $x \in
Z(\BNC_a)$ in the $(\und{\ell},\und{\ell'})$ entry and zeros elsewhere to
$\sigma_{\und{\ell}}x\lambda_{\und{\ell'}}$.

%
\section{Brief review of calculus for categorified sl(2)}
%

 \subsection{The algebra $\U$}

The quantum group $\Uq$ is the associative algebra (with unit) over
$\Q(q)$ with generators $E$, $F$, $K$, $K^{-1}$ and relations
\begin{eqnarray}
  KK^{-1}=&1&=K^{-1}K, \label{eq_UqI}\\
  KE &=& q^2EK, \\
  KF&=&q^{-2}FK, \\
  EF-FE&=&\frac{K-K^{-1}}{q-q^{-1}}. \label{eq_UqIV}
\end{eqnarray}
For simplicity the algebra $\Uq$ is written ${\bf U}$. Divided powers are defined by
\begin{equation}
  E^{(a)}:= \frac{E^a}{[a]!}, \qquad F^{(b)}:= \frac{F^b}{[b]!}.
\end{equation}

The BLM idempotented algebra $\U=\U(\mf{sl}_2)$ is the $\Q(q)$-algebra
obtained from ${\bf U}$ by substituting the unit element by a collection of orthogonal
idempotents $1_n$ for $n \in \Z$
\begin{equation} \label{eq_orthog_idemp}
  1_n1_m=\delta_{n,m}1_n,
\end{equation}
indexed by the elements of the weight lattice of $\mathfrak{sl}_2$, such that
\begin{equation}
K1_n =1_nK= q^n 1_n, \quad E1_n = 1_{n+2}E, \quad F1_n = 1_{n-2}F.
\end{equation}
The algebra $\UA$ is the $\Z[q,q^{-1}]$-subalgebra of $\U$ generated by products of divided powers $E^{(a)}1_{n}$ and $F^{(a)}1_{n}$.
There are direct sum decompositions of $\Q(q)$ and $\Z[q,q^{-1}]$-modules, respectively
\[
 \U = \bigoplus_{n,m \in \Z}1_m\U1_n ,\qquad \qquad \UA = \bigoplus_{n,m \in
 \Z}1_m(\UA)1_n.
\]
The $\Z[q,q^{-1}]$-submodule $1_m(\UA)1_n$ can be shown to be spanned by
$1_mE^{(a)}F^{(b)}1_n$ or $1_mF^{(b)}E^{(a)}1_n$ for $a,b \in \Z_+$ such that $m-n=2a-2b$. Throughout the paper we use a shorthand notation when writing products of generators.  For instance, $E^{(a)}E^{(b)}1_n$ means $E^{(a)}1_{n+2b}E^{(b)}1_n=1_{n+2a+2b}E^{(a)}E^{(b)}1_n$.

The algebra $\U$ does not have a unit since the infinite sum
$\sum_{n\in \Z}1_n$ is not an element in $\U$; however, the system
of idempotents $\{1_n | n \in \Z \}$ in a sense serves as a
substitute for a unit. A complete set of defining relations in $\UA$ is collected below:
\begin{eqnarray} \label{eq_AUrel1}
  1_n 1_m &=& \delta_{n,m} 1_n,
\end{eqnarray}
\begin{equation} \label{eq_AUrel2}
   E^{(a)}1_n = 1_{n+2a}E^{(a)}1_n =
  1_{n+2a}E^{(a)},
  \qquad
  F^{(b)}1_n = 1_{n-2b}F^{(b)}1_n =
  1_{n-2b}F^{(b)},
\end{equation}
\begin{eqnarray}
\label{eq_EaEb}
 E^{(a)}E^{(b)}1_n &=& \qbin{a+b}{a}E^{(a+b)}1_n, \\
 F^{(a)}F^{(b)}1_n &=& \qbin{a+b}{a}F^{(a+b)}1_n, \label{eq_EaEb2}
\\ \label{eq_FbEa1}
 E^{(a)}F^{(b)}1_{n}&=&
\sum_{j=0}^{\min(a,b)}\qbin{a-b+n}{j}F^{(b-j)}E^{(a-j)}1_{n}, \\
\label{eq_EaFb1}
F^{(b)}E^{(a)}1_n&=&
\sum_{j=0}^{\min(a,b)}\qbin{b-a-n}{j}E^{(a-j)}F^{(b-j)}1_n.
\end{eqnarray}
Relations \eqref{eq_AUrel1}--\eqref{eq_EaEb2}, \eqref{eq_FbEa1} for $n \geq b-a$, and \eqref{eq_EaFb1} for $n \leq b-a$ already constitute a set of defining relations in $\UA$.

Lusztig's canonical basis $\B$ of $\U$ consists of the elements
\begin{enumerate}[(i)]
     \item $E^{(a)}F^{(b)}1_{n} \quad $ for $a$,$b\in \Z_+$,
     $n\in\Z$, $n\leq b-a$,
     \item $F^{(b)}E^{(a)}1_{n} \quad$ for $a$,$b\in\Z_+$, $n\in\Z$,
     $n\geq
     b-a$,
\end{enumerate}
where $E^{(a)}F^{(b)}1_{b-a}=F^{(b)}E^{(a)}1_{b-a}$.  Taking $a=0$ or $b=0$ in the above two equations shows that the elements $E^{(a)}1_n$ and $F^{(b)}1_n$ are canonical for all $a,b \geq 0$ and $n \in \Z$.

Denote by ${}_m\B_n$ the set of elements of $\B$ that belong to $1_m\U1_n$.  The set ${}_m\B_n$ is a basis of $1_m\U1_n$ as a $\Q(q)$-vector space, and a basis of $1_m(\UA)1_n$ as a free $\Z[q,q^{-1}]$-module.  The importance of this
basis is that the structure constants are in $\Z_{+}[q,q^{-1}]$.  In particular, for
$x, y \in \B$
\[
 x y = \sum_{z \in \B}m_{x,y}^z z
\]
with $m_{x,y}^z \in \Z_{+}[q,q^{-1}]$.

\begin{lem} \label{lem_triple_canonical}
Let $a,b,c >0$.
\begin{enumerate}[a)]
  \item \label{part_triple1}In the product
\begin{equation}
  E^{(a)}F^{(b)}E^{(c)}1_n
\end{equation}
the subproducts $E^{(a)}F^{(b)}1_{n+2c}$ and $F^{(b)}E^{(c)}1_n$ are
never both canonical basis elements.
  \item \label{part_triple2} In the product
\begin{equation}
  F^{(a)}E^{(b)}F^{(c)}1_n
\end{equation}
the subproducts $F^{(a)}E^{(b)}1_{n-2c}$ and $E^{(b)}F^{(c)}1_n$ are
never both canonical basis elements.
\end{enumerate}
\end{lem}

\begin{proof}
We prove \eqref{part_triple1} leaving \eqref{part_triple2} as an exercise for the reader.  When $a,b,c>0$ the subproduct $E^{(a)}F^{(b)}1_{n+2c}$ is
canonical if and only if $n+2c \leq b-a$, and $F^{(b)}E^{(c)}1_n$ is canonical
if and only if $n \geq b-c$.  If both are canonical then $b-c \leq n \leq b-a-2c$, so that $0 \leq
-(a+c)$ which is impossible.
\end{proof}

 \subsection{The 2-category $\Ucat$}

For an introductory reference on 2-categories see \cite[Chapter
7]{Bor}.

Let $\Bbbk$ be a field.

\begin{defn}[\cite{Lau1}, see also \cite{KL3}] \label{def_Ucat}
The 2-category $\Ucat$ is the additive $\Bbbk$-linear 2-category consisting of
\begin{itemize}
  \item objects: $n$ for $n \in \Z$.
\end{itemize}
The hom $\Ucat(n,n')$ between two objects $n$, $n'$ is an additive $\Bbbk$-linear category:
\begin{itemize}
  \item objects of $\Ucat(n,n')$: for a sequence of signs $\ep = \epsilon_1\epsilon_2\dots\epsilon_m$, where $\epsilon_1, \dots, \epsilon_m \in \{ +,-\}$, define
 $$\cal{E}_{\ep} := \cal{E}_{\epsilon_1} \cal{E}_{\epsilon_2}\dots \cal{E}_{\epsilon_m}$$
where $\cal{E}_{+}:= \cal{E}$ and $\cal{E}_{-}:= \cal{F}$.  An object of $\Ucat(n,n')$, called a 1-morphism
in $\Ucat$, is a formal finite direct sum of 1-morphisms
  \[
 \cal{E}_{\ep} \onen\{t\} =\onenp \cal{E}_{\ep} \onen\{t\}
  \]
for $t\in \Z$ and sequences $\ep$ such that $n'=n+2\sum_{j=1}^m
\epsilon_j1$.

  \item morphisms of $\Ucat(n,n')$: for 1-morphisms $\cal{E}_{\ep} \onen\{t\}$ and $\cal{E}_{\ep'} \onen\{t'\}$ in $\Ucat$, hom
sets $\Ucat(\cal{E}_{\ep} \onen\{t\},\cal{E}_{\ep'} \onen\{t'\})$ in
$\Ucat(n,n')$ are $\Bbbk$-vector spaces given by linear
combinations of degree $t-t'$ diagrams, modulo certain relations,
built from the following generating 2-morphisms:
\begin{enumerate}[i)]
  \item  Degree zero identity 2-morphism $\Id_x$ for each 1-morphism $x$ in
$\Ucat$. Identity 2-morphisms $\Id_{\cal{E} \onen}\{t\}$ and
$\Id_{\cal{F} \onen}\{t\}$  are represented graphically by
\[
\begin{array}{ccc}
  \Id_{\cal{E} \onen\{t\}} &\quad  & \Id_{\cal{F} \onen\{t\}} \\ \\
    \xy
 (0,8);(0,-8); **\dir{-} ?(.5)*\dir{>}+(2.3,0)*{\scriptstyle{}};
 (0,11)*{};
 (6,2)*{ n};
 (-8,2)*{ n +2};
 (-10,0)*{};(10,0)*{};
 \endxy
 & &
 \;\;   \xy
 (0,8);(0,-8); **\dir{-} ?(.5)*\dir{<}+(2.3,0)*{\scriptstyle{}};
 (0,-11)*{};(0,11)*{};
 (6,2)*{ n};
 (-8,2)*{ n -2};
 (-12,0)*{};(12,0)*{};
 \endxy
\end{array}
\]
and, more generally,  the identity $\Id_{\cal{E}_{\ep} \onen\{t\}}$
2-morphism for a sequence $\ep = \epsilon_1\epsilon_2
\dots \epsilon_m$, is depicted as
\begin{equation*}
\begin{array}{ccc}
  \xy
 (-12,8);(-12,-8); **\dir{-};
 (-4,8);(-4,-8); **\dir{-};
 (4,0)*{\cdots};
 (12,8);(12,-8); **\dir{-};
  (-12,-11)*{\cal{E}_{\epsilon_1}}; (-2,-11)*{\cal{E}_{\epsilon_2}};(14,-11)*{\cal{E}_{\epsilon_m}};
 (18,2)*{ n}; (-20,2)*{ n'};
 \endxy
\end{array}
\end{equation*}
The strand labelled $\cal{E}_{\epsilon_{\alpha}}$ is oriented
up if $\epsilon_{\alpha}=+$ and down if $\epsilon_{\alpha}=-$.

  \item For each $n \in \Z$ the 2-morphisms
\[
\begin{tabular}{|l|c|c|c|c|}
 \hline
 {\bf 2-morphism:} &   \xy
 (0,7);(0,-7); **\dir{-} ?(.75)*\dir{>};
 (0,-2)*{\txt\large{$\bullet$}};
 (6,4)*{n}; (-8,4)*{n +2}; (-10,0)*{};(10,0)*{};
 \endxy
 &
     \xy
 (0,7);(0,-7); **\dir{-} ?(.75)*\dir{<};
 (0,-2)*{\txt\large{$\bullet$}};
 (-6,4)*{n}; (8,4)*{n+2}; (-10,0)*{};(10,9)*{};
 \endxy
 &
   \xy
  (0,0)*{\xybox{
    (-4,-4)*{};(4,4)*{} **\crv{(-4,-1) & (4,1)}?(1)*\dir{>} ;
    (4,-4)*{};(-4,4)*{} **\crv{(4,-1) & (-4,1)}?(1)*\dir{>};
     (8,1)*{n};     (-12,0)*{};(12,0)*{};     }};
  \endxy
 &
   \xy
  (0,0)*{\xybox{
    (-4,4)*{};(4,-4)*{} **\crv{(-4,1) & (4,-1)}?(1)*\dir{>} ;
    (4,4)*{};(-4,-4)*{} **\crv{(4,1) & (-4,-1)}?(1)*\dir{>};
     (8,1)*{ n};     (-12,0)*{};(12,0)*{};     }};
  \endxy
\\ & & & &\\
\hline
 {\bf Degree:} & \;\;\text{  2 }\;\;
 &\;\;\text{  2}\;\;& \;\;\text{ -2}\;\;
 & \;\;\text{  -2}\;\; \\
 \hline
\end{tabular}
\]
\[
\begin{tabular}{|l|c|c|c|c|}
 \hline
  {\bf 2-morphism:} &   \xy
    (0,-3)*{\bbpef{}};
    (8,-5)*{n};    (-12,4)*{};(12,0)*{};
    \endxy
  &\xy
    (0,-3)*{\bbpfe{}};
    (8,-5)*{n};    (-12,4)*{};(12,0)*{};
    \endxy
  & \xy
    (0,-5)*{\bbcef{}};
    (8,-3)*{n};     (-12,4)*{};(12,0)*{};
    \endxy
  & \xy
    (0,-5)*{\bbcfe{}};
    (8,-3)*{n};    (-12,4)*{};(12,0)*{};
    \endxy\\& & &  &\\ \hline
 {\bf Degree:} & \;\;\text{  $1+n$}\;\;
 & \;\;\text{ $1-n$}\;\;
 & \;\;\text{ $1+n$}\;\;
 & \;\;\text{  $1-n$}\;\;
 \\
 \hline
\end{tabular}
\]
\end{enumerate}
such that the following identities hold:

\item  cups and caps are biadjointness morphisms up to grading shifts:
\begin{equation} \label{eq_biadjoint1}
  \xy   0;/r.18pc/:
    (-8,0)*{}="1";
    (0,0)*{}="2";
    (8,0)*{}="3";
    (-8,-10);"1" **\dir{-};
    "1";"2" **\crv{(-8,8) & (0,8)} ?(0)*\dir{>} ?(1)*\dir{>};
    "2";"3" **\crv{(0,-8) & (8,-8)}?(1)*\dir{>};
    "3"; (8,10) **\dir{-};
    (12,-9)*{n};
    (-6,9)*{n+2};
    \endxy
    \; =
    \;
\xy   0;/r.18pc/:
    (-8,0)*{}="1";
    (0,0)*{}="2";
    (8,0)*{}="3";
    (0,-10);(0,10)**\dir{-} ?(.5)*\dir{>};
    (5,8)*{n};
    (-9,8)*{n+2};
    \endxy
\qquad \quad  \xy   0;/r.18pc/:
    (-8,0)*{}="1";
    (0,0)*{}="2";
    (8,0)*{}="3";
    (-8,-10);"1" **\dir{-};
    "1";"2" **\crv{(-8,8) & (0,8)} ?(0)*\dir{<} ?(1)*\dir{<};
    "2";"3" **\crv{(0,-8) & (8,-8)}?(1)*\dir{<};
    "3"; (8,10) **\dir{-};
    (12,-9)*{n+2};
    (-6,9)*{ n};
    \endxy
    \; =
    \;
\xy   0;/r.18pc/:
    (-8,0)*{}="1";
    (0,0)*{}="2";
    (8,0)*{}="3";
    (0,-10);(0,10)**\dir{-} ?(.5)*\dir{<};
   (9,8)*{n+2};
    (-6,8)*{ n};
    \endxy
\end{equation}

\begin{equation}
 \xy   0;/r.18pc/:
    (8,0)*{}="1";
    (0,0)*{}="2";
    (-8,0)*{}="3";
    (8,-10);"1" **\dir{-};
    "1";"2" **\crv{(8,8) & (0,8)} ?(0)*\dir{>} ?(1)*\dir{>};
    "2";"3" **\crv{(0,-8) & (-8,-8)}?(1)*\dir{>};
    "3"; (-8,10) **\dir{-};
    (12,9)*{n};
    (-5,-9)*{n+2};
    \endxy
    \; =
    \;
      \xy 0;/r.18pc/:
    (8,0)*{}="1";
    (0,0)*{}="2";
    (-8,0)*{}="3";
    (0,-10);(0,10)**\dir{-} ?(.5)*\dir{>};
    (5,-8)*{n};
    (-9,-8)*{n+2};
    \endxy
\qquad \quad \xy  0;/r.18pc/:
    (8,0)*{}="1";
    (0,0)*{}="2";
    (-8,0)*{}="3";
    (8,-10);"1" **\dir{-};
    "1";"2" **\crv{(8,8) & (0,8)} ?(0)*\dir{<} ?(1)*\dir{<};
    "2";"3" **\crv{(0,-8) & (-8,-8)}?(1)*\dir{<};
    "3"; (-8,10) **\dir{-};
    (12,9)*{n+2};
    (-6,-9)*{ n};
    \endxy
    \; =
    \;
\xy  0;/r.18pc/:
    (8,0)*{}="1";
    (0,0)*{}="2";
    (-8,0)*{}="3";
    (0,-10);(0,10)**\dir{-} ?(.5)*\dir{<};
    (9,-8)*{n+2};
    (-6,-8)*{ n};
    \endxy
\end{equation}

\item NilHecke relations hold:
 \begin{equation}
  \vcenter{\xy 0;/r.18pc/:
    (-4,-4)*{};(4,4)*{} **\crv{(-4,-1) & (4,1)}?(1)*\dir{>};
    (4,-4)*{};(-4,4)*{} **\crv{(4,-1) & (-4,1)}?(1)*\dir{>};
    (-4,4)*{};(4,12)*{} **\crv{(-4,7) & (4,9)}?(1)*\dir{>};
    (4,4)*{};(-4,12)*{} **\crv{(4,7) & (-4,9)}?(1)*\dir{>};
 \endxy}
 =0, \qquad \quad
 \vcenter{
 \xy 0;/r.18pc/:
    (-4,-4)*{};(4,4)*{} **\crv{(-4,-1) & (4,1)}?(1)*\dir{>};
    (4,-4)*{};(-4,4)*{} **\crv{(4,-1) & (-4,1)}?(1)*\dir{>};
    (4,4)*{};(12,12)*{} **\crv{(4,7) & (12,9)}?(1)*\dir{>};
    (12,4)*{};(4,12)*{} **\crv{(12,7) & (4,9)}?(1)*\dir{>};
    (-4,12)*{};(4,20)*{} **\crv{(-4,15) & (4,17)}?(1)*\dir{>};
    (4,12)*{};(-4,20)*{} **\crv{(4,15) & (-4,17)}?(1)*\dir{>};
    (-4,4)*{}; (-4,12) **\dir{-};
    (12,-4)*{}; (12,4) **\dir{-};
    (12,12)*{}; (12,20) **\dir{-};
  (18,8)*{n};
\endxy}
 \;\; =\;\;
 \vcenter{
 \xy 0;/r.18pc/:
    (4,-4)*{};(-4,4)*{} **\crv{(4,-1) & (-4,1)}?(1)*\dir{>};
    (-4,-4)*{};(4,4)*{} **\crv{(-4,-1) & (4,1)}?(1)*\dir{>};
    (-4,4)*{};(-12,12)*{} **\crv{(-4,7) & (-12,9)}?(1)*\dir{>};
    (-12,4)*{};(-4,12)*{} **\crv{(-12,7) & (-4,9)}?(1)*\dir{>};
    (4,12)*{};(-4,20)*{} **\crv{(4,15) & (-4,17)}?(1)*\dir{>};
    (-4,12)*{};(4,20)*{} **\crv{(-4,15) & (4,17)}?(1)*\dir{>};
    (4,4)*{}; (4,12) **\dir{-};
    (-12,-4)*{}; (-12,4) **\dir{-};
    (-12,12)*{}; (-12,20) **\dir{-};
  (10,8)*{n};
\endxy} \label{eq_nil_rels}
  \end{equation}
\begin{eqnarray}
  \xy
  (4,4);(4,-4) **\dir{-}?(0)*\dir{<}+(2.3,0)*{};
  (-4,4);(-4,-4) **\dir{-}?(0)*\dir{<}+(2.3,0)*{};
  (9,2)*{n};
 \endxy
 \quad =
\xy
  (0,0)*{\xybox{
    (-4,-4)*{};(4,4)*{} **\crv{(-4,-1) & (4,1)}?(1)*\dir{>}?(.25)*{\bullet};
    (4,-4)*{};(-4,4)*{} **\crv{(4,-1) & (-4,1)}?(1)*\dir{>};
     (8,1)*{ n};
     (-10,0)*{};(10,0)*{};
     }};
  \endxy
 \;\; -
 \xy
  (0,0)*{\xybox{
    (-4,-4)*{};(4,4)*{} **\crv{(-4,-1) & (4,1)}?(1)*\dir{>}?(.75)*{\bullet};
    (4,-4)*{};(-4,4)*{} **\crv{(4,-1) & (-4,1)}?(1)*\dir{>};
     (8,1)*{ n};
     (-10,0)*{};(10,0)*{};
     }};
  \endxy
 \;\; =
\xy
  (0,0)*{\xybox{
    (-4,-4)*{};(4,4)*{} **\crv{(-4,-1) & (4,1)}?(1)*\dir{>};
    (4,-4)*{};(-4,4)*{} **\crv{(4,-1) & (-4,1)}?(1)*\dir{>}?(.75)*{\bullet};
     (8,1)*{ n};
     (-10,0)*{};(10,0)*{};
     }};
  \endxy
 \;\; -
  \xy
  (0,0)*{\xybox{
    (-4,-4)*{};(4,4)*{} **\crv{(-4,-1) & (4,1)}?(1)*\dir{>} ;
    (4,-4)*{};(-4,4)*{} **\crv{(4,-1) & (-4,1)}?(1)*\dir{>}?(.25)*{\bullet};
     (8,1)*{ n};
     (-10,0)*{};(10,0)*{};
     }};
  \endxy \nn \\ \label{eq_nil_dotslide}
\end{eqnarray}

  \item All 2-morphisms are cyclic (see \cite{Lau1}) with respect to the above biadjoint structure.  Cyclicity is described by the relations:
\begin{equation} \label{eq_cyclic_dot}
    \xy
    (-8,2)*{}="1";
    (0,2)*{}="2";
    (0,-2)*{}="2'";
    (8,-2)*{}="3";
    (-8,-10);"1" **\dir{-};
    "2";"2'" **\dir{-} ?(.5)*\dir{<};
    "1";"2" **\crv{(-8,8) & (0,8)} ?(0)*\dir{<};
    "2'";"3" **\crv{(0,-8) & (8,-8)}?(1)*\dir{<};
    "3"; (8,10) **\dir{-};
    (15,-9)*{ n+2};
    (-12,9)*{n};
    (0,2)*{\txt\large{$\bullet$}};
    \endxy
    \quad = \quad
      \xy
 (0,10);(0,-10); **\dir{-} ?(.75)*\dir{<}+(2.3,0)*{\scriptstyle{}}
 ?(.1)*\dir{ }+(2,0)*{\scs };
 (0,0)*{\txt\large{$\bullet$}};
 (-6,5)*{ n};
 (8,5)*{n +2};
 (-10,0)*{};(10,0)*{};
 \endxy
    \quad = \quad
    \xy
    (8,2)*{}="1";
    (0,2)*{}="2";
    (0,-2)*{}="2'";
    (-8,-2)*{}="3";
    (8,-10);"1" **\dir{-};
    "2";"2'" **\dir{-} ?(.5)*\dir{<};
    "1";"2" **\crv{(8,8) & (0,8)} ?(0)*\dir{<};
    "2'";"3" **\crv{(0,-8) & (-8,-8)}?(1)*\dir{<};
    "3"; (-8,10) **\dir{-};
    (15,9)*{n+2};
    (-12,-9)*{n};
    (0,2)*{\txt\large{$\bullet$}};
    \endxy
\end{equation}
\begin{equation} \label{eq_cyclic_cross-gen}
\xy 0;/r.19pc/:
  (0,0)*{\xybox{
    (-4,-4)*{};(4,4)*{} **\crv{(-4,-1) & (4,1)}?(1)*\dir{>};
    (4,-4)*{};(-4,4)*{} **\crv{(4,-1) & (-4,1)};
     (4,4)*{};(-18,4)*{} **\crv{(4,16) & (-18,16)} ?(1)*\dir{>};
     (-4,-4)*{};(18,-4)*{} **\crv{(-4,-16) & (18,-16)} ?(1)*\dir{<}?(0)*\dir{<};
     (18,-4);(18,12) **\dir{-};(12,-4);(12,12) **\dir{-};
     (-18,4);(-18,-12) **\dir{-};(-12,4);(-12,-12) **\dir{-};
     (8,1)*{ n};
     (-10,0)*{};(10,0)*{};
      (4,-4)*{};(12,-4)*{} **\crv{(4,-10) & (12,-10)}?(1)*\dir{<}?(0)*\dir{<};
      (-4,4)*{};(-12,4)*{} **\crv{(-4,10) & (-12,10)}?(1)*\dir{>}?(0)*\dir{>};
     }};
  \endxy
\quad =  \quad \xy
  (0,0)*{\xybox{
    (-4,-4)*{};(4,4)*{} **\crv{(-4,-1) & (4,1)}?(0)*\dir{<} ;
    (4,-4)*{};(-4,4)*{} **\crv{(4,-1) & (-4,1)}?(0)*\dir{<};
     (-8,0)*{ n};
     (-12,0)*{};(12,0)*{};
     }};
  \endxy \quad =  \quad
 \xy 0;/r.19pc/:
  (0,0)*{\xybox{
    (4,-4)*{};(-4,4)*{} **\crv{(4,-1) & (-4,1)}?(1)*\dir{>};
    (-4,-4)*{};(4,4)*{} **\crv{(-4,-1) & (4,1)};
     (-4,4)*{};(18,4)*{} **\crv{(-4,16) & (18,16)} ?(1)*\dir{>};
     (4,-4)*{};(-18,-4)*{} **\crv{(4,-16) & (-18,-16)} ?(1)*\dir{<}?(0)*\dir{<};
     (-18,-4);(-18,12) **\dir{-};(-12,-4);(-12,12) **\dir{-};
     (18,4);(18,-12) **\dir{-};(12,4);(12,-12) **\dir{-};
     (8,1)*{ n};
     (-10,0)*{};(10,0)*{};
     (-4,-4)*{};(-12,-4)*{} **\crv{(-4,-10) & (-12,-10)}?(1)*\dir{<}?(0)*\dir{<};
      (4,4)*{};(12,4)*{} **\crv{(4,10) & (12,10)}?(1)*\dir{>}?(0)*\dir{>};
     }};
  \endxy
\end{equation}
These relations imply that isotopic diagrams represent
the same 2-morphism in $\Ucat$.

It is convenient to define degree zero 2-morphisms:
\begin{equation} \label{eq_crossl-gen}
  \xy
  (0,0)*{\xybox{
    (-4,-4)*{};(4,4)*{} **\crv{(-4,-1) & (4,1)}?(1)*\dir{>} ;
    (4,-4)*{};(-4,4)*{} **\crv{(4,-1) & (-4,1)}?(0)*\dir{<};
     (8,2)*{ n};
     (-12,0)*{};(12,0)*{};
     }};
  \endxy
:=
 \xy 0;/r.19pc/:
  (0,0)*{\xybox{
    (4,-4)*{};(-4,4)*{} **\crv{(4,-1) & (-4,1)}?(1)*\dir{>};
    (-4,-4)*{};(4,4)*{} **\crv{(-4,-1) & (4,1)};
     (-4,4);(-4,12) **\dir{-};
     (-12,-4);(-12,12) **\dir{-};
     (4,-4);(4,-12) **\dir{-};(12,4);(12,-12) **\dir{-};
     (16,1)*{n};
     (-10,0)*{};(10,0)*{};
     (-4,-4)*{};(-12,-4)*{} **\crv{(-4,-10) & (-12,-10)}?(1)*\dir{<}?(0)*\dir{<};
      (4,4)*{};(12,4)*{} **\crv{(4,10) & (12,10)}?(1)*\dir{>}?(0)*\dir{>};
     }};
  \endxy
  \quad = \quad
  \xy 0;/r.19pc/:
  (0,0)*{\xybox{
    (-4,-4)*{};(4,4)*{} **\crv{(-4,-1) & (4,1)}?(1)*\dir{<};
    (4,-4)*{};(-4,4)*{} **\crv{(4,-1) & (-4,1)};
     (4,4);(4,12) **\dir{-};
     (12,-4);(12,12) **\dir{-};
     (-4,-4);(-4,-12) **\dir{-};(-12,4);(-12,-12) **\dir{-};
     (16,1)*{n};
     (10,0)*{};(-10,0)*{};
     (4,-4)*{};(12,-4)*{} **\crv{(4,-10) & (12,-10)}?(1)*\dir{>}?(0)*\dir{>};
      (-4,4)*{};(-12,4)*{} **\crv{(-4,10) & (-12,10)}?(1)*\dir{<}?(0)*\dir{<};
     }};
  \endxy
\end{equation}
\begin{equation} \label{eq_crossr-gen}
  \xy
  (0,0)*{\xybox{
    (-4,-4)*{};(4,4)*{} **\crv{(-4,-1) & (4,1)}?(0)*\dir{<} ;
    (4,-4)*{};(-4,4)*{} **\crv{(4,-1) & (-4,1)}?(1)*\dir{>};
     (-8,2)*{ n};
     (-12,0)*{};(12,0)*{};
     }};
  \endxy
:=
 \xy 0;/r.19pc/:
  (0,0)*{\xybox{
    (-4,-4)*{};(4,4)*{} **\crv{(-4,-1) & (4,1)}?(1)*\dir{>};
    (4,-4)*{};(-4,4)*{} **\crv{(4,-1) & (-4,1)};
     (4,4);(4,12) **\dir{-};
     (12,-4);(12,12) **\dir{-};
     (-4,-4);(-4,-12) **\dir{-};(-12,4);(-12,-12) **\dir{-};
     (-16,1)*{n};
     (10,0)*{};(-10,0)*{};
     (4,-4)*{};(12,-4)*{} **\crv{(4,-10) & (12,-10)}?(1)*\dir{<}?(0)*\dir{<};
      (-4,4)*{};(-12,4)*{} **\crv{(-4,10) & (-12,10)}?(1)*\dir{>}?(0)*\dir{>};
     }};
  \endxy
  \quad = \quad
  \xy 0;/r.19pc/:
  (0,0)*{\xybox{
    (4,-4)*{};(-4,4)*{} **\crv{(4,-1) & (-4,1)}?(1)*\dir{<};
    (-4,-4)*{};(4,4)*{} **\crv{(-4,-1) & (4,1)};
     (-4,4);(-4,12) **\dir{-};
     (-12,-4);(-12,12) **\dir{-};
     (4,-4);(4,-12) **\dir{-};(12,4);(12,-12) **\dir{-};
     (-16,1)*{n};
     (-10,0)*{};(10,0)*{};
     (-4,-4)*{};(-12,-4)*{} **\crv{(-4,-10) & (-12,-10)}?(1)*\dir{>}?(0)*\dir{>};
      (4,4)*{};(12,4)*{} **\crv{(4,10) & (12,10)}?(1)*\dir{<}?(0)*\dir{<};
     }};
  \endxy
\end{equation}
where the second equalities in \eqref{eq_crossl-gen} and \eqref{eq_crossr-gen}
follow from \eqref{eq_cyclic_cross-gen}.  We also write
\[
  \xy
 (0,7);(0,-7); **\dir{-} ?(.75)*\dir{>}+(2.3,0)*{\scriptstyle{}};
 (0.1,-2)*{\txt\large{$\bullet$}};
 (6,4)*{ n};(-3,-1)*{\alpha};
 (-10,0)*{};(10,0)*{};(0,-10)*{};(0,10)*{};
 \endxy \quad := \quad \left(
  \xy
 (0,7);(0,-7); **\dir{-} ?(.75)*\dir{>}+(2.3,0)*{\scriptstyle{}};
 (0.1,-2)*{\txt\large{$\bullet$}};
 (6,4)*{ n};(-3,-1)*{};
 (-10,0)*{};(10,0)*{};(0,-10)*{};(0,10)*{};
 \endxy
 \right)^{\alpha}
\]
to denote the $\alpha$-fold vertical composite of a dot with itself.

\item  All dotted bubbles of negative degree are zero. That is,
\begin{equation} \label{eq_positivity_bubbles}
 \xy
 (-12,0)*{\cbub{\alpha}{}};
 (-8,8)*{n};
 \endxy
  = 0
 \qquad
  \text{if $\alpha<n-1$,} \qquad
 \xy
 (-12,0)*{\ccbub{\alpha}{}};
 (-8,8)*{n};
 \endxy = 0\quad
  \text{if $\alpha< -n-1$, \;\; $\alpha \in \Z_+$.}
\end{equation}
A dotted bubble of degree zero equals 1:
\[
\xy 0;/r.18pc/:
 (0,0)*{\cbub{n-1}{}};
  (4,8)*{n};
 \endxy
  = 1 \quad \text{for $n \geq 1$,}
  \qquad \quad
  \xy 0;/r.18pc/:
 (0,0)*{\ccbub{-n-1}{}};
  (4,8)*{n};
 \endxy
  = 1 \quad \text{for $n \leq -1$.}
\]
We will often express dotted bubbles via a notation that emphasizes the degree:
\[
\xy 0;/r.18pc/:
 (0,0)*{\cbub{\spadesuit+\alpha}{}};
  (4,8)*{n};
 \endxy\quad := \quad
\xy 0;/r.18pc/:
 (0,0)*{\cbub{(n-1)+\alpha}{}};
  (4,8)*{n};
 \endxy
   \qquad \quad
   \qquad
   \xy 0;/r.18pc/:
 (0,0)*{\ccbub{\spadesuit+\alpha}{}};
  (4,8)*{n};
 \endxy \quad := \quad
   \xy 0;/r.18pc/:
 (0,0)*{\ccbub{(-n-1)+\alpha}{}};
  (4,8)*{n};
 \endxy
\]
so that
\[
 \deg\left(\xy 0;/r.18pc/:
 (0,0)*{\cbub{\spadesuit+\alpha}{}};
  (4,8)*{n};
 \endxy \right) \;\; = \;\; 2\alpha
 \qquad \qquad
  \deg\left(    \xy 0;/r.18pc/:
 (0,0)*{\ccbub{\spadesuit+\alpha}{}};
  (4,8)*{n};
 \endxy\right) \;\; = \;\; 2\alpha.
\]
The value of $\spadesuit$ depends on the orientation, $\spadesuit =
n-1$ for clockwise oriented bubbles and $\spadesuit = -n-1$ for
counter-clockwise oriented bubbles.  For some values of
$n$ the number $\spadesuit +\alpha$ might be negative even though
$\alpha \geq 0$.  While vertically composing a generator with itself
a negative number of times makes no sense, having these symbols
around greatly simplifies the calculus.  For  $\spadesuit + \alpha
<0$ such that
\[
 \text{either} \qquad \deg\left(\xy 0;/r.18pc/:
 (0,0)*{\cbub{\spadesuit+\alpha}{}};
  (4,8)*{n};
 \endxy \right) \geq 0
 \qquad \text{or} \qquad
  \deg\left(    \xy 0;/r.18pc/:
 (0,0)*{\ccbub{\spadesuit+\alpha}{}};
  (4,8)*{n};
 \endxy\right) \geq 0
\]
we introduce formal symbols, called {\em fake bubbles},  inductively defined by the homogeneous terms of the equation
\begin{center}
\begin{eqnarray}
 \makebox[0pt]{ $
\left( \xy 0;/r.15pc/:
 (0,0)*{\ccbub{\spadesuit}{}};
  (4,8)*{n};
 \endxy
 +
 \xy 0;/r.15pc/:
 (0,0)*{\ccbub{\spadesuit+1}{}};
  (4,8)*{n};
 \endxy t
 + \cdots +
 \xy 0;/r.15pc/:
 (0,0)*{\ccbub{\spadesuit+\alpha}{}};
  (4,8)*{n};
 \endxy t^{\alpha}
 + \cdots
\right)
%
\left( \xy 0;/r.15pc/:
 (0,0)*{\cbub{\spadesuit}{}};
  (4,8)*{n};
 \endxy
 + \cdots +
 \xy 0;/r.15pc/:
 (0,0)*{\cbub{\spadesuit+\alpha}{}};
 (4,8)*{n};
 \endxy t^{\alpha}
 + \cdots
\right) =1$ } \nn \\ \label{eq_infinite_Grass}
\end{eqnarray}
\end{center}
and the additional condition
\[
\xy 0;/r.18pc/:
 (0,0)*{\cbub{\spadesuit}{}};
  (4,8)*{n};
 \endxy
 \quad = \quad
  \xy 0;/r.18pc/:
 (0,0)*{\ccbub{\spadesuit}{}};
  (4,8)*{n};
 \endxy
  \quad = \quad 1.
\]
Equation~\eqref{eq_infinite_Grass} is called the infinite
Grassmannian relation.  It remains valid even in high degrees when
none of the bubbles involved is fake,  see \cite[Proposition 5.5]{Lau1}.

\item For the following relations we employ the convention that all summations
are increasing, so that $\sum_{f=0}^{\alpha}$ is zero if $\alpha < 0$.
\begin{equation} \label{eq_reduction}
  \text{$\xy 0;/r.18pc/:
  (14,8)*{n};
  (-3,-8)*{};(3,8)*{} **\crv{(-3,-1) & (3,1)}?(1)*\dir{>};?(0)*\dir{>};
    (3,-8)*{};(-3,8)*{} **\crv{(3,-1) & (-3,1)}?(1)*\dir{>};
  (-3,-12)*{\bbsid};  (-3,8)*{\bbsid};
  (3,8)*{}="t1";  (9,8)*{}="t2";
  (3,-8)*{}="t1'";  (9,-8)*{}="t2'";
   "t1";"t2" **\crv{(3,14) & (9, 14)};
   "t1'";"t2'" **\crv{(3,-14) & (9, -14)};
   "t2'";"t2" **\dir{-} ?(.5)*\dir{<};
   (9,0)*{}; (-6,-8)*{\scs };
 \endxy$} \;\; = \;\; -\sum_{f_1+f_2=-n}
   \xy
  (19,4)*{n};
  (0,0)*{\bbe{}};(-2,-8)*{\scs };
  (12,-2)*{\cbub{\spadesuit+f_2}{}};
  (0,6)*{\bullet}+(3,-1)*{\scs f_1};
 \endxy
\qquad \qquad
  \text{$ \xy 0;/r.18pc/:
  (-12,8)*{n};
   (-3,-8)*{};(3,8)*{} **\crv{(-3,-1) & (3,1)}?(1)*\dir{>};?(0)*\dir{>};
    (3,-8)*{};(-3,8)*{} **\crv{(3,-1) & (-3,1)}?(1)*\dir{>};
  (3,-12)*{\bbsid};
  (3,8)*{\bbsid}; (6,-8)*{\scs };
  (-9,8)*{}="t1";
  (-3,8)*{}="t2";
  (-9,-8)*{}="t1'";
  (-3,-8)*{}="t2'";
   "t1";"t2" **\crv{(-9,14) & (-3, 14)};
   "t1'";"t2'" **\crv{(-9,-14) & (-3, -14)};
  "t1'";"t1" **\dir{-} ?(.5)*\dir{<};
 \endxy$} \;\; = \;\;
 \sum_{g_1+g_2=n}^{}
   \xy
  (-12,8)*{n};
  (0,0)*{\bbe{}};(2,-8)*{\scs};
  (-12,-2)*{\ccbub{\spadesuit+g_2}{}};
  (0,6)*{\bullet}+(3,-1)*{\scs g_1};
 \endxy
\end{equation}
\begin{eqnarray}
 \vcenter{\xy 0;/r.18pc/:
  (-8,0)*{};
  (8,0)*{};
  (-4,10)*{}="t1";
  (4,10)*{}="t2";
  (-4,-10)*{}="b1";
  (4,-10)*{}="b2";(-6,-8)*{\scs };(6,-8)*{\scs };
  "t1";"b1" **\dir{-} ?(.5)*\dir{<};
  "t2";"b2" **\dir{-} ?(.5)*\dir{>};
  (10,2)*{n};
  (-10,2)*{n};
  \endxy}
&\quad = \quad&
 -\;\;
 \vcenter{   \xy 0;/r.18pc/:
    (-4,-4)*{};(4,4)*{} **\crv{(-4,-1) & (4,1)}?(1)*\dir{>};
    (4,-4)*{};(-4,4)*{} **\crv{(4,-1) & (-4,1)}?(1)*\dir{<};?(0)*\dir{<};
    (-4,4)*{};(4,12)*{} **\crv{(-4,7) & (4,9)};
    (4,4)*{};(-4,12)*{} **\crv{(4,7) & (-4,9)}?(1)*\dir{>};
  (8,8)*{n};(-6,-3)*{\scs };
     (6.5,-3)*{\scs };
 \endxy}
  \quad + \quad
   \sum_{ \xy  (0,3)*{\scs f_1+f_2+f_3}; (0,0)*{\scs =n-1};\endxy}
    \vcenter{\xy 0;/r.18pc/:
    (-10,10)*{n};
    (-8,0)*{};
  (8,0)*{};
  (-4,-15)*{}="b1";
  (4,-15)*{}="b2";
  "b2";"b1" **\crv{(5,-8) & (-5,-8)}; ?(.05)*\dir{<} ?(.93)*\dir{<}
  ?(.8)*\dir{}+(0,-.1)*{\bullet}+(-3,2)*{\scs f_3};
  (-4,15)*{}="t1";
  (4,15)*{}="t2";
  "t2";"t1" **\crv{(5,8) & (-5,8)}; ?(.15)*\dir{>} ?(.95)*\dir{>}
  ?(.4)*\dir{}+(0,-.2)*{\bullet}+(3,-2)*{\scs \; f_1};
  (0,0)*{\ccbub{\scs \quad \spadesuit+f_2}{}};
  \endxy} \nn
 \\  \; \nn \\
 \vcenter{\xy 0;/r.18pc/:
  (-8,0)*{};(-6,-8)*{\scs };(6,-8)*{\scs };
  (8,0)*{};
  (-4,10)*{}="t1";
  (4,10)*{}="t2";
  (-4,-10)*{}="b1";
  (4,-10)*{}="b2";
  "t1";"b1" **\dir{-} ?(.5)*\dir{>};
  "t2";"b2" **\dir{-} ?(.5)*\dir{<};
  (10,2)*{n};
  (-10,2)*{n};
  \endxy}
&\quad = \quad&
 -\;\;
   \vcenter{\xy 0;/r.18pc/:
    (-4,-4)*{};(4,4)*{} **\crv{(-4,-1) & (4,1)}?(1)*\dir{<};?(0)*\dir{<};
    (4,-4)*{};(-4,4)*{} **\crv{(4,-1) & (-4,1)}?(1)*\dir{>};
    (-4,4)*{};(4,12)*{} **\crv{(-4,7) & (4,9)}?(1)*\dir{>};
    (4,4)*{};(-4,12)*{} **\crv{(4,7) & (-4,9)};
  (8,8)*{n};(-6,-3)*{\scs };  (6,-3)*{\scs };
 \endxy}
  \quad + \quad
    \sum_{ \xy  (0,3)*{\scs g_1+g_2+g_3}; (0,0)*{\scs =-n-1};\endxy}
    \vcenter{\xy 0;/r.18pc/:
    (-8,0)*{};
  (8,0)*{};
  (-4,-15)*{}="b1";
  (4,-15)*{}="b2";
  "b2";"b1" **\crv{(5,-8) & (-5,-8)}; ?(.1)*\dir{>} ?(.95)*\dir{>}
  ?(.8)*\dir{}+(0,-.1)*{\bullet}+(-3,2)*{\scs g_3};
  (-4,15)*{}="t1";
  (4,15)*{}="t2";
  "t2";"t1" **\crv{(5,8) & (-5,8)}; ?(.15)*\dir{<} ?(.9)*\dir{<}
  ?(.4)*\dir{}+(0,-.2)*{\bullet}+(3,-2)*{\scs g_1};
  (0,0)*{\cbub{\scs \quad\; \spadesuit + g_2}{}};
  (-10,10)*{n};
  \endxy} \label{eq_ident_decomp}
\end{eqnarray}
for all $n\in \Z$.  In equations \eqref{eq_reduction} and
\eqref{eq_ident_decomp} whenever the summations are nonzero they
involve fake bubbles.

\item the additive $\Bbbk$-linear composition functor $\Ucat(n,n')
 \times \Ucat(n',n'') \to\Ucat(n,n'') $ is given on
 1-morphisms of $\Ucat$ by
\begin{equation}
  \cal{E}_{\ep'}\mathbf{1}_{n'}\{t'\} \times \cal{E}_{\ep}\onen\{t\} \mapsto
  \cal{E}_{\ep'\ep}\onen\{t+t'\} ,
\end{equation}
and on 2-morphisms of $\Ucat$ by juxtaposition of diagrams
\[
\left(\;\;\vcenter{\xy 0;/r.16pc/:
 (-4,-15)*{}; (-20,25) **\crv{(-3,-6) & (-20,4)}?(0)*\dir{<}?(.6)*\dir{}+(0,0)*{\bullet};
 (-12,-15)*{}; (-4,25) **\crv{(-12,-6) & (-4,0)}?(0)*\dir{<}?(.6)*\dir{}+(.2,0)*{\bullet};
 ?(0)*\dir{<}?(.75)*\dir{}+(.2,0)*{\bullet};?(0)*\dir{<}?(.9)*\dir{}+(0,0)*{\bullet};
 (-28,25)*{}; (-12,25) **\crv{(-28,10) & (-12,10)}?(0)*\dir{<};
  ?(.2)*\dir{}+(0,0)*{\bullet}?(.35)*\dir{}+(0,0)*{\bullet};
 (-36,-15)*{}; (-36,25) **\crv{(-34,-6) & (-35,4)}?(1)*\dir{>};
 (-28,-15)*{}; (-42,25) **\crv{(-28,-6) & (-42,4)}?(1)*\dir{>};
 (-42,-15)*{}; (-20,-15) **\crv{(-42,-5) & (-20,-5)}?(1)*\dir{>};
 (6,10)*{\cbub{}{}};
 (-23,0)*{\cbub{}{}};
 (8,-4)*{n'};(-44,-4)*{n''};
 \endxy}\;\;\right) \;\; \times \;\;
\left(\;\;\vcenter{ \xy 0;/r.18pc/: (-14,8)*{\xybox{
 (0,-10)*{}; (-16,10)*{} **\crv{(0,-6) & (-16,6)}?(.5)*\dir{};
 (-16,-10)*{}; (-8,10)*{} **\crv{(-16,-6) & (-8,6)}?(1)*\dir{}+(.1,0)*{\bullet};
  (-8,-10)*{}; (0,10)*{} **\crv{(-8,-6) & (-0,6)}?(.6)*\dir{}+(.2,0)*{\bullet}?
  (1)*\dir{}+(.1,0)*{\bullet};
  (0,10)*{}; (-16,30)*{} **\crv{(0,14) & (-16,26)}?(1)*\dir{>};
 (-16,10)*{}; (-8,30)*{} **\crv{(-16,14) & (-8,26)}?(1)*\dir{>};
  (-8,10)*{}; (0,30)*{} **\crv{(-8,14) & (-0,26)}?(1)*\dir{>}?(.6)*\dir{}+(.25,0)*{\bullet};
   }};
 (-2,-4)*{n}; (-26,-4)*{n'};
 \endxy} \;\;\right)
 \;\;\mapsto \;\;
\vcenter{\xy 0;/r.16pc/:
 (-4,-15)*{}; (-20,25) **\crv{(-3,-6) & (-20,4)}?(0)*\dir{<}?(.6)*\dir{}+(0,0)*{\bullet};
 (-12,-15)*{}; (-4,25) **\crv{(-12,-6) & (-4,0)}?(0)*\dir{<}?(.6)*\dir{}+(.2,0)*{\bullet};
 ?(0)*\dir{<}?(.75)*\dir{}+(.2,0)*{\bullet};?(0)*\dir{<}?(.9)*\dir{}+(0,0)*{\bullet};
 (-28,25)*{}; (-12,25) **\crv{(-28,10) & (-12,10)}?(0)*\dir{<};
  ?(.2)*\dir{}+(0,0)*{\bullet}?(.35)*\dir{}+(0,0)*{\bullet};
 (-36,-15)*{}; (-36,25) **\crv{(-34,-6) & (-35,4)}?(1)*\dir{>};
 (-28,-15)*{}; (-42,25) **\crv{(-28,-6) & (-42,4)}?(1)*\dir{>};
 (-42,-15)*{}; (-20,-15) **\crv{(-42,-5) & (-20,-5)}?(1)*\dir{>};
 (6,10)*{\cbub{}{}};
 (-23,0)*{\cbub{}{}};
 \endxy}
 \vcenter{ \xy 0;/r.16pc/: (-14,8)*{\xybox{
 (0,-10)*{}; (-16,10)*{} **\crv{(0,-6) & (-16,6)}?(.5)*\dir{};
 (-16,-10)*{}; (-8,10)*{} **\crv{(-16,-6) & (-8,6)}?(1)*\dir{}+(.1,0)*{\bullet};
  (-8,-10)*{}; (0,10)*{} **\crv{(-8,-6) & (-0,6)}?(.6)*\dir{}+(.2,0)*{\bullet}?
  (1)*\dir{}+(.1,0)*{\bullet};
  (0,10)*{}; (-16,30)*{} **\crv{(0,14) & (-16,26)}?(1)*\dir{>};
 (-16,10)*{}; (-8,30)*{} **\crv{(-16,14) & (-8,26)}?(1)*\dir{>};
  (-8,10)*{}; (0,30)*{} **\crv{(-8,14) & (-0,26)}?(1)*\dir{>}?(.6)*\dir{}+(.25,0)*{\bullet};
   }};
 (0,-5)*{n};
 \endxy}
\]
\end{itemize}
\end{defn}

We record here some additional relations that hold in $\Ucat$.  See \cite{Lau1} for more details.

\begin{equation} \label{eq_ind_dotslide}
\xy
  (0,0)*{\xybox{
    (-4,-4)*{};(4,4)*{} **\crv{(-4,-1) & (4,1)}?(1)*\dir{>}?(.25)*{\bullet}+(-2.5,1)*{\alpha};
    (4,-4)*{};(-4,4)*{} **\crv{(4,-1) & (-4,1)}?(1)*\dir{>};
     (8,-4)*{n};
     (-10,0)*{};(10,0)*{};
     }};
  \endxy
 \;\; -
 \xy
  (0,0)*{\xybox{
    (-4,-4)*{};(4,4)*{} **\crv{(-4,-1) & (4,1)}?(1)*\dir{>}?(.75)*{\bullet}+(2.5,-1)*{\alpha};
    (4,-4)*{};(-4,4)*{} **\crv{(4,-1) & (-4,1)}?(1)*\dir{>};
     (8,-4)*{n};
     (-10,0)*{};(10,0)*{};
     }};
  \endxy
 \;\; =
\xy
  (0,0)*{\xybox{
    (-4,-4)*{};(4,4)*{} **\crv{(-4,-1) & (4,1)}?(1)*\dir{>};
    (4,-4)*{};(-4,4)*{} **\crv{(4,-1) & (-4,1)}?(1)*\dir{>}?(.75)*{\bullet}+(-2.5,-1)*{\alpha};
     (8,3)*{ n};
     (-10,0)*{};(10,0)*{};
     }};
  \endxy
 \;\; -
  \xy
  (0,0)*{\xybox{
    (-4,-4)*{};(4,4)*{} **\crv{(-4,-1) & (4,1)}?(1)*\dir{>} ;
    (4,-4)*{};(-4,4)*{} **\crv{(4,-1) & (-4,1)}?(1)*\dir{>}?(.25)*{\bullet}+(2.5,1)*{\alpha};
     (8,3)*{n};
     (-10,0)*{};(10,0)*{};
     }};
  \endxy
  \;\; = \;\;
  \sum_{f_1 + f_2 = \alpha-1}
  \xy
  (3,4);(3,-4) **\dir{-}?(0)*\dir{<} ?(.5)*\dir{}+(0,0)*{\bullet}+(2.5,1)*{\scs f_2};
  (-3,4);(-3,-4) **\dir{-}?(0)*\dir{<}?(.5)*\dir{}+(0,0)*{\bullet}+(-2.5,1)*{\scs f_1};;
  (9,-4)*{n};
 \endxy
\end{equation}

\begin{equation} \label{eq_reduction-dots}
  \text{$\xy 0;/r.18pc/:
  (14,8)*{n};
  (-3,-8)*{};(3,8)*{} **\crv{(-3,-1) & (3,1)}?(1)*\dir{>};?(0)*\dir{>};
    (3,-8)*{};(-3,8)*{} **\crv{(3,-1) & (-3,1)}?(1)*\dir{>};
  (-3,-12)*{\bbsid};  (-3,8)*{\bbsid};
  (3,8)*{}="t1";  (9,8)*{}="t2";
  (3,-8)*{}="t1'";  (9,-8)*{}="t2'";
   "t1";"t2" **\crv{(3,14) & (9, 14)};
   "t1'";"t2'" **\crv{(3,-14) & (9, -14)};
   (9,-4)*{\bullet}+(3,-1)*{\scs x};
   "t2'";"t2" **\dir{-} ?(.5)*\dir{<};
   (9,0)*{}; (-6,-8)*{\scs };
 \endxy$} \;\; = \;\; -\sum_{f_1+f_2=x-n}
   \xy
  (19,4)*{n};
  (0,0)*{\bbe{}};(-2,-8)*{\scs };
  (12,-2)*{\cbub{\spadesuit+f_2}{}};
  (0,6)*{\bullet}+(3,-1)*{\scs f_1};
 \endxy
\qquad \qquad
  \text{$ \xy 0;/r.18pc/:
  (-12,8)*{n};
   (-3,-8)*{};(3,8)*{} **\crv{(-3,-1) & (3,1)}?(1)*\dir{>};?(0)*\dir{>};
    (3,-8)*{};(-3,8)*{} **\crv{(3,-1) & (-3,1)}?(1)*\dir{>};
  (3,-12)*{\bbsid};
  (3,8)*{\bbsid}; (6,-8)*{\scs };
  (-9,8)*{}="t1";
  (-3,8)*{}="t2";
  (-9,-8)*{}="t1'";
  (-3,-8)*{}="t2'";
   "t1";"t2" **\crv{(-9,14) & (-3, 14)};
   "t1'";"t2'" **\crv{(-9,-14) & (-3, -14)};
  "t1'";"t1" **\dir{-} ?(.5)*\dir{<};  (-9,-4)*{\bullet}+(-3,-1)*{\scs x};
 \endxy$} \;\; = \;\;
 \sum_{g_1+g_2=x+n}^{}
   \xy
  (-12,8)*{n};
  (0,0)*{\bbe{}};(2,-8)*{\scs};
  (-12,-2)*{\ccbub{\spadesuit+g_2}{}};
  (0,6)*{\bullet}+(3,-1)*{\scs g_1};
 \endxy
\end{equation}
\begin{eqnarray}
 \vcenter{\xy 0;/r.18pc/:
  (-8,0)*{};
  (8,0)*{};
  (-4,10)*{}="t1";
  (4,10)*{}="t2";
  (-4,-10)*{}="b1";
  (4,-10)*{}="b2";(-6,-8)*{\scs };(6,-8)*{\scs };
  "t1";"b1" **\dir{-} ?(.5)*\dir{<};
  "t2";"b2" **\dir{-} ?(.5)*\dir{>};
   (-4,-4)*{\bullet}+(-3,-1)*{\scs x};
   (4,-4)*{\bullet}+(3,-1)*{\scs y};
  (10,2)*{n};
  (-10,2)*{n};
  \endxy}
&\quad = \quad&
 -\;\;
 \vcenter{   \xy 0;/r.18pc/:
    (-4,-4)*{};(4,4)*{} **\crv{(-4,-1) & (4,1)};
    (4,-4)*{};(-4,4)*{} **\crv{(4,-1) & (-4,1)}?(0)*\dir{<};
    (-4,4)*{};(4,12)*{} **\crv{(-4,7) & (4,9)};
    (4,4)*{};(-4,12)*{} **\crv{(4,7) & (-4,9)}?(1)*\dir{>};
  (8,8)*{n};(-6,-3)*{\scs };
     (6.5,-3)*{\scs };
   (-4,4)*{\bullet}+(-3,-1)*{\scs y};
   (4,4)*{\bullet}+(3,-1)*{\scs x};
 \endxy}
  \quad + \quad
   \sum_{ \xy  (0,3)*{\scs f_1+f_2+f_3}; (0,0)*{\scs =x+y+n-1};\endxy}
    \vcenter{\xy 0;/r.18pc/:
    (-10,10)*{n};
    (-8,0)*{};
  (8,0)*{};
  (-4,-15)*{}="b1";
  (4,-15)*{}="b2";
  "b2";"b1" **\crv{(5,-8) & (-5,-8)}; ?(.05)*\dir{<} ?(.93)*\dir{<}
  ?(.8)*\dir{}+(0,-.1)*{\bullet}+(-3,2)*{\scs f_3};
  (-4,15)*{}="t1";
  (4,15)*{}="t2";
  "t2";"t1" **\crv{(5,8) & (-5,8)}; ?(.15)*\dir{>} ?(.95)*\dir{>}
  ?(.4)*\dir{}+(0,-.2)*{\bullet}+(3,-2)*{\scs \; f_1};
  (0,0)*{\ccbub{\scs \quad \spadesuit+f_2}{}};
  \endxy} \nn
 \\  \; \nn \\
 \vcenter{\xy 0;/r.18pc/:
  (-8,0)*{};(-6,-8)*{\scs };(6,-8)*{\scs };
  (8,0)*{};
  (-4,10)*{}="t1";
  (4,10)*{}="t2";
  (-4,-10)*{}="b1";
  (4,-10)*{}="b2";
  "t1";"b1" **\dir{-} ?(.5)*\dir{>};
  "t2";"b2" **\dir{-} ?(.5)*\dir{<};
  (-4,-4)*{\bullet}+(-3,-1)*{\scs x};
   (4,-4)*{\bullet}+(3,-1)*{\scs y};
  (10,2)*{n};
  (-10,2)*{n};
  \endxy}
&\quad = \quad&
 -\;\;
   \vcenter{\xy 0;/r.18pc/:
    (-4,-4)*{};(4,4)*{} **\crv{(-4,-1) & (4,1)}?(0)*\dir{<};
    (4,-4)*{};(-4,4)*{} **\crv{(4,-1) & (-4,1)};
    (-4,4)*{};(4,12)*{} **\crv{(-4,7) & (4,9)}?(1)*\dir{>};
    (4,4)*{};(-4,12)*{} **\crv{(4,7) & (-4,9)};
  (8,8)*{n};(-6,-3)*{\scs };  (6,-3)*{\scs };
   (-4,4)*{\bullet}+(-3,-1)*{\scs y};
   (4,4)*{\bullet}+(3,-1)*{\scs x};
 \endxy}
  \quad + \quad
    \sum_{ \xy  (0,3)*{\scs g_1+g_2+g_3}; (0,0)*{\scs =x+y-n-1};\endxy}
    \vcenter{\xy 0;/r.18pc/:
    (-8,0)*{};
  (8,0)*{};
  (-4,-15)*{}="b1";
  (4,-15)*{}="b2";
  "b2";"b1" **\crv{(5,-8) & (-5,-8)}; ?(.1)*\dir{>} ?(.95)*\dir{>}
  ?(.8)*\dir{}+(0,-.1)*{\bullet}+(-3,2)*{\scs g_3};
  (-4,15)*{}="t1";
  (4,15)*{}="t2";
  "t2";"t1" **\crv{(5,8) & (-5,8)}; ?(.15)*\dir{<} ?(.9)*\dir{<}
  ?(.4)*\dir{}+(0,-.2)*{\bullet}+(3,-2)*{\scs g_1};
  (0,0)*{\cbub{\scs \quad\; \spadesuit + g_2}{}};
  (-10,10)*{n};
  \endxy} \label{eq_ident_decomp-dots}
\end{eqnarray}

The 2-category $\Ucat$ has an enriched hom that associates a graded $\Bbbk$-vector space
\begin{equation}
\HOM_{\Ucat}(x,y) := \bigoplus_{t\in \Z}\Hom_{\Ucat}(x\{t\},y)
\end{equation}
to each pair of 1-morphisms $x,y$ in $\Ucat$.  We write $\END_{\Ucat}(x):= \HOM_{\Ucat}(x,x)$.

\begin{rem}
When defining the 2-category $\Ucat$ we can take $\Bbbk$ to be any commutative ring rather than a field.  In that case, by an additive $\Bbbk$-linear 2-category we mean a 2-category enriched in $\Bbbk-\cat{Mod}-\cat{Cat}$, see \cite[Section 5.1]{Lau1} and references therein.  In particular, the set $\Ucat(\cal{E}_{\ep}\onen,\cal{E}_{\ep'}\onen)$ has the structure of a $\Bbbk$-module and horizonal composition is given by $\Bbbk$-module homomorphisms.
\end{rem}

 \subsection{Box notation for $\Ucat$}

The basic box style diagrammatic notation introduced in
Section~\ref{sec_box-line-split} can be immediately applied to
represent 2-morphisms in $\Ucat$ from $\cal{E}^{a}\onen\{t\}$ to
$\cal{E}^{a}\onen\{t'\}$ for various $a$ and $n$ and appropriate
grading shifts $t$ and $t'$.  To interpret a box diagram as a
2-morphism in $\Ucat$ we simply label the regions of such a box
diagram by weights $n\in \Z$.  Since our conventions for upward and
downward oriented arrows determine the labels of all regions of a
diagram from the label of any particular region, we label only one
region of such diagrams.   For example, we have
\begin{equation}
   \xy
 (0,0)*{\includegraphics[scale=0.5]{figs/c1-1.eps}};
 (0,0)*{D_4}; (12,3)*{n};
  \endxy
 \quad = \quad
   \xy
 (0,0)*{\includegraphics[scale=0.5]{figs/c1-2.eps}};
 (12,3)*{n};
  \endxy \maps \cal{E}^{4}\onen\{t-12\} \to \cal{E}^{4} \onen
  \{t\}.
\end{equation}

We can use the same labelled boxes from
Section~\ref{sec_box-line-split} with downward orientations by
applying the composite 2-functor $\tilde{\psi}\tilde{\omega}\tilde{\sigma}$ from \cite[Section 5.6]{Lau1}. For simplicity we write $\tau'$ for the composite $\tilde{\psi}\tilde{\omega}\tilde{\sigma}$ since it differs from the 2-functor $\tilde{\tau}$ by a grading shift.   Recall from \cite[Section
5.6]{Lau1} that the 2-functor $\tau' \maps \Ucat \to \Ucat^{{\rm
coop}}$ can be interpreted as a symmetry of the graphical
calculus which rotates a diagram by 180 degrees, in particular, it
is contravariant on 1-morphisms and 2-morphisms.  The 2-functor
$\tau'$ fixes objects of $\Ucat$, but maps a morphism $\onenn{m}\cal{E}_{\ep}\onen\{t\}$ to $\onen \cal{E}_{-\ep} \onenn{m}\{-t\}$.

We will simplify the depiction of rotated boxes by omitting the reference to $\tau'$.

Some of the most important downward oriented boxes are collected
below:
\begin{equation} 
  \xy
 (0,0)*{\includegraphics[scale=0.5,angle=180]{figs/single-up.eps}};
 (-3,3)*{a}; (6,-4)*{n};
  \endxy
\;\; = \;\;
 \tau'\left( \;\;\xy
 (0,0)*{\includegraphics[scale=0.5]{figs/single-up.eps}};
 (-3,3)*{a}; (-6,-4)*{n};
  \endxy \;\;\right)
  \;\; = \;\;
  \xy
 (0,0)*{\includegraphics[scale=0.5,angle=180]{figs/multi-up.eps}};
 (0,-11)*{\underbrace{\hspace{0.7in}}};  (0,-14)*{a}; (12,3)*{n};
  \endxy
\end{equation}
\begin{equation}  
    \xy
 (0,0)*{\includegraphics[scale=0.5,angle=180]{figs/box-up.eps}};
 (0,0)*{\delta_a}; (8,5)*{n};
  \endxy
 \;\; = \;\;
 \tau'\left(\;\;     \xy
 (0,0)*{\includegraphics[scale=0.5]{figs/box-up.eps}};
 (0,0)*{\delta_a}; (-8,5)*{n};
  \endxy \;\; \right)
 \;\; = \;\;
 \xy
 (2.6,0)*{\includegraphics[scale=0.5, angle=180]{figs/multi-up-wide.eps}};
 (17.7,3)*{\bullet}+(4,1)*{\scs a-1};
 (7.7,3)*{\bullet}+(4,1)*{\scs a-2};
 (3,-2)*{\cdots};
 (-2.3,3)*{\bullet}+(2,1)*{\scs 2};
 (-7.3,3)*{\bullet}; (26,-4)*{n};
  \endxy
 \end{equation}
\begin{equation} 
   \xy
 (0,0)*{\includegraphics[scale=0.5,angle=180]{figs/c1-1.eps}};
 (0,0)*{D_a};(0,-11)*{\underbrace{\hspace{0.7in}}};
 (0,-14)*{a}; (12,5)*{n};
  \endxy
 \quad = \quad
 \tau'\left( \;\;   \xy
 (0,0)*{\includegraphics[scale=0.5]{figs/c1-1.eps}};
 (0,0)*{D_a};(0,-11)*{\underbrace{\hspace{0.7in}}};
 (0,-14)*{a}; (-12,5)*{n};
  \endxy \;\;\right)
 \quad = \quad
   \xy
 (0,0)*{\includegraphics[scale=0.5,angle=180]{figs/c1-2.eps}};
 (12,5)*{n};
  \endxy
\end{equation}
\begin{equation} 
    \xy
 (0,0)*{\includegraphics[scale=0.5,angle=180]{figs/c2-1.eps}};
 (0,1.5)*{e_a}; (16,5)*{n};
  \endxy
 \quad = \;\;
 \xy
 (0,0)*{\includegraphics[scale=0.5,angle=180]{figs/c2-2.eps}};
 (17.7,-5)*{\bullet}+(4,-1)*{\scs a-1};
 (7.7,-5)*{\bullet}+(4,-1)*{\scs a-2};
 (3,-4)*{\cdots};(3,7)*{\cdots};
 (-2.3,-5)*{\bullet}+(2,-1)*{\scs 2};
 (-9.9,-5)*{\bullet};(0,1.5)*{D_a}; (24,5)*{n};
  \endxy \qquad \quad
   \xy
 (0,0)*{\includegraphics[angle=180,scale=0.5]{figs/box-up.eps}};
 (0,0)*{e_a}; (8,5)*{n};
  \endxy
  \quad = \quad
  \xy
 (0,0)*{\includegraphics[angle=180,scale=0.5]{figs/two-box-up.eps}};
 (0,5.5)*{D_a};(0,-4.5)*{\delta_a}; (8,5)*{n};
  \endxy
\end{equation}

 \subsection{Karoubi envelope}

The Karoubi envelope $Kar(\cal{C})$ of a category $\cal{C}$ is a minimal
enlargement of $\cal{C}$ in which all idempotents split
(see \cite[Section 9]{Lau1} and references therein). An idempotent
$e \maps b\to b$ in a category $\cal{C}$  is said to split if there
exist morphisms
\[
 \xymatrix{ b \ar[r]^g & b' \ar[r]^h &b}
\]
such that $e=hg$ and $g h = \Id_{b'}$. More precisely, the
Karoubi envelope $Kar(\cal{C})$ is a category whose objects are pairs $(b,e)$
where $e \maps b \to b$ is an idempotent of $\cal{C}$ and whose
morphisms are triples of the form
\[
 (e,f,e') \maps (b,e) \to (b',e')
\]
where $f \maps b \to b'$ in $\cal{C}$ such that $f =e'fe$. Note that this implies $f=e'f=fe=e'fe$. Composition is induced from the composition in $\cal{C}$,
and the identity morphism is $(e,e,e) \maps (b,e) \to (b,e)$.  When
$\cal{C}$ is an additive category, the splitting of idempotents
allows us to write $(b,e)\in Kar(\cal{C})$ as $\im e$, and $b
\cong \im e \oplus \im (\Id_b-e)$ in $Kar(\cal{C})$.

The identity map $\Id_b\maps b \to b$ is an idempotent and this defines a fully faithful functor $\cal{C} \to Kar(\cal{C})$. In
$Kar(\cal{C})$ all idempotents of $\cal{C}$ are split and this
functor is universal among functors which split
idempotents in $\cal{C}$.  When
$\cal{C}$ is additive the inclusion $\cal{C} \to Kar(\cal{C})$ is an
additive functor.

\begin{defn}
Define the additive $\Bbbk$-linear 2-category $\UcatD$ to have the
same objects as $\Ucat$ and hom additive $\Bbbk$-linear categories
given by $\UcatD(\lambda,\lambda') =
Kar\left(\Ucat(\lambda,\lambda')\right)$. The fully-faithful
additive $\Bbbk$-linear functors $\Ucat(\lambda,\lambda') \to
\UcatD(\lambda,\lambda')$ combine to form an additive $\Bbbk$-linear
2-functor $\Ucat \to \UcatD$ universal with respect to splitting
idempotents in the hom categories $\UcatD(\lambda,\lambda')$.  The
composition functor $\UcatD(\lambda,\lambda') \times
\UcatD(\lambda',\lambda'') \to \UcatD(\lambda,\lambda'')$ is induced
by the universal property of the Karoubi envelope from the
composition functor for $\Ucat$. The 2-category $\UcatD$ has graded
2-homs given by
\begin{equation}
\HOM_{\UcatD}(x,y) := \bigoplus_{t\in \Z}\Hom_{\UcatD}(x\{t\},y).
\end{equation}
\end{defn}

We define 1-morphisms $\cal{E}^{(a)}\onen := (\cal{E}^a\onen \{
\frac{-a(a-1)}{2}\}, e_a)$ and $\cal{F}^{(a)}\onen :=
(\cal{F}^a\onen \{ \frac{a(a-1)}{2}\}, \tau'(e_a))$ in $\UcatD$. Note
that under the fully-faithful 2-functor $\Ucat \to \UcatD$ the
1-morphisms $\cal{E}^a\onen$ and $\cal{F}^a\onen$ map to
$(\cal{E}^a\onen,\Id_{\cal{E}^a\onen})$ and
$(\cal{F}^a\onen,\Id_{\cal{F}^a\onen})$.

 \section{Thick calculus and $\UcatD$}

We can extend the box notation of section~\ref{sec_box-line-split}
to the 2-category $\Ucat$.

\subsection{Thick lines oriented up or down}

A 2-morphism in $\UcatD$ between summands of $\cal{E}_{\ep}\onen\{t\}$ and $\cal{E}_{\ep'}\onen\{t'\}$ is a triple $(e,f,e')\maps
(\cal{E}_{\ep}\onen\{t\},e) \to (\cal{E}_{\ep'}\onen\{t'\},e')$
where $e$ and $e'$ are idempotent 2-morphisms in $\Ucat$ and $f$ is
a degree $(t-t')$ 2-morphism in $\Ucat$ with the property that
\begin{equation} \label{eq_efep}
  \xy
 (0,-2)*{\includegraphics[scale=0.5]{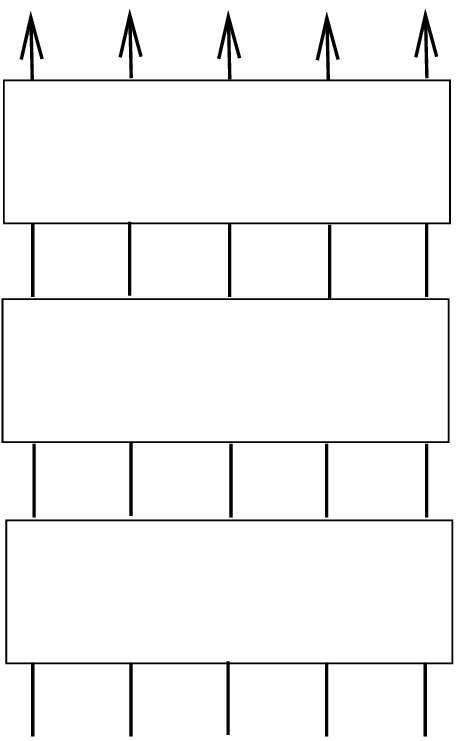}};
 (0,-2)*{f};(0,9)*{e'};(0,-13.5)*{e};
  \endxy
 \quad
 =
 \quad
   \xy
 (0,-2)*{\includegraphics[scale=0.5]{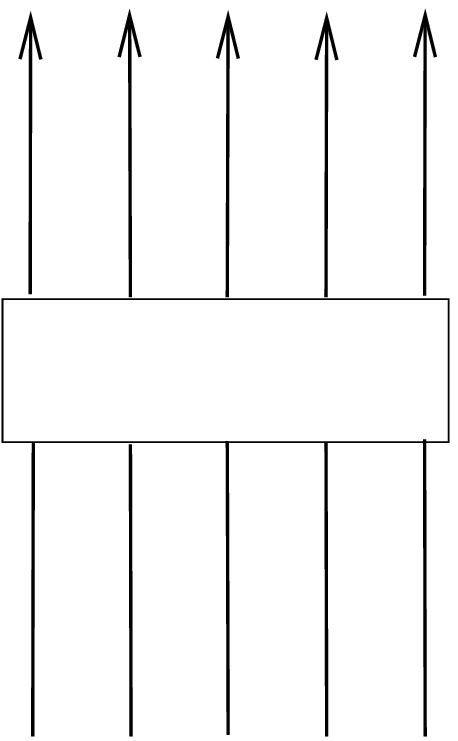}};
 (0,-2)*{f};
  \endxy
\end{equation}
In practice we often omit the idempotents from the top and bottom of
such a 2-morphism when it is clear that $f$ has the above form.

We will be primarily concerned with 2-morphisms between products of
divided powers $\cal{E}^{(a)}\onen := (\cal{E}^a\onen \{
\frac{-a(a-1)}{2}\}, e_a)$ and $\cal{F}^{(a)}\onen :=
(\cal{F}^a\onen \{ \frac{a(a-1)}{2}\}, \tau'(e_a))$.  We will use the
thick calculus to simplify the presentation of such 2-morphisms. For
example, the identity 2-morphism for $\cal{E}^{(a)}\onen$ is given
by the triple
\begin{equation}
 (e_a,e_a,e_a)
 \quad
 =
 \quad
 \left( e_a,\;
 \xy
 (0,0)*{\includegraphics[scale=0.5]{figs/c1-1.eps}};
 (0,0)*{e_a}; (11,-5)*{n};
  \endxy
 , e_a \right)
  \;\; =: \;\;     \xy
 (0,0)*{\includegraphics[scale=0.5]{figs/single-tup.eps}};
 (-2.5,-6)*{a}; (5,5)*{n};
  \endxy
\end{equation}
An upward-oriented thick line of thickness $a$
compactly describes the triple.  We represent the identity
2-morphism $(\tau'(e_b),\tau'(e_b),\tau'(e_b)) \maps \cal{F}^{(b)}\onen
\To \cal{F}^{(b)}\onen$ using a downward-oriented thick line
labelled $b$
\[ 
  (\tau'(e_b),\tau'(e_b),\tau'(e_b))  \;\; =\;\;
 \left(\tau'(e_b),\;
 \xy
 (0,0)*{\includegraphics[scale=0.5,angle=180]{figs/c1-1.eps}};
 (0,0)*{e_b}; (11,-5)*{n};
  \endxy,
\tau'(e_b)\right) \;\; =: \;\;
    \xy
 (0,0)*{\includegraphics[scale=0.5,angle=180]{figs/single-tup.eps}};
 (-2.5,6)*{b}; (5,-5)*{n};
  \endxy
\]

Products of divided powers are represented in the thick calculus by
placing appropriate upward and downward oriented thick lines side by
side.  For example, the identity 2-morphism of the product of
divided powers $\cal{E}^{(a)}\cal{F}^{(b)}\cal{F}^{(c)}\onen$ is
depicted as
\[
    \xy
(-20,0)*{\includegraphics[scale=0.5]{figs/single-tup.eps}};
 (-10,0)*{\includegraphics[scale=0.5,angle=180]{figs/single-tup.eps}};
 (0,0)*{\includegraphics[scale=0.5,angle=180]{figs/single-tup.eps}};
 (-2.5,6)*{c}; (-12.5,6)*{b};(-22.5,6)*{a}; (5,-5)*{n};
  \endxy
 \;\;  := \;\;
 \left( e_a \tau'(e_b) \tau'(e_c),e_a \tau'(e_b) \tau'(e_c), e_a \tau'(e_b) \tau'(e_c)
 \right),
\]
where $e_a \tau'(e_b) \tau'(e_c)$ denotes the horizontal composite of
2-morphisms in $\Ucat$.

\begin{rem}
Just as 2-morphisms in $\Ucat$ between objects
$\cal{E}_{\ep}\onen\{t\}$ and $\cal{E}_{\ep'}\onen\{t'\}$ are given by
linear combinations of diagrams of degree $t-t'$, the diagrams in
the thick calculus can be used to represent morphisms between
degree-shifted products of divided powers $\cal{E}^{(a)}\onen$ and
$\cal{F}^{(b)}\onen$.
\end{rem}

 \subsection{Splitters as diagrams for the inclusion of a summand}
The 1-morphism $\cal{E}^{(a)}\onen$ was defined up to a degree shift as a summand of
$\cal{E}^{a}\onen$ corresponding to the idempotent $e_a$ (see
\cite{Lau1} where it was shown that $\cal{E}^a\onen \cong
\bigoplus_{[a]!} \cal{E}^{(a)}\onen$). We depict the inclusion of
$\cal{E}^{(a)}\onen$ into the lowest degree summand of
$\cal{E}^{a}\onen$ in $\UcatD$ using the triple $(e_a,D_a,
\Id_{\cal{E}^a\onen})$, which we represent in the graphical calculus
as
\begin{equation} \label{eq_incEa}
  \xy
 (0,0)*{\includegraphics[scale=0.5]{figs/ufork-u.eps}};
 (-3,-5)*{a};
  \endxy
  \quad = \quad \left( e_a, \;
  \xy
 (0,0)*{\includegraphics[scale=0.5]{figs/c1-1.eps}};
 (0,0)*{D_a};
  \endxy\;,\Id_{\cal{E}^a\onen}\right).
\end{equation}
Notice that the 2-morphism $D_a$ of $\Ucat$ is invariant under
composing with the idempotent $e_a$ on the bottom by \eqref{eq_partial-en} and the idempotent $\Id_{\cal{E}^a\onen}$ on the top. Taking into account the degree
shift in the definition of $\cal{E}^{(a)}\onen$, the inclusion given by \eqref{eq_incEa} has degree $a(a-1)/2$.

Similarly, the projection of $\cal{E}^a\onen$ onto a highest
degree copy of $\cal{E}^{(a)}\onen$ is given by the triple
$(\Id_{\cal{E}^a\onen},e_a,e_a)$.  We represent this map in the
graphical calculus by
\begin{equation}
   \xy
 (0,0)*{\includegraphics[scale=0.5,angle=180]{figs/ufork-d.eps}};
 (-3,5)*{a};
  \endxy
  \quad = \quad \left( \Id_{\cal{E}^a\onen}, \;
  \xy
 (0,0)*{\includegraphics[scale=0.5]{figs/c1-1.eps}};
 (0,0)*{e_a};
  \endxy\;, e_a \right)
\end{equation}
This diagram also has degree $a(a-1)/2$. Similar inclusions and
projections can be defined for $\cal{F}^a\onen$ by applying the
2-functor $\tau'$ which rotates all diagrams by 180 degrees.

More generally we can represent other natural inclusions and
projections of divided powers as follows:
\begin{eqnarray}
    \xy
 (0,0)*{\includegraphics[scale=0.5]{figs/tsplit.eps}};
 (-5,-10)*{a+b};(-8,4)*{a};(8,4)*{b}; (7,-6)*{n};
  \endxy
& :=&
     \left(e_{a+b}, \;\xy
 (0,0)*{\includegraphics[scale=0.5]{figs/def-tsplitu.eps}};
 (-7,-10)*{b};(7,-10)*{a};(-6,4)*{e_a};(6,4)*{e_b}; (14,-2)*{n};
  \endxy \;, e_ae_b \right) \maps \cal{E}^{(a+b)}\onen \{t-ab\} \to
  \cal{E}^{(a)}\cal{E}^{(b)}\onen\{t\} \nn
  \\
\xy
 (0,0)*{\includegraphics[scale=0.5,angle=180]{figs/tsplitd.eps}};
 (-5,10)*{a+b};(-8,-4)*{a};(8,-4)*{b}; (7,6)*{n};
  \endxy
  &:= &
  \left( e_a e_b,\;
     \xy
 (0,0)*{\includegraphics[scale=0.5,angle=180]{figs/def-tsplitd.eps}};
 (-8,-7)*{a};(8,-7)*{b};(0,1)*{e_{a+b}};(-5,10)*{a+b}; (11,6)*{n};
  \endxy\;, e_{a+b} \right) \maps \cal{E}^{(a)}\cal{E}^{(b)}\onen \{t-ab\} \to
  \cal{E}^{(a+b)}\onen\{t\} \nn
\\
   \xy
 (0,0)*{\includegraphics[scale=0.5]{figs/tsplitd.eps}};
 (-5,-10)*{a+b};(-8,4)*{a};(8,4)*{b}; (7,-6)*{n};
  \endxy
  & := &
    \left( \tau'(e_{a+b}), \; \xy
 (0,0)*{\includegraphics[scale=0.5]{figs/def-tsplitd.eps}};
 (-8,7)*{a};(8,7)*{b};(0,-1)*{e_{a+b}};(5,-10)*{a+b};
 (11,-6)*{n};
  \endxy\; , \tau'(e_a)\tau'(e_b)\right) \maps \cal{F}^{(a+b)}\onen
  \{t-ab\} \to \cal{F}^{(a)}\cal{F}^{(b)} \onen\{t\} \nn
  \\
      \xy
 (0,0)*{\includegraphics[scale=0.5, angle=180]{figs/tsplit.eps}};
 (-5,10)*{a+b};(-8,-4)*{a};(8,-4)*{b}; (7,6)*{n};
  \endxy
  & =&
     \left(\tau'(e_a)\tau'(e_b), \;  \xy
 (0,0)*{\includegraphics[scale=0.5, angle=180]{figs/def-tsplitu.eps}};
 (-7,10)*{b};(7,10)*{a};(-6,-4)*{e_a};(6,-4)*{e_b}; (-11,6)*{n};
  \endxy\; ,\tau'(e_{a+b}) \right) \maps \cal{F}^{(a)}\cal{F}^{(b)} \onen
  \{t-ab\} \to \cal{F}^{(a+b)}\onen\{t\} \nn
\end{eqnarray}
In particular, all of the above diagrams in the thick calculus have
degree $-ab$. For example, in the first diagram the degree shift in
the divided power for $\cal{E}^{(a+b)}\onen$ is
$-\frac{(a+b)(a+b-1)}{2}$, while the degree shift in the composite
$\cal{E}^{(a)}\cal{E}^{(b)}\onen$ is
$-\frac{a(a-1)}{2}-\frac{b(b-1)}{2}$, so that the net difference is
$-\frac{2ab}{2} = -ab$.

\begin{rem}
To simplify computations in $\UcatD$ we will henceforth use the
thick calculus as defined above.  We will often omit explicit
reference to the idempotents $e$ and $e'$ in the triples $(e,f,e')$
and work with the 2-morphism $f$ in $\Ucat$ with the property in
equation \eqref{eq_efep}. Using the thick calculus the idempotents
are easily recovered via the convention that each upward oriented
thick line of thickness $a$ correspond to the idempotent $e_a$ and
each downward oriented thick line of thickness $b$ corresponds to
the idempotent $\tau'(e_b)$.

In Section~\ref{sec_nil_thick} we developed the thick calculus of
lines oriented up for computations in the nilHecke algebra.
Identities obtained there give identities on 2-morphisms between
compositions of $\cal{E}^{(a)}\onen$ for various $a$ and $n$.
\end{rem}

\subsection{Adding isotopies via thick caps and cups}

In this section we augment the thick calculus to include thick cups
and caps giving the almost biadjointness data for divided powers
$\cal{E}^{(a)}\onen$ and $\cal{F}^{(a)}\onen$.

\begin{eqnarray} 
 \xy
 (0,0)*{\includegraphics[scale=0.5, angle=180]{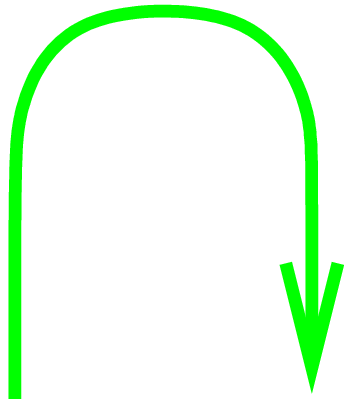}};
 (-12,-5)*{a}; (12,6)*{n};
  \endxy
  \;\; := \;\;
  \left(e_a\tau'(e_a), \; \xy
 (0,0)*{\reflectbox{\includegraphics[scale=0.5, angle=180]{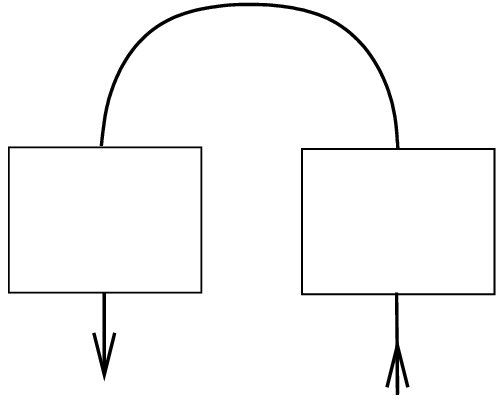}}};
 (-7,-1)*{e_a}; (7,-1)*{e_a};(12,8)*{n};
  \endxy\;,\; \Id_{\onen}
  \right) \maps \cal{E}^{(a)}\cal{F}^{(a)}\onen\{t+a(a-n)\} \to \onen
  \{t\} \nn
   \\
   \xy
 (0,0)*{\reflectbox{\includegraphics[scale=0.5, angle=180]{figs/tcup-ef.eps}}};
 (-12,5)*{a}; (12,6)*{n};
  \endxy
  \;\; := \;\;
  \left(\tau'(e_a)e_a, \;   \xy
 (0,0)*{\includegraphics[scale=0.5, angle=180]{figs/cup-twobox-fe.eps}};
 (-7,-1)*{e_a}; (7,-1)*{e_a};(12,8)*{n};
  \endxy\;,\; \Id_{\onen}
  \right) \maps \cal{F}^{(a)}\cal{E}^{(a)}\onen\{t+a(a+n)\} \to \onen
  \{t\} \nn
  \\ 
  \xy
 (0,0)*{\includegraphics[scale=0.5]{figs/tcup-ef.eps}};
 (12,5)*{a};(12,-6)*{n};
  \endxy
  \;\; := \;\;
  \left( \Id_{\onen}, \;    \xy
   (0,0)*{\reflectbox{\includegraphics[scale=0.5]{figs/cup-twobox-fe.eps}}};
 (-7,1)*{e_a}; (7,1)*{e_a}; (12,-8)*{n};
  \endxy\;,\; e_a\tau'(e_a)
  \right) \maps \onen\{t+a(a-n)\} \to \cal{F}^{(a)}\cal{E}^{(a)}\onen
  \{t\} \nn
   \\
   \xy
 (0,0)*{\reflectbox{\includegraphics[scale=0.5]{figs/tcup-ef.eps}}};
 (-12,5)*{a}; (12,-6)*{n};
  \endxy
  \;\; := \;\;
  \left( \Id_{\onen}, \;   \xy
 (0,0)*{\includegraphics[scale=0.5]{figs/cup-twobox-fe.eps}};
 (-7,1)*{e_a}; (7,1)*{e_a}; (12,-8)*{n};
  \endxy\;,\; \tau'(e_a)e_a
  \right) \maps \onen\{t+a(a-n)\} \to \cal{E}^{(a)}\cal{F}^{(a)}\onen
  \{t\} \nn
\end{eqnarray}
Note that clockwise oriented caps and cups of thickness $a$ have
degree $a(a-n)$, while counter-clockwise oriented caps and cups of
thickness $a$ have degree $a(a+n)$.

\begin{prop}[Biadjointness]
\[
   \xy
 (0,0)*{\includegraphics[scale=0.5]{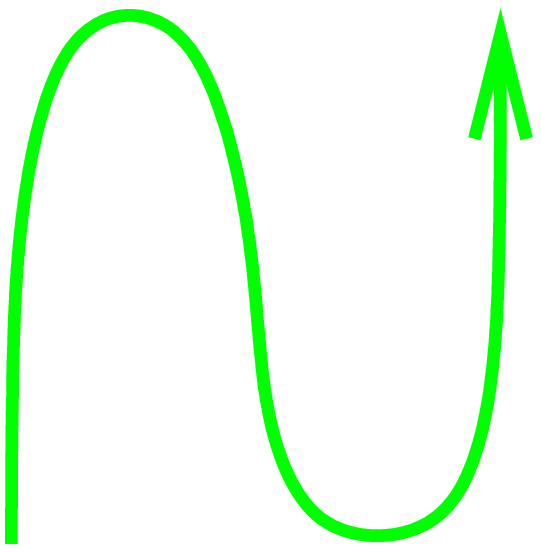}};
 (-15,-7)*{a};
  \endxy
 \quad = \qquad
    \xy
 (0,0)*{\includegraphics[scale=0.5]{figs/tlong-up.eps}};
 (-3,-7)*{a};
  \endxy
 \qquad = \quad
    \xy
 (0,0)*{\reflectbox{\includegraphics[scale=0.5]{figs/zigzag1.eps}}};
 (-15,-7)*{a};
  \endxy
\]
\[
   \xy
 (0,0)*{\includegraphics[scale=0.5, angle=180]{figs/zigzag1.eps}};
 (-15,-7)*{a};
  \endxy
 \quad = \qquad
    \xy
 (0,0)*{\includegraphics[scale=0.5, angle=180]{figs/tlong-up.eps}};
 (-3,-7)*{a};
  \endxy
 \qquad = \quad
    \xy
 (0,0)*{\reflectbox{\includegraphics[scale=0.5, angle=180]{figs/zigzag1.eps}}};
 (-15,-7)*{a};
  \endxy
\]
\end{prop}

\begin{proof}
Since all 2-morphisms in $\Ucat$ are cyclic with respect to the biadjoint structure, idempotents can be slid around caps and cups. The proposition then follows using that $e_a^2=e_a$ and equation \eqref{eq_partial-en}.
\end{proof}

Using the biadjoint structure any relation involving upward pointing
arrows can be translated into a relation involving downward pointing
arrows by rotating the relation by 180 degrees using $\tau'$.

\begin{prop}[Cyclic Structure]
\begin{equation} \label{eq_lem_cyclic1}
   \xy
 (0,0)*{\includegraphics[scale=0.5]{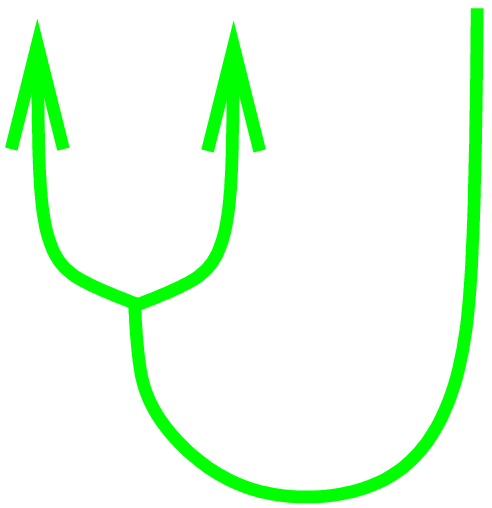}};
 (-13,3)*{a};(2,3)*{b};
  \endxy
 \quad = \quad
    \xy
 (0,0)*{\includegraphics[scale=0.5]{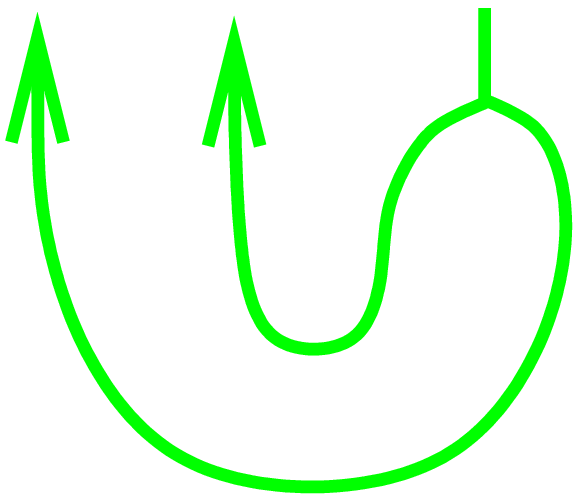}};
 (-15,3)*{a};(0,3)*{b};
  \endxy
\qquad \qquad
   \xy
 (0,0)*{\reflectbox{\includegraphics[scale=0.5]{figs/cyclic1.eps}}};
 (13,3)*{b};(-2,3)*{a};
  \endxy
 \quad = \quad
    \xy
 (0,0)*{\reflectbox{\includegraphics[scale=0.5]{figs/cyclic3.eps}}};
 (15,3)*{b};(0,3)*{a};
  \endxy
\end{equation}
\begin{equation} \label{eq_lem_cyclic2}
   \xy
 (0,0)*{\reflectbox{ \includegraphics[scale=0.5, angle=180]{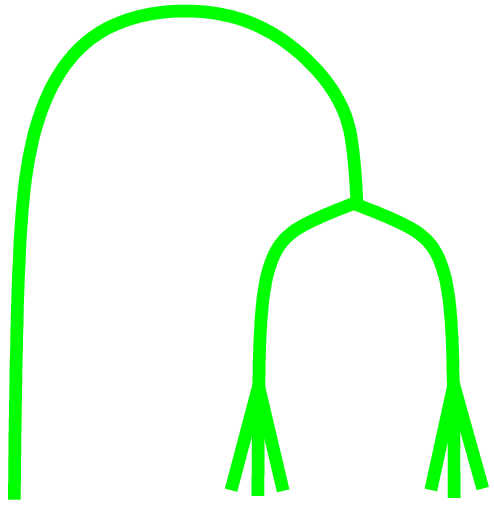}}};
 (-13,-5)*{a};(1,-5)*{b};
  \endxy
 \quad = \quad
    \xy
 (0,0)*{\reflectbox{\includegraphics[scale=0.5, angle=180]{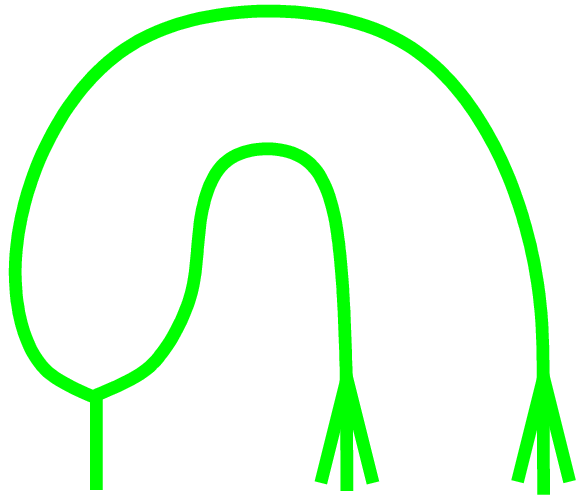}}};
  (-15,-5)*{a};(-1,-5)*{b};
  \endxy
 \qquad \qquad
   \xy
 (0,0)*{\includegraphics[scale=0.5, angle=180]{figs/cyclic4.eps}};
 (13,-5)*{b};(-1,-5)*{a};
  \endxy
 \quad = \quad
    \xy
 (0,0)*{\includegraphics[scale=0.5, angle=180]{figs/cyclic2.eps}};
   (15,-5)*{b};(1,-5)*{a};
  \endxy
\end{equation}
\end{prop}

\begin{proof}
For equation \eqref{eq_lem_cyclic1} use \eqref{eq_various_partial}.
For \eqref{eq_lem_cyclic2} use $e_a^2 = e_a$ together with
\eqref{eq_projector_absorb}.
\end{proof}

 \subsection{Thin bubble slides}

Recall the following thin bubble slide moves from \cite{Lau1}
\begin{eqnarray} \label{eq_bubslide1}
 \xy
  (-5,8)*{n};
  (0,0)*{\bbe{}};
  (12,-2)*{\cbub{\spadesuit+j}{}};
 \endxy
  & =&
     \xy
  (0,8)*{n};
  (12,0)*{\bbe{}};
  (0,-2)*{\cbub{\spadesuit+(j-2)}{}};
  (12,6)*{\bullet}+(3,-1)*{\scs 2};
 \endxy
   -2 \;
         \xy
  (0,8)*{n};
  (12,0)*{\bbe{}};
  (0,-2)*{\cbub{\spadesuit+(j-1)}{}};
  (12,6)*{\bullet}+(8,-1)*{\scs };
 \endxy
 + \;\;
     \xy
  (0,8)*{n};
  (12,0)*{\bbe{}};
  (0,-2)*{\cbub{\spadesuit+j}{}};
  (12,6)*{}+(8,-1)*{\scs };
 \endxy
 \\ \nn \\ \nn \\
  \xy
  (17,8)*{n};
  (12,0)*{\bbe{}};
  (0,-2)*{\ccbub{\spadesuit+j}{}};
  (12,6)*{}+(8,-1)*{\scs };
 \endxy
  &=&
    \xy
  (15,8)*{n};
  (0,0)*{\bbe{}};
  (12,-2)*{\ccbub{\quad\spadesuit+(j-2)}{}};
  (0,6)*{\bullet }+(3,1)*{\scs 2};
 \endxy
  -2 \;
      \xy
  (15,8)*{n};
  (0,0)*{\bbe{}};
  (12,-2)*{\ccbub{\quad\spadesuit+(j-1)}{}};
  (0,6)*{\bullet }+(5,-1)*{\scs };
 \endxy
 + \;\;
      \xy
  (15,8)*{n};
  (0,0)*{\bbe{}};
  (12,-2)*{\ccbub{\spadesuit+j}{}};
 \endxy
 \end{eqnarray}

\begin{eqnarray}
\xy
  (14,8)*{n};
  (0,0)*{\bbe{}};
  (12,-2)*{\ccbub{\spadesuit+j}{}};
  (0,6)*{ }+(7,-1)*{\scs  };
 \endxy
 & = &
  \xsum{f=0}{j}(j+1-f)
   \xy
  (0,8)*{n+2};
  (12,0)*{\bbe{}};
  (0,-2)*{\ccbub{\spadesuit+f}{}};
  (12,6)*{\bullet}+(5,-1)*{\scs j-f};
 \endxy \\ \nn \\  \nn \\
     \xy
  (15,8)*{n};
  (11,0)*{\bbe{}};
  (0,-2)*{\cbub{\spadesuit+j\quad }{}};
 \endxy
   &= &
     \xsum{f=0}{j}(j+1-f)
     \xy
  (18,8)*{n};
  (0,0)*{\bbe{}};
  (14,-4)*{\cbub{\spadesuit+f}{}};
  (0,6)*{\bullet }+(5,-1)*{\scs j-f};
 \endxy \label{eq_bubslide2}
\end{eqnarray}

A simple inductive argument shows that these bubble slide moves are valid for fake bubbles as well as real bubbles.

\begin{prop}  \label{prop_thinslide} For any $n \in \Z$ and $j \geq 0$ the equalities
\begin{equation}
  \xy
 (0,0)*{\includegraphics[scale=0.5]{figs/tlong-up.eps}};
 (-2.5,-11)*{a};
 (10,-2)*{\cbub{\spadesuit+j}{}}; (9,9)*{n};
  \endxy \quad = \quad
  \sum_{a \geq p \geq q \geq 0} (-1)^{p+q}(p-q+1)
   \xy
 (0,0)*{\includegraphics[scale=0.5]{figs/tlong-up.eps}};
 (-2.5,-11)*{a};(0,-2)*{\bigb{\pi_{\overline{p,q}}}};
 (-12,-2)*{\cbub{\spadesuit+j-(p+q)}{}}; (9,9)*{n};
  \endxy
\end{equation}
\begin{equation}
  \xy
 (0,0)*{\includegraphics[scale=0.5]{figs/tlong-up.eps}};
 (-2.5,-11)*{a};
 (-12,-2)*{\ccbub{\spadesuit+j}{}}; (9,9)*{n};
  \endxy \quad = \quad
  \sum_{a \geq p \geq q \geq 0} (-1)^{p+q}(p-q+1) \;\;
   \xy
 (0,0)*{\includegraphics[scale=0.5]{figs/tlong-up.eps}};
 (-2.5,-11)*{a};(0,-2)*{\bigb{\pi_{\overline{p,q}}}};
 (14,-2)*{\ccbub{\spadesuit+j-(p+q)}{}}; (9,9)*{n};
  \endxy
\end{equation}
hold in $\UcatD$.
\end{prop}

\begin{proof}
The proof is by induction on the thickness $a$.  The case $a=1$
follows from the thin calculus \eqref{eq_bubslide1}.  Assume the
result holds for thickness $a$, we will show that it holds for
thickness $a+1$.
\begin{equation}
  \xy
 (0,0)*{\includegraphics[scale=0.5]{figs/tlong-up.eps}};
 (-5,-11)*{a+1};
 (9,-2)*{\xy 0;/r.19pc/: (-2,0)*{\cbub{\spadesuit+j}{}}; \endxy}; (9,9)*{n};
  \endxy \;\; \refequal{\eqref{eq_schur_left}} \;\;
   \xy
 (0,0)*{\includegraphics[scale=0.5]{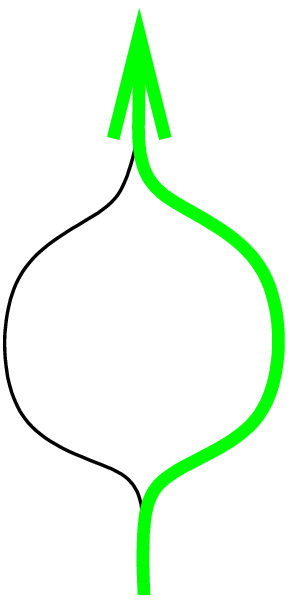}};
 (-5,-13)*{a+1}; (8,-8)*{a}; (9,9)*{n}; (-6,2)*{\bullet}+(-2,1)*{\scs a};
 (15,0)*{\xy 0;/r.18pc/: (-2,0)*{\cbub{\spadesuit+j}{}}; \endxy};
  \endxy
   \refequal{{\rm Induction}}
  \sum_{a \geq p \geq q \geq 0} (-1)^{p+q}(p-q+1) \;\;
   \xy
 (0,0)*{\includegraphics[scale=0.5]{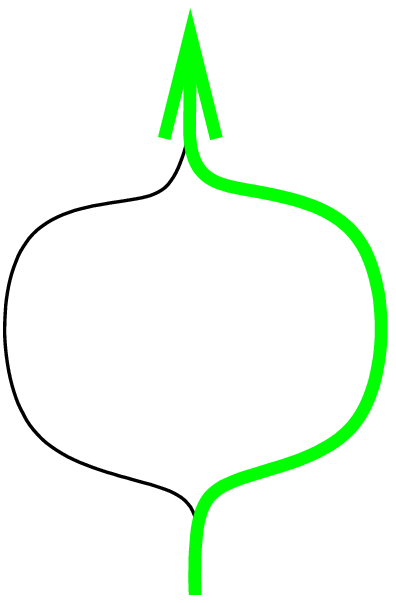}};
 (-5,-13)*{a+1}; (8,-10)*{a};(10,0)*{\bigb{\pi_{\overline{p,q}}}};
 (0,-1)*{\xy 0;/r.19pc/: (-2,0)*{\cbub{\spadesuit+j-(p+q)}{}}; \endxy}; (9,9)*{n};
 (-9,2)*{\bullet}+(-2,1)*{\scs a};
  \endxy
\end{equation}
\begin{equation}
=
    \sum_{a \geq p \geq q \geq 0} (-1)^{p+q}(p-q+1) \left[\;\;
   \xy
 (0,0)*{\includegraphics[scale=0.5]{figs/split-thinthick.eps}};
 (-5,-13)*{a+1}; (8,-8)*{a}; (9,9)*{n}; (-6,2)*{\bullet}+(-4,2)*{\scs a+2};
 (-17,-4)*{\xy 0;/r.18pc/: (-2,-2)*{\cbub{\spadesuit+j-(p+q)-2}{}}; \endxy};
  (8,-2)*{\bigb{\pi_{\overline{p,q}}}};
  \endxy
 \right.
 \end{equation}
\begin{equation}
\hspace{1in} \left.
   \;\; -2\;\;
     \xy
 (0,0)*{\includegraphics[scale=0.5]{figs/split-thinthick.eps}};
 (-5,-13)*{a+1}; (8,-8)*{a}; (9,9)*{n}; (-6,2)*{\bullet}+(-4,2)*{\scs a+1};
 (-17,-4)*{\xy 0;/r.18pc/: (-2,-2)*{\cbub{\spadesuit+j-(p+q)-1}{}}; \endxy};
 (8,-2)*{\bigb{\pi_{\overline{p,q}}}};
  \endxy \;\; +\;\;
     \xy
 (0,0)*{\includegraphics[scale=0.5]{figs/split-thinthick.eps}};
 (-5,-13)*{a+1}; (8,-8)*{a}; (9,9)*{n}; (-6,2)*{\bullet}+(-2,2)*{\scs a};
 (-15,-4)*{\xy 0;/r.18pc/: (-2,-2)*{\cbub{\spadesuit+j-(p+q)}{}}; \endxy};
  (8,-2)*{\bigb{\pi_{\overline{p,q}}}};
  \endxy\quad
  \right]
\end{equation}
\begin{equation}
\refequal{\eqref{eq_schur_left}} \;\; \sum_{a \geq p \geq q \geq 0} (-1)^{p+q}(p-q+1)
  \xy
 (8,0)*{\includegraphics[scale=0.5]{figs/tlong-up.eps}};
 (13,-11)*{a+1};
 (-4,-4)*{\xy 0;/r.19pc/: (-2,0)*{\cbub{\spadesuit+j-(p+q)-2}{}}; \endxy}; (12,9)*{n};
 (8,-2)*{\bigb{\pi_{\overline{p+1,q+1}}}};
  \endxy
 \quad - \;\; 2\sum_{a \geq p \geq 0} (-1)^{p}(p+1)
   \xy
 (8,0)*{\includegraphics[scale=0.5]{figs/tlong-up.eps}};
 (13,-11)*{a+1};
 (-4,-4)*{\xy 0;/r.19pc/: (-2,0)*{\cbub{\spadesuit+j-p-1}{}}; \endxy}; (12,9)*{n};
 (8,-2)*{\bigb{\pi_{\overline{p+1,0}}}};
  \endxy
\end{equation}
\begin{equation}
  \;\; - \;\;  \sum_{a \geq p \geq 1} (-1)^{p+1} p
     \xy
 (8,0)*{\includegraphics[scale=0.5]{figs/tlong-up.eps}};
 (13,-11)*{a+1};
 (-4,-4)*{\xy 0;/r.19pc/: (-2,0)*{\cbub{\spadesuit+j-p-1}{}}; \endxy}; (12,9)*{n};
 (8,-2)*{\bigb{\pi_{\overline{p+1,0}}}};
  \endxy \quad + \;\;
       \xy
 (8,0)*{\includegraphics[scale=0.5]{figs/tlong-up.eps}};
 (13,-11)*{a+1};
 (-4,-4)*{\xy 0;/r.19pc/: (-2,0)*{\cbub{\spadesuit+j}{}}; \endxy}; (12,9)*{n};
 (8,-2)*{\bigb{\pi_{\overline{0,0}}}};
  \endxy
\end{equation}
Shifting the indices in the first term and simplifying the remaining
terms gives:
\begin{equation}
= \;\; \sum_{a+1 \geq p \geq q \geq 1} (-1)^{p+q}(p-q+1)
     \xy
 (8,0)*{\includegraphics[scale=0.5]{figs/tlong-up.eps}};
 (13,-11)*{a+1};
 (-4,-4)*{\xy 0;/r.19pc/: (-2,0)*{\cbub{\spadesuit+j-(p+q)}{}}; \endxy}; (12,9)*{n};
 (8,-2)*{\bigb{\pi_{\overline{p,q}}}};
  \endxy\;\; + \;\;
  \sum_{a+1 \geq p  \geq 0} (-1)^{p}(p+1)
       \xy
 (8,0)*{\includegraphics[scale=0.5]{figs/tlong-up.eps}};
 (13,-11)*{a+1};
 (-4,-4)*{\xy 0;/r.19pc/: (-2,0)*{\cbub{\spadesuit+j-p}{}}; \endxy}; (12,9)*{n};
 (8,-2)*{\bigb{\pi_{\overline{p,0}}}};
  \endxy
\end{equation}
as desired.  The second identity is proved similarly.
\end{proof}

\begin{prop} \label{prop_Schur-elem}
\begin{equation}
\sum_{a \geq p \geq q \geq 0} (-1)^{p+q}(p-q+1)\pi_{\overline{(p,q)}} \quad
= \quad \sum_{m=0}^{2a}(-1)^m  \sum_{ \xy  (0,3)*{\scs x,y \leq a}; (0,0)*{\scs x+y=m};\endxy} \varepsilon_x \varepsilon_{y}
\end{equation}
where $\varepsilon_x$ is the $x$th elementary symmetric polynomial in $a$ variables.
\end{prop}

\begin{proof}
Using the (dual) Giambelli formula
\begin{equation}
\pi_{\overline{(p,q)}} = \varepsilon_p\varepsilon_q - \varepsilon_{p+1}\varepsilon_{q-1}.
\end{equation}
Therefore,
\begin{eqnarray}
\sum_{a \geq p \geq q \geq 0} (-1)^{p+q}(p-q+1) \left( \varepsilon_p\varepsilon_q - \varepsilon_{p+1}\varepsilon_{q-1} \right)\hspace{2in} \nn\\ = \sum_{m=0}^{2a}(-1)^m
\left[
 \sum_{ \xy  (0,3)*{\scs a \geq p \geq q \geq 0}; (0,0)*{\scs p+q=m};\endxy}
 (p-q+1) \varepsilon_p \varepsilon_q -
 \sum_{ \xy  (0,3)*{\scs a \geq p \geq q +2\geq 2}; (0,0)*{\scs p+q=m};\endxy}
 (p-q-1) \varepsilon_p \varepsilon_q
\right] \\
= \sum_{m=0}^{2a}(-1)^m
 \left[
 \varepsilon^2_{\frac{m}{2}} + 2 \sum_{ \xy  (0,3)*{\scs a \geq p > q \geq 0}; (0,0)*{\scs p+q=m};\endxy}\varepsilon_p \varepsilon_q
 \right] \;\;=\;\; \sum_{m=0}^{2a}(-1)^m \sum_{ \xy  (0,3)*{\scs x,y \leq a}; (0,0)*{\scs x+y=m};\endxy} \varepsilon_x \varepsilon_y.
\end{eqnarray}
\end{proof}

\begin{cor} \label{cor_thinslide}
\begin{equation}
  \xy
 (0,0)*{\includegraphics[scale=0.5]{figs/tlong-up.eps}};
 (-2.5,-11)*{a};
 (10,-2)*{\cbub{\spadesuit+j}{}}; (9,9)*{n};
  \endxy \quad = \quad
  \sum_{m=0}^{2a}(-1)^m \sum_{ \xy  (0,3)*{\scs x,y \leq a}; (0,0)*{\scs x+y=m};\endxy}
   \xy
 (0,0)*{\includegraphics[scale=0.5]{figs/tlong-up.eps}};
 (-2.5,-11)*{a};(0,-5)*{\bigb{\varepsilon_y}};(0,3)*{\bigb{\varepsilon_x}};
 (-12,-2)*{\cbub{\spadesuit+j-(x+y)}{}}; (9,9)*{n};
  \endxy
\end{equation}
\begin{equation}
  \xy
 (0,0)*{\includegraphics[scale=0.5]{figs/tlong-up.eps}};
 (-2.5,-11)*{a};
 (-12,-2)*{\ccbub{\spadesuit+j}{}}; (9,9)*{n};
  \endxy \quad = \quad
\sum_{m=0}^{2a}(-1)^m \sum_{ \xy  (0,3)*{\scs x,y \leq a}; (0,0)*{\scs x+y=m};\endxy}
   \xy
 (0,0)*{\includegraphics[scale=0.5]{figs/tlong-up.eps}};
 (-2.5,-11)*{a};(0,-5)*{\bigb{\varepsilon_y}};(0,3)*{\bigb{\varepsilon_x}};
 (14,-2)*{\ccbub{\spadesuit+j-(x+y)}{}}; (9,9)*{n};
  \endxy
\end{equation}
\end{cor}

\begin{prop} \label{prop_thinslide2}
\begin{equation}
  \xy
 (0,0)*{\includegraphics[scale=0.5]{figs/tlong-up.eps}};
 (-2.5,-11)*{a};
 (10,-2)*{\ccbub{\spadesuit+j}{}}; (9,9)*{n};
  \endxy \quad = \quad
  \sum_{j\geq p \geq q \geq 0} (p-q+1)
  \xy
 (0,0)*{\includegraphics[scale=0.5]{figs/tlong-up.eps}};
 (-2.5,-11)*{a};(0,0)*{\bigb{\pi_{(p,q)}}};
 (-12,-2)*{\ccbub{\spadesuit+j-(p+q)}{}}; (9,9)*{n};
  \endxy
  \quad = \quad
  \sum_{m=0}^{j} \sum_{ \xy  (0,3)*{\scs x,y \leq a}; (0,0)*{\scs x+y=m};\endxy}
   \xy
 (0,0)*{\includegraphics[scale=0.5]{figs/tlong-up.eps}};
 (-2.5,-11)*{a};(0,-5)*{\bigb{h_y}};(0,3)*{\bigb{h_x}};
 (-12,-2)*{\ccbub{\spadesuit+j-m}{}}; (9,9)*{n};
  \endxy
\end{equation}
\begin{equation}
  \xy
 (0,0)*{\includegraphics[scale=0.5]{figs/tlong-up.eps}};
 (-2.5,-11)*{a};
 (-12,-2)*{\cbub{\spadesuit+j}{}}; (9,9)*{n};
  \endxy \quad = \quad
  \sum_{j\geq p \geq q \geq 0} (p-q+1) \qquad
     \xy
 (0,0)*{\includegraphics[scale=0.5]{figs/tlong-up.eps}};
 (-2.5,-11)*{a};(0,2)*{\bigb{\pi_{(p,q)}}};
 (14,-3)*{\cbub{\spadesuit+j-m}{}}; (9,9)*{n};
  \endxy
  \quad = \quad
\sum_{m=0}^{j} \sum_{ \xy  (0,3)*{\scs x,y \leq a}; (0,0)*{\scs
x+y=m};\endxy}
   \xy
 (0,0)*{\includegraphics[scale=0.5]{figs/tlong-up.eps}};
 (-2.5,-11)*{a};(0,-5)*{\bigb{h_y}};(0,3)*{\bigb{h_x}};
 (14,-2)*{\cbub{\spadesuit+j-m}{}}; (9,9)*{n};
  \endxy
\end{equation}
\end{prop}

\begin{proof}
The proof is similar to the proof of Proposition~\ref{prop_thinslide} and Corollary~\ref{cor_thinslide}.
\end{proof}

 \subsection{Thick bubbles}

Using the equations above, a thick bubble can always be written as a
closed diagram in the thin calculus.
\begin{equation}
 \xy
  (0,0)*{\reflectbox{\includegraphics[scale=0.5]{figs/ccbub.eps}}};
  (-5,9)*{n}; (3,-8)*{a};
 \endxy
\quad =  \quad \xy
 (0,0)*{\includegraphics[scale=0.35]{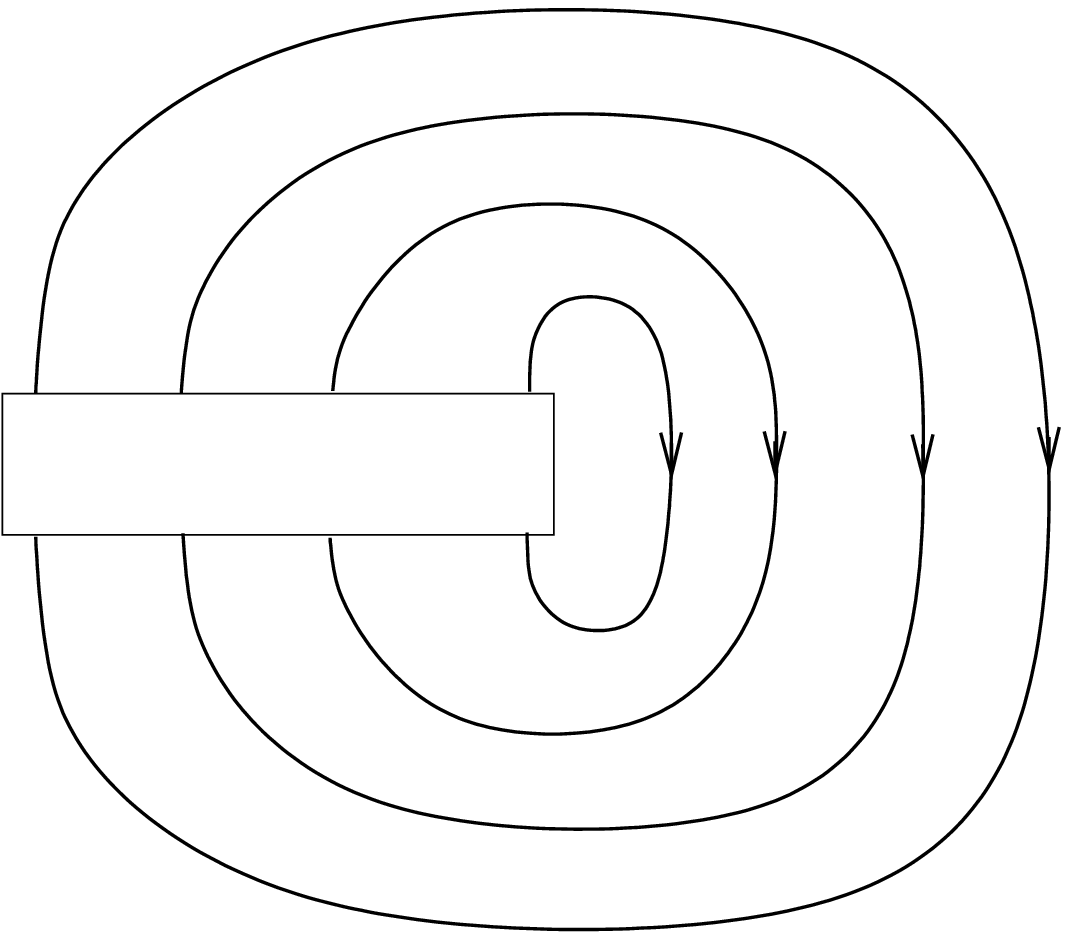}};
 (-9,0)*{e_a};
 \endxy
\qquad
 \xy
  (0,0)*{\includegraphics[scale=0.5]{figs/ccbub.eps}};
  (-5,9)*{n}; (3,-8)*{a};
 \endxy
\quad =  \quad \xy
 (0,0)*{\reflectbox{\includegraphics[scale=0.35]{figs/box-loop.eps}}};
 (9,0)*{e_a};
 \endxy
\end{equation}
where we have cancelled redundant copies of the idempotent $e_a$
using that $e_a^2 = e_a$.  Thick bubbles have the following degrees
\begin{equation}
\deg \left(\;
  \xy
 (0,0)*{\includegraphics[scale=0.5]{figs/ccbub.eps}};
(-5,9)*{n}; (3,-8)*{a};
 \endxy \right) \;\; =\;\; 2a(a+n)\qquad \qquad \quad
 \deg \left(\;
  \xy
 (0,0)*{\reflectbox{\includegraphics[scale=0.5]{figs/ccbub.eps}}};
(-5,9)*{n}; (3,-8)*{a};
 \endxy \right) \;\; =\;\; \quad 2a(a-n).
 \end{equation}
By the positivity of bubbles axioms, the degree of a
nonvanishing thick bubble must be nonnegative.   We use the $\spadesuit$-notation
when we want to emphasize the degree of a  bubble.

Recall our conventions for partitions from \eqref{eq_def_partition_add}. For $\alpha,\beta \in P(a)$ with $\alpha_a \geq n+a$ and $\beta_a \geq a-n$ define
\begin{equation}
 \xy
 (0,0)*{\stccbub{a}{\alpha}};
 (0,9)*{n};
 \endxy \quad := \qquad \quad
   \xy
 (0,0)*{\includegraphics[scale=0.5]{figs/ccbub.eps}};
(-8,9)*{n}; (3,-8)*{a}; (-8,0)*{\bigb{\pi_{\alpha-(n+a)}}};
 \endxy
   \qquad \qquad \quad
 \xy
  (0,0)*{\stcbub{a}{\beta}};
 (0,9)*{n};
 \endxy \qquad := \qquad \quad
   \xy
 (0,0)*{\reflectbox{\includegraphics[scale=0.5]{figs/ccbub.eps}}};
(8,9)*{n}; (3,-8)*{a}; (8,0)*{\bigb{\pi_{\beta+(n-a)}}};
 \endxy
 \end{equation}
so that
\begin{equation}
\deg \left(\quad
 \xy
 (0,0)*{\stccbub{a}{\alpha}};
 (0,9)*{n};
 \endxy \right) \;\; =\;\; \deg\left( \pi_{\alpha}\right),\qquad \qquad
 \deg \left(\;
 \xy
  (0,0)*{\stcbub{a}{\beta}};
 (0,9)*{n};
 \endxy \quad \right) \;\; =\;\; \deg \left( \pi_{\beta}\right).
 \end{equation}
In terms of the thin calculus we have
\begin{equation}
 \xy
  (0,0)*{\stcbub{a}{\beta}};
 (0,9)*{n};
 \endxy
\quad \refequal{\eqref{eq_schur_thin}}  \quad \xy
 (0,0)*{\includegraphics[scale=0.45]{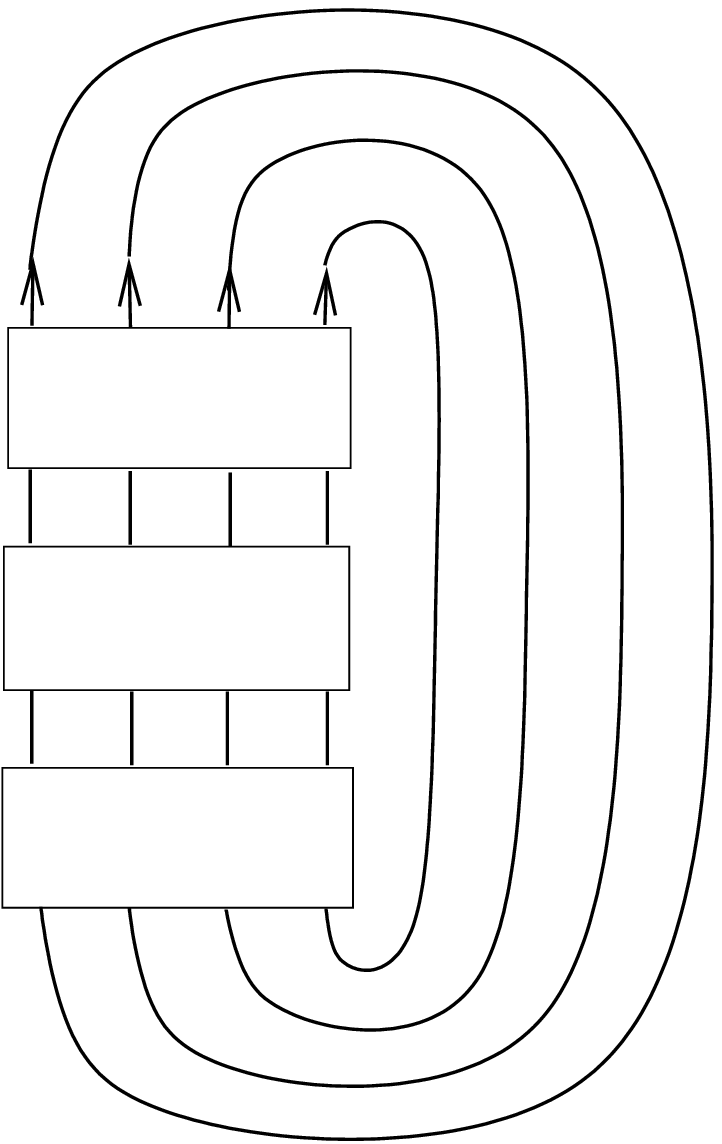}};
 (-8,8)*{e_a};(-8,-12)*{e_a};(-8,-2)*{x^{\beta+(n-a)}};
 \endxy
\quad =  \quad \xy
 (0,0)*{\includegraphics[scale=0.45]{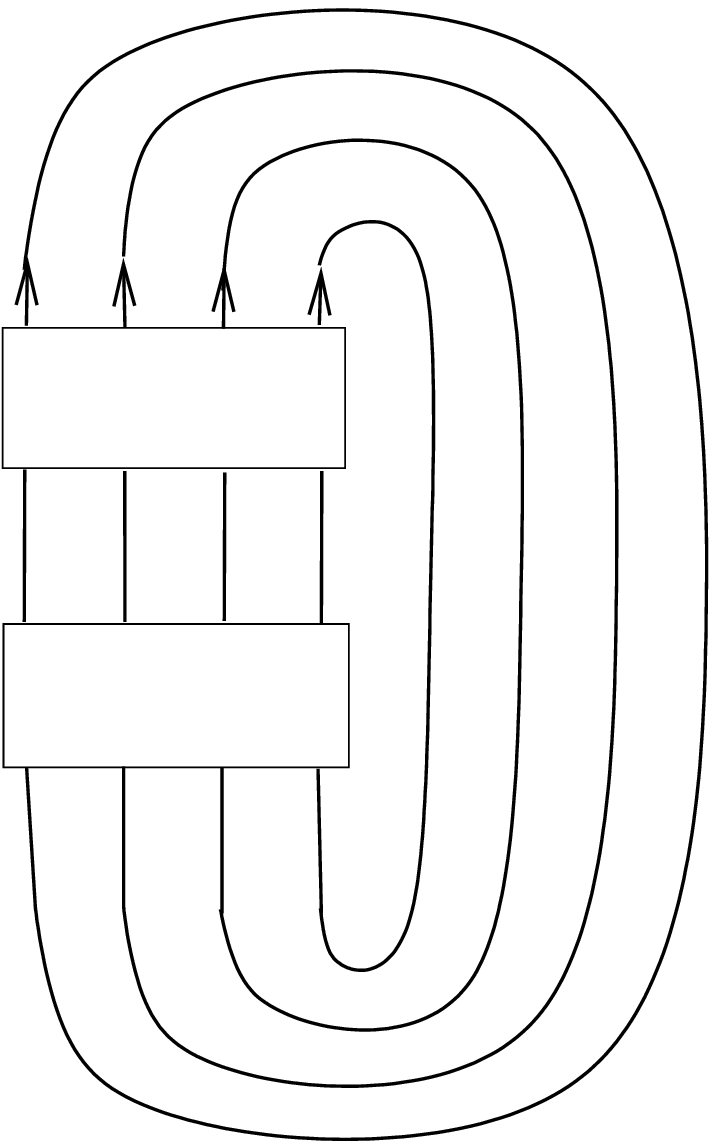}};
 (-8,-6)*{e_a};(-8,8)*{x^{\beta+(n-a)}};
 \endxy
\end{equation}
\begin{equation}
\quad =  \quad \xy
 (0,0)*{\includegraphics[scale=0.45]{figs/box-loop.eps}};
 (-22.5,5)*{\bullet}+(-2.5,1)*{\scs \beta_1'};
 (-16,5)*{\bullet}+(-2.5,1)*{\scs \beta_2'};
 (-9,5)*{\bullet}+(-2.5,1)*{\scs \beta_3'};(-4.5,-5)*{\cdots};
 (-0,5)*{\bullet}+(-2.5,1)*{\scs \beta_a'};
 (-12,0)*{D_a};
 \endxy
\end{equation}
where $\beta_j' := \beta_j+n-j$.  In the second equality above the
top idempotent $e_a$ is slid around the loop and then the identity
$e_a^2=e_a$ is used.  For the last equality the dots from the
idempotent $e_a$ are combined with the dots appearing in the
monomial $x^{\beta+(n-a)}$.  We have
a similar expression for the counter-clockwise thick dotted bubble.

For partitions $\alpha,\beta \in P(a)$ with $\alpha_a < n+a$ or $\beta_a < a-n$ some of the entries in $\alpha^{\spadesuit}$ or $\beta^{\spadesuit}$ will be negative.  By analogy with the fake bubbles in the thin calculus, when this happens, the thick bubbles labelled by a spaded Schur polynomial are regarded as formal symbols defined as follows:
\begin{equation} \label{eq_ccbub_det}
 \xy
 (0,0)*{\stccbub{a}{\alpha}};
 (0,9)*{n};
 \endxy \quad := \quad (-1)^{\frac{a(a-1)}{2}}
\left|
\begin{array}{ccccc}
  \xy 0;/r.18pc/: (0,0)*{\ccbub{\spadesuit+\alpha_1}{}};(5,7)*{n}; \endxy &
  \xy 0;/r.18pc/: (0,0)*{\ccbub{\spadesuit+\alpha_1+1}{}};(5,7)*{n}; \endxy &
  \xy 0;/r.18pc/: (0,0)*{\ccbub{\spadesuit+\alpha_1+2}{}};(5,7)*{n}; \endxy & \cdots &
  \xy 0;/r.18pc/: (0,0)*{\ccbub{\spadesuit+\alpha_1+(a-1)}{}};(5,7)*{n}; \endxy \\ \\
 \xy 0;/r.18pc/: (0,0)*{\ccbub{\spadesuit+\alpha_2-1}{}};(5,7)*{n}; \endxy &
  \xy 0;/r.18pc/: (0,0)*{\ccbub{\spadesuit+\alpha_2}{}};(5,7)*{n}; \endxy &
  \xy 0;/r.18pc/: (0,0)*{\ccbub{\spadesuit+\alpha_2+1}{}};(5,7)*{n}; \endxy & \cdots &
  \xy 0;/r.18pc/: (0,0)*{\ccbub{\spadesuit+\alpha_2+(a-2)}{}};(5,7)*{n}; \endxy \\ \\
  \cdots & \cdots & \cdots & \cdots & \cdots \\ \\
 \xy 0;/r.18pc/: (0,0)*{\ccbub{\spadesuit+\alpha_a-a+1}{}};(5,7)*{n}; \endxy &
  \xy 0;/r.18pc/: (0,0)*{\ccbub{\spadesuit+\alpha_a-a+2}{}};(5,7)*{n}; \endxy &
  \xy 0;/r.18pc/: (0,0)*{\ccbub{\spadesuit+\alpha_a-a+3}{}};(5,7)*{n}; \endxy & \cdots &
  \xy 0;/r.18pc/: (0,0)*{\ccbub{\spadesuit+\alpha_a}{}};(5,7)*{n}; \endxy
\end{array}
\right|
\end{equation}

\begin{equation} \label{eq_cbub_det}
 \xy
 (0,0)*{\stcbub{a}{\beta}};
 (0,9)*{n};
 \endxy \quad := \quad (-1)^{\frac{a(a-1)}{2}}
\left|
\begin{array}{ccccc}
  \xy 0;/r.18pc/: (0,0)*{\cbub{\spadesuit+\beta_1}{}};(5,7)*{n}; \endxy &
  \xy 0;/r.18pc/: (0,0)*{\cbub{\spadesuit+\beta_1+1}{}};(5,7)*{n}; \endxy &
  \xy 0;/r.18pc/: (0,0)*{\cbub{\spadesuit+\beta_1+2}{}};(5,7)*{n}; \endxy & \cdots &
  \xy 0;/r.18pc/: (0,0)*{\cbub{\spadesuit+\beta_1+(a-1)}{}};(5,7)*{n}; \endxy \\ \\
 \xy 0;/r.18pc/: (0,0)*{\cbub{\spadesuit+\beta_2-1}{}};(5,7)*{n}; \endxy &
  \xy 0;/r.18pc/: (0,0)*{\cbub{\spadesuit+\beta_2}{}};(5,7)*{n}; \endxy &
  \xy 0;/r.18pc/: (0,0)*{\cbub{\spadesuit+\beta_2+1}{}};(5,7)*{n}; \endxy & \cdots &
  \xy 0;/r.18pc/: (0,0)*{\cbub{\spadesuit+\beta_2+(a-2)}{}};(5,7)*{n}; \endxy \\ \\
  \cdots & \cdots & \cdots & \cdots & \cdots \\ \\
 \xy 0;/r.18pc/: (0,0)*{\cbub{\spadesuit+\beta_a-a+1}{}};(5,7)*{n}; \endxy &
  \xy 0;/r.18pc/: (0,0)*{\cbub{\spadesuit+\beta_a-a+2}{}};(5,7)*{n}; \endxy &
  \xy 0;/r.18pc/: (0,0)*{\cbub{\spadesuit+\beta_a-a+3}{}};(5,7)*{n}; \endxy & \cdots &
  \xy 0;/r.18pc/: (0,0)*{\cbub{\spadesuit+\beta_a}{}};(5,7)*{n}; \endxy
\end{array}
\right|
\end{equation}

The following Theorem shows that the above formulas remain valid even when the thick bubbles are not fake:

\begin{thm} \label{thm_thickbub_det}
For $\alpha,\beta \in P(a)$ with $\alpha_a \geq n+a$ and $\beta_a \geq a-n$ equations \eqref{eq_ccbub_det} and \eqref{eq_cbub_det} hold.  In other words, equations \eqref{eq_ccbub_det} and \eqref{eq_cbub_det} hold for all $\alpha$ and $\beta$.
\end{thm}

\begin{proof}
The proof is by induction on the thickness $a$ of the bubble.  The
base case is trivial.  Assume the result holds for bubbles of
thickness $a$, we will show the result holds for thickness $a+1$.
Let $\alpha=(x,\alpha_1,\dots,\alpha_a) \in P(a+1)$ be a partition with $\alpha_a\geq n+a+1$, so that
$\pi_{\alpha}^{\spadesuit} = \pi_{\alpha-n-a-1}$ were $\alpha-n-a-1$ is the partition $(x-n-a-1, \alpha_1-n-a-1, \dots ,
\alpha_a-n-a-1)$. First split the thickness $a+1$ line into a
thickness $a$ line and a thickness 1 line using
\eqref{eq_schur_left}:
\begin{equation}
   \xy
 (0,0)*{\includegraphics[scale=0.5]{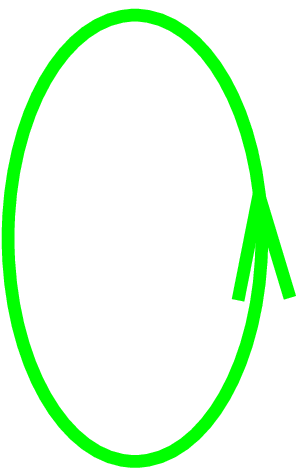}};
 (-5,-13)*{a+1};(4,-7)*{\bigb{\pi_{\alpha}^{\spadesuit}}}; (9,9)*{n};
  \endxy
 \;\; \refequal{\eqref{eq_schur_left}}\;\;
     \xy
 (0,0)*{\includegraphics[scale=0.5]{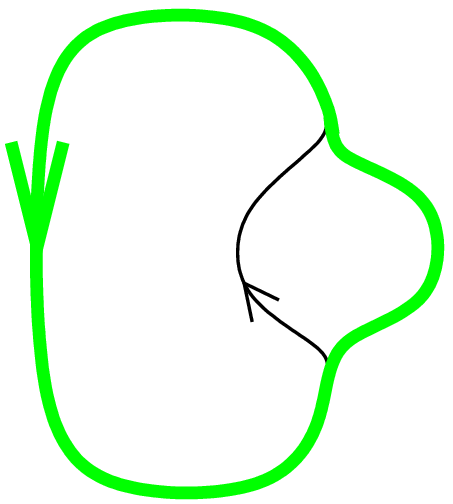}};
 (-9,-13)*{a+1};(12,0)*{\bigb{\pi_{\beta}}}; (9,9)*{n};
 (.5,1)*{\bullet}+(-3,3)*{\scs x-n-1};(9,-5)*{a};
  \endxy
 \qquad  =  \;  \xy
 (0,0)*{\includegraphics[scale=0.5]{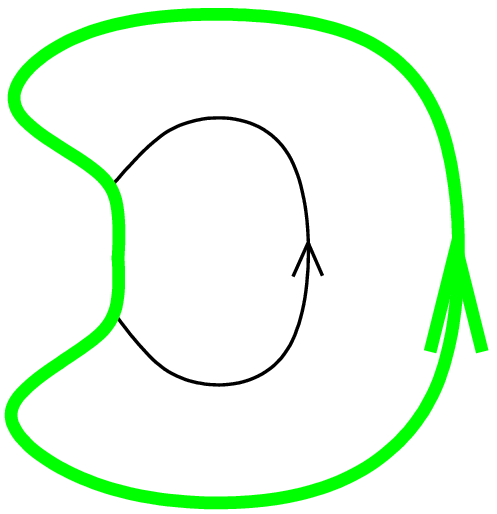}};
 (11,4)*{\bigb{\pi_{\beta}}}; (15,9)*{n};
 (2,5)*{\bullet}+(1,3)*{\scs x-n-1};(7,-12)*{a};
  \endxy
   \quad \refequal{\eqref{eq_triangle}}\quad
   \xy
 (0,0)*{\includegraphics[scale=0.5]{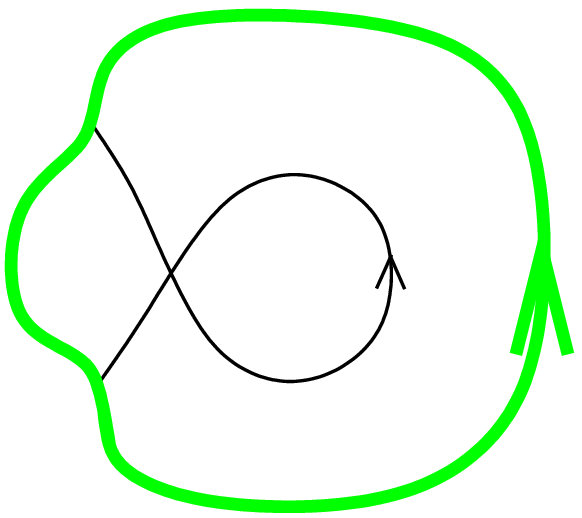}};
 (12,4)*{\bigb{\pi_{\beta}}}; (15,9)*{n};
 (2,4)*{\bullet}+(1,3)*{\scs x-n-1};(9,-12)*{a};
  \endxy
\end{equation}
where $\beta= (\alpha_1-n-a-1, \dots,
\alpha_a-n-a-1)$.
\begin{equation}
  \quad =\quad
  \sum_{g_1+g_2=2a+x-1}
   \xy
 (0,0)*{\includegraphics[scale=0.5]{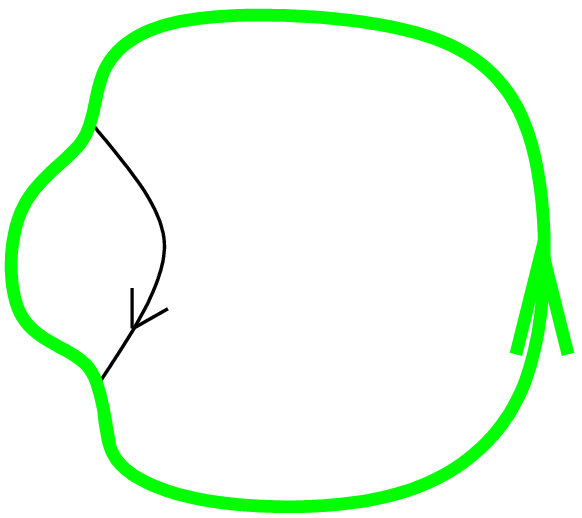}};
 (12,4)*{\bigb{\pi_{\beta}}}; (15,9)*{n};
 (-7,3)*{\bullet}+(1,3)*{\scs g_1};(9,-12)*{a};
 (2,0)*{\ccbub{\spadesuit+g_2}{}};
  \endxy
      \quad =\quad
  \sum_{g=0}^{x+a} \qquad \xy
 (0,0)*{\includegraphics[scale=0.5]{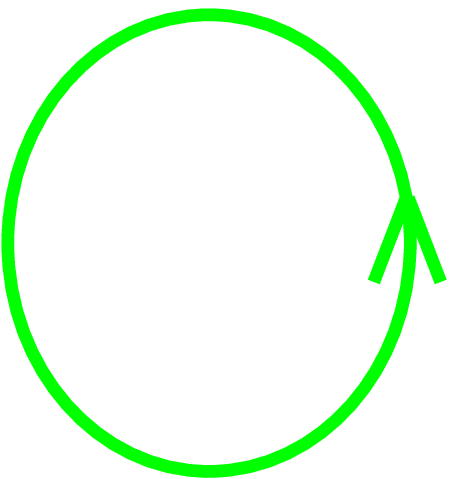}};
 (8,7)*{\bigb{\pi^{\spadesuit}_{\alpha'}}}; (15,9)*{n};
 (-12,-2)*{\bigb{h_{g}}};(9,-12)*{a};
 (-1,0)*{\ccbub{\spadesuit+a+x-g}{}};
  \endxy
\end{equation}
where $\alpha' = (\alpha_1-1, \dots , \alpha_a-1)$, so that the thick bubble labelled $\pi_{\alpha'}^{\spadesuit}$ is not fake.  In the first summation the diagram vanishes whenever $g_1 < a-1$. Removing those vanishing terms, re-indexing the summation to include only the nonzero terms, and using that $\pi_{(x+a-g, 0, \dots, 0)}=h_{x+a-g}$ gives the last equation.

Now slide the inner thin bubble outside the large thick bubble using
Corollary~\ref{cor_thinslide}.
\begin{equation}
  \;\; = \;\;
  \sum_{g=0}^{x+a} \sum_{m=0}^{2a}(-1)^m
  \sum_{ \xy  (0,3)*{\scs p,q \leq a}; (0,0)*{\scs p+q=m};\endxy}
   \xy
   (-24,0)*{\ccbub{\spadesuit+a+x-g-(p+q)}{}};
 (0,0)*{\includegraphics[scale=0.5]{figs/tbigloop.eps}};
 (-1,-13)*{a};(4,-7)*{\bigb{\pi_{\alpha'}^{\spadesuit}}}; (9,9)*{n};
 (-6,-8)*{\bigb{h_{g}}}; (-8,0)*{\bigb{\varepsilon_{q}}};  (-6,8)*{\bigb{\varepsilon_{p}}};
  \endxy
\end{equation}
\begin{equation}
  \;\; = \;\;
  \sum_{g=0}^{x+a} \sum_{m=0}^{\min(2a,x+a-g)}\sum_{p=0}^{\min(a,m)}(-1)^m
   \xy
   (-24,0)*{\ccbub{\spadesuit+a+x-g-m}{}};
 (0,0)*{\includegraphics[scale=0.5]{figs/tbigloop.eps}};
 (-1,-13)*{a};(4,-7)*{\bigb{\pi_{\alpha'}^{\spadesuit}}}; (9,9)*{n};
 (-6,-8)*{\bigb{h_{g}}}; (-8,0)*{\bigb{\varepsilon_{m-p}}};  (-6,8)*{\bigb{\varepsilon_{p}}};
  \endxy
\end{equation}
After switching the indices by letting $s=g+m$ the above
summation becomes
\begin{equation}
  \;\; = \;\;
  \sum_{s=0}^{x+a} \sum_{p=0}^{a}\sum_{m=p}^{\min(2a,s)}(-1)^m
   \xy
   (-24,0)*{\ccbub{\spadesuit+a+x-s}{}};
 (0,0)*{\includegraphics[scale=0.5]{figs/tbigloop.eps}};
 (-1,-13)*{a};(4,-7)*{\bigb{\pi_{\alpha'}^{\spadesuit}}}; (9,9)*{n};
 (-6,-8)*{\bigb{h_{s-m}}}; (-8,0)*{\bigb{\varepsilon_{m-p}}};  (-6,8)*{\bigb{\varepsilon_{p}}};
  \endxy
\end{equation}
Note that
\begin{eqnarray}
  \sum_{p=0}^a \sum_{m=p}^{\min(2a,s)}(-1)^m \varepsilon_p \varepsilon_{m-p}h_{s-m} &=&
  \sum_{p=0}^a(-1)^p \varepsilon_p \sum_{m=0}^{\min(s-p,2a-p)}(-1)^m \varepsilon_m
  h_{s-p-m} \nn
 \\  &=&
  \sum_{a=0}^p(-1)^p \varepsilon_p \delta_{s-p,0} = (-1)^s \varepsilon_s \nn
\end{eqnarray}
and since we are working with elementary symmetric functions in $a$
variables we have that the $\varepsilon_s=0$ for $s>a$.  Hence, we have shown
\begin{equation}
     \xy
 (0,0)*{\includegraphics[scale=0.5]{figs/tbigloop.eps}};
 (-5,-13)*{a+1};(4,-7)*{\bigb{\pi_{\alpha}^{\spadesuit}}}; (9,9)*{n};
  \endxy \;\; = \;\;
  \sum_{s=0}^{a}(-1)^s
   \xy
   (-24,0)*{\ccbub{\spadesuit+a+x-s}{}};
 (0,0)*{\includegraphics[scale=0.5]{figs/tbigloop.eps}};
 (-1,-13)*{a};(4,-7)*{\bigb{\pi_{\alpha'}^{\spadesuit}}}; (9,9)*{n};
 (-8,0)*{\bigb{\varepsilon_{s}}};
  \endxy
\end{equation}
and using the induction hypothesis, the above summation is the
expansion of the required determinant along the top row.
\end{proof}

Let $\Lambda(\underline{x})$ denote the graded ring of symmetric functions in infinitely many variables $\underline{x}=(x_1,x_2, \dots)$. This ring can be identified with the polynomial ring $\Z[h_1,h_2,\dots]$ and with the polynomial ring $\Z[\varepsilon_1,\varepsilon_2,\dots]$.  It was shown in \cite[Proposition 8.2]{Lau1} that there is an isomorphism
\begin{eqnarray}
 \phi^n \maps \Lambda(\underline{x}) &\to& \HOM_{\Ucat}(\onen,\onen),
  \nn \\ \label{eq_isom_h}
  h_r(\und{x}) & \mapsto & \xy
 (0,0)*{\cbub{\spadesuit+r}{}};
 (6,8)*{n};
 \endxy,
\end{eqnarray}
that we denote here by $\phi^n$.  This isomorphism can alternatively be described using counter-clockwise oriented bubbles
\begin{eqnarray}
 \phi^n \maps \Lambda(\underline{x}) &\to& \HOM_{\Ucat}(\onen,\onen),
 \nn \\ \label{eq_isom_e}
 (-1)^r\varepsilon_r(\und{x}) & \mapsto &   \xy
 (0,0)*{\ccbub{\spadesuit+r}{}};
 (6,8)*{n};
 \endxy,
\end{eqnarray}
so that the homogeneous term of degree $r$ in the infinite
Grassmannian relation \eqref{eq_infinite_Grass} becomes \eqref{eq_eh_rel}.

\begin{prop} \label{prop_image_thickbub}
The equalities
  \begin{equation}
 (\phi^n)^{-1}\left( \; \xy
  (0,0)*{\stcbub{a}{\alpha}};
 (0,9)*{n};
 \endxy\quad\right) \;\; = \;\;
(-1)^{\frac{a(a-1)}{2}}\pi_{\alpha}(\underline{x}),
   \qquad \qquad
 (\phi^n)^{-1}\left( \quad\xy
 (0,0)*{\stccbub{a}{\beta}};
 (0,9)*{n};
 \endxy \;
 \right)
 \;\; = \;\; (-1)^{\frac{a(a-1)}{2}+|\beta|}
 \pi_{\overline{\beta}}(\underline{x}).
 \end{equation}
hold provided $\alpha,\beta\in P(a)$.
\end{prop}

\begin{proof}
The first equality follows immediately from the reduction of thick bubbles to thin bubbles given in Theorem~\ref{thm_thickbub_det} using the isomorphism $\phi^n$. For the second equality observe that
\begin{equation}
  (\phi^n)^{-1}\left(
   \xy
  (0,0)*{\ccbub{\spadesuit+x}{}};
 \endxy
  \right) \;\; =\;\; (-1)^x \varepsilon_x(\underline{z}).
\end{equation}
Reducing the thick to thin bubbles via \eqref{eq_ccbub_det} the second equality follows since
 \begin{eqnarray}
 (-1)^{\frac{a(a-1)}{2}}\left|
 \begin{array}{cccc}
   (-1)^{\beta_1}\varepsilon_{\beta_1} & (-1)^{\beta_1+1}\varepsilon_{\beta_1+1} & \cdots & (-1)^{\beta_1+i-1}\varepsilon_{\beta_1+i-1} \\
   \vdots &  &  & \vdots \\
   (-1)^{\beta_a-a+1}\varepsilon_{\beta_a-a+1} & (-1)^{\beta_a-a+2}\varepsilon_{\beta_a-a+2} & \cdots & (-1)^{\beta_a}\varepsilon_{\beta_a}
 \end{array}
 \right| \nn \hspace{0.8in} \\
 = (-1)^{\frac{a(a-1)}{2}}\det\left[
 (-1)^{\beta_s+t-s} \varepsilon_{\beta_s+t-s}
 \right]_{s,t=1}^{a} = (-1)^{\frac{a(a-1)}{2}}(-1)^{|\gamma|}\det\left[ (-1)^{t-s}\varepsilon_{\beta_s+t-s}\right]_{s,t=1}^{a} \nn
 \end{eqnarray}
and multiplying every even row and column by $(-1)$ in the last equality gives the determinant formula \eqref{eq_Schur_e} for the Schur polynomial $\pi_{\overline{\beta}}$.
\end{proof}

\begin{thm}[Higher infinite Grassmannian relations]
For $\alpha\in P(a)$ the equality
\begin{equation}
 \sum_{\beta,\gamma\in P(a)}
 c_{\beta,\gamma}^{\alpha}\;\;\;\;\;
 \xy
 (-6,0)*{\stccbub{a}{\beta}};
 (7,0)*{\stcbub{a}{\gamma}};
 (0,9)*{n};
 \endxy \qquad = \quad \delta_{\alpha,0}
 \end{equation}
 holds in $\UcatD$.
\end{thm}

\begin{proof}
Applying the isomorphism $\phi^n$ the theorem becomes the statement
\begin{equation}
  \sum_{\beta,\gamma} (-1)^{|\beta|} c_{\beta,\gamma}^{\alpha} \pi_{\bar{\beta}}(\und{x})\pi_{\gamma}(\und{x}) = \delta_{\alpha,0}.
\end{equation}
This is proven in the following proposition.
\end{proof}

\begin{prop}
For every partition $\alpha$ we have
\begin{equation}\label{hgr}
\sum_{\beta,\gamma} {(-1)^{|\beta|} c_{\beta,\gamma}^{\alpha} \pi_{\bar{\beta}}(\underline{x}) \pi_{\gamma}(\underline{x})} = \delta_{\alpha,\emptyset}.
\end{equation}
\end{prop}

\begin{proof}
The multiplication in the ring of symmetric functions $\Lambda(\und{x})$ can be enhanced to a Hopf algebra structure~\cite{Gei,Zel,FJ}.  The comultiplication is given by taking a function $f(\und{x}) \in \Lambda(\und{x})$ and writing it as a function of two sets of variables $x \otimes 1$ and $1 \otimes x$:
\begin{equation}
  \Delta(f) = f(\und{x} \otimes 1, 1 \otimes \und{x}).
\end{equation}
For instance,
\begin{equation}
  \Delta(\varepsilon_n) = \sum_{i=0}^n \varepsilon_i \otimes \varepsilon_{n-i}.
\end{equation}
In the basis of Schur functions the comultiplication is given by
\begin{equation}
  \Delta(\pi_{\alpha}(\und{x})) =
  \sum c_{\beta \gamma}^{\alpha} \pi_{\beta}(\und{x}) \otimes \pi_{\gamma}(\und{x}).
\end{equation}

The ring of symmetric functions has a natural symmetric bilinear form
\begin{equation}
  \langle , \rangle \maps \Lambda(\und{x}) \otimes \Lambda(\und{x}) \to \Bbbk
\end{equation}
in which the basis of Schur functions is orthonormal
\begin{equation}
  \langle \pi_{\alpha}(\und{x}), \pi_{\beta}(\und{x}) \rangle
  = \delta_{\alpha,\beta}.
\end{equation}
The multiplication and comultiplication are adjoint relative to this bilinear form
\begin{equation}
  \langle f, gh \rangle =
  \langle \Delta(f), g\otimes h \rangle.
\end{equation}
The counit of $\Lambda(\und{x})$ is
\begin{equation}
  \varepsilon(\pi_{\alpha}(\und{x})) = \delta_{\alpha, \emptyset}.
\end{equation}
The counit is nonzero only on the Schur function $\pi_{\emptyset}(\und{x})=1$ of the empty partition.  The unit map is
\begin{equation}
  \iota \maps \Bbbk \to \Lambda(\und{x}), \quad 1 \mapsto 1.
\end{equation}
The antipode is given by
\begin{equation}
  S(\pi_{\alpha}(\und{x})) = (-1)^{|\alpha|}\pi_{\bar{\alpha}}(\und{x}),
\end{equation}
where $\bar{\alpha}$ is the conjugate partition.

One of the Hopf algebra axioms is
\begin{equation}
  M(S \otimes I) \Delta = \iota \varepsilon
\end{equation}
(see \cite{Kup1,Kup2} for a diagrammatic interpretation of Hopf algebra axioms with applications to 3-manifold invariants).
Applying this to $\pi_{\alpha}(\und{x})$ we get that
\begin{equation}
  \xy
  (-57,0)*+{\pi_{\alpha}(\und{x})}="1";
  (-30,0)*+{\sum c_{\beta,\gamma}^{\alpha}\pi_{\beta}(\und{x})\otimes \pi_{\gamma}(\und{x})}="2";
  (17,0)*+{\sum(-1)^{|\beta|} c_{\beta,\gamma}^{\alpha}\pi_{\bar{\beta}}(\und{x})\otimes \pi_{\gamma}(\und{x})}="3";
  (65,0)*+{\sum(-1)^{|\beta|} c_{\beta,\gamma}^{\alpha}\pi_{\bar{\beta}}(\und{x}) \pi_{\gamma}(\und{x})}="4";
  {\ar^-{\Delta} "1";"2"};
  {\ar^-{S \otimes I} "2";"3"};
  {\ar^-{M} "3";"4"};
  \endxy
\end{equation}
must be equal to
\begin{equation}
  \xy
  (-20,0)*+{\pi_{\alpha}(\und{x})}="1";
  (-0,0)*+{\delta_{\alpha,\emptyset}}="2";
  (20,0)*+{\delta_{\alpha,\emptyset} \cdot 1}="3";
  {\ar^-{\varepsilon} "1";"2"};
  {\ar^-{\iota} "2";"3"};
  \endxy
\end{equation}
finishing the proof.
\end{proof}

Under the isomorphism $\Lambda(\und{x})\cong \END_{\UcatD}(\onen)$ the antipode $S$ of the algebra of symmetric functions corresponds to the reversal of orientation in the closed diagrams: compare the equation $S(h_r)(\und{x})=(-1)^r\varepsilon_r(\und{x})$, with \eqref{eq_isom_h} and \eqref{eq_isom_e}: and, more generally, see Proposition~\ref{prop_image_thickbub}.

 \subsection{Thick bubble slides and some key lemmas}

The following lemmas are needed for sliding thick bubbles through thick lines and they are also used in Section~\ref{sec_decomp_EaFb} for the most difficult decomposition result.

\begin{lem}
\begin{eqnarray} \label{eq_lem1}
\sum_{i+j+k=M} \xy
  (15,8)*{n};
  (0,0)*{\bbe{}};
  (12,-2)*{\cbub{\spadesuit+k}{}};
  (0,6)*{\bullet}+(5,-1)*{\scs i+j+x};
 \endxy
  & =&
     \xy
  (0,8)*{n+2};
  (12,0)*{\bbe{}};
  (12,-12)*{\scs j};
  (0,-2)*{\cbub{\spadesuit+M}{}};
  (12,6)*{\bullet}+(3,-1)*{\scs x};
 \endxy
 \end{eqnarray}
\end{lem}

\begin{proof}
This is just a restatement of \eqref{eq_bubslide2}.  The coefficient appearing in \eqref{eq_bubslide2} arises from the additional summation indices.
\end{proof}

\begin{lem}[Central Element Lemma] \label{lem_ladder_slide1}
\begin{equation}
\sum_{i+j+k=M}
 \xy
 (0,0)*{\includegraphics[scale=0.5]{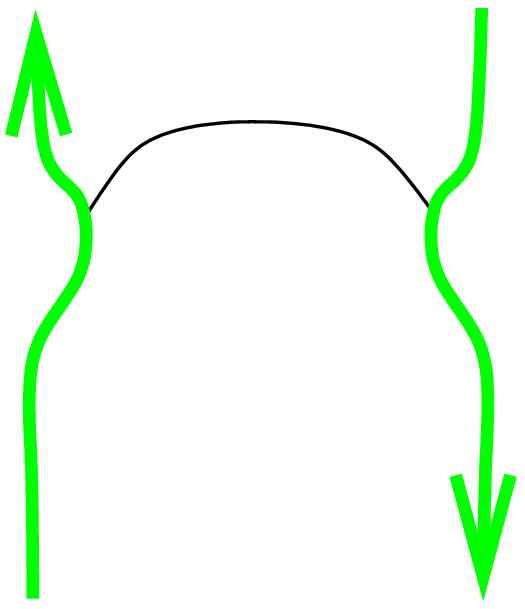}};
 (-12,-4)*{\bigb{h_i}}; (12,-4)*{\bigb{h_j}}; (0,-6)*{\cbub{\spadesuit+k}{}};
 (-14,-12)*{a};(15,-12)*{b};(-18,12)*{a-1};(18,12)*{b-1};
  (0,9)*{\bullet}+(0,3)*{x}; (19,0)*{n};
  \endxy
 \quad = \quad
 \sum_{i+j+k=M}
 \xy
 (0,0)*{\includegraphics[scale=0.5]{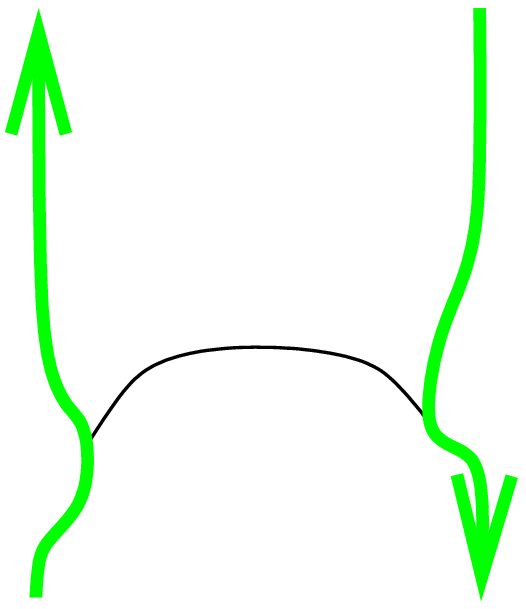}};
 (-12,4)*{\bigb{h_i}}; (12,4)*{\bigb{h_j}}; (-2,8)*{\cbub{\spadesuit+k}{}};
(-14,-12)*{a};(15,-12)*{b};(-18,12)*{a-1};(18,12)*{b-1};
  (2,-2.5)*{\bullet}+(0,-3)*{x}; (19,0)*{n};
  \endxy
\end{equation}
Recall that $h_i=\pi_{i,0,\dots,0}$ is the $i$th complete symmetric
polynomial.
\end{lem}

\begin{proof}
We use
\begin{equation}
   \xy
 (0,0)*{\includegraphics[scale=0.5]{figs/tonesplit.eps}};
 (-3,-11)*{a};(-11,8)*{a-1};(8,8)*{};
  (0,-7)*{\bigb{h_i}};
  \endxy
  \quad =
 \quad
 \sum_{\ell+\ell'=i} \;
    \xy
 (0,0)*{\includegraphics[scale=0.5]{figs/tonesplit.eps}};
 (-3,-11)*{a};(-11,8)*{a-1};(8,8)*{};
 (-5,2)*{\bigb{h_{\ell}}};(5,2)*{\bullet}+(3,0)*{\ell'};
  \endxy
\end{equation}
and a similar statement for oppositely oriented lines of thickness $b$, so that

\begin{eqnarray}
  \sum_{i+j+k=M}
 \xy
 (0,0)*{\includegraphics[scale=0.5]{figs/marko-lem2-1}};
 (-12,-4)*{\bigb{h_i}}; (12,-4)*{\bigb{h_j}}; (0,-6)*{\cbub{\spadesuit+k}{}};
 (-14,-12)*{a};(15,-12)*{b};(-18,12)*{a-1};(18,12)*{b-1};
  (0,9)*{\bullet}+(0,3)*{x}; (19,0)*{n};
  \endxy & = &
  \sum_{\ell+\ell'+m+m'+k=M}
 \xy
 (0,0)*{\includegraphics[scale=0.5]{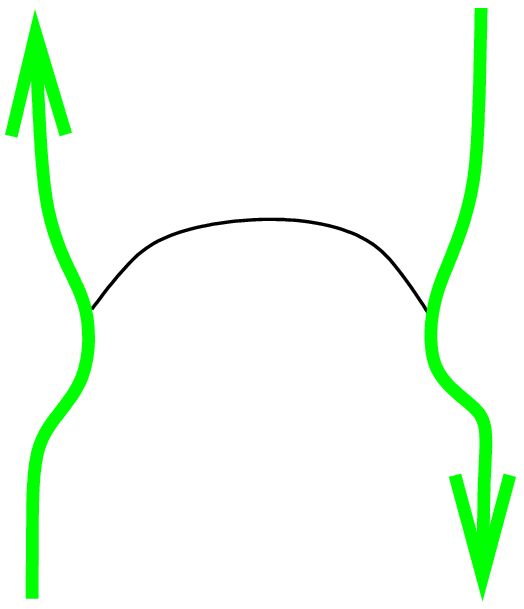}};
 (-12.5,5)*{\bigb{h_{\ell'}}}; (12.5,5)*{\bigb{h_{m'}}}; (0,-6)*{\cbub{\spadesuit+k}{}};
 (-14,-12)*{a};(15,-12)*{b};(-18,12)*{a-1};(18,12)*{b-1};
  (0,4.3)*{\bullet}+(0,3)*{x+\ell+m}; (19,0)*{n};
  \endxy \nn
\end{eqnarray}
\begin{eqnarray}
 & = &
  \sum_{M'=0}^M \sum_{\ell'+m'=M'}\sum_{\ell+m+k=M-M'}
 \xy
 (0,0)*{\includegraphics[scale=0.5]{figs/marko-lem2-1p}};
 (-12.5,5)*{\bigb{h_{\ell'}}}; (12.5,5)*{\bigb{h_{m'}}}; (0,-6)*{\cbub{\spadesuit+k}{}};
 (-14,-12)*{a};(15,-12)*{b};(-18,12)*{a-1};(18,12)*{b-1};
  (0,4.3)*{\bullet}+(0,3)*{x+\ell+m}; (19,0)*{n};
  \endxy \nn \\
  & \refequal{\eqref{eq_lem1}} & \sum_{M'=0}^M \sum_{\ell'+m'=M'}
   \xy
 (0,0)*{\includegraphics[scale=0.5]{figs/marko-lem2-2}};
 (-12,4)*{\bigb{h_{\ell'}}}; (12,4)*{\bigb{h_{m'}}}; (-1,8)*{\cbub{\spadesuit+M-M'}{}};
(-14,-12)*{a};(15,-12)*{b};(-18,12)*{a-1};(18,12)*{b-1};
  (2,-2.5)*{\bullet}+(0,-3)*{x}; (19,0)*{n};
  \endxy
\end{eqnarray}
which after reindexing completes the proof.
\end{proof}

We can restate Lemma~\ref{lem_ladder_slide1} in alternative form
using the injective homomorphism
\begin{eqnarray} \label{eq_EcFdgamma_isom}
  \phi_{a,b,c}^n \maps
  \Z[x_1, \dots x_a, y_1, \dots, y_{b}]^{S_a \times S_b} \otimes \Lambda(\und{z})
   & \to& \HOM_{\UcatD}(\cal{E}^{(a)} \cal{F}^{(b)}\onen,\cal{E}^{(a)} \cal{F}^{(b)}\onen)
   \nn \\
  \pi_{\alpha}(\underline{x})
  \pi_{\beta}(\underline{y}) \pi_{\gamma}(\underline{z})& \mapsto &
  (-1)^{\frac{c(c-1)}{2}} \xy
 (-12,0)*{\includegraphics[scale=0.5]{figs/tlong-up.eps}};
 (-15,-11)*{a};(-12,-2)*{\bigb{\pi_{\alpha}}};
 (-2,0)*{\stcbub{c}{\gamma}};
 (12,0)*{\includegraphics[angle=180, scale=0.5]{figs/tlong-up.eps}};
 (15,-11)*{b};(12,-2)*{\bigb{\pi_{\beta}}};  (19,8)*{n};
  \endxy
\end{eqnarray}
where $\gamma\in P(c)$.   For all $a,b,c \geq 0$ we write
\begin{eqnarray} \label{eq_phiofh}
  \phi_{a,b,c}^n(h_i(\underline{x},\underline{y},\underline{z}))
  \quad := \quad (-1)^{\frac{c(c-1)}{2}}\xy
 (-12,0)*{\includegraphics[scale=0.5]{figs/tlong-up.eps}};
  (12,0)*{\includegraphics[angle=180, scale=0.5]{figs/tlong-up.eps}};
 (15,-11)*{b}; (-15,-11)*{a};
   (0,0)*{\bigb{\hspace{0.3in}h_i(a,b,c)\hspace{0.3in}}};
   (19,8)*{n};
  \endxy
\end{eqnarray}
where $\underline{x}=(x_1,x_2,\ldots,x_a)$,
$\underline{y}=(y_1,y_2,\ldots,y_b)$, and
$\underline{z}=(z_1,z_2,\ldots)$.  Then
Lemma~\ref{lem_ladder_slide1} can be restated as the equality in
$\HOM_{\UcatD}(\cal{E}^{(a)}\cal{F}^{(b)}\onen,
\cal{E}^{(a-1)}\cal{F}^{(b-1)}\onen)$
\begin{equation} \label{eq_alternate_central_element}
 \xy
 (0,0)*{\includegraphics[scale=0.5]{figs/marko-lem2-1}};
 (-14,-12)*{a};(15,-12)*{b};(-18,12)*{a-1};(18,12)*{b-1};
  (0,-4)*{\bigb{\hspace{0.3in}h_i(a,b,1)\hspace{0.3in}}};
  (0,9)*{\bullet}+(0,3)*{x}; (19,0)*{n};
  \endxy
 \quad = \quad
 \xy
 (0,0)*{\includegraphics[scale=0.5]{figs/marko-lem2-2}};
 (0,4)*{\bigb{\hspace{0.1in}h_i(a-1,b-1,1)\hspace{0.1in}}};
(-14,-12)*{a};(15,-12)*{b};(-18,12)*{a-1};(18,12)*{b-1};
  (2,-2.5)*{\bullet}+(0,-3)*{x}; (19,0)*{n};
  \endxy
\end{equation}

A function $p=p(\underline{x},\underline{y}, \underline{z}) \in
\Z[x_1, \dots x_a, y_1, \dots, y_{b}]^{S_a\times S_b} \otimes \Lambda(\und{z})$ that is symmetric in all variables $\und{x}, \und{y}, \und{z}$ can be depicted as in \eqref{eq_phiofh}
\begin{eqnarray}
  \phi_{a,b,c}^n(p(\underline{x},\underline{y}, \underline{z}))
  \quad := \quad (-1)^{\frac{c(c-1)}{2}} \xy
 (-12,0)*{\includegraphics[scale=0.5]{figs/tlong-up.eps}};
  (12,0)*{\includegraphics[angle=180, scale=0.5]{figs/tlong-up.eps}};
 (15,-11)*{b}; (-15,-11)*{a};
   (0,0)*{\bigb{\hspace{0.3in}p(a,b,c)\hspace{0.3in}}};
   (19,8)*{n};
  \endxy
\end{eqnarray}
For simplicity we sometimes label the diagram for the image of the symmetric function $p$ as  $p(c)$
\begin{eqnarray}
\xy
 (-12,0)*{\includegraphics[scale=0.5]{figs/tlong-up.eps}};
  (12,0)*{\includegraphics[angle=180, scale=0.5]{figs/tlong-up.eps}};
 (15,-11)*{b}; (-15,-11)*{a};
   (0,0)*{\bigb{\hspace{0.3in}p(a,b,c)\hspace{0.3in}}};
   (19,8)*{n};
  \endxy
 \;\; =\;\;  \xy
 (-12,0)*{\includegraphics[scale=0.5]{figs/tlong-up.eps}};
  (12,0)*{\includegraphics[angle=180, scale=0.5]{figs/tlong-up.eps}};
 (15,-11)*{b}; (-15,-11)*{a};
   (0,0)*{\bigb{\hspace{0.45in}p(c)\hspace{0.45in}}};
   (19,8)*{n};
  \endxy \nn \\
\end{eqnarray}
since the indices $a$ and $b$ can be determined from the thick lines intersecting the coupon.  Any symmetric function $p \in \Lambda(\und{x},\und{y},\und{z})$ can be expressed via generators $h_i(\underline{x},\underline{y},\underline{z})$.  The following lemma is immediate.

\begin{lem}[General Central Element Lemma] \label{lem_gen_central}
If $p=p(\underline{x},\underline{y}, \underline{z})$ is symmetric in
all variables, then
\begin{equation}
 \xy
 (0,0)*{\includegraphics[scale=0.5]{figs/marko-lem2-1}};
 (-14,-12)*{a};(15,-12)*{b};(-18,12)*{a-1};(18,12)*{b-1};
  (0,-4)*{\bigb{\hspace{0.5in}p(c)\hspace{0.5in}}};
  (0,9)*{\bullet}+(0,3)*{x}; (19,0)*{n};
  \endxy
 \quad = \quad
 \xy
 (0,0)*{\includegraphics[scale=0.5]{figs/marko-lem2-2}};
 (0,4)*{\bigb{\hspace{0.5in}p'(c)\hspace{0.5in}}};
(-14,-12)*{a};(15,-12)*{b};(-18,12)*{a-1};(18,12)*{b-1};
  (2,-2.5)*{\bullet}+(0,-3)*{x}; (19,0)*{n};
  \endxy
\end{equation}
where $p'(c)$ is the polynomial $p(\und{x}',\und{y}',\und{z})$ in variables $\und{x}'=(x_1,\dots,x_{a-1})$, $\und{y}'=(y_1,\dots,y_{b-1})$, $\und{z}=(z_1,\dots)$ obtained from $p(\und{x},\und{y},\und{z})$ by setting $x_a=y_b=0$.  In what follows we will write $p(c)$ rather than $p'(c)$ for simplicity.
\end{lem}

\begin{lem}[Square Flop] \label{lem_square_flop}
\begin{eqnarray}
 \xy
 (0,0)*{\includegraphics[scale=0.45]{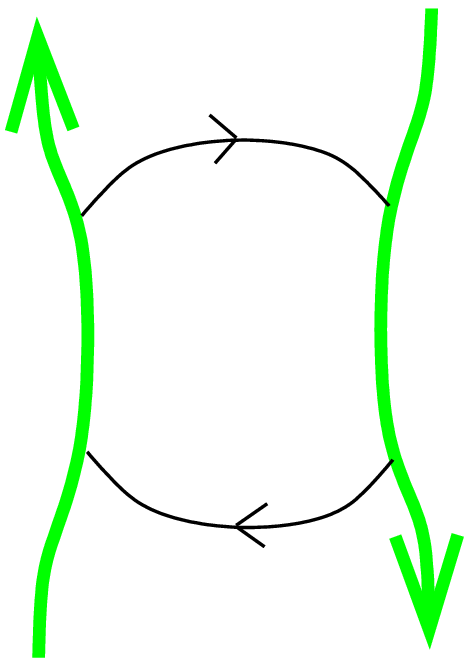}};
 (-13,-12)*{a};(14,-12)*{b};(-13,12)*{a};(14,12)*{b};
 (-13,2)*{a+1};(13,2)*{b+1};
  (3,8.5)*{\bullet}+(0,3)*{y};
  (3,-9)*{\bullet}+(0,-3)*{x}; (19,-8)*{n};
  \endxy
  &= & - \;\;
     \xy
 (0,0)*{\includegraphics[scale=0.45]{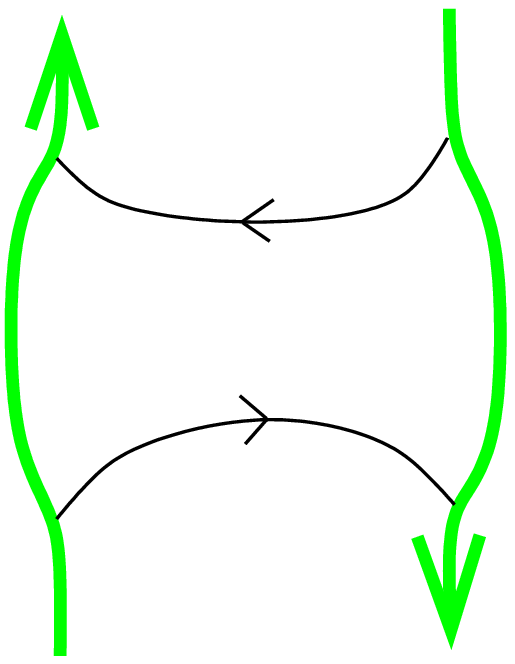}};
 (-13,-13)*{a};(14,-12)*{b};(-13,12)*{a};(13,13)*{b};
  (3,5)*{\bullet}+(0,3)*{x};
  (3,-4.5)*{\bullet}+(0,-3)*{y}; (19,-8)*{n};
  \endxy \hspace{-0.1in}+(-1)^{a-b}\hspace{-0.2in}\sum_{  \xy
  (0,-1)*{\scs p+q+r =};
  (0,-4)*{\scs x+y+b-a+1-n};
  \endxy}
  \xy
 (-10,0)*{\includegraphics[scale=0.5]{figs/tlong-up.eps}};
 (-12.5,-11)*{a};(-10,-2)*{\bigb{h_p}};
  (10,0)*{\includegraphics[angle=180, scale=0.5]{figs/tlong-up.eps}};
 (12.5,-11)*{b};(10,-2)*{\bigb{h_q}};
  (0,-2)*{\cbub{\spadesuit +r}{}}; (16,8)*{n};
  \endxy \nn \\ &=&
  - \;\;
     \xy
 (0,0)*{\includegraphics[scale=0.45]{figs/marko-lem3-2}};
 (-13,-13)*{a};(14,-12)*{b};(-13,12)*{a};(13,13)*{b};
  (3,5)*{\bullet}+(0,3)*{x};
  (3,-4.5)*{\bullet}+(0,-3)*{y}; (19,-8)*{n};
  \endxy
  \;\;+(-1)^{a-b}\;\;
\xy
 (-10,0)*{\includegraphics[scale=0.5]{figs/tlong-up.eps}};
  (10,0)*{\includegraphics[angle=180, scale=0.5]{figs/tlong-up.eps}};
 (15,-11)*{b}; (-15,-11)*{a}; (18,8)*{n};
   (0,0)*{\bigb{h_{x+y+b-a+1-n}(a,b,1)}};
  \endxy \nn
\end{eqnarray}
\end{lem}

\begin{proof} If $a=0$ and $b$ is arbitrary the lemma follows from \eqref{eq_triangle} and the thin reduction to bubbles axiom \eqref{eq_reduction}.  Similarly, when $b=0$ and $a$ is arbitrary. Otherwise, it suffices to assume that $a,b >0$ so that
\begin{eqnarray}
   \xy
 (0,0)*{\includegraphics[scale=0.4]{figs/marko-lem3-1}};
 (-12,-10)*{a};(13,-10)*{b};(-12,10)*{a};(13,10)*{b};
 (-12,2)*{a+1};(12,2)*{b+1};
  (3,7.5)*{\bullet}+(0,3)*{y};
  (3,-8)*{\bullet}+(0,-3)*{x}; (19,-8)*{n};
  \endxy
  & \refequal{\eqref{eq_triangle}} &
     \xy
 (0,0)*{\includegraphics[scale=0.4]{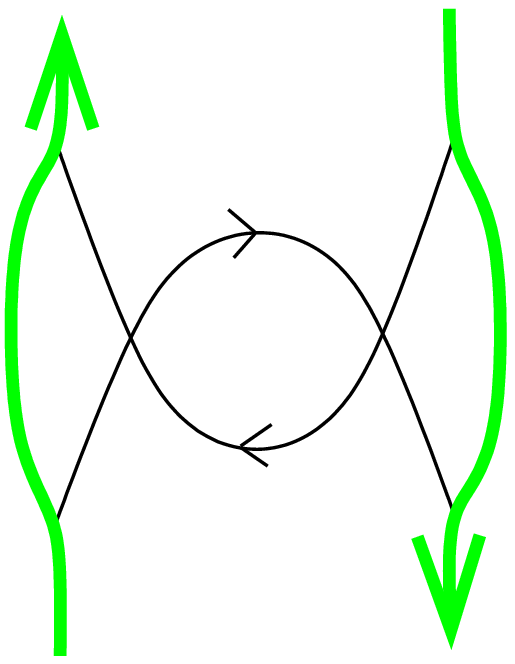}};
 (-13,-13)*{a};(14,-12)*{b};(-13,12)*{a};(13,13)*{b};
 (-17,0)*{a-1};(17,0)*{b-1};
  (3,3)*{\bullet}+(0,3)*{x};
  (3,-4)*{\bullet}+(0,-3)*{y};
  \endxy \nn \hspace{2in}
\end{eqnarray}
\begin{eqnarray}
 & \refequal{\eqref{eq_ident_decomp-dots}} &
  - \;\;
     \xy
 (0,0)*{\includegraphics[scale=0.4]{figs/marko-lem3-2}};
 (-13,-13)*{a};(14,-12)*{b};(-13,12)*{a};(13,13)*{b};
  (3,4.5)*{\bullet}+(0,3)*{x};
  (3,-4)*{\bullet}+(0,-3)*{y}; (19,-8)*{n};
  \endxy
  + \sum_{  \xy
  (0,-1)*{\scs p+q+r =};
  (0,-4)*{\scs x+y-n+2b-1};
  \endxy}
    \xy
 (-14,0)*{\reflectbox{\includegraphics[scale=0.45]{figs/split-thinthick2.eps}}};
 (-16.5,-11)*{a};(-9.5,0)*{\bullet}+(2,1)*{ p};
  (14,0)*{\reflectbox{\includegraphics[angle=180,scale=0.45]{figs/split-thinthick2.eps}}};
 (16.5,-11)*{b};(9.5,0)*{\bullet}+(-2,1)*{ q};
  (0,-4)*{\cbub{\spadesuit +r}{}}; (19,8)*{n};
  \endxy \nn
\end{eqnarray}
The Lemma follows from (\ref{eq_schur_left}).
\end{proof}

\begin{lem} \label{lem_bottom_zero}
If $x < n+a-b$, then
\begin{equation}
   \xy
 (0,0)*{\includegraphics[scale=0.4]{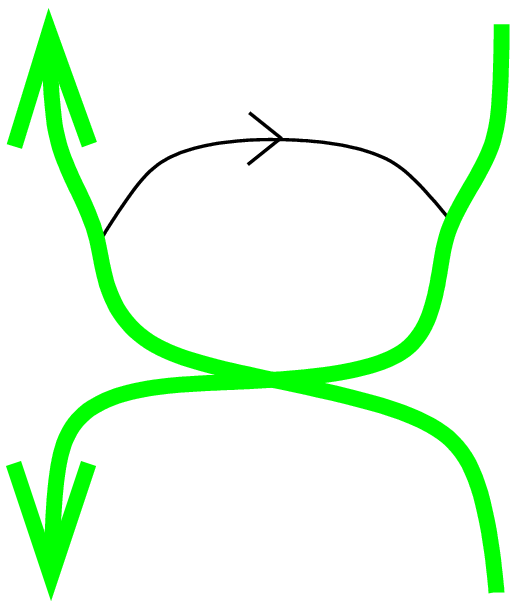}};
 (-16,-9)*{b-i};(16,-10)*{a-i};(-18,10)*{a-i-1};(18,9)*{b-i-1};
  (3,6.5)*{\bullet}+(0,3)*{x};
  (15,0)*{n};
  \endxy = 0.
\end{equation}
for any $0 \leq i \leq \min(a-1,b-1)$.
\end{lem}

\begin{proof}
Using thick calculus relations we have
  \begin{equation}
 \xy
 (0,0)*{\includegraphics[scale=0.4]{figs/marko-lem6-1}};
 (-16,-9)*{b-i};(16,-9)*{a-i};(-18,9)*{a-i-1};(18,9)*{b-i-1};
  (3,6.5)*{\bullet}+(0,3)*{x};
  (15,0)*{n};
  \endxy
  \;\;\refequal{\eqref{eq_defn_thick_cross}} \;\;
 \xy
 (0,0)*{\includegraphics[scale=0.4]{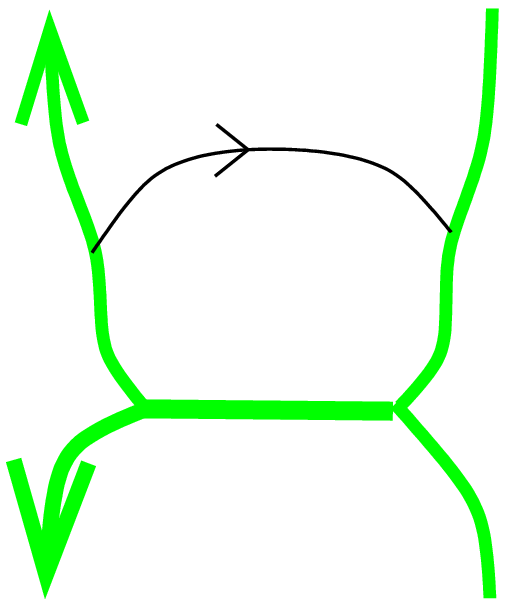}};
 (-16,-9)*{b-i};(16,-9)*{a-i};(-18,9)*{a-i-1};(18,9)*{b-i-1};
  (3,6.5)*{\bullet}+(0,3)*{x};
  (15,0)*{n};
  \endxy  \hspace{1in}
 \end{equation}
 \begin{equation}
   \;\; \refequal{\eqref{eq_split_assoc}} \;\;
 \xy
 (0,0)*{\includegraphics[scale=0.4]{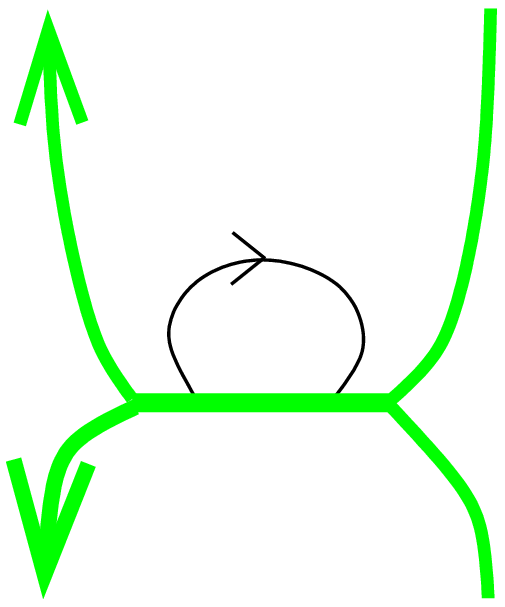}};
  (-16,-9)*{b-i};(16,-9)*{a-i};(-18,9)*{a-i-1};(18,9)*{b-i-1};
  (3,1)*{\bullet}+(0,3)*{x};
  (15,0)*{n};
  \endxy
    \;\; \refequal{\eqref{eq_triangle}} \;\;
 \xy
 (0,0)*{\includegraphics[scale=0.4]{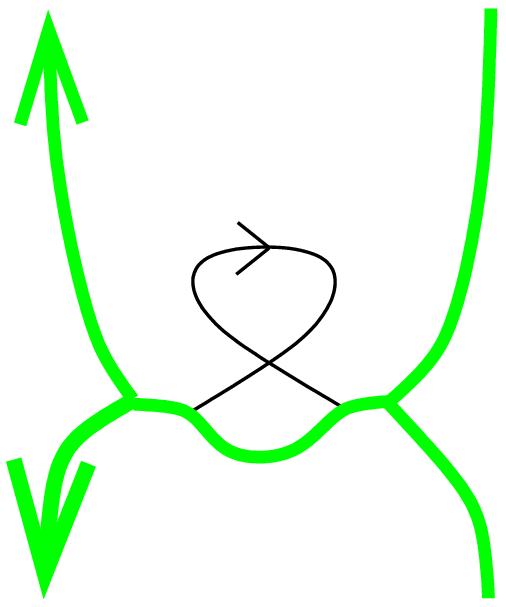}};
  (-16,-9)*{b-i};(16,-9)*{a-i};(-18,9)*{a-i-1};(18,9)*{b-i-1};
  (3,2)*{\bullet}+(0,3)*{x};
  (15,0)*{n};
  \endxy
\end{equation}
and the last diagram is zero by \eqref{eq_reduction-dots} and \eqref{eq_split_right_thin} since $x<n+a-b$.
\end{proof}

\begin{prop}[Thick Bubble Slides] For $\alpha \in P(b)$
\begin{eqnarray}
  \xy
 (0,0)*{\includegraphics[scale=0.5]{figs/tlong-up.eps}};
 (-2.5,-11)*{a};
 (12,0)*{\stccbub{b}{\alpha}};
 (9,9)*{n};
  \endxy
  &=&
   \sum_{\beta,x,y} c_{\beta,x,y}^{\alpha}\qquad
   \xy
 (0,0)*{\includegraphics[scale=0.5]{figs/tlong-up.eps}};
 (-2.5,-11)*{a};(0,-5)*{\bigb{\pi_y}};(0,3)*{\bigb{\pi_x}};
 (-12,-2)*{\stccbub{b}{\beta}}; (9,9)*{n};
  \endxy
\\
  \xy
 (1,0)*{\includegraphics[scale=0.5]{figs/tlong-up.eps}};
 (-1.5,-11)*{a};
 (10,0)*{\stcbub{b}{\alpha}};
 (9,9)*{n};
  \endxy \quad&=&
   \sum_{\beta,x,y}
     (-1)^{|x|+|y|}c_{\beta,x,y}^{\alpha}\qquad
   \xy
 (0,0)*{\includegraphics[scale=0.5]{figs/tlong-up.eps}};
 (-2.5,-11)*{a};(0,-5)*{\bigb{\pi_{\overline{y}}}};(0,3)*{\bigb{\pi_{\overline{x}}}};
 (-12,-2)*{\stccbub{b}{\beta}}; (9,9)*{n};
  \endxy
\end{eqnarray}
\end{prop}

\begin{proof}
The proof follows by taking $a=1$ in the General Central Element Lemma \ref{lem_gen_central}.  This Lemma can be used to slide an arbitrary thick bubble through a thin line.  By exploding an arbitrary thick line and iterating the application of this Lemma the Proposition follows.
\end{proof}

%
\section{Decompositions of functors and other applications}
%

We continue to use the conventions for bounds on summation indices discussed in Remark~\ref{rem_diagram_conventions}.

\subsection{Decomposition of $\cal{E}^{(a)}\cal{E}^{(b)}\onen$}

In this section we categorify the relation
\begin{eqnarray}
 E^{(a)}E^{(b)}1_n &=& \qbin{a+b}{a}E^{(a+b)}1_n
\end{eqnarray}
in $\UA$.

For $\alpha\in P(a,b)$ the diagrams $\sigma_{\alpha}$ and $\lambda_{\alpha}$ introduced in \eqref{eq_def_sigmalambda_alpha} can now be viewed as maps
\begin{eqnarray}
  \sigma_{\alpha}:= \xy
 (0,0)*{\includegraphics[scale=0.5]{figs/tsplit.eps}};
 (-5,-11)*{a+b};(-8,8)*{a};(8,8)*{b};
 (-5,2)*{\bigb{\pi_{\alpha}}};
  \endxy
  &\maps & \cal{E}^{(a+b)}\onen \{2|\alpha|-ab\} \longrightarrow \cal{E}^{(a)}\cal{E}^{(b)}\onen, \\ \nn \\
 \lambda_{\alpha}:= (-1)^{|\hat{\alpha}|}
     \xy
 (0,0)*{\includegraphics[scale=0.5,angle=180]{figs/tsplitd.eps}};
 (-5,11)*{a+b};(-8,-8)*{a};(8,-8)*{b};
 (5,-2)*{\bigb{\pi_{\hat{\alpha}}}};
  \endxy
& \maps& \cal{E}^{(a)}\cal{E}^{(b)}\onen\longrightarrow
\cal{E}^{(a+b)}\onen \{2|\alpha|-ab\}.
\end{eqnarray}

\begin{thm} \label{thm_cal_EaEb}
The maps
\begin{eqnarray}
 \sum_{\alpha \in P(a,b)} \sigma_{\alpha} &\maps &
 \bigoplus_{\alpha \in P(a,b)}  \cal{E}^{(a+b)}\onen \{2|\alpha|-ab\} \longrightarrow \cal{E}^{(a)}\cal{E}^{(b)}\onen
\end{eqnarray}
and
\begin{eqnarray}
 \sum_{\alpha \in P(a,b)} \lambda_{\alpha} & \maps &  \cal{E}^{(a)}\cal{E}^{(b)}\onen \longrightarrow \bigoplus_{\alpha \in P(a,b)}  \cal{E}^{(a+b)}\onen \{2|\alpha|-ab\}
\end{eqnarray}
are mutually-inverse isomorphisms, giving a canonical isomorphism
\begin{eqnarray}
\cal{E}^{(a)}\cal{E}^{(b)}\onen \cong \bigoplus_{\qbin{a+b}{a}}\cal{E}^{(a+b)}\onen
\refequal{\eqref{eq_Pab_card}} \bigoplus_{\alpha \in P(a,b)}\cal{E}^{(a+b)}\onen \{ 2|\alpha|-ab\}
\end{eqnarray}
in $\UcatD$ for all $n \in \Z$.
\end{thm}

\begin{proof}
The theorem follows from formulas in Corollary~\ref{cor_oval_small} and Theorem~\ref{thm_nil-eaeb} and from the results in Section~\ref{sec_split_exp_idem}.
\end{proof}

 \subsection{Decomposition of $\cal{E}^{(a)}\cal{F}^{(b)}\onen$} \label{sec_decomp_EaFb}

Our goal in this section is to categorify the relations in $\UA$
\begin{align} \label{eq_defining1}
 E^{(a)}F^{(b)}1_{n}&=
\sum_{j=0}^{\min(a,b)}\qbin{a-b+n}{j}F^{(b-j)}E^{(a-j)}1_{n}, \quad
\text{if $n\ge b-a$}, \\ \label{eq_defining2}
F^{(b)}E^{(a)}1_n&=
\sum_{j=0}^{\min(a,b)}\qbin{b-a-n}{j}E^{(a-j)}F^{(b-j)}1_n, \quad \text{if $n\le b-a$},
\end{align}
that decompose products $E^{(a)}F^{(b)}1_n$ and $F^{(b)}E^{(a)}1_n$ when they are not canonical basis vectors into positive linear combinations of canonical basis vectors. Relations \eqref{eq_defining1} and \eqref{eq_defining2} hold without restrictions on the values of $n,a,b$, but under these restrictions the coefficients in the right hand side lie in $\Z_{+}[q,q^{-1}]$, allowing for categorification.  These relations turn into explicit decompositions of $\cal{E}^{(a)}\cal{F}^{(b)}\onen$ and $\cal{F}^{(b)}\cal{E}^{(a)}\onen$ when these 1-morphisms are decomposable into a direct sum of indecomposable 1-morphisms:
\begin{align}
 \cal{E}^{(a)}\cal{F}^{(b)}\onen&\cong
\bigoplus_{j=0}^{\min(a,b)}\bigoplus_{\qbin{a-b+n}{j}}\cal{F}^{(b-j)}
\cal{E}^{(a-j)}\onen, \nn \\ \nn &\refequal{\eqref{eq_Pab_card}}  \bigoplus_{j=0}^{\min(a,b)} \bigoplus_{\alpha \in P(j,n+a-b-j)}\cal{F}^{(b-j)}\cal{E}^{(a-j)}\onen\{2|\alpha|-j(a-b+n)\} \quad
\text{if $n\ge b-a$}, \\
 \cal{F}^{(b)}\cal{E}^{(a)}\onen&\cong
\bigoplus_{j=0}^{\min(a,b)}\bigoplus_{\qbin{b-a-n}{j}}\cal{E}^{(a-j)}\cal{F}^{(b-j)}\onen, \nn \\ \nn &\refequal{\eqref{eq_Pab_card}}  \bigoplus_{j=0}^{\min(a,b)} \bigoplus_{\alpha \in P(j,-n+b-a-j)}\cal{E}^{(a-j)}\cal{F}^{(b-j)}\onen\{2|\alpha|-j(b-a-n)\}
\quad
\text{if $n\leq b-a$.}\\ \label{eq_tired}
\end{align}
in $\UcatD$.  The goal of this section is to show that these decompositions hold over any ground commutative ring $\Bbbk$, see Theorem~\ref{eq_cat_EaFb}. Existence of these decompositions for $\Bbbk$ a field was shown implicitly in \cite{Lau1}.

Let $i \in \{ 0,1, \dots, \min(a,b)\}$ be fixed.  For every
partition $\alpha\in P(i,n+a-b-i)$
we define the following maps
\begin{eqnarray}
 \lambda_{\alpha}^i &:=&   \xy
  (0,0)*{\includegraphics[scale=0.5]{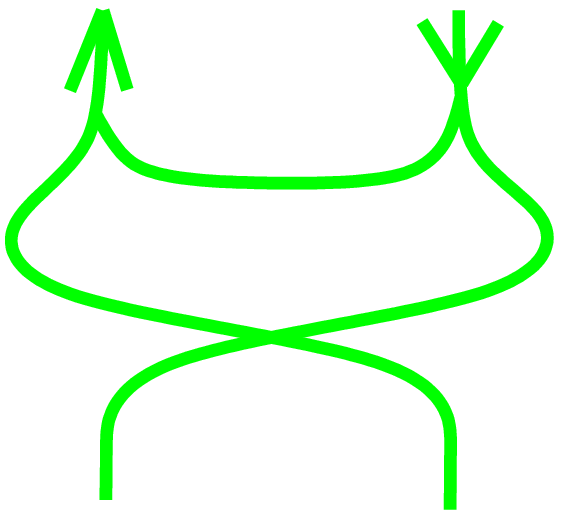}};
 (-2,4)*{\bigb{\pi_{\alpha} }};  (4,2)*{i};
  (15,0)*{}; (18,7)*{n};  (14,-12)*{a-i};(-15,-12)*{b-i};(13,12)*{b};(-13,12)*{a};
 \endxy \maps \F{b-i}\E{a-i}\onen \{ 2|\alpha|-i(n+a-b-i) \} \longrightarrow \E{a}\F{b} \onen \nn \\ \\
\sigma_{\alpha}^i  &:=&   (-1)^{ab}\sum_{\beta, \gamma, x, y}(-1)^{\frac{i(i+1)}{2}+|x|+|y|}
c_{\alpha,\beta,\gamma,x, y}^{K_i}\xy
  (0,0)*{\includegraphics[scale=0.5]{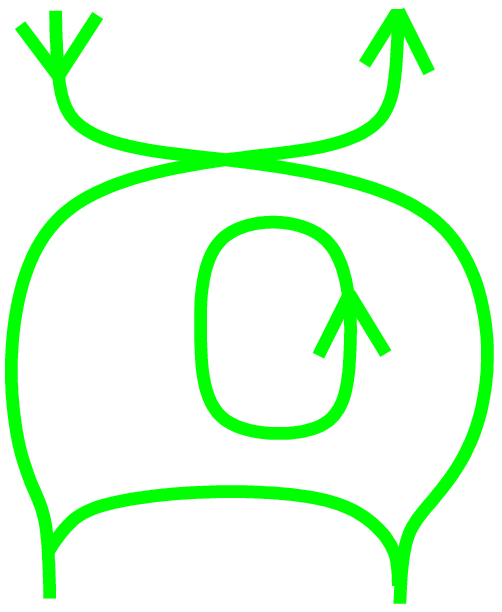}};
   (-2,0)*{\bigb{\pi_{\gamma}^{\spadesuit} }};(0,-10)*{\bigb{\pi_{\beta} }}; (-12,0)*{\bigb{\pi_{\overline{x}} }};(13,0)*{\bigb{\pi_{\overline{y}} }}; (-7,-7)*{i};  (7,-5)*{i};
  (15,0)*{}; (18,6)*{n};  (14,12)*{a-i};(-16,12)*{b-i};(10,-13)*{b};(-12,-13)*{a};
 \endxy \nn \\
 & & \hspace{1in} \maps \E{a}\F{b}\onen \longrightarrow \F{b-i}\E{a-i}\onen\{2|\alpha|-i(n+a-b-i) \},
\end{eqnarray}
where $K_0=\emptyset$ and $K_i=((n+a-b-i)^i)$.  Recall from \eqref{eq_Pab_card} that there are $|P(i,n+a-b-i)|={n+a-b \choose i}$ such partitions $\alpha$ and that
\begin{equation}
 \sum_{\alpha \in P(i,n+a-b-i)}q^{2|\alpha|-i(n+a-b-i)} = \qbin{n+a-b}{i}.
\end{equation}

Let \begin{equation}
  e_{\alpha}^i = \lambda_{\alpha}^i \sigma_{\alpha}^i \maps \E{a}\F{b}\onen \longrightarrow
  \E{a}\F{b} \onen.
\end{equation}

\begin{lem} \label{lem_orthogonal-pre}
Given $0\le i,j \le \min(a,b)$ and
partitions $\alpha\in P(i,n+a-b-i)$ and $\alpha'\in P(j,n+a-b-j)$ then
\begin{equation} \label{eq_orthog_siglam}
\sigma_{\alpha'}^{j}\lambda_{\alpha}^i = \delta_{i,j}\delta_{\alpha,\alpha'}
  (-1)^{(a-i)(b-i)} \xy
 (0,0)*{\includegraphics[scale=0.5]{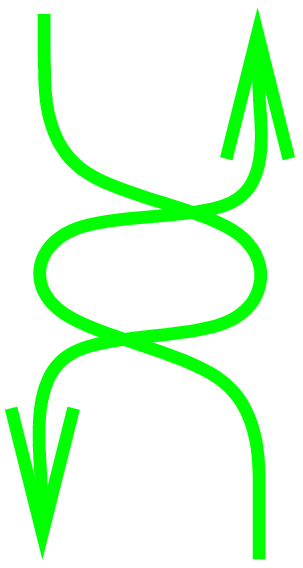}};
 (10 ,-12)*{a-i}; (-11 ,-12)*{b-i}; (12,8)*{n};
 \endxy
\end{equation}
\end{lem}

\begin{proof}
The composite $\sigma_{\alpha'}^{j}\lambda_{\alpha}^i$ has the following form
\begin{equation} \label{eq_sigma_lambda}
(-1)^{ab}\sum_{\beta, \gamma, x, y}(-1)^{\frac{i(i+1)}{2}+|x|+|y|}
c_{\alpha',\beta,\gamma,x, y}^{K_i}
 \xy
  (0,9)*{\includegraphics[scale=0.5]{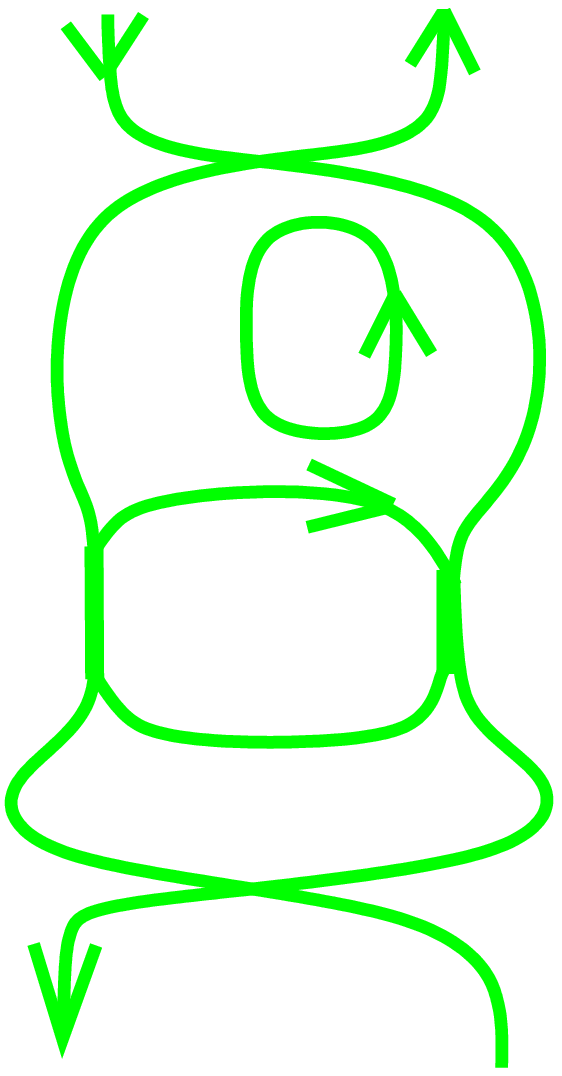}};
 (-4,0)*{\bigb{\pi_{\alpha} }};(-4,12)*{\bigb{\pi_{\beta} }};
  (6,-3)*{i};(7,13)*{j};(7,24)*{j};
  (15,0)*{}; (22,8)*{n};
  (13,5)*{b};(-13,5)*{a};
  (16,-14)*{a-i};(-17,-14)*{b-i};
  (16,32)*{a-j};(-17,32)*{b-j};
   (-2,20)*{\bigb{\pi_{\gamma}^{\spadesuit} }};
   (-12,20)*{\bigb{\pi_{\overline{x}} }};
   (13,20)*{\bigb{\pi_{\overline{y}} }};
 \endxy
\end{equation}
Since $K_i=((n+a-b-i)^i)$ is a rectangular partition, we have
\begin{equation}
  c_{\alpha',\beta,\gamma,x,y}^{K_i} = \sum_{\theta}
  c_{\alpha',\beta}^{\theta} c_{\gamma,x,y}^{K_i-\theta}.
\end{equation}
Furthermore, using Proposition~\ref{prop_image_thickbub} we have
 \begin{eqnarray}
   (\phi_{a-j,b-j,j}^n)^{-1}\left(
   \xy
  (-12,0)*{\includegraphics[scale=0.5]{figs/single-tup.eps}};
  (10,0)*{\includegraphics[scale=0.5,angle=180]{figs/single-tup.eps}};
  (0,0)*{\stccbub{j}{\gamma}};
  (-17,-6)*{a-j};(16,-6)*{b-j};
 \endxy
  \right)  \;\; =\;\;  (-1)^{\frac{j(j-1)}{2}+|\gamma|}\pi_{\overline{\gamma}}(\underline{z}).
 \end{eqnarray}

\begin{equation}
  p_{\theta} := \phi_{a-j,b-j,j}^n\left((-1)^{\frac{j(j-1)}{2}+\frac{i(i+1)}{2}}\sum_{\gamma, x, y} c_{\gamma,x,y}^{K_i-\theta} (-1)^{|x| +
  |y|+|\gamma|} \pi_{\overline{x}}(\underline{x})
  \pi_{\overline{y}}(\underline{y}) \pi_{\bar{\gamma}}(\underline{z})\right).
\end{equation}
To simplify notation we will identify $p_{\theta}$ with its pre-image in $\Z[\und{x},\und{y},\und{z}]$. Note that the coefficients $c_{\gamma,x,y}^{K_i-\theta}$ are zero unless $|x|+|y|+|\gamma|=|K_i - \theta|$.  Hence, we can write
\begin{eqnarray}
  p_{\theta} &= & (-1)^{\frac{j(j-1)}{2}+\frac{i(i+1)}{2}}(-1)^{|K_i-\theta|} \sum_{\gamma, x, y} c_{\gamma,x,y}^{K_i-\theta} \pi_{\overline{x}}(\underline{x})
  \pi_{\overline{y}}(\underline{y}) \pi_{\bar{\gamma}}(\underline{z})  \\
  &= & (-1)^{\frac{j(j-1)}{2}+\frac{i(i+1)}{2}}(-1)^{|K_i-\theta|} \pi_{\overline{K_i-\theta}}(\underline{x},\und{y},\und{z}), \label{eq_simple_ptheta}
\end{eqnarray}
where the last equality used that $c_{\gamma,x,y}^{K_i-\theta} = c_{\overline{\gamma},\overline{x},\overline{y}}^{\overline{K_i-\theta}}$ together with \eqref{eq_prod_schur_twovar}.

Observe that $p_{\theta}$ is symmetric in all variables
$\underline{x}$, $\underline{y}$, $\underline{z}$.  Thus, we can
write \eqref{eq_sigma_lambda} as
\begin{equation} \label{eq_sigma_lambda2}
  (-1)^{ab}\sum_{\beta,\theta} c_{\alpha',\beta}^{\theta}
 \xy
  (0,9)*{\includegraphics[scale=0.5]{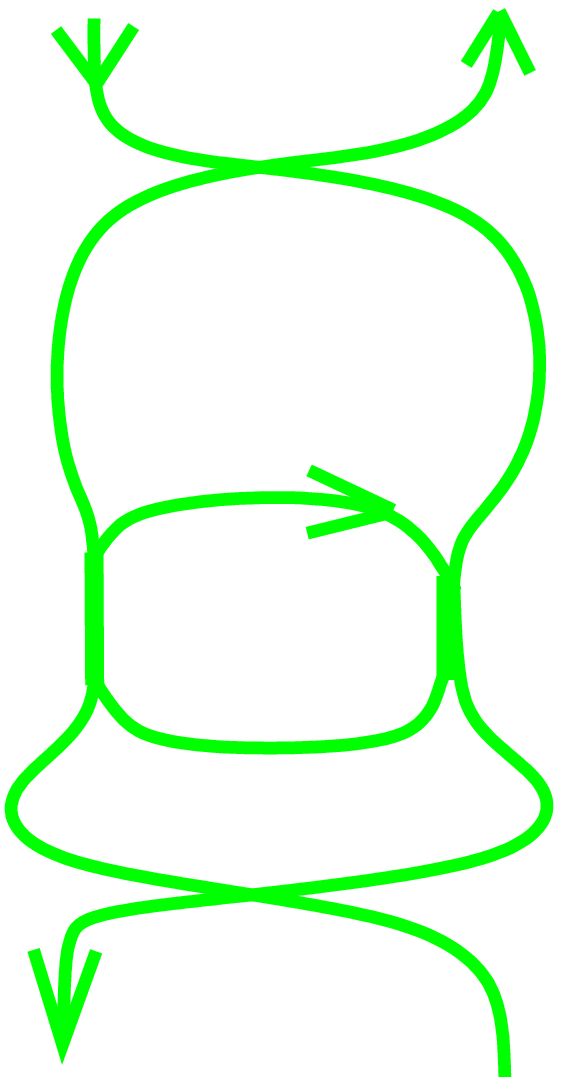}};
 (-4,0)*{\bigb{\pi_{\alpha} }};(-4,12)*{\bigb{\pi_{\beta} }};
  (6,-3)*{i};(7,13)*{j};
  (15,0)*{}; (22,8)*{n};
  (13,5)*{b};(-13,5)*{a};
  (16,-14)*{a-i};(-17,-14)*{b-i};
  (16,32)*{a-j};(-17,32)*{b-j};
   (0,20)*{\bigb{\hspace{0.5in}p_{\theta}\hspace{0.5in}}};
 \endxy
\end{equation}

Next use \eqref{eq_schur_thin} to split the Schur polynomials
$\pi_{\alpha}$ and $\pi_{\beta}$ into $i$, respectively $j$, dotted
thin lines. Repeatedly using associativity and coassociativity of
splitters
\begin{eqnarray}
 \xy
  (0,0)*{\includegraphics[scale=0.4]{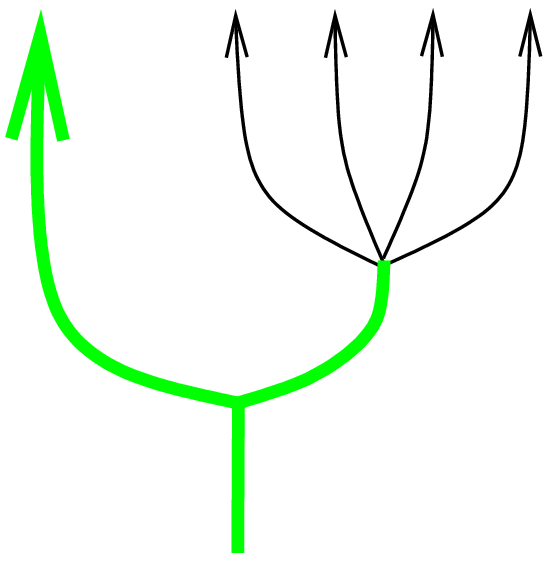}};
   (6,-4)*{j};
  (15,0)*{}; (10,-9)*{n};
  (-16,2)*{a-j};
  (-4,-9)*{a};
 \endxy
 \;\; = \;\;
  \xy
  (0,0)*{\includegraphics[scale=0.4]{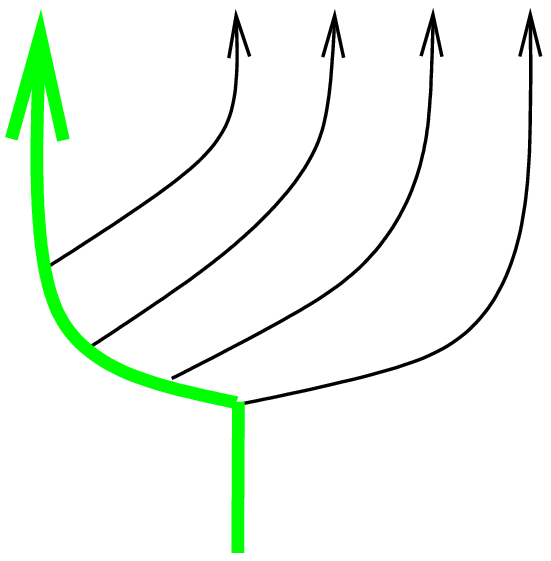}};
  (15,0)*{}; (10,-9)*{n};
  (-16,2)*{a-j};
  (-4,-9)*{a};
 \endxy
 \\ \nn
\\
 \xy
  (0,0)*{\includegraphics[scale=0.4]{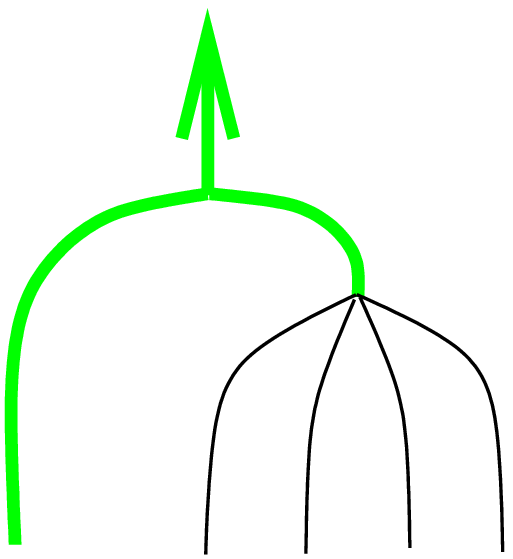}};
   (6,4)*{i};
  (15,0)*{}; (10,9)*{n};
  (-16,-9)*{a-i};
  (-6,9)*{a};
 \endxy
\;\; = \;\;
  \xy
  (0,0)*{\includegraphics[scale=0.4]{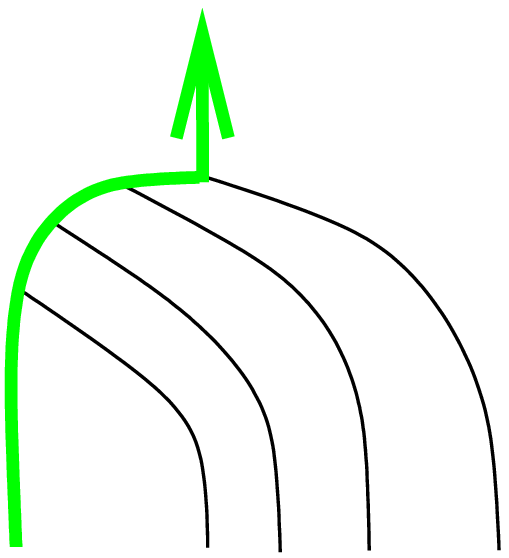}};
  (15,0)*{}; (10,9)*{n};
 (-16,-9)*{a-i};
  (-6,9)*{a};
 \endxy
\end{eqnarray}
\eqref{eq_sigma_lambda2} can be rewritten as the following `ladder
diagram'
\begin{equation} \label{eq_big_ladder}
 \sum_{\beta,\theta} c_{\alpha',\beta}^{\theta}\vcenter{\xy
  (0,8.2)*{\includegraphics[scale=0.5]{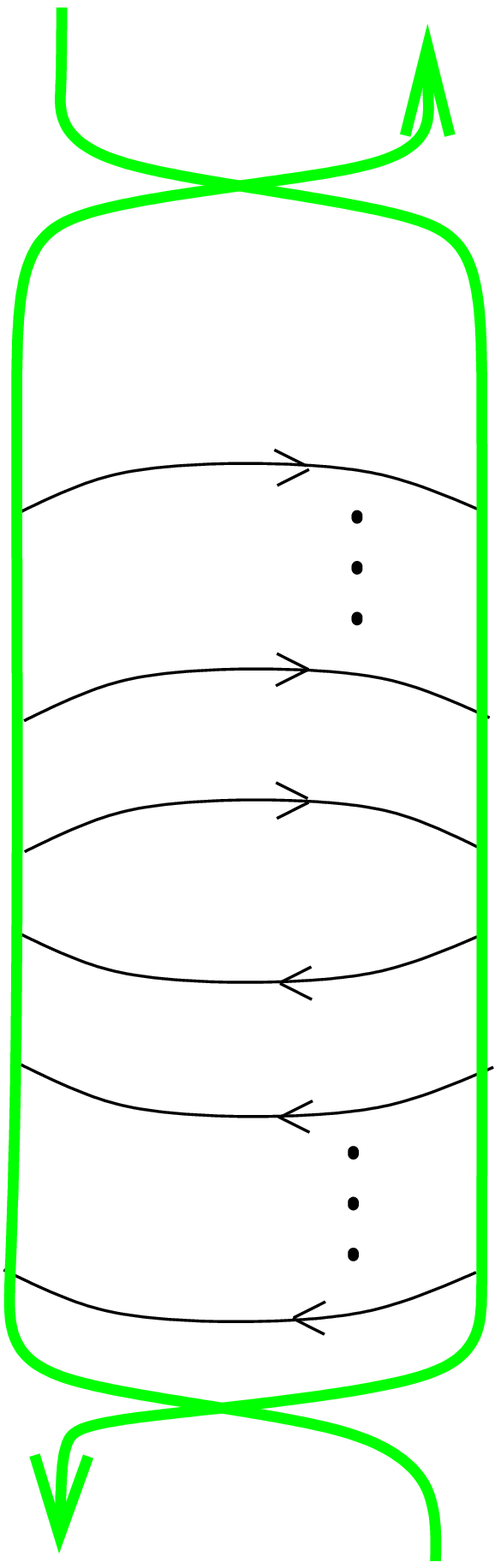}};
  (-4,27)*{\bullet}+(-1,3)*{\scs \beta_1+j-1};
  (-4,15)*{\bullet}+(-1,3)*{\scs \beta_{j-1}+1};
  (-4,7)*{\bullet}+(-1,3)*{\scs \beta_j};
  (-4,-3.5)*{\bullet}+(-1,3)*{\scs \alpha_i};
  (-4,-11.5)*{\bullet}+(-1,3)*{\scs \alpha_{i-1}+1};
  (-4,-23.5)*{\bullet}+(-1,3)*{\scs \alpha_1+i-1};
  (15,0)*{}; (21,9)*{n};
  (17,2)*{b};(-16,2)*{a};
  (19,-32)*{a-i};(-19,-32)*{b-i};
  (19,50)*{a-j};(-19,50)*{b-j};
  (0,36)*{\bigb{\hspace{0.5in}p_{\theta}\hspace{0.5in}}};
 \endxy}
 \quad = \quad \sum_{\beta,\theta} c_{\alpha',\beta}^{\theta}\vcenter{\xy
  (0,8.2)*{\includegraphics[scale=0.5]{figs/eafb5.eps}};
  (-4,27)*{\bullet}+(-1,3)*{\scs \overline{\beta_1}};
  (-4,15)*{\bullet}+(-1,3)*{\scs \overline{\beta_{j-1}}};
  (-4,7)*{\bullet}+(-1,3)*{\scs \overline{\beta_j}};
  (-4,-3.5)*{\bullet}+(-1,3)*{\scs \overline{\alpha_i}};
  (-4,-11.5)*{\bullet}+(-1,3)*{\scs \overline{\alpha_{i-1}}};
  (-4,-23.5)*{\bullet}+(-1,3)*{\scs \overline{\alpha_1}};
  (15,0)*{}; (21,9)*{n};
  (17,2)*{b};(-16,2)*{a};
  (19,-32)*{a-i};(-19,-32)*{b-i};
  (19,50)*{a-j};(-19,50)*{b-j};
  (0,36)*{\bigb{\hspace{0.5in}p_{\theta}\hspace{0.5in}}};
 \endxy}
\end{equation}
where we ease notation by writing $\overline{\alpha_{\ell}}:=\alpha_{\ell}+i-\ell$ and $\overline{\beta_{\ell'}}:=\beta_{\ell'}+i-\ell'$ for $1\leq \ell \leq i$ and $1 \leq \ell' \leq j$.

Now we look at a single term in the summation over $\beta$ and $\theta$. Apply the Square Flop Lemma~\ref{lem_square_flop} to simplify the square in the center of the diagram.
\begin{equation}
 \vcenter{\xy0;/r.19pc/:
  (0,8.2)*{\includegraphics[scale=0.4]{figs/eafb5.eps}};
  (-4,27)*{\bullet}+(-1,3)*{\scs \overline{\beta_1}};
  (-4,15)*{\bullet}+(-1,3)*{\scs \overline{\beta_{j-1}}};
  (-4,7)*{\bullet}+(-1,3)*{\scs \overline{\beta_j}};
  (-4,-3.5)*{\bullet}+(-1,3)*{\scs \overline{\alpha_i}};
  (-4,-11.5)*{\bullet}+(-1,3)*{\scs \overline{\alpha_{i-1}}};
  (-4,-23.5)*{\bullet}+(-1,3)*{\scs \overline{\alpha_1}};
  (15,0)*{}; (21,20)*{n};
  (17,2)*{b};(-16,2)*{a};
  (19,-32)*{a-i};(-19,-32)*{b-i};
  (19,50)*{a-j};(-19,50)*{b-j};
  (0,36)*{\bigb{\hspace{0.5in}p_{\theta}\hspace{0.5in}}};
 \endxy}
\;\; = \;\;- \hspace{-0.1in}
 \vcenter{\xy0;/r.19pc/:
  (0,8.2)*{\includegraphics[scale=0.4]{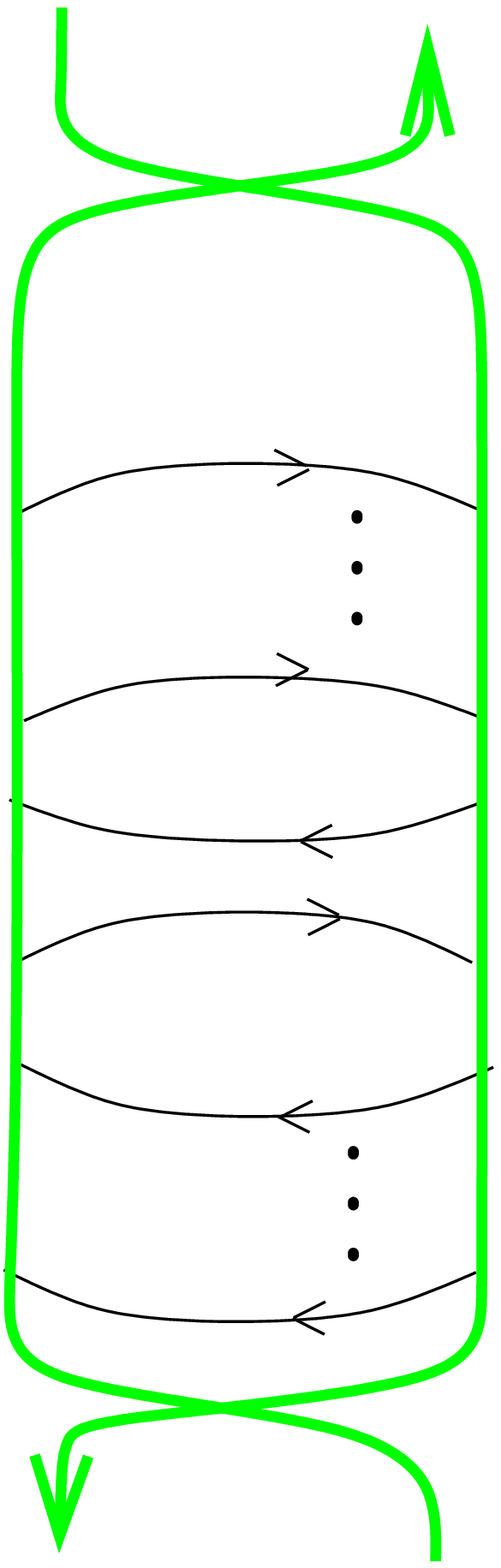}};
  (-4,27)*{\bullet}+(-1,3)*{\scs \overline{\beta_1}};
  (-4,15)*{\bullet}+(-1,3)*{\scs \overline{\beta_{j-1}}};
  (-4,5)*{\bullet}+(-1,3)*{\scs \overline{\alpha_i}};
  (-4,0.5)*{\bullet}+(-1,-3.5)*{\scs \overline{\beta_j}};
  (-4,-11.5)*{\bullet}+(-1,3)*{\scs \overline{\alpha_{i-1}}};
  (-4,-23.5)*{\bullet}+(-1,3)*{\scs \overline{\alpha_1}};
  (15,0)*{}; (21,20)*{n};
  (19,-32)*{a-i};(-19,-32)*{b-i};
  (19,50)*{a-j};(-19,50)*{b-j};
  (0,36)*{\bigb{\hspace{0.5in}p_{\theta}\hspace{0.5in}}};
 \endxy}
 \hspace{-0.1in} + \hspace{-0.1in}
 \sum_{  \xy
  (0,-1)*{\scs p_i+q_i+r_i};
  (0,-4)*{\scs =\alpha_i+\beta_j};
  (0,-7)*{\scs -n+b-a+1};
  \endxy}
\hspace{-0.1in}(-1)^{a-b} \hspace{-0.1in}
 \vcenter{\xy0;/r.19pc/:
  (0,8.2)*{\includegraphics[scale=0.4]{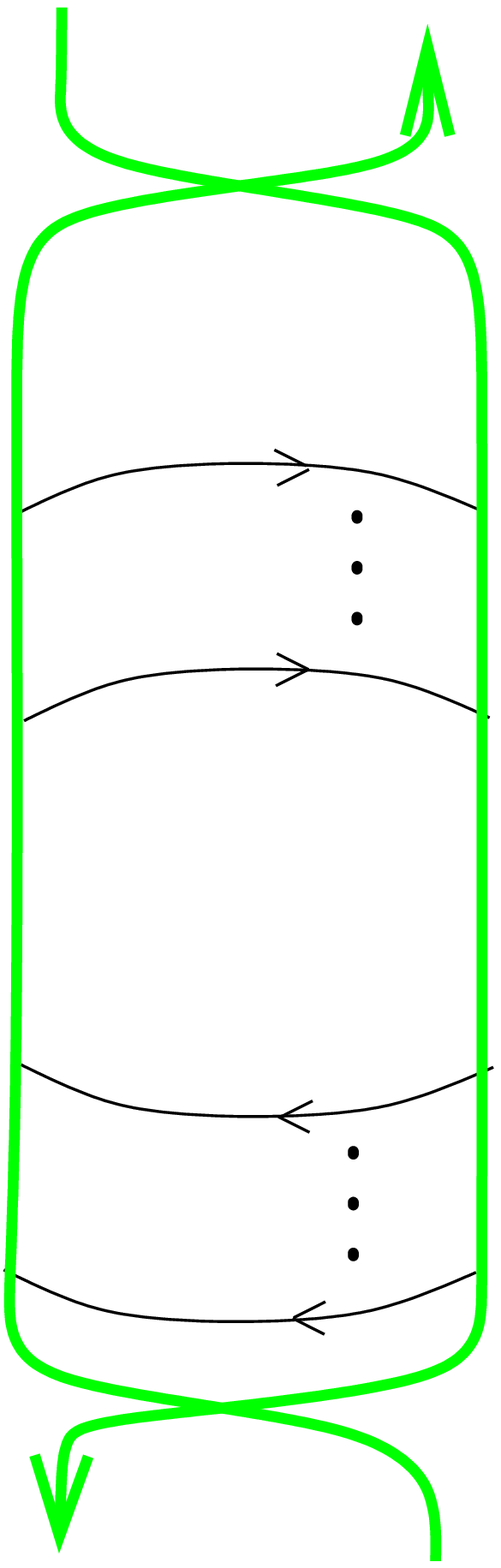}};
  (-4,27)*{\bullet}+(-1,3)*{\scs \overline{\beta_1}};
  (-4,15)*{\bullet}+(-1,3)*{\scs \overline{\beta_{j-1}}};
  (-14,4)*{\bigb{h_{p_i}}};
  (14,4)*{\bigb{h_{q_i}}};
  (0,2)*{\cbub{\spadesuit +r_i}{}};
  (-4,-11.5)*{\bullet}+(-1,3)*{\scs \overline{\alpha_{i-1}}};
  (-4,-23.5)*{\bullet}+(-1,3)*{\scs \overline{\alpha_1}};
  (15,0)*{}; (24,20)*{n};
  (17,-4)*{b};(-16,-4)*{a};
  (19,-32)*{a-i};(-19,-32)*{b-i};
  (19,50)*{a-j};(-19,50)*{b-j};
  (0,36)*{\bigb{\hspace{0.5in}p_{\theta}\hspace{0.5in}}};
 \endxy}
\end{equation}
or using the notation introduced in \eqref{eq_phiofh}
\begin{equation}
 \vcenter{\xy0;/r.19pc/:
  (0,8.2)*{\includegraphics[scale=0.4]{figs/eafb5.eps}};
  (-4,27)*{\bullet}+(-1,3)*{\scs \overline{\beta_1}};
  (-4,15)*{\bullet}+(-1,3)*{\scs \overline{\beta_{j-1}}};
  (-4,7)*{\bullet}+(-1,3)*{\scs \overline{\beta_j}};
  (-4,-3.5)*{\bullet}+(-1,3)*{\scs \overline{\alpha_i}};
  (-4,-11.5)*{\bullet}+(-1,3)*{\scs \overline{\alpha_{i-1}}};
  (-4,-23.5)*{\bullet}+(-1,3)*{\scs \overline{\alpha_1}};
  (15,0)*{}; (21,9)*{n};
  (17,2)*{b};(-16,2)*{a};
  (19,-32)*{a-i};(-19,-32)*{b-i};
  (19,50)*{a-j};(-19,50)*{b-j};
  (0,36)*{\bigb{\hspace{0.5in}p_{\theta}\hspace{0.5in}}};
 \endxy}
\;\; = \;\;-\;
 \vcenter{\xy0;/r.19pc/:
  (0,8.2)*{\includegraphics[scale=0.4]{figs/eafb6.eps}};
  (-4,27)*{\bullet}+(-1,3)*{\scs \overline{\beta_1}};
  (-4,15)*{\bullet}+(-1,3)*{\scs \overline{\beta_{j-1}}};
  (-4,5)*{\bullet}+(-1,3)*{\scs \overline{\alpha_i}};
  (-4,0.5)*{\bullet}+(-1,-3.5)*{\scs \overline{\beta_j}};
  (-4,-11.5)*{\bullet}+(-1,3)*{\scs \overline{\alpha_{i-1}}};
  (-4,-23.5)*{\bullet}+(-1,3)*{\scs \overline{\alpha_1}};
  (15,0)*{}; (21,9)*{n};
  (19,-32)*{a-i};(-19,-32)*{b-i};
  (19,50)*{a-j};(-19,50)*{b-j};
  (0,36)*{\bigb{\hspace{0.5in}p_{\theta}\hspace{0.5in}}};
 \endxy}
 \quad +
(-1)^{a-b}
 \vcenter{\xy0;/r.19pc/:
  (0,8.2)*{\includegraphics[scale=0.4]{figs/eafb7.eps}};
  (-4,27)*{\bullet}+(-1,3)*{\scs \overline{\beta_1}};
  (-4,15)*{\bullet}+(-1,3)*{\scs \overline{\beta_{j-1}}};
  (0,4)*{\bigb{h_{\alpha_i+\beta_j-n+b-a+1}(1)}};
  (-4,-11.5)*{\bullet}+(-1,3)*{\scs \overline{\alpha_{i-1}}};
  (-4,-23.5)*{\bullet}+(-1,3)*{\scs \overline{\alpha_1}};
  (15,0)*{}; (24,13)*{n};
  (17,-4)*{b};(-16,-4)*{a};
  (19,-32)*{a-i};(-19,-32)*{b-i};
  (19,50)*{a-j};(-19,50)*{b-j};
  (0,36)*{\bigb{\hspace{0.5in}p_{\theta}\hspace{0.5in}}};
 \endxy}
\end{equation}
where we use the simplified notation $h_{\alpha_i+\beta_j-n+b-a+1}(1)$ for $h_{\alpha_i+\beta_j-n+b-a+1}(a,b,1)$ since the $a$ and $b$ labels can be recovered from the diagram.

Now we would like to iterate the application of the Square Flop
Lemma~\ref{lem_square_flop} moving the rightward oriented thin
lines down to the bottom of the ladder diagram, and the leftward oriented thin lines toward the top of the diagram, commuting the thin lines past the central elements using the Central Element Lemma~\ref{lem_ladder_slide1}, see \eqref{eq_alternate_central_element}.  If a rightward oriented thin line makes it to the bottom of the ladder diagram, or a leftward oriented thin line makes it to the top, then the entire ladder diagram is zero by Lemma~\ref{lem_bottom_zero} since
$x < n+a-b$ so that both diagrams
\begin{equation}
   \xy
 (0,0)*{\includegraphics[scale=0.4]{figs/marko-lem6-1}};
 (-16,-9)*{b-i};(16,-9)*{a-i};(-18,9)*{a-i-1};(18,9)*{b-i-1};
  (3,6.5)*{\bullet}+(0,3)*{x};
  (8,0)*{};
  \endxy  \qquad
     \xy
 (0,0)*{\includegraphics[angle=180,scale=0.4]{figs/marko-lem6-1}};
 (-16,9)*{b-i};(16,9)*{a-i};(-18,-9)*{a-i-1};(18,-9)*{b-i-1};
  (3,-6.5)*{\bullet}+(0,-3)*{x};
  (8,0)*{};
  \endxy
\end{equation}
are equal to zero.  Hence, the ladder diagram \eqref{eq_big_ladder} can only be nonzero if all the rightward oriented lines cancel with the leftward oriented lines.  This is only possible when $i=j$.

This implies that the idempotents $e_{\alpha}^i$ and $e_{\beta}^j$ are orthogonal when $i\neq j$.  Thus, we only need to consider the case when  $i=j$.  In this case we will have nonzero terms arising from terms where the rightward oriented thin lines cancel with the leftward oriented thin lines.
As we slide all the leftward oriented thin lines to the top of the ladder diagram and the rightward oriented thin lines to the bottom,
we must apply the Square Flop Lemma to a square where the leftward oriented thin line carries $\overline{\alpha_{\ell}}$ dots and the thin rightward oriented line carries $\overline{\beta_{\sigma(\ell)}}$ dots for some $\sigma\in S_i$. Keeping track of the signs appearing in the Square Flop Lemma, the overall sign arising from the permutation $\sigma$ is $\sgn(\sigma)$. Thus, the ladder diagram can be reduced to
\[
\eqref{eq_big_ladder} \;\; \refequal{i=j}\;\; (-1)^{i(a-b)}\sum_{\sigma\in S_i} \sgn(\sigma)\sum_{  \xy
  (0,-1)*{\scs p_1+q_1+r_1 =};
  (0,-4)*{\scs \overline{\alpha_1}+\overline{\beta_{\sigma(1)}}-n};
  (0,-7)*{\scs +b-a+1};
  \endxy} \dots
  \sum_{  \xy
  (0,-1)*{\scs p_i+q_i+r_i =};
  (0,-4)*{\scs \overline{\alpha_i}+\overline{\beta_{\sigma(i)}}-n};
  (0,-7)*{\scs +b-a+1};
  \endxy}
 \vcenter{\xy0;/r.19pc/:
  (0,8.2)*{\includegraphics[scale=0.4]{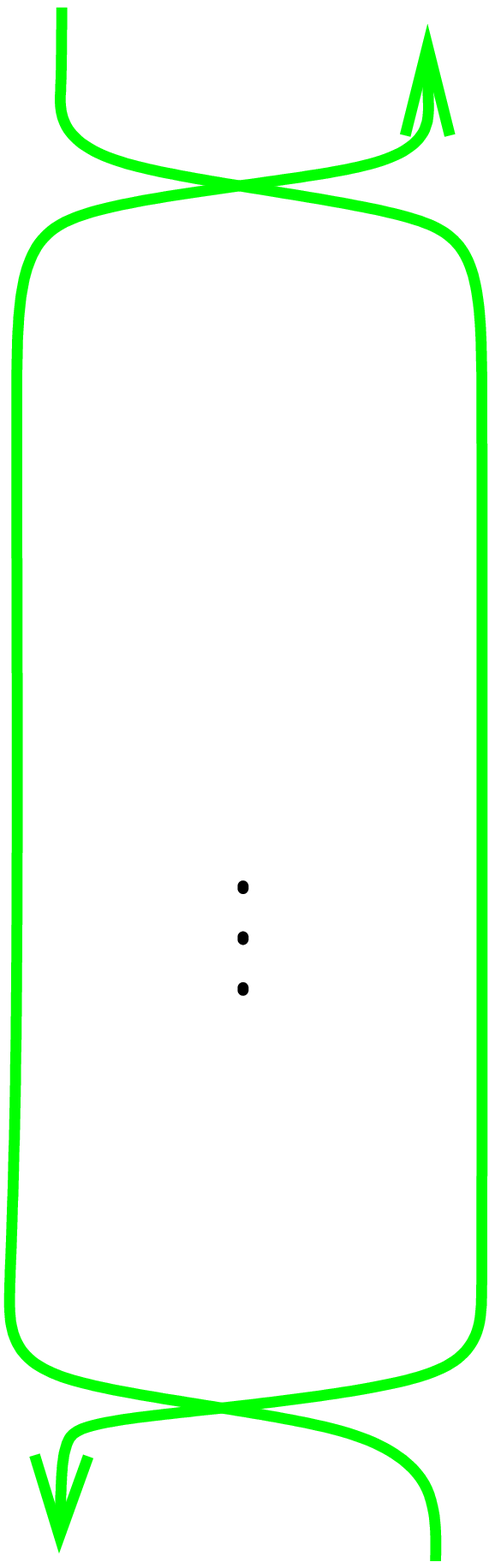}};
  (-14,-12)*{\bigb{h_{p_i}}};
  (14,-12)*{\bigb{h_{q_i}}};
  (0,-14)*{\cbub{\spadesuit +r_i}{}};
  (-14,20)*{\bigb{h_{p_1}}};
  (14,20)*{\bigb{h_{q_1}}};
  (0,18)*{\cbub{\spadesuit +r_1}{}};
  (15,0)*{}; (24,9)*{n};
  (19,-32)*{a-i};(-19,-32)*{b-i};
  (19,50)*{a-i};(-19,50)*{b-i};
  (0,36)*{\bigb{\hspace{0.5in}p_{\theta}\hspace{0.5in}}};
 \endxy}
 \]
which we can rewrite as
\begin{equation} \label{eq_last_ladder}
\eqref{eq_big_ladder} \;\; =\;\;(-1)^{i(a-b)}\sum_{\sigma\in S_i}\sgn(\sigma) \qquad \;\;\;
 \vcenter{\xy0;/r.19pc/:
  (0,8.2)*{\includegraphics[scale=0.4]{figs/eafb8.eps}};
  (0,-12)*{\bigb{
  h_{\overline{\alpha_i}+\overline{\beta_{\sigma(i)}}-n+b-a+1}(a-i,b-i,1)}};
  (0,20)*{\bigb{
  h_{\overline{\alpha_1}+\overline{\beta_{\sigma(1)}}-n+b-a+1}(a-i,b-i,1)}};
  (15,0)*{}; (24,9)*{n};
  (19,-32)*{a-i};(-19,-32)*{b-i};
  (19,50)*{a-i};(-19,50)*{b-i};
  (0,36)*{\bigb{\hspace{0.5in}p_{\theta}\hspace{0.5in}}};
 \endxy}
 \end{equation}

To complete the proof we utilize the injective homomorphism $\phi_{a-i,b-i,1}^n$  from \eqref{eq_EcFdgamma_isom} to pass from diagrams to the ring of symmetric polynomials. In particular,
\begin{equation}
  (\phi_{a-i,b-i,1}^n)^{-1}(h_{i}(a-i,b-i,1)) \;\; =\;\; h_{i}(\underline{x},\underline{y},\underline{z})
\end{equation}
where $\underline{x}=(x_1,\dots, x_{a-i})$, $\underline{y}=(y_1,\dots,y_{b-i})$, and $\underline{z}=(z_1,z_2,\dots)$.   Then we can write
\begin{eqnarray} \label{eq_phiinv_det}
 (\phi_{a-i,b-i,1}^n)^{-1}\left( \;\;
 \sum_{\sigma\in S_i}\sgn(\sigma)  \xy
 (-12,0)*{\includegraphics[scale=0.6]{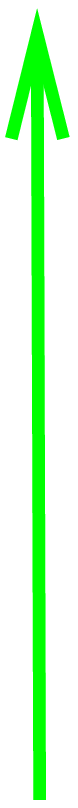}};
  (12,0)*{\includegraphics[angle=180, scale=0.6]{figs/tlonger-up.eps}};
 (18,-21)*{b-i}; (-18,-21)*{a-i};
   (0,-10)*{\bigb{ h_{\overline{\alpha_i}+\overline{\beta_{\sigma(i)}}-n+b-a+1}(1)}};
  (0,10)*{\bigb{
  h_{\overline{\alpha_1}+\overline{\beta_{\sigma(1)}}-n+b-a+1}(1)}};
   (21,0)*{n};(0,2)*{\vdots};
  \endxy \;\;\right)  =\hspace{1in}\end{eqnarray}
 \begin{eqnarray}
 \hspace{0.3in}& = &
 \sum_{\sigma\in S_i}\sgn(\sigma) h_{\overline{\alpha_1}+\overline{\beta_{\sigma(1)}}-n+b-a+1} \dots h_{\overline{\alpha_i}+\overline{\beta_{\sigma(i)}}-n+b-a+1} \nn \\
    & =& \left|
 \begin{array}{cccc}
   h_{\overline{\alpha_1}+\overline{\beta_{1}}-n+b-a+1} &
   h_{\overline{\alpha_1}+\overline{\beta_{2}}-n+b-a+1} & \cdots &
   h_{\overline{\alpha_1}+\overline{\beta_{i}}-n+b-a+1} \\
   h_{\overline{\alpha_2}+\overline{\beta_{1}}-n+b-a+1} &
   h_{\overline{\alpha_2}+\overline{\beta_{2}}-n+b-a+1} & \cdots &
   h_{\overline{\alpha_2}+\overline{\beta_{i}}-n+b-a+1} \\
   \vdots &   &  \ddots & \vdots \\
   h_{\overline{\alpha_i}+\overline{\beta_{1}}-n+b-a+1} &
   h_{\overline{\alpha_i}+\overline{\beta_{2}}-n+b-a+1} & \cdots &
   h_{\overline{\alpha_i}+\overline{\beta_{i}}-n+b-a+1} \\
 \end{array}
 \right|  \nn
 \end{eqnarray}
By permuting the columns and replacing the values of $\overline{\alpha_{\ell}}$ and $\overline{\beta_{\ell'}}$, \eqref{eq_phiinv_det} becomes
 \begin{equation}
    (-1)^{\frac{i(i-1)}{2}}\left|
 \begin{array}{cccc}
   h_{\beta_{1}+\alpha_i-(n+a-b-i)} &
   h_{\beta_{1}+\alpha_{i-1}+1-(n+a-b-i)} & \cdots &
   h_{\beta_{1}+\alpha_{i-1}+(i-1)-(n+a-b-i)} \\
   h_{\beta_{2}+\alpha_i-1-(n+a-b-i)} &
   h_{\beta_{2}+\alpha_{i-1}-(n+a-b-i)} & \cdots &
   h_{\beta_{2}+\alpha_{i-1}+(i-2)-(n+a-b-i)} \\
   \vdots &   &  \ddots & \vdots \\
   h_{\beta_{i}+\alpha_i-(i-1)-(n+a-b-i)} &
   h_{\beta_{i}+\alpha_{i-1}-(i-2)-(n+a-b-i)} & \cdots &
   h_{\beta_{i}+\alpha_{i-1}-(n+a-b-i)} \\
 \end{array}
 \right|  \nn
\end{equation}
\begin{equation} \label{eq_weird_det}
  =  (-1)^{\frac{i(i-1)}{2}} \det\left[
  h_{\beta_s+\alpha_{i+1-t} +t-s-(n+a-b-i)}\right]^i_{s,t=1}
\end{equation}
The determinant in formula \eqref{eq_weird_det} describes the Skew Schur polynomial $\pi_{\beta/{(K_i-\alpha)}}(\underline{x},\underline{y},\underline{z})$,
by \cite[Formula (5.4), page 70]{Mac}.

Moreover we have
\begin{equation}\label{eq_skew_schur}
\pi_{\beta/{(K_i-\alpha)}}(\underline{x},\underline{y},
\underline{z})=
\sum_{\chi}c_{K_i-\alpha,\chi}^\beta\pi_{\chi}(\underline{x},\underline{y},\underline{z})
\end{equation}
Furthermore, by \eqref{eq_simple_ptheta}
\begin{eqnarray}
 \sum_{\theta,\beta} c_{\alpha',\beta}^{\theta} p_{\vartheta}
 &=& (-1)^{\frac{i(i-1)}{2}+\frac{i(i+1)}{2}}\sum_{\theta,\beta} c_{\alpha',\beta}^{\theta}(-1)^{|K_i-\theta|} \pi_{\overline{K_i-\theta}}(\und{x},\und{y},\und{z})
  \\
& =& (-1)^{i}\sum_{\beta,\theta} (-1)^{|\theta|}c_{\alpha',\beta,\theta}^{K_i}
 \pi_{\bar{\theta}}(\underline{x},\underline{y},\underline{z}). \label{eq_pe_theta}
\end{eqnarray}

By replacing \eqref{eq_weird_det}, \eqref{eq_skew_schur} and \eqref{eq_pe_theta} into the general formula (see
\eqref{eq_sigma_lambda2}), and using that the composite homomorphism $\phi_{a-i,b-i,i}^n(\phi_{a-i,b-i,1}^n)^{-1}$ is just multiplication by $(-1)^{i(i-1)/2}$, we have
\begin{eqnarray}
 (\phi_{a-i,b-i,i}^n)^{-1}\left( \;\;
 (-1)^{ab}(-1)^{i(a-b)} \sum_{\theta,\beta} c_{\alpha',\beta}^{\theta} \sum_{\sigma\in S_i}\sgn(\sigma)  \xy
 (-12,0)*{\includegraphics[scale=0.6]{figs/tlonger-up.eps}};
  (12,0)*{\includegraphics[angle=180, scale=0.6]{figs/tlonger-up.eps}};
 (18,-21)*{b-i}; (-18,-21)*{a-i};
 (0,-12)*{\bigb{h_{\overline{\alpha_i}+\overline{\beta_{\sigma(i)}}-n+b-a+1}(1)}};
 (0,4)*{\bigb{  h_{\overline{\alpha_1}+\overline{\beta_{\sigma(1)}}-n+b-a+1}(1)}};
   (0,12)*{\bigb{\hspace{0.5in}p_{\theta}\hspace{0.5in}} };
   (21,0)*{n};(0,-3)*{\vdots};
  \endxy \;\;\right)  =\hspace{1in}
\end{eqnarray}
\begin{equation}
(-1)^{ab}(-1)^{i(a-b)}(-1)^{\frac{i(i-1)}{2}} \sum_{\chi,\theta,\beta}(-1)^{|\theta|}
c_{\alpha',\beta,\theta}^{K_i}c_{K_i-\alpha,\chi}^\beta
\pi_{\bar{\theta}}(\underline{x},\underline{y},\underline{z})
\pi_{\chi}(\underline{x},\underline{y},\underline{z}).\label{1}
\end{equation}
Using that
\begin{equation}
  \sum_{\beta}c_{\alpha',\beta,\theta}^{K_i}c_{K_i-\alpha,\chi}^\beta=
  c_{\alpha',\theta,K_i-\alpha,\chi}^{K_i}
=c_{\alpha',\theta,\chi}^{K_i-(K_i-\alpha)}=\sum_{A}c_{\alpha',A}^{\alpha}c_{\theta,\chi}^A
\end{equation}
equation \eqref{1} becomes
$$(-1)^{(a-i)(b-i)}\sum_A c_{\alpha',A}^{\alpha}\left(\sum_{\theta,\chi}(-1)^{|\theta|}c_{\theta,\chi}^A\pi_{\bar{\theta}}(\underline{x},\underline{y},\underline{z})\pi_{\chi}(\underline{x},\underline{y},\underline{z})\right)=$$
$$=(-1)^{(a-i)(b-i)}\sum_A c_{\alpha',A}^{\alpha} \delta_{A,0}=(-1)^{(a-i)(b-i)}\delta_{\alpha,\alpha'}$$ by thick Grassmannian relation, implying the lemma.

This concludes the proof of equation \eqref{eq_orthog_siglam} showing that $e_{\alpha}^{i}$ and $e_{\alpha'}^j$ are orthogonal unless $i=j$ and $\alpha=\alpha'$, in which case the ladder diagram in \eqref{eq_last_ladder} simplifies to the diagram
\begin{equation}
  (-1)^{(a-i)(b-i)} \xy
 (0,0)*{\includegraphics[scale=0.5]{figs/tr2.eps}};
 (10 ,-12)*{a-i}; (-11 ,-12)*{b-i}; (12,8)*{n};
 \endxy
\end{equation}
\end{proof}

\begin{lem} \label{eq_easy_pullapart}
If $n \geq d-c$, then
\begin{equation} \label{eq_thick_r2}
 \xy
 (0,0)*{\includegraphics[scale=0.5]{figs/tr2.eps}};
 (8 ,-12)*{c}; (-9 ,-12)*{d}; (12,8)*{n};
 \endxy
 \quad = \quad
 (-1)^{cd} \;\;
 \xy
  (-4,-0.5)*{\includegraphics[scale=0.5,angle=180]{figs/tlong-up.eps}};
 (4,0)*{\includegraphics[scale=0.5]{figs/tlong-up.eps}};
 (8 ,-12)*{c}; (-8 ,-12)*{d}; (12,8)*{n};
 \endxy
\end{equation}
and if $n \leq d-c$
\begin{equation} \label{eq_thick_r2p}
 \xy
 (0,0)*{\reflectbox{\includegraphics[scale=0.5]{figs/tr2.eps}}};
 (8 ,-12)*{d}; (-9 ,-12)*{c}; (12,8)*{n};
 \endxy
 \quad = \quad
 (-1)^{cd} \;\;
 \xy
  (-4,0)*{\includegraphics[scale=0.5]{figs/tlong-up.eps}};
 (4,-0.5)*{\includegraphics[scale=0.5,angle=180]{figs/tlong-up.eps}};
 (8 ,-12)*{d}; (-8 ,-12)*{c}; (12,8)*{n};
 \endxy
\end{equation}
\end{lem}

This lemma gives us canonical (up to minus sign) direct summand inclusions and projections
\begin{equation}
  \xymatrix{\cal{F}^{(d)}\cal{E}^{(c)}\onen \ar@<-1ex>@{^{(}->}[r] & \ar@<-1ex>@{->>}[l] \cal{E}^{(c)}\cal{F}^{(d)}\onen } \qquad n \geq d-c,
\end{equation}
and
\begin{equation}
  \xymatrix{\cal{E}^{(c)}\cal{F}^{(d)}\onen \ar@<-1ex>@{^{(}->}[r] & \ar@<-1ex>@{->>}[l] \cal{F}^{(d)}\cal{E}^{(c)}\onen } \qquad n \leq d-c.
\end{equation}

\begin{proof}
We explode the left-hand-side of \eqref{eq_thick_r2}
\begin{eqnarray}
 \xy
 (0,0)*{\includegraphics[scale=0.5]{figs/tr2.eps}};
 (8 ,-12)*{c}; (-9 ,-12)*{d}; (12,8)*{n};
 \endxy \;\; := \;\;
 \xy
 (0,0)*{\includegraphics[scale=0.5]{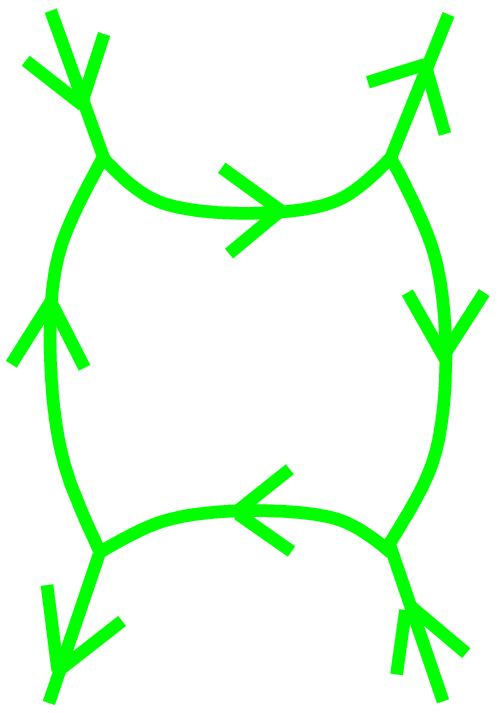}};
 (-14,-16)*{d};(14,-16)*{c};
 (-14,16)*{d};(14,16)*{c};
 (-14,-4)*{c};(14,-4)*{d}; (0,-12)*{c+d}; (0,3)*{c+d};
 (15,8)*{n};
  \endxy
  \;\; = \;\;
 \xy
 (0,0)*{\includegraphics[scale=0.5]{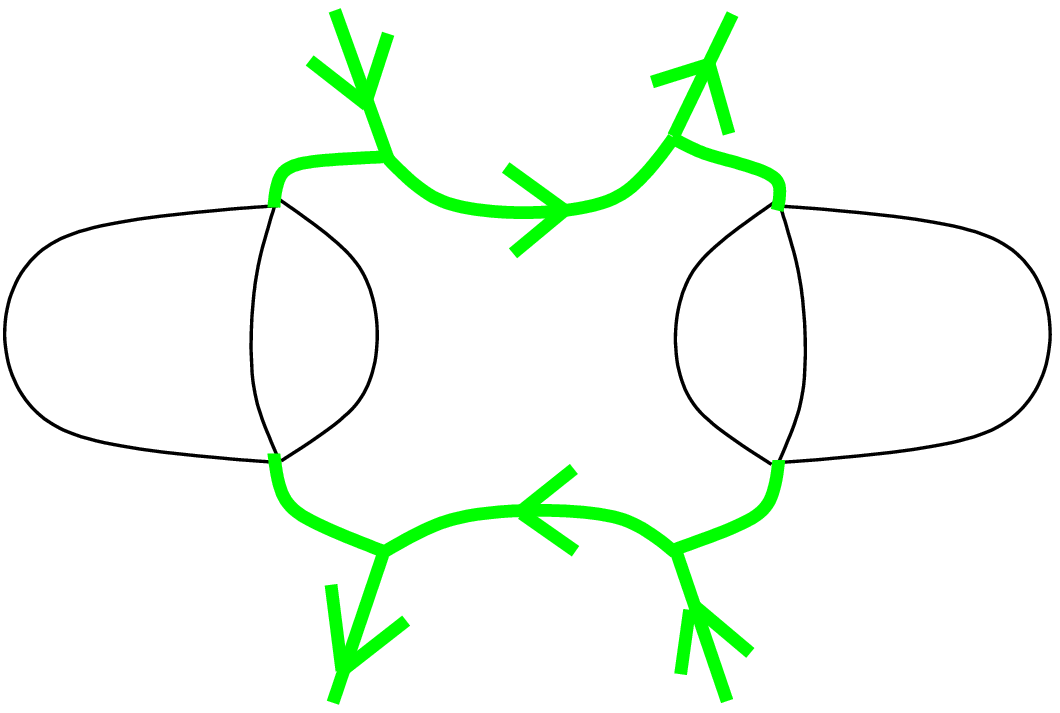}};
 (-14,-16)*{d};(14,-16)*{c};
 (-14,16)*{d};(14,16)*{c};(0,-12)*{c+d}; (0,3)*{c+d};
 (-26.5,1)*{\bullet}+(-4,2)*{\scs c-1};
 (-14,1)*{\bullet}+(-3,2)*{\scs 1};
 (26.5,1)*{\bullet}+(4,2)*{\scs d-1};
 (14,1)*{\bullet}+(3,2)*{\scs 1};
 (-20,-1)*{\cdots};(20,-1)*{\cdots};
 (22,-11)*{n};
  \endxy
\end{eqnarray}
Repeatedly using associativity and coassociativity of splitters this is equal to the ladder diagram
\begin{equation} \label{eq_sideways_ladder}
 = \;\;  \xy
 (0,0)*{\includegraphics[scale=0.6]{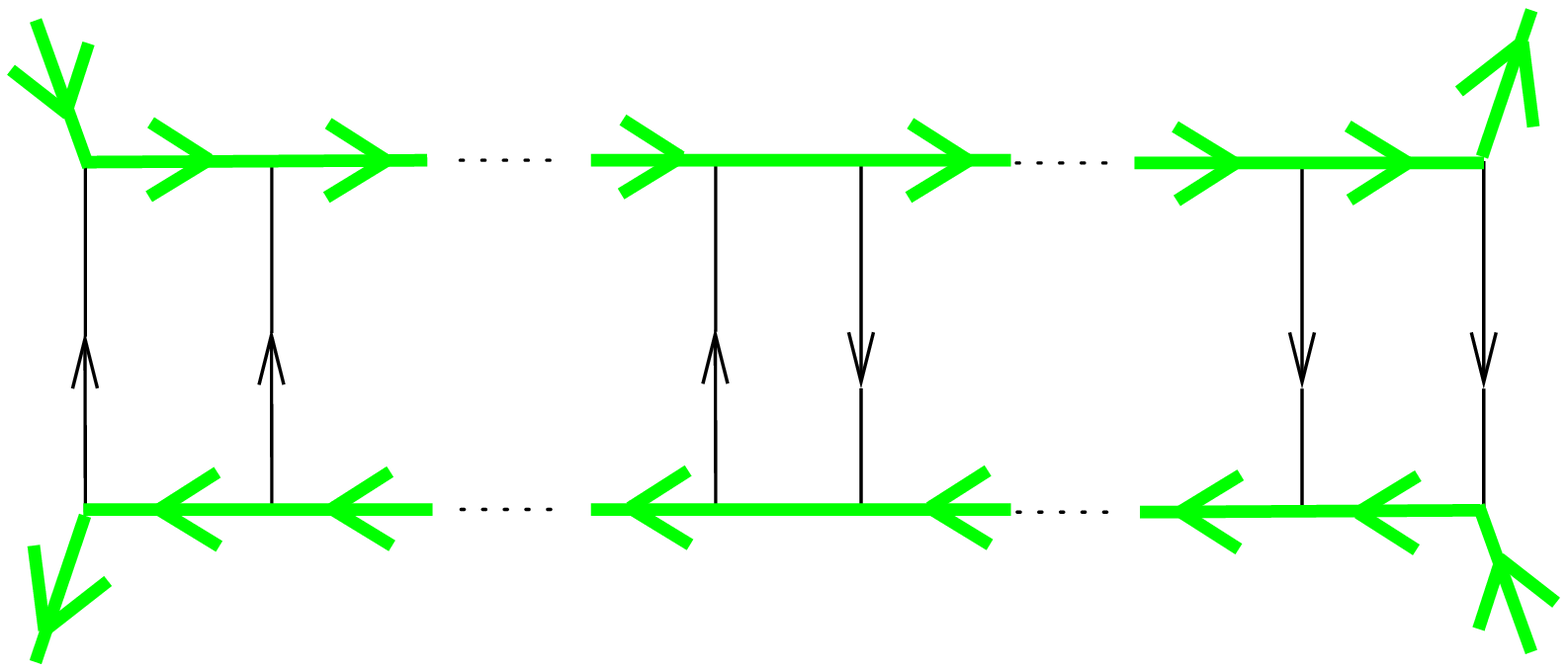}};
 (-50,-16)*{d};(50,-16)*{c};
 (-50,16)*{d};(50,16)*{c};(0,-15)*{c+d}; (0,15)*{c+d};
 (37,15)*{c+1};(37,-15)*{c+1};(-37,15)*{d+1};(-37,-15)*{d+1};
 (-43,3)*{\bullet}+(-4,2)*{\scs c-1};
 (-32,3)*{\bullet}+(-3,2)*{\scs c-2};
 (43,3)*{\bullet}+(4,2)*{\scs d-1};
 (32,3)*{\bullet}+(3,2)*{\scs d-2};
 (-20,-1)*{\cdots};(20,-1)*{\cdots};
 (56,-5)*{n};
  \endxy
\end{equation}
Again we will simplify the innermost square diagram using the Square Flop Lemma.  In this case the Square Flop Lemma simplifies since for $t\geq 1$ and $x$, $y$ and $m$ such that $x+y\leq m-2$ the equality
\begin{equation}
   \xy
 (0,0)*{\includegraphics[scale=0.5]{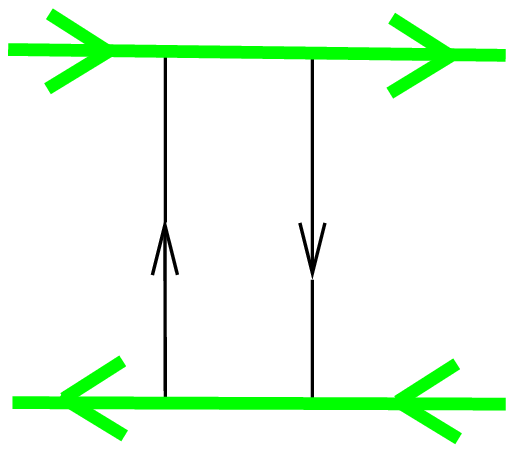}};
 (-13,-13)*{t};(13,-13)*{t};
 (-13,13)*{t};(13,13)*{t};(0,-12)*{t+1}; (0,12)*{t+1};
 (-5,3)*{\bullet}+(-3,2)*{\scs x};
 (3,3)*{\bullet}+(3,2)*{\scs y};
 (4,18)*{m};
  \endxy
  \quad = \quad -\;\;
   \xy
 (0,0)*{\includegraphics[scale=0.5]{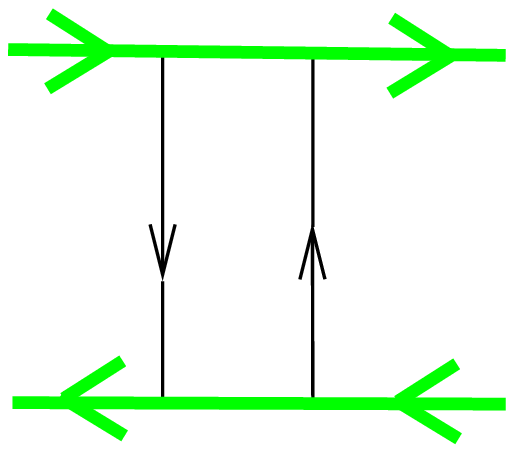}};
 (-13,-13)*{t};(13,-13)*{t};
 (-13,13)*{t};(13,13)*{t};(0,-12)*{t-1}; (0,12)*{t-1};
 (-5,3)*{\bullet}+(-3,2)*{\scs y};
 (3,3)*{\bullet}+(3,2)*{\scs x};
 (4,18)*{m};
  \endxy
\end{equation}
holds for the additional terms arising in the Square Flop Lemma are always zero by our conventions for $h_r$ with $r<0$.  In \eqref{eq_sideways_ladder} we have $x \leq c-1$, $y \leq d-1$ and $m=n+2c$ and since $n \geq d-c$ by assumption, the condition $x+y \leq m-2$ is satisfied.  Hence, repeatedly applying the Square Flop Lemma in this simplified setting produces the ladder diagram
\begin{equation}
 =  (-1)^{cd} \xy
 (0,0)*{\includegraphics[scale=0.6]{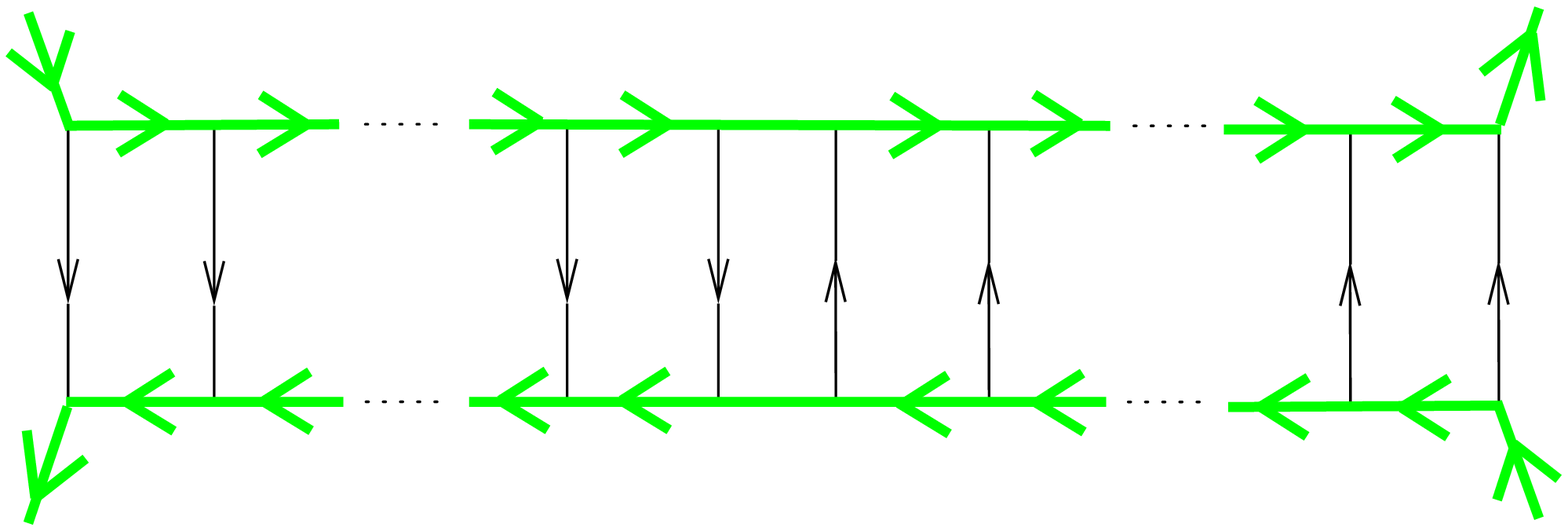}};
 (-63,-16)*{d};(63,-16)*{c};
 (-63,16)*{d};(63,16)*{c};
 (0,-15)*{0}; (0,15)*{0};
 (-13,-15)*{1}; (-13,15)*{1};
 (10,15)*{1}; (10,-15)*{1};
 (48,15)*{c-1};(48,-15)*{c-1};
 (-49,15)*{d-1};(-49,-15)*{d-1};
 (4,3)*{\bullet}+(4,2)*{\scs c-1};
 (16,3)*{\bullet}+(4,2)*{\scs c-2};
 (-5,3)*{\bullet}+(-4,2)*{\scs d-1};
 (-17,3)*{\bullet}+(-4,2)*{\scs d-2};
 (-44.5,3)*{\bullet}+(-3,2)*{\scs 1};
 (44,3)*{\bullet}+(3,2)*{\scs 1};
 (62,-5)*{n};
  \endxy \nn
\end{equation}
\begin{equation}
  = (-1)^{cd} \;\;
    \xy
 (0,-1)*{\includegraphics[scale=0.5,angle=180]{figs/texplode.eps}};
 (-4.5,0)*{\bullet}+(-3,1)*{\scs 1};
 (4.5,0)*{\bullet}+(3.5,1)*{\scs d-2};
 (14,0)*{\bullet}+(3.5,1)*{\scs d-1};
 (0,-1)*{\cdots};  (-3,-12)*{d};
  \endxy\quad
  \xy
 (0,0)*{\includegraphics[scale=0.5]{figs/texplode.eps}};
 (-14,0)*{\bullet}+(-3,1)*{\scs c-1};
 (-4.5,0)*{\bullet}+(-3,1)*{\scs c-2};
 (4.5,0)*{\bullet}+(2.5,1)*{\scs 1};
 (0,-2)*{\cdots};  (-3,-12)*{c};
  \endxy
  \quad = \quad
 (-1)^{cd} \;\;
 \xy
  (-4,-0.5)*{\includegraphics[scale=0.5,angle=180]{figs/tlong-up.eps}};
 (4,0)*{\includegraphics[scale=0.5]{figs/tlong-up.eps}};
 (8 ,-12)*{c}; (-9 ,-12)*{d}; (12,8)*{n};
 \endxy
\end{equation}
completing the proof of the first claim.  The second claim is proven similarly.
\end{proof}

\begin{cor} \label{lem_orthogonal}
Given $0\le i,j \le \min(a,b)$ and
partitions $\alpha\in P(i,n+a-b-i)$ and $\alpha'\in P(j,n+a-b-j)$, then
\begin{equation}
\sigma_{\alpha'}^{j}\lambda_{\alpha}^i =
\delta_{\alpha,\alpha'}\delta_{i,j}\Id_{\cal{F}^{(b-i)}\cal{E}^{(a-i)}\onen}.
\end{equation}
and
\begin{equation} \label{orte}
  e_{\alpha}^i e_{\alpha'}^{j} = \delta_{i,j} \delta_{\alpha,\alpha'}e_{\alpha}^i,
\end{equation}
so that the 2-morphisms $e_{\alpha}^i$ are orthogonal
idempotents.
\end{cor}

\begin{cor} \label{cor_n_b-a}
There is a canonical (up to sign) isomorphism
$\cal{E}^{(a)}\cal{F}^{(b)}\onenn{b-a} \cong
\cal{F}^{(b)}\cal{E}^{(a)}\onenn{b-a}$.
\end{cor}

\begin{proof}
The isomorphism is given by the following maps
\begin{equation}
  \xy
   (-20,0)*{\cal{E}^{(a)}\cal{F}^{(b)}\onenn{b-a}}="1";
   (20,0)*{\cal{F}^{(b)}\cal{E}^{(a)}\onenn{b-a}}="2";
  {\ar^{(-1)^{ab}
     \xy (0,0)*{\reflectbox{\includegraphics[scale=0.35]{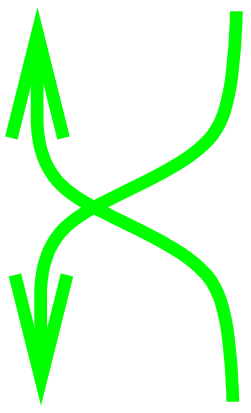}}};
   (-5,-5)*{\scs a}; (5,-5)*{\scs b}; \endxy} "1";"2"};
    \endxy
  \qquad
  \qquad
    \xy
   (20,0)*{\cal{E}^{(a)}\cal{F}^{(b)}\onenn{b-a}}="1";
   (-20,0)*{\cal{F}^{(b)}\cal{E}^{(a)}\onenn{b-a}}="2";
   {\ar^{\xy (0,0)*{\includegraphics[scale=0.35]{figs/tcross-side.eps}};
   (-5,-5)*{\scs b}; (5,-5)*{\scs a}; \endxy} "2";"1"};
  \endxy
\end{equation}
which are mutually-inverse isomorphisms by the previous lemma.
\end{proof}

\begin{thm}[Sto\v si\'c Formula] \label{thm_EaFb} There is an
equality
\begin{equation} \label{decef}
 \xy
 (-4,0)*{\includegraphics[scale=0.5]{figs/tlong-up.eps}};
  (4,-1)*{\includegraphics[angle=180, scale=0.5]{figs/tlong-up.eps}};
 (-6 ,-5)*{a}; (6 ,-5)*{b}; (10,10)*{n};
 \endxy\quad = \quad
 (-1)^{ab}\sum_{i=0}^{\min(a,b)}
  \sum_{\alpha,\beta,\gamma,x,y} (-1)^{\frac{i(i+1)}{2}+|x|+|y|} c_{\alpha,\beta,\gamma,x,y}^{K_i}\;\;
 \xy
  (0,0)*{\includegraphics[scale=0.5]{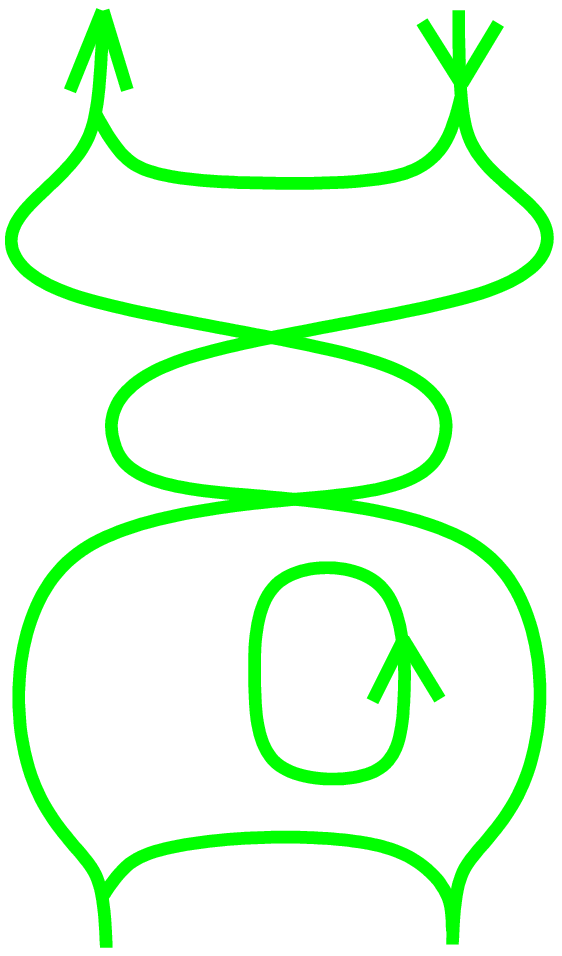}};
  (-2,-10)*{\bigb{\pi_{\gamma}^{\spadesuit} }};(-2,-19)*{\bigb{\pi_{\beta} }}; (-12,-10)*{\bigb{\pi_{\overline{x}} }};(13,-10)*{\bigb{\pi_{\overline{y}} }};
  (-2,16)*{\bigb{\pi_{\alpha} }};
  (-7,-17)*{i};  (7,-15)*{i}; (4,13)*{i};
  (18,16)*{n};  (13,2)*{a-i};(-14,2)*{b-i};(12,-23)*{b};(-12,-23)*{a};
  (13,21)*{b};(-13,21)*{a};
 \endxy
\end{equation}
where the sum is over all partitions $\alpha,\beta,\gamma \in P(i)$,
$x\in P(i,a-i)$,  $y \in P(i,b-i)$;  $K_0=\emptyset$, and $K_i=((n+a-b-i)^i)$ for $1\leq i \leq \min(a,b)$.
\end{thm}

\begin{rem}
Recall that $P(0)$ is the set of all partitions with at most $0$ parts.  In particular, $P(0)$ contains only the empty partition 0. Then
$c_{\alpha,\beta}^{\emptyset}=\delta_{\alpha,\emptyset}\delta_{\beta,\emptyset}$.  Hence, the $i=0$ term of \eqref{decef} is just
\begin{equation}
  (-1)^{ab}  \xy
 (0,0)*{\reflectbox{\includegraphics[scale=0.5]{figs/tr2.eps}}};
 (8 ,-12)*{b}; (-9 ,-12)*{a}; (12,8)*{n};
 \endxy
\end{equation}
When $n+a-b \leq 0$ then $E^{(a)}F^{(b)}1_n$ is a canonical basis vector and equation \eqref{decef} reduces to
\begin{equation}
 \xy
 (-4,0)*{\includegraphics[scale=0.5]{figs/tlong-up.eps}};
  (4,-1)*{\includegraphics[angle=180, scale=0.5]{figs/tlong-up.eps}};
 (-6 ,-5)*{a}; (6 ,-5)*{b}; (10,10)*{n};
 \endxy\quad = \quad
  (-1)^{ab}  \xy
 (0,0)*{\reflectbox{\includegraphics[scale=0.5]{figs/tr2.eps}}};
 (8 ,-12)*{b}; (-9 ,-12)*{a}; (12,8)*{n};
 \endxy
\end{equation}
since if $n+a-b-i<0$ and $i > 0$ the corresponding term should be omitted from the sum ($K_i$ would not make sense).
\end{rem}

\begin{proof}
Lemma \ref{lem_orthogonal} gives a collection of mutually orthogonal idempotents
\begin{equation}
  e_{\alpha}^i = \lambda_{\alpha}^i \sigma_{\alpha}^i \maps \E{a}\F{b}\onen \longrightarrow
  \E{a}\F{b} \onen
\end{equation}
for $0\le i \le \min(a,b)$ and $\alpha\in P(i,n+a-b-i)$ which are projections onto direct summands isomorphic to $\cal{F}^{(b-i)}\cal{E}^{(a-i)}\onen\{2|\alpha|-i(a-b+n)\}$.

When the ground field $\Bbbk=\Q$, it follows from \cite{Lau1} that $\cal{E}^{(a)}\cal{F}^{(b)}\onen$ is isomorphic to the direct sum of these summands over all values of the parameters $i$ and $\alpha$.  Therefore, equality \eqref{decef} holds over $\Q$, and, then, over $\Z$, since all coefficients are integers.
\end{proof}

For $i \in \{ 0,1, \dots, \min(a,b)\}$ and $\alpha\in P(i,-n+a-b-i)$
define maps
\begin{eqnarray}
 \bar{\lambda}_{\alpha}^i &:=&   \xy
  (0,0)*{\reflectbox{\includegraphics[scale=0.5]{figs/lambda-alpha.eps}}};
 (-2,4)*{\bigb{\pi_{\alpha} }};  (4,2)*{i};
  (15,0)*{}; (18,7)*{n};  (-14,-12)*{a-i};(15,-12)*{-b-i};(-13,12)*{b};(13,12)*{a};
 \endxy \maps \E{a-i}\F{b-i}\onen \{ 2|\alpha|-i(-n+a-b-i) \} \longrightarrow \F{b}\E{a} \onen \nn \\ \\
\bar{\sigma}_{\alpha}^i  &:=&   (-1)^{ab+i(a+b)}\sum_{\beta, \gamma, x, y}(-1)^{\frac{i(i+1)}{2}+|x|+|y|}
c_{\alpha,\beta,\gamma,x, y}^{K_i}\xy
  (-2,0)*{\reflectbox{\includegraphics[scale=0.5]{figs/idemp1.eps}}};
   (1,-2)*{\bigb{\pi_{\gamma}^{\spadesuit} }};(0,-10)*{\bigb{\pi_{\beta} }}; (10,-2)*{\bigb{\pi_{\overline{x}} }};
   (-15,-2)*{\bigb{\pi_{\overline{y}} }}; (-8,-7)*{i};  (8,-7)*{i};
  (15,0)*{}; (18,6)*{n};  (-16,12)*{a-i};(14,12)*{b-i};(-12,-13)*{b};(10,-13)*{a};
 \endxy \nn \\
 & & \hspace{1in} \maps \F{b}\E{a}\onen \longrightarrow \E{a-i}\F{b-i}\onen\{2|\alpha|-i(-n+a-b-i) \}.
\end{eqnarray}
where $K_0=\emptyset$ and $K_i=((-n+a-b-i)^i)$.  One can show analogously to Lemma~\ref{lem_orthogonal-pre} that
\begin{equation}
\bar{\sigma}_{\alpha'}^{j}\bar{\lambda}_{\alpha}^i =
\delta_{\alpha,\alpha'}\delta_{i,j}\Id_{\cal{E}^{(a-i)}\cal{F}^{(b-i)}\onen},
\end{equation}
so that if we define
\begin{equation}
  \bar{e}_{\alpha}^i = \bar{\lambda}_{\alpha}^i \bar{\sigma}_{\alpha}^i \maps \E{a}\F{b}\onen \longrightarrow
  \E{a}\F{b} \onen,
\end{equation}
then we have mutually orthogonal idempotents
\begin{equation} \label{orte2}
  \bar{e}_{\alpha}^i \bar{e}_{\alpha'}^{j} = \delta_{i,j} \delta_{\alpha,\alpha'}\bar{e}_{\alpha}^i.
\end{equation}

In \cite[Section 5.6]{Lau1} symmetries of the 2-category $\Ucat$ were defined that were shown to extend to invertible 2-functors on $\UcatD$.  We refer the reader there for details, but briefly recall the symmetry 2-functor $\tilde{\sigma}\maps \Ucat \to \Ucat$ that acts on diagrams by rescaling the crossing $\Ucross
\mapsto -\Ucross$ for all $n \in \Z$, reflecting a diagram across the vertical axis, and sending $n$ to $-n$.  Notice that
\begin{equation}
\tilde{\sigma} \left( \; \;\;
  \xy
 (0,0)*{\includegraphics[scale=0.5]{figs/c2-1.eps}};
 (0,-1.5)*{e_{a}};  (16,2)*{-n};
  \endxy \; \right)
\quad := \quad
 (-1)^{\frac{a(a-1)}{2}}\xy
 (0,0)*{\reflectbox{\includegraphics[scale=0.5]{figs/c2-2.eps}}};
 (17.7,5)*{\bullet}+(4,1)*{\scs a-1};
 (7.7,5)*{\bullet}+(4,1)*{\scs a-2};
 (3,4)*{\cdots};(3,-7)*{\cdots};
 (-2.3,5)*{\bullet}+(2,1)*{\scs 2};
 (-9.9,5)*{\bullet};(0,-1.5)*{D_a}; (25,2)*{n};
  \endxy
\end{equation}
and that the idempotent $e_a\onenn{-n}$ is equivalent to the idempotent $\tilde{\sigma}(e_a\onen)$ since $\tilde{\sigma}(e_a\onen) e_a\onenn{-n} = \tilde{\sigma}(e_a\onen)$ and $e_a\onenn{-n} \tilde{\sigma}(e_a\onen)=e_a\onenn{-n}$.  Since these idempotents are equivalent they give rise to isomorphic 1-morphisms in $\UcatD$.

By writing Theorem~\ref{thm_EaFb} at $-n$ using the definition of the thick calculus we can apply the 2-functor $\tilde{\sigma}$ and obtain a new equality for $n$.  The diagrams in the resulting equation are not immediately related to the the thick calculus.  However, by composing both sides of the equality with the 2-morphism $\tau'(e_b)\Id_{\cal{E}^a\onen}$ on the bottom and $\Id_{\cal{F}^b\onenn{n+2a}}e_a$ on the top of both sides of the equality, one can then use \eqref{eq_split_bs} to flip the order of the dots appearing in the middle of the diagrams on the right hand side so that the resulting equation can be expressed in terms of the thick calculus.  Keeping careful track of the signs one obtains the following corollary.

\begin{cor}
  \begin{equation} \label{decfe}
 \xy
 (4,0)*{\includegraphics[scale=0.5]{figs/tlong-up.eps}};
  (-4,-1)*{\includegraphics[angle=180, scale=0.5]{figs/tlong-up.eps}};
 (6 ,-5)*{a}; (-6 ,-5)*{b}; (10,10)*{n};
 \endxy =
 (-1)^{ab}\sum_{i=0}^{\min(a,b)}(-1)^{i(a+b)+\frac{i(i+1)}{2}}
  \sum_{\alpha,\beta,\gamma,x,y} (-1)^{|x|+|y|} c_{\alpha,\beta,\gamma,x,y}^{K}\;\;
 \xy
  (0,0)*{\reflectbox{\includegraphics[scale=0.5]{figs/EaFb.eps}}};
  (2,-10)*{\bigb{\pi_{\gamma}^{\spadesuit} }};(2,-19)*{\bigb{\pi_{\beta} }}; (12,-10)*{\bigb{\pi_{\overline{x}} }};(-13,-10)*{\bigb{\pi_{\overline{y}} }};
  (2,16)*{\bigb{\pi_{\alpha} }};
  (7,-17)*{i};  (-7,-15)*{i}; (-4,13)*{i};
  (18,16)*{n};  (-13,2)*{a-i};(14,2)*{b-i};(-12,-23)*{b};(12,-23)*{a};
  (-13,21)*{b};(13,21)*{a};
 \endxy
\end{equation}
where the sum is over all partitions $\alpha,\beta,\gamma \in P(i)$,
$x\in P(i,a-i)$,  $y \in P(i,b-i)$, $K_0=\emptyset$, and $K_i=((-n+a-b-i)^i)$ for $1\leq i \leq \min(a,b)$.
\end{cor}

\begin{thm} \label{eq_cat_EaFb}
For $n \geq b-a$ maps
\begin{eqnarray}
 \sum_{j=0}^{\min(a,b)}\sum_{\alpha} \sigma_{\alpha}^j &\maps &
 \bigoplus_{j=0}^{\min(a,b)} \bigoplus_{\alpha \in P(j,n+a-b-j)}\cal{F}^{(b-j)}\cal{E}^{(a-j)}\onen\{2|\alpha|-j(a-b+n)\} \longrightarrow \cal{E}^{(a)}\cal{F}^{(b)}\onen  \nn \\
 \nn \\
 \sum_{j=0}^{\min(a,b)}\sum_{\alpha} \lambda_{\alpha}^j & \maps &  \cal{E}^{(a)}\cal{F}^{(b)}\onen \longrightarrow \bigoplus_{j=0}^{\min(a,b)} \bigoplus_{\alpha \in P(j,n+a-b-j)}\cal{F}^{(b-j)}\cal{E}^{(a-j)}\onen\{2|\alpha|-j(a-b+n)\} \nn
\end{eqnarray}
are mutually-inverse isomorphisms, giving a canonical isomorphism
\begin{equation}
   \E{a}\F{b}\onen \cong \bigoplus_{j=0}^{\min(a,b)} \bigoplus_{\alpha \in P(j,n+a-b-i)}\cal{F}^{(b-j)}\cal{E}^{(a-j)}\onen\{2|\alpha|-j(n+a-b)\}
\end{equation}
in $\UcatD$.  Likewise, for $n \leq b-a$ maps
\begin{eqnarray}
 \sum_{j=0}^{\min(a,b)}\sum_{\alpha} \bar{\sigma}_{\alpha}^j &\maps &
 \bigoplus_{j=0}^{\min(a,b)} \bigoplus_{\alpha \in P(j,-n+a-b-j)}\cal{E}^{(a-j)}\cal{F}^{(b-j)}\onenn{n}\{2|\alpha|-j(a-b-n)\} \longrightarrow \cal{F}^{(b)}\cal{E}^{(a)}\onenn{n}  \nn \\
 \nn \\
 \sum_{j=0}^{\min(a,b)}\sum_{\alpha}  \bar{\lambda}_{\alpha}^j & \maps &  \cal{F}^{(b)}\cal{E}^{(a)}\onenn{n} \longrightarrow \bigoplus_{j=0}^{\min(a,b)} \bigoplus_{\alpha \in P(j,-n+a-b-j)}\cal{E}^{(a-j)}\cal{F}^{(b-j)}\onenn{n}\{2|\alpha|-j(a-b-n)\} \nn
\end{eqnarray}
are mutually-inverse isomorphisms, giving a canonical isomorphism
\begin{equation}
 \F{b} \E{a}\onen \cong \bigoplus_{j=0}^{\min(a,b)} \bigoplus_{\alpha \in P(j,-n+b-a-j)}\cal{E}^{(a-j)}\cal{F}^{(b-j)}\onen\{2|\alpha|-j(b-a-n)\}.
\end{equation}
\end{thm}

We end this section with a useful relation for simplifying curl diagrams in $\UcatD$.

\begin{prop}[Higher reduction to bubbles]
For $\beta\in P(b)$, $n \in \Z$ we
have
\begin{equation} \label{eq_thick_red}
 \xy
  (0,0)*{\includegraphics[scale=0.5]{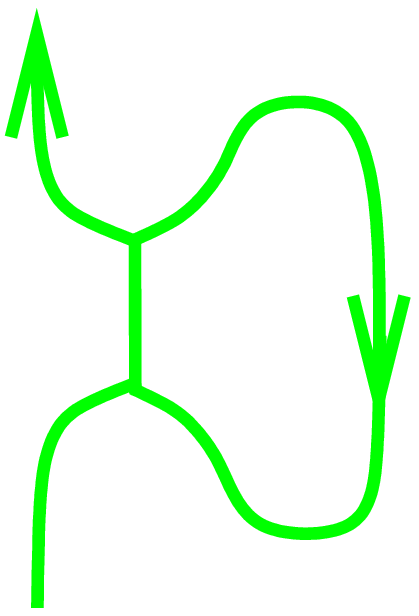}};
  (-11,-11)*{a}; (0,-11)*{b};  (-9,0)*{a+b};  (15,0)*{}; (15,11)*{n};
  (8,5)*{\bigb{\pi_{\beta} }};
 \endxy \quad = \quad (-1)^{ab}
\sum_{\gamma,\delta} c_{\gamma,\delta}^{\beta-(n+a-b)} \;\;\;\;
 \xy
 (-4,0)*{\includegraphics[scale=0.5]{figs/tlong-up.eps}};
 (6,0)*{\stcbub{b}{\delta}};
 (-6 ,-9)*{a}; (10,10)*{n}; (-4,0)*{\bigb{\pi_{\gamma} }};
 \endxy
\end{equation}
\begin{equation}
 \xy
  (0,0)*{\reflectbox{\includegraphics[scale=0.5]{figs/treduction.eps}}};
   (11,-11)*{a}; (0,-11)*{b};  (9,0)*{a+b};  (-15,0)*{}; (-15,11)*{n};
  (-8,5)*{\bigb{\pi_{\beta} }};
 \endxy \quad = \quad
\sum_{\gamma,\delta} c_{\gamma,\delta}^{\beta-(-n+a-b)} \;\;\;\;
 \xy
 (4,0)*{\includegraphics[scale=0.5]{figs/tlong-up.eps}};
 (-6,0)*{\stccbub{b}{\delta}};
 (6 ,-9)*{a}; (-10,10)*{n}; (4,0)*{\bigb{\pi_{\gamma} }};
 \endxy
\end{equation}
\end{prop}

\begin{proof}
The proof is very similar to the proof of Theorem~\ref{thm_EaFb}.
Again, using Lemmas \ref{lem_ladder_slide1},
\ref{lem_square_flop}, and \ref{lem_bottom_zero} the result follows
by a degenerate version of the argument in Theorem~\ref{thm_EaFb}.
\end{proof}

\begin{rem}
In \eqref{eq_thick_red} if $\beta-(n+a-b)$ is not defined, there are no terms on the right hand side, and the left hand side is equal to zero.
\end{rem}

To each $x \in \B$ we associate a 1-morphism in $\UcatD$:
\begin{equation}
  x \mapsto \cal{E}(x) := \left\{
\begin{array}{cl}
  \cal{E}^{(a)}\cal{F}^{(b)}\onen & \text{if $x=E^{(a)}F^{(b)}1_n$,} \\
  \cal{F}^{(b)}\cal{E}^{(a)}\onen & \text{if $x=F^{(b)}E^{(a)}1_n$.}
\end{array}
  \right.
\end{equation}
By Corollary~\ref{cor_n_b-a}, when $n=b-a$  we have a canonical (up to sign) isomorphism $\cal{E}^{(a)}\cal{F}^{(b)}\onenn{b-a} \cong \cal{F}^{(b)}\cal{E}^{(a)}\onenn{b-a}$ allowing us to switch between the 1-morphisms in this case.

Defining relations \eqref{eq_AUrel1}--\eqref{eq_EaEb2}, \eqref{eq_FbEa1} for $n \geq b-a$, and \eqref{eq_EaFb1} for $n \leq b-a$ have all been categorified, see Theorems~\ref{thm_cal_EaEb} and \ref{eq_cat_EaFb}. Therefore, we have a homomorphism of $\Z[q,q^{-1}]$-modules
\begin{eqnarray} \label{def_gamma}
\gamma\maps  \UA & \longrightarrow& K_0(\UcatD) \\
   x & \mapsto & [\cal{E}(x)],
\end{eqnarray}
where $\UcatD$ is defined over a commutative ring $\Bbbk$.

 \subsection{Indecomposables over $\Z$}

Let $A$ be a graded ring, $A=\bigoplus_{i\in \Z} A^i$,
with the grading bounded from below, equipped with a
decomposition of $1$ into a sum of  degree $0$ orthogonal idempotents:
$$ 1 = \sum_{\alpha} 1_{\alpha}, \quad \alpha\in I, \quad \quad 1_{\alpha}1_{\beta}=
\delta_{\alpha,\beta} 1_{\alpha}, \quad \deg(1_{\alpha})=0,$$
for a finite index set $I$. We can decompose $A$ as a graded abelian group
\begin{equation} \label{A_decomp}
  A = \bigoplus_{\alpha,\beta} {}_{\alpha} A_{\beta}
\end{equation}
where ${}_{\alpha} A_{\beta} = 1_{\alpha} A1_{\beta}$.
Assume furthermore that ${}_{\alpha} A^0_{\alpha}=\Z$,
 ${}_{\alpha} A^{<0}_{\alpha}=0$ and
$$ {}_{\alpha}A_{\beta} \cdot {}_{\beta}A_{\alpha} \subset {}_{\alpha} A^{>0}_{\alpha} $$
for all $\alpha\not= \beta$. Here $A^{>0}= \bigoplus_{i>0}A^i$, etc.
Let
\begin{equation}
  A_{\alpha} = \bigoplus_{\beta} {}_{\beta}A_{\alpha}.
\end{equation}
Then $A_{\alpha}$
is an indecomposable graded projective left $A$-module.

\begin{prop} \label{prop_MK1}
Under these assumptions, any finitely-generated indecomposable graded projective
left $A$-module $P$ is isomorphic to a direct sum of modules $A_{\alpha}$ with multiplicities
being Laurent polynomials in $q$:
 $$ P \cong \bigoplus_{\alpha\in I} A_{\alpha}^{f_{\alpha}}, \quad \quad f_{\alpha} \in \Z_+[q,q^{-1}].$$
The multiplicities $f_{\alpha}$ are invariants of $P$.
\end{prop}

\begin{proof}
To prove this result, let
\begin{equation}
  J(A) = (\bigoplus_{\alpha \not= \beta} {}_{\alpha} A_{\beta}) \oplus
 (\bigoplus_{\alpha} {}_{\alpha} A_{\alpha}^{>0}).
\end{equation}
Then $J(A)$ is the graded Jacobson radical of $A$, and it is locally nilpotent:
for any $m\in \Z$ there exist $N\in \Z_+$ such that
 $J(A)^N \cap A^m = 0$. The quotient $A/J(A)$ is the graded ring
$\Z\times \Z \times \dots \times \Z = \Z^I$. Uniqueness of decomposition of
graded projectives and invariance of multiplicities for this ring is clear, and
that for $A$ follows by standard arguments, as in~\cite[Chapter 1]{Benson}.
\end{proof}

An analogous result holds when the index set $I$ is infinite. Then $A$ is a nonunital
idempotented ring, with the decomposition \eqref{A_decomp}, and finitely-generated graded projective
modules are defined as direct summands of finite direct sums of $A_{\alpha}\{i\}$
for various $\alpha$ and $i$.

We now apply this result for infinite $I$ to the classification of 1-morphisms of $\UcatD$ when
the ground ring $\Bbbk$ is $\Z$ rather than a field. When $\Bbbk$ is a field,
it was shown in~\cite[Proposition 9.10]{Lau1} that each category $_m\UcatD_n:=\UcatD(n,m)$ has the unique decomposition
property, with isomorphism classes of indecomposable objects represented by
$\cal{E}(x)\{i\}$, for $b$ in the Lusztig canonical basis $_{m}\B_n$ of $_m\U_n$ and $\cal{E}(x)$ the object of $_m\UcatD_n$ associated to $x$ via \eqref{def_gamma}.

\begin{prop} \label{prop_unique_decomp}
Let $\Bbbk=\Z$. Then categories $_m\UcatD_n$ possess the unique
decomposition property, with isomorphism classes of indecomposable objects
represented by $\cal{E}(x)\{i\}$, over $x\in \B$ and $i\in \Z$.
\end{prop}

\begin{proof}
Fix $n,m$ and consider the category $_m\UcatD_n$. It is the Karoubi
envelope of the category $_m\Ucat_n:=\Ucat(n,m)$, whose morphisms are finite direct sums of products $\cal{E}_{\ep} \onen$. Repeatedly using decompositions in Theorems~\ref{thm_cal_EaEb} and \ref{thm_EaFb}, which hold over $\Z$, and applying Lemma~\ref{lem_triple_canonical},  we can realize $\cal{E}_{\ep} \onen$
as a direct sum of $\cal{E}(x)\{i\}$ for $x\in {}_m\B_n$ and $i \in \Z$. Define
\begin{equation}
  \cal{E}_{m,n} := \bigoplus_{x \in {}_m\B_n} \cal{E}(x).
\end{equation}
Then the category
$_m\UcatD_n$ is graded Morita equivalent to the category of graded
finitely-generated projective modules over the endomorphism ring
\begin{equation}
  R_{m,n}:= \END_{\UcatD} \left( \cal{E}_{m,n}\right).
\end{equation}
This ring is nonunital, with a family of orthogonal idempotents $1_{\cal{E}(x)}$,
over $x\in {}_m\B_n$. That this ring satisfies the conditions of Proposition~\ref{prop_MK1} above (more precisely, of its generalization to nonunital rings) was shown
in \cite[Propositions 9.9, 9.10]{Lau1} in the case of $\Bbbk$ being a field. Since $\Z\subset \Q$,
the ring $R_{m,n}$ for $\Bbbk=\Z$ is a subring of $R_{m,n}$ for $\Bbbk=\Q$,
as implied by the nondegeneracy of the graphical calculus, see \cite{Lau1}. Thus,
the assumptions of Proposition~\ref{prop_MK1} (for infinite $I$) hold for it as well. Therefore, any
indecomposable finitely-generated graded projective module over $R_{m,n}$
is isomorphic to $R_{m,n}1_{\cal{E}(x)}\{i\}$ for a unique $x$ and $i$, and any indecomposable object of $_m\UcatD_n$ is isomorphic to
$\cal{E}(x) \{i\}$ for a unique $x\in {}_m\B_n$ and $i\in \Z$.
\end{proof}

\begin{cor} \label{cor_overZ}
The homomorphism $\gamma \maps \UA \to K_0(\UcatD)$ in \eqref{def_gamma} is an isomorphism when $\Bbbk=\Z$.
\end{cor}

Previously this was shown for $\Bbbk$ a field in \cite[Theorem 9.13]{Lau1}. Proposition~\ref{prop_unique_decomp} and Corollary~\ref{cor_overZ} hold more generally, for any Noetherian commutative ring $\Bbbk$ such that any finitely-generated projective $\Bbbk$-module is free.

\subsection{Bases for HOMs between some 1-morphisms} \label{sec_basis}

\begin{prop} Let $a,b,\delta \in \Z_{+}$ and $n\in \Z$.  Then the graded abelian group
$$\HOM_{\UcatD}(\cal{E}^{(a)}
\cal{F}^{(b)}\onen,\cal{E}^{(a+\delta)}\cal{F}^{(b+\delta)}\onen)$$
is a free module over $\END_{\UcatD}(\onen)$ with the basis
\begin{equation} \label{eq_basisX}
 \left\{\;\; \xy
 (0,0)*{\includegraphics[scale=0.6]{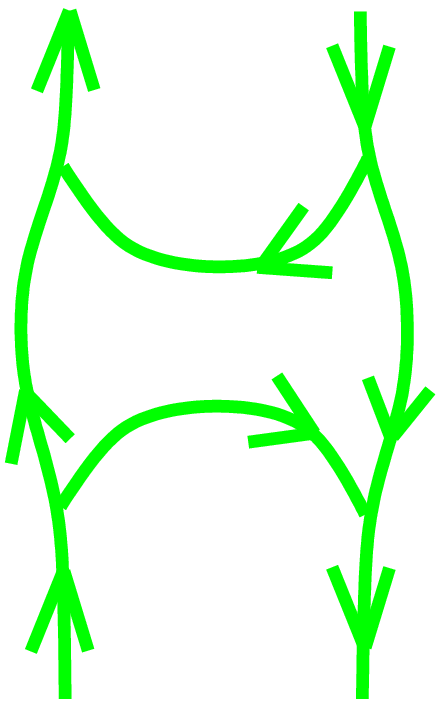}};
 (-13,2)*{\bigb{\pi_{\alpha} }};
 (13,2)*{\bigb{\pi_{\beta} }};
 (-3,6)*{\bigb{\pi_{\gamma} }};
 (-3,-3)*{\bigb{\pi_{\sigma} }};
 (-14,-16)*{a};(14,-16)*{b};
 (-17,16)*{a+\delta};(17,16)*{b+\delta};
 (-17,-4)*{a-j};(17,-4)*{b-j}; (4,-8)*{j}; (5,3)*{\delta+j};
  \endxy \;\; \right\}
\end{equation}
for $0 \leq j \leq \min(a,b)$ and all $\alpha \in P(a-j)$, $\beta\in
P(b-j)$, $\gamma \in P(\delta+j)$, $\sigma\in P(j)$.
\end{prop}

\begin{proof}
First we show that the set of 2-morphisms in \eqref{eq_basisX} spans $\HOM_{\UcatD}(\cal{E}^{(a)}
\cal{F}^{(b)}\onen,\cal{E}^{(a+\delta)}\cal{F}^{(b+\delta)}\onen)$ as a module over $\END_{\UcatD}(\onen)$.  Consider an arbitrary 2-morphism $f \maps \cal{E}^{(a)}\cal{F}^{(b)}\onen \to \cal{E}^{(a+\delta)}\cal{F}^{(b+\delta)}\onen\{t\}$
\begin{equation}
 \xy
 (0,0)*{\includegraphics[scale=0.5, angle=180]{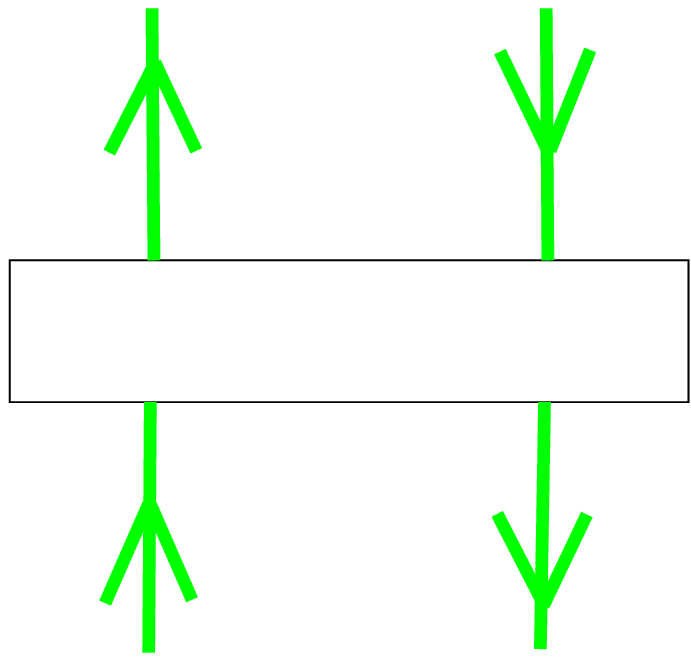}};
 (-14,-15)*{a};(-17,15)*{a+\delta}; (12,-15)*{b};(17,15)*{b+\delta};(22,8)*{n};
 (0,0)*{f};
  \endxy
\end{equation}
in $\UcatD$. By exploding the thick lines we can write it as
\begin{equation} \label{eq_explode_general}
\xy
 (0,0)*{\includegraphics[scale=0.5]{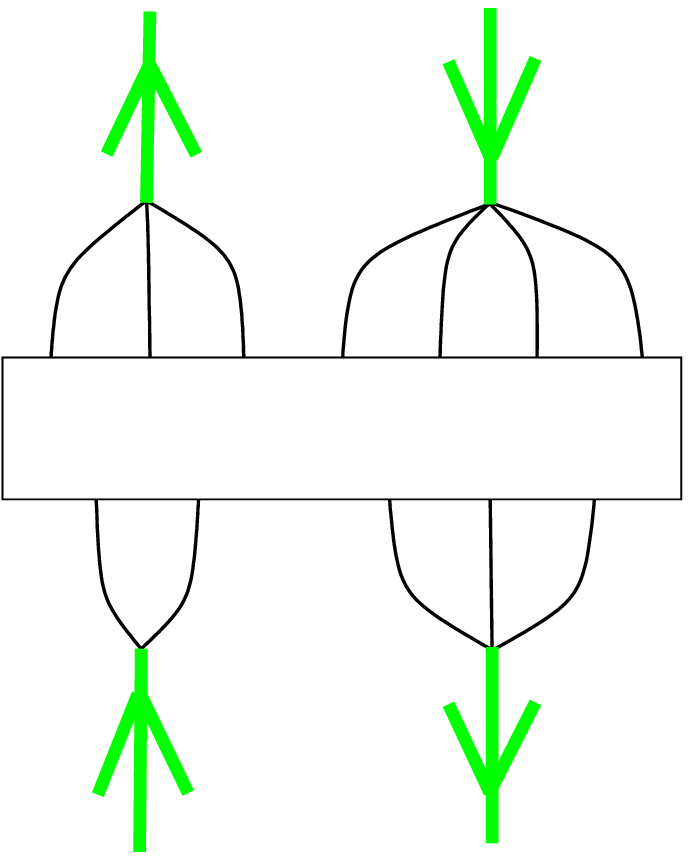}};
 (-14,-20)*{a};(-17,20)*{a+\delta}; (12,-20)*{b};(17,20)*{b+\delta};(22,8)*{n};
 (0,0)*{f'};
  \endxy
\end{equation}
for some $f' \maps \cal{E}^a\cal{F}^b\onen \to \cal{E}^{a+\delta}\cal{F}^{b+\delta}\onen$.

A basis for the $\HOM_{\UcatD}(\cal{E}^a\cal{F}^b\onen, \cal{E}^{a+\delta}\cal{F}^{b+\delta}\onen)$ is given by diagrams where strands have no self intersections, no two strands intersect more than once, all dots are confined to an interval on each arc, and all closed diagrams have been reduced to dotted bubbles with the same orientation that have been moved to the far right of a diagram, see \cite[Section 3.2]{KL3}. An example is shown below:
\[
\xy 0;/r.18pc/:
 (0,15)*{}; (20,-25) **\crv{(1,-6) & (20,-4)}?(0)*\dir{<}?(.6)*\dir{}+(0,0)*{\bullet};
 (8,15)*{}; (4,-25) **\crv{(8,6) & (4,0)}?(0)*\dir{<}?(.6)*\dir{}+(.2,0)*{\bullet};
 ?(0)*\dir{<}?(.75)*\dir{}+(.2,0)*{\bullet};?(0)*\dir{<}?(.9)*\dir{}+(0,0)*{\bullet};
 (28,-25)*{}; (12,-25) **\crv{(28,-10) & (12,-10)}?(0)*\dir{<};
  ?(.2)*\dir{}+(0,0)*{\bullet}?(.35)*\dir{}+(0,0)*{\bullet};
  (54,15)*{}; (50,-25) **\crv{(51,-10) & (49,-10)}?(1)*\dir{>}?(.35)*\dir{}+(.2,0)*{\bullet};;
 (36,15)*{}; (36,-25) **\crv{(34,6) & (35,-4)}?(1)*\dir{>};
 (28,15)*{}; (42,-25) **\crv{(28,6) & (42,-4)}?(1)*\dir{>};
 (42,15)*{}; (20,15) **\crv{(42,5) & (20,5)}?(1)*\dir{>};
 (48,15)*{}; (14,15) **\crv{(46,-18) & (20,-10)}?(1)*\dir{>};
 (64,5)*{\cbub{n-1+\alpha_2}{}};
 (86,5)*{\cbub{n-1+\alpha_4}{}};
 (64,-15)*{\cbub{n-1+\alpha_1}{}};
 (87,-15)*{\cbub{n-1+\alpha_3}{}};
 (105,-8)*{\cbub{n-1+\alpha_5}{}};
 (72,18)*{n};
 (41,18)*{\overbrace{\hspace{.85in}}};(41,22)*{\scs b+\delta};
 (10,18)*{\overbrace{\hspace{.65in}}};(10,22)*{\scs a+\delta};
 (39,-28)*{\underbrace{\hspace{.75in}}};(39,-32)*{\scs b};
 (12,-28)*{\underbrace{\hspace{.55in}}};(12,-32)*{\scs a};
 \endxy
\]
However, the splitters in \eqref{eq_explode_general} imply that \eqref{eq_explode_general} reduces to diagrams of the form
\begin{equation}
\xy
 (0,0)*{\includegraphics[scale=0.6]{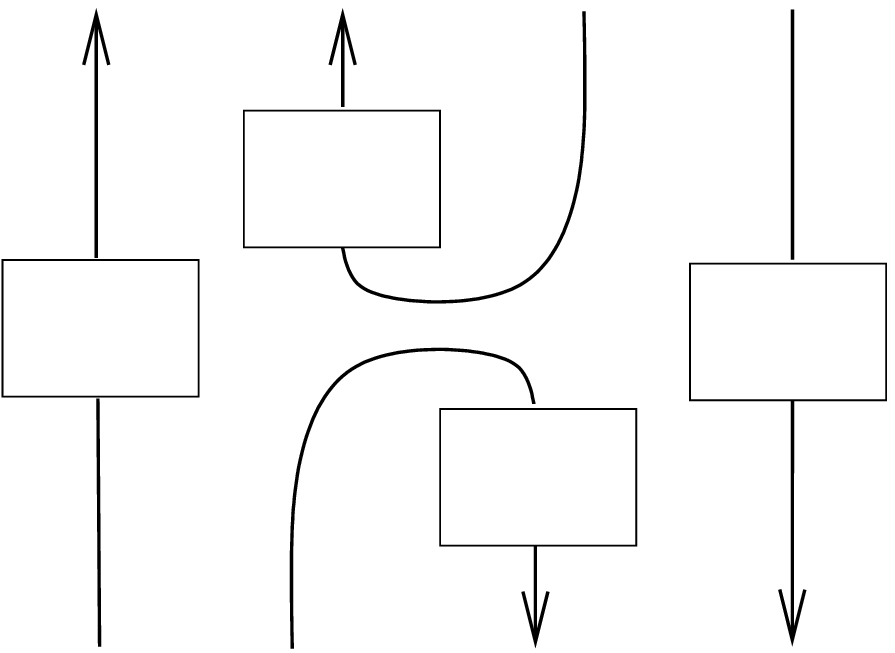}};
 (-28,-17)*{a-j};(28,17)*{b-j};
 (-12,-17)*{j};(29,8)*{n};
 (14,17)*{\delta+j};
 (-6,9)*{p_2}; (6,-9)*{p_3};(-21,0)*{p_1}; (21,0)*{p_4};
  (62,-15)*{\underbrace{\hspace{1.9in}}};
 (62,-20)*{\text{bubble monomial}};
  (45,-5)*{\cbub{n-1+\alpha_2}{}};
 (65,-5)*{\cbub{n-1+\alpha_4}{}};
 (45,15)*{\cbub{n-1+\alpha_1}{}};
 (65,15)*{\cbub{n-1+\alpha_3}{}};
 (80,8)*{\cbub{n-1+\alpha_5}{}};
  \endxy
\end{equation}
where $p_1 \in \Z[x_1,\dots, x_{a-j}]$, $p_2 \in \Z[x_1,\dots,x_{\delta+j}]$,  $p_3 \in \Z[x_1,\dots,x_{j}]$, and  $p_4 \in \Z[x_1,\dots,x_{b-j}]$.  Therefore, using the associativity of splitters and the results of Section~\ref{sec_nilHecke} any 2-morphism $f \maps \cal{E}^{(a)}\cal{F}^{(b)}\onen \to \cal{E}^{(a+\delta)}\cal{F}^{(b+\delta)}\onen$ can be written as a linear combination of diagrams of the form \eqref{eq_basisX} as a module over $\END_{\UcatD}(\onen)$.

To see that the spanning set in \eqref{eq_basisX} is a basis, observe that
 \begin{equation}
 \deg \left( \;\; \xy
 (0,0)*{\includegraphics[scale=0.6]{figs/basis.eps}};
 (-14,-16)*{a};(14,-16)*{b};
 (-17,16)*{a+\delta};(17,16)*{b+\delta};
 (-17,-4)*{a-j};(17,-4)*{b-j}; (4,-8)*{j}; (5,3)*{\delta+j};
  \endxy \;\; \right) \quad = \quad
2j(j+b-a-n)+\delta^2+\delta(b-a-n+2j).
\end{equation}
But using the formula for the semilinear form given in
\cite[Proposition 2.8]{Lau1} we have
\begin{eqnarray}
 \rkq \;\HOM_{\UcatD}(\cal{E}^{(a)}
\cal{F}^{(b)}\onen,\cal{E}^{(a+\delta)}\cal{F}^{(b+\delta)} \onen)=
\rkq\; \HOM_{\UcatD}(\cal{E}^{(a+\delta)}\cal{F}^{(b+\delta)}\onen,\cal{E}^{(a)}
\cal{F}^{(b)} \onen) \nn \\ =\sum_{j=0}^{\min(a,b)}
 q^{2j(j+b-a-n)+\delta^2+\delta(b-a-n+2j)}
 g(a-j)g(b-j)g(\delta+j)g(j),
 \nn
 \end{eqnarray}
where $g(x) := \prod_{j=1}^x\frac{1}{ (1-q^{2j})}
 =\sum_{\alpha\in P(x)}q^{\deg(\pi_{\alpha})}$. Summands in the formula match degrees of diagrams in \eqref{eq_basisX}.
\end{proof}

Reflecting the basis in \eqref{eq_basisX} across the horizontal axis and inverting the orientation gives a basis $\HOM_{\UcatD}(\cal{E}^{(a+\delta)}\cal{F}^{(b+\delta)}\onen,\cal{E}^{(a)}
\cal{F}^{(b)}\onen)$. Similarly, by reflecting across the vertical axis one obtains bases for  graded vector spaces $$
\HOM_{\UcatD}(\cal{F}^{(b+\delta)}\cal{E}^{(a+\delta)}\onen,
\cal{F}^{(b)}\cal{E}^{(a)}\onen) \quad  \text{and} \quad \HOM_{\UcatD}(\cal{F}^{(b)}\cal{E}^{(a)}
\onen,\cal{F}^{(b+\delta)}\cal{E}^{(a+\delta)}\onen).$$

\begin{lem}
For $x,y \in {}_m\B_{n}$ either $\cal{E}(x) = \cal{E}^{(a)}\cal{F}^{(b)}\onen$ and $\cal{E}(y)=\cal{E}^{(a+\delta)}\cal{F}^{(b+\delta)}\onen$,  or $\cal{E}(x) = \cal{F}^{(b)}\cal{E}^{(a)}\onen$ and $\cal{E}(y)=\cal{F}^{(b+\delta)}\cal{E}^{(a+\delta)}\onen$ for $a,b \in \N$ and $\delta \in \Z$.
\end{lem}

\begin{proof}
The proof is by a direct computation, see \cite[Lemma 2.6]{Lau1}.
\end{proof}

In view of the lemma we get a basis of
\begin{equation}
  \HOM_{\UcatD}\left(\cal{E}(x),\cal{E}(y)\right)
\end{equation}
for any canonical basis vectors $x,y \in {}_m\B_{n}$ as a free $\END_{\UcatD}(\onen)$-module via diagrams in \eqref{eq_basisX} or diagrams obtained from these by suitable symmetries (reflections or rotations).

%

\bigskip


%


\vspace{0.1in}

\noindent M.K.: { \sl \small Department of Mathematics, Columbia University, New
York, NY 10027} \newline \noindent {\tt \small email: khovanov@math.columbia.edu}

\vspace{0.1in}

\noindent A.L.:  { \sl \small Department of Mathematics, Columbia University, New
York, NY 10027} \newline \noindent
  {\tt \small email: lauda@math.columbia.edu}

\vspace{0.1in}

\noindent M.M.: { \sl \small Departamento de Matem\'atica, Universidade do Algarve, Campus de Gambelas, 8005-139 Faro, Portugal and CAMGSD, Instituto Superior T\'ecnico, Av. Rovisco Pais, 1049-001 Lisboa, Portugal} \newline
{\sl \small }\noindent {\tt \small email: mmackaay@ualg.pt}

\vspace{0.1in}

\noindent M.S.:  { \sl \small Instituto de Sistemas e Rob\'otica and CAMGSD, Instituto Superior T\'ecnico, Av. Rovisco Pais, 1049-001 Lisboa, Portugal} \newline \noindent
  {\tt \small email: mstosic@math.ist.utl.pt}

\end{document}